# Data Sampling Strategies in Stochastic Algorithms for Empirical Risk Minimization

*Dominik Csiba*

Doctor of Philosophy
University of Edinburgh
2017

# Declaration

I declare that this thesis was composed by myself and that the work contained therein is my own, except where explicitly stated otherwise in the text.

<div style="text-align: right;">(*Dominik Csiba*)</div>

# Acknowledgements

First and foremost, I would like to express my gratitude towards my first supervisor Peter Richtárik. His genuine interest in my progress throughout my academic years has pushed me further than I would ever expect. On top of his guidance in research, he gave me a lot of teaching advice, which I appreciate a lot.

I would like to thank my second supervisor Charles Sutton and my PCDS scholarship mentors Mike Davies and Ilias Diakonikolas for numerous discussions on topics related to my research field, especially for the connections and insights to machine learning, compressed sensing, and statistics.

For all the research and non-research discussions I would like to thank to the past and current members of our research group, namely: Nicolas Loizou, Filip Hanzely, Robert Gower, Ademir Ribeiro, Zheng Qu – for her collaboration on my very first paper – and last but not least Jakub Konečný - without whom I would not pursue my PhD degree in the first place.

I am very grateful to the School of Mathematics and all its members for making my stay so smooth. Especially I would like to thank Gill Law for her administrative magic, and Joan Simon and Iain Gordon for their patience and cooperation regarding various concessions I had to ask throughout my study – it really changed my life.

I would like to say a huge thank you to all the people from Edinburgh who made my stay more pleasant. Notably, my office mates Adria Caballe Mestres and Hattice Kursungecmez, my lunch mates Wenyi Qin, Ivet Galabova, Xavier Cabezas, Tom Byrne, Saranthorn Phusinga, the people from chess society, especially Zorbey Turkalp, Motiejus Gudenas, Julius Schwartz, Adam Bremner, and Lilli Hahn, the members of the Czech-Slovak society, in particular Ivan Bartoš, Matúš Falis, Juraj Labant, Jano Horváth and Benino Číčel, and my flatmates Soňa Galovičová and Jozef Mokrý.

Further, I would like to thank to the researchers who deepened my knowledge of various areas of optimization and mathematics in general. Notably, I would like to thank Julien Mairal and his group in INRIA Grenoble for inviting me for a research visit. Also, I am very grateful for my mentors during my internship at Amazon Berlin, namely Jim Huang, Cedric Archambeau, and especially Rodolphe Jenatton. For other random encounters and discussions I would like to thank Aurelien Lucchi, Martin Jaggi, Ohad Shamir, Nathan Srebro, Shai Shalev-Shwartz, Mark Schmidt, Joseph Salmon, Martin Takáč, and Petros Drineas.

Finally, I would like to thank my family Jarmila, Viktor and Peter, and my fiancée Veronika for all their love and support.

# Abstract

Gradient descent methods and especially their stochastic variants have become highly popular in the last decade due to their efficiency on big data optimization problems. In this thesis we present the development of data sampling strategies for these methods. In the first four chapters we focus on four views on the sampling for convex problems, developing and analyzing new state-of-the-art methods using non-standard data sampling strategies. Finally, in the last chapter we present a more flexible framework, which generalizes to more problems as well as more sampling rules.

In the first chapter we propose an adaptive variant of stochastic dual coordinate ascent (SDCA) for solving the regularized empirical risk minimization (ERM) problem. Our modification consists in allowing the method to adaptively change the probability distribution over the dual variables throughout the iterative process. AdaSDCA achieves a provably better complexity bound than SDCA with the best fixed probability distribution, known as importance sampling. However, it is of a theoretical character as it is expensive to implement. We also propose AdaSDCA+: a practical variant which in our experiments outperforms existing non-adaptive methods.

In the second chapter we extend the dual-free analysis of SDCA, to arbitrary mini-batching schemes. Our method is able to better utilize the information in the data defining the ERM problem. For convex loss functions, our complexity results match those of QUARTZ, which is a primal-dual method also allowing for arbitrary mini-batching schemes. The advantage of a dual-free analysis comes from the fact that it guarantees convergence even for non-convex loss functions, as long as the average loss is convex. We illustrate through experiments the utility of being able to design arbitrary mini-batching schemes.

In the third chapter we study importance sampling of minibatches. Minibatching is a well studied and highly popular technique in supervised learning, used by practitioners due to its ability to accelerate training through better utilization of parallel processing power and reduction of stochastic variance. Another popular technique is importance sampling – a strategy for preferential sampling of more important examples also capable of accelerating the training process. However, despite considerable effort by the community in these areas, and due to the inherent technical difficulty of the problem, there is no existing work combining the power of importance sampling with the strength of minibatching. In this chapter we propose the first importance sampling for minibatches and give simple and rigorous complexity analysis of its performance. We illustrate on synthetic problems that for training data of certain properties, our sampling can lead to several orders of magnitude improvement in training time. We then test the new sampling on several popular datasets, and show that the improvement can reach an order of magnitude.

In the fourth chapter we ask whether randomized coordinate descent (RCD) methods should be applied to the ERM problem or rather to its dual. When the number of examples ($n$) is much larger than the number of features ($d$), a common strategy is to apply RCD to the dual problem. On the other hand, when the number of features is much larger than the number of examples, it makes sense to apply RCD directly to the primal problem. In this paper we provide the first joint study of these two approaches when applied to L2-regularized ERM. First, we show through a rigorous analysis that for dense data, the above intuition is precisely correct. However, we find that for sparse and structured data, primal RCD can significantly outperform dual RCD even if $d \ll n$, and vice versa, dual RCD can be much faster than primal RCD even if $n \gg d$. Moreover, we show that, surprisingly, a single sampling strategy minimizes both the (bound on the) number of iterations and the overall expected complexity of RCD. Note that the latter complexity measure also takes into account the average cost of the iterations, which depends on the structure and sparsity of the data, and on the sampling strategy employed. We confirm our theoretical predictions using extensive experiments with both synthetic and real data sets.

In the last chapter we introduce two novel generalizations of the theory for gradient descent type methods in the proximal setting. Firstly, we introduce the proportion function, which we further use to analyze all the known block-selection rules for coordinate descent methods under a single framework. This framework includes randomized methods with uniform, non-uniform or even adaptive sampling strategies, as well as deterministic methods with batch, greedy or cyclic selection rules. We additionally introduce a novel block selection technique called greedy minibatches, for which we provide competitive convergence guarantees. Secondly, the whole theory of strongly-convex optimization was recently generalized to a specific class of non-convex functions satisfying the so-called Polyak-Łojasiewicz condition. To mirror this generalization in the weakly convex case, we introduce the Weak Polyak-Łojasiewicz condition, using which we give global convergence guarantees for a class of non-convex functions previously not considered in theory. Additionally, we give local convergence guarantees for an even larger class of non-convex functions satisfying only a certain smoothness assumption. By combining the two abovementioned generalizations we recover the state-of-the-art convergence guarantees for a large class of previously known methods and setups as special cases of our framework. Also, we provide new guarantees for many previously not considered combinations of methods and setups, as well as a huge class of novel non-convex objectives. The flexibility of our approach offers a lot of potential for future research, as any new block selection procedure will have a convergence guarantee for all objectives considered in our framework, while any new objective analyzed under our approach will have a whole fleet of block selection rules with convergence guarantees readily available.

# Contents

# Chapter 1

# Introduction

Stochastic optimization methods have been around since the famous paper of Robbins and Monro in 1951 [61] and have enjoyed a lot of success ever since. The above work introduces the well-known method called **S**tochastic **G**radient **D**escent (SGD). The inherent simplicity of SGD caused that it is successfully applied in many areas of science and it is still one of the most important baselines for large-scale optimization problems, even after more than 60 years.

Although SGD was sufficient to solve a lot of challenging problems, with the increasing size of the processed datasets, it has become too slow. Under these circumstances, the optimization community was posed with a new challenge – how to devise methods able to accurately solve problems over arbitrarily large datasets in a reasonably short time. With the arrival of the so-called "big data" era, this challenge has become one of the top priorities of the community.

In the following we will focus our attention to one of the fundamental problems of machine learning – supervised learning. The goal of supervised learning is to learn a predictor which is going to accurately predict the labels of given examples, e.g., identify the word a person said. Formally, the problem can be defined in the following way. Let $\mathcal{X}$ be a space of examples and $\mathcal{Y}$ the corresponding space of labels. Our goal is to find a predictor $h^* : \mathcal{X} \to \mathcal{Y}$, such that for any example-label pair drawn from the example-label space $(\mathbf{x}, y) \sim \mathcal{X} \times \mathcal{Y}$, we will have $h^*(\mathbf{x}) = y$. To write this down as a learning problem, we can formulate the following

$$h^* = \arg\min_{h:\mathcal{X}\to\mathcal{Y}} \left\{ \mathbf{E}_{(\mathbf{x},y)\sim\mathcal{X}\times\mathcal{Y}}[\mathbf{1}_{[h(\mathbf{x})=y]}] \right\},$$

where the function $\mathbf{1}_{[z]}$ is equal to 1 if $z$ is true and it is 0 otherwise. In words, we want to find a predictor $h^*$ such that it has the lowest *prediction error* on average for a randomly drawn example-label pair. Functions as $\mathbf{1}$ are usually referred to as loss/risk functions, as they penalize the predictor for mispredicting.

The above formulation has one key drawback – we do not know the underlying probability distribution of the example-label pairs. Indeed, we do not know what is the probability of seeing a specific image or hearing a specfic sound. To overcome this issue, the theory suggests to simulate the above quantity by randomly drawing $n$ samples $\{(\mathbf{x}^j, y_j)\}_{j=1}^n$ from $\mathcal{X} \times \mathcal{Y}$ in an i.i.d. fashion and approximating the above expectation by an empirical average

$$\mathbf{E}_{(\mathbf{x},y)\sim\mathcal{X}\times\mathcal{Y}}[\mathbf{1}_{[h(\mathbf{x})=y]}] \quad \approx \quad \frac{1}{n} \sum_{i=1}^{n} \mathbf{1}_{[h(\mathbf{x}^j)=y_j]}.$$

If the number of samples $n$ is large enough, we can expect the empirical average to approximate the expectation very closely. If that is the case, we can hope that the predictor $\hat{h}$ found by minimizing the empirical average

$$\hat{h} = \arg\min_{h:\mathcal{X}\to\mathcal{Y}} \left\{ \frac{1}{n} \sum_{j=1}^{n} \mathbf{1}_{[h(\mathbf{x}^j)=y_j]} \right\}$$

will be approximately the same as $h^*$. This idea is known as the problem of **E**mpirical **R**isk

**M**inimization (ERM) and it is one of the key optimization problems arising in machine learning [81]. However, the above problem has several caveats, which are well described in e.g. [68].

First issue is the function being minimized. In the current formulation, it is non-continuous and therefore difficult to minimize. Also, in a case that the label is a scalar value, e.g., temperature, it seems harsh to penalize a prediction off by 1 by the same amount as a prediction off by 100. The standard solution is to introduce a loss/risk function $\ell : \mathcal{Y} \times \mathcal{Y} \to \mathbb{R}_+$, which is going to measure, how good our current predictor $h$ is on a given example as the value $\ell(h(\mathbf{x}^j), y_j)$. A basic example of a loss function is the squared loss, which is simply $\ell(h(\mathbf{x}^j), y_j) = (h(\mathbf{x}^j) - y_j)^2$. The main property of a loss function is, that it should output zero, if the prediction is correct, and output a positive penalty based on how far off our prediction is in the case of an incorrect prediction. The squared loss satisfies these requirements exactly. More losses are going to be introduced later in this chapter in Section 1.2.3.

By introducing a loss function, we generalize our risk minimization problem to

$$\hat{h} = \arg\min_{h:\mathcal{X}\to\mathcal{Y}} \left\{ \frac{1}{n} \sum_{j=1}^{n} \ell(h(\mathbf{x}^j), y_j) \right\}.$$

It is indeed a generalization as by choosing $\ell(h(\mathbf{x}), y) = \mathbf{1}_{[h(\mathbf{x})=y_j]}$ we recover the original formulation. However, solving the above problem still does not guarantee as a good predictor and so we need to look at the next caveat.

To tackle the second issue, we need to restrict our potential predictors $h : \mathcal{X} \to \mathcal{Y}$ to a class $\mathcal{H}$. Indeed, what is known as the "No Free Lunch Theorem" [82] proves that we are unable to find a good predictor using the empirical risk, if we consider all the possible predictors from $\mathcal{X}$ to $\mathcal{Y}$. This makes intuitive sense, as there are many predictors which will correctly predict the labels for our finitely many examples, but only some of them will be able to predict correct labels for unseen examples as well. The formal proof in the above work uses the fact, that we are unable to distinguish between the two cases mentioned above, unless we restrict the class of predictors. Examples of such classes of predictors can be linear models or deep learning.

One important approach for class restriction is the so-called "regularization". The main idea behind the regularization is, that we prefer simple predictors over complex, as complex predictors are usually more likely to overfit to the training examples. This approach usually leads to predictors, which generalize better to unseen data. Formally, regularization is a function $R : \mathcal{H} \to \mathbb{R}_+$, which measures the complexity of a given predictor $h$ in some way. An example can be a regularizer, which outputs the number of characters needed to code up the predictor. This way, shorter (simpler) codes will be penalized by less than long and complex predictors. More examples of regularizers are introduced in Section 1.2.4.

Putting everything together, the empirical risk minimization problem can be written as

$$\hat{h} = \arg\min_{h\in\mathcal{H}} \left\{ \frac{1}{n} \sum_{j=1}^{n} \ell(h(\mathbf{x}^j), y_j) + R(h) \right\}.$$

In words, we are looking for the predictor $h^*$, which has the smallest penalty for misprediction combined with the penalty for complexity.

Usually, the classes of predictors $\mathcal{H}$ can be parametrized by a weight vector $\mathbf{w} \in \mathcal{W} \subseteq \mathbb{R}^d$, which transforms the problem from a functional space to a euclidean space

$$\mathbf{w}^* = \arg\min_{\mathbf{w}\in\mathcal{W}} \left\{ \frac{1}{n} \sum_{j=1}^{n} f_j(\mathbf{w}) + r(\mathbf{w}) \right\}, \qquad (1.1)$$

where $f_j(\mathbf{w}) := \ell(h(\mathbf{x}^j), y_j)$ and $r(\mathbf{w}) := R(h)$. The above problem (1.1) is often regarded as the ERM problem.

One important special case of the ERM problem is the linear ERM problem, i.e., the ERM problem specified to linear models. For linear models, we assume that the examples $\mathbf{x}^j \in \mathbb{R}^d$ are vectors and our class of predictors is given by

$$\mathcal{H}_{\text{linear}} \quad := \quad \left\{ h(\mathbf{x}) = \langle \mathbf{x}, \mathbf{w} \rangle + c, \quad \mathbf{w} \in \mathcal{W} \subseteq \mathbb{R}^d, c \in \mathbb{R} \right\}.$$

For easier notation, the above class is often redefined to contain the constant $c$ inside the vector $\mathbf{w}$ by adding an additional entry equal to 1 to all examples $\mathbf{x}^j$. Using this, the linear ERM problem takes the form

$$\mathbf{w}^* \quad = \quad \arg\min_{\mathbf{w} \in \mathcal{W}} \left\{ \frac{1}{n} \sum_{j=1}^n \phi_j(\langle \mathbf{x}^j, \mathbf{w} \rangle) + r(\mathbf{w}) \right\}, \tag{1.2}$$

where $\phi_j(s) := \ell(s, y_j)$.

Both the general ERM (1.1) and linear ERM (1.2) problems have their specific applications and it is important to study the ERM problem on both of these levels – which we do in this work. In Chapters 2, 4, and 5 we study the problem of linear ERM (1.2), while in Chapter 3 we propose a method for the general ERM (1.1). In Chapter 6 we go even further in generality, and we propose methods for the framework

$$\mathbf{x}^* \quad = \quad \arg\min_{\mathbf{x} \in X} \left\{ f(\mathbf{x}) + g(\mathbf{x}) \right\}, \tag{1.3}$$

where our only assumptions are, that $f : \mathbb{R}^d \to \mathbb{R}$ is smooth and $g : \mathbb{R}^d \to \mathbb{R}$ is separable. This framework can be specified to both general and linear ERM, as well as many other problems appearing outside of machine learning.

In this work we consider stochastic iterative methods for solving (1.1), (1.2), and (1.3). The methods we consider are stochastic for one of two reasons: i) we update only one random coordinate of our objective vector on each iteration ii) we randomly sample a single example and we update our objective vector only using information from this single example. In both of these cases, we randomly sample an index, either of a coordinate, or an example. One of the main goals of this work is to answer the question: how should we sample this index? In this work we show, that the most natural sampling – the uniform sampling, where we consider each index with equal probability – might be far from optimal. In each chapter we consider a separate approach to this question. These approaches are discussed in more detail in Section 1.4

In the rest of this chapter we will present the standard theory developed prior to this work, which is assumed to be known by the reader throughout the rest of the chapters. We also introduce the empirical risk minimization problem in more details together with specific forms of the losses and regularizers. Finally, we present notational conventions used in the rest of this work.

## 1.1 Technical glossary

We will introduce the basic technical glossary in our most general setup. We are interested in minimization problems $\min F(\mathbf{x})$, where the function $F : \mathbb{R}^d \to \bar{\mathbb{R}} := \mathbb{R} \cup \{+\infty\}$ can be decomposed as

$$F(\mathbf{x}) \quad = \quad f(\mathbf{x}) + g(\mathbf{x}), \tag{1.4}$$

where function $f : \mathbb{R}^d \to \mathbb{R}$ is assumed to be smooth and the function $g : \mathbb{R}^d \to \bar{\mathbb{R}}$ is usually assumed to be non-smooth and simple, sometimes even separable. Throughout this work we impose additional assumptions on the functions $f$ and $g$, which are further specified in the following.

### 1.1.1 Smooth functions

We say that a function $f : \mathbb{R}^d \to \mathbb{R}$ is $\mathbf{M}$-smooth, if there exists a positive definite matrix $\mathbf{M} \in \mathbb{R}^{d \times d}$, such that

$$f(\mathbf{x} + \mathbf{h}) \quad \leq \quad f(\mathbf{x}) + \langle \nabla f(\mathbf{x}), \mathbf{h} \rangle + \frac{1}{2} \mathbf{h}^\top \mathbf{M} \mathbf{h}, \tag{1.5}$$

for every $\mathbf{x}, \mathbf{h} \in \mathbb{R}^d$. If a function $f$ is $\beta \cdot \mathbf{I}$-smooth, where $\mathbf{I}$ is the identity matrix, we also say that the function is $\beta$-smooth. By standard theory, this is also equivalent to

$$\|\nabla f(\mathbf{x}) - \nabla f(\mathbf{y})\|_2 \leq \beta \|\mathbf{x} - \mathbf{y}\|_2, \tag{1.6}$$

for all $\mathbf{x}, \mathbf{y} \in \mathbb{R}^d$ and means that the gradient of $f$ is $\beta$-Lipschitz.

### 1.1.2 Convex and strongly-convex functions

Let $\lambda \geq 0$. A function $F : \mathbb{R}^d \to \bar{\mathbb{R}}$ is said to be $\lambda$-*strongly convex*, if for all $\mathbf{x}, \mathbf{y} \in \mathbb{R}^d$ and all $\beta \in [0, 1]$ we have

$$F(\beta \mathbf{x} + (1-\beta)\mathbf{y}) \leq \beta F(\mathbf{x}) + (1-\beta) F(\mathbf{y}) - \frac{\lambda \beta(1-\beta)}{2} \|\mathbf{x} - \mathbf{y}\|_2^2. \tag{1.7}$$

A function $F$ which is 0-strongly convex, i.e.,

$$F(\beta \mathbf{x} + (1-\beta)\mathbf{y}) \leq \beta F(\mathbf{x}) + (1-\beta) F(\mathbf{y}) \tag{1.8}$$

is usually referred to as *convex*.

Consider a smooth function $f : \mathbb{R}^d \to \mathbb{R}$. If $f$ is $\lambda$-strongly convex, it follows by standard theory that $\forall \mathbf{x}, \mathbf{h} \in \mathbb{R}^d$,

$$f(\mathbf{x} + \mathbf{h}) \geq f(\mathbf{x}) + \langle \nabla f(\mathbf{x}), \mathbf{h} \rangle + \frac{\lambda}{2} \|\mathbf{h}\|_2^2. \tag{1.9}$$

Similarly, a convex and smooth function $f$ satisfies $\forall \mathbf{x}, \mathbf{h} \in \mathbb{R}^d$

$$f(\mathbf{x} + \mathbf{h}) \geq f(\mathbf{x}) + \langle \nabla f(\mathbf{x}), \mathbf{h} \rangle. \tag{1.10}$$

### 1.1.3 Separable functions

We say that a function $g : \mathbb{R}^d \to \mathbb{R}$ is separable, if there exist $d$ scalar functions $\{g_1, \ldots, g_d\}$, such that

$$g(\mathbf{x}) = \sum_{i=1}^{d} g_i(x_i). \tag{1.11}$$

### 1.1.4 Fenchel conjugation

For a general function $\rho : \mathbb{R}^k \to \mathbb{R}$ we define the Fenchel (also called convex) conjugate $\rho^* : \mathbb{R}^k \to \mathbb{R}$ as the function

$$\rho^*(\mathbf{s}) := \sup_{\mathbf{u} \in \mathbb{R}^k} \{\langle \mathbf{s}, \mathbf{u} \rangle - \rho(\mathbf{u})\}. \tag{1.12}$$

It also follows from standard theory [28], that if the function $\rho$ is $1/\gamma$-smooth and convex, then the conjugate function $\rho^*$ is $\gamma$-strongly convex and vice-versa – if the function $\rho$ is $\gamma$-strongly convex, then the conjugate function $\rho^*$ is $1/\gamma$-strongly convex. For more details on Fenchel conjugation we refer the reader to the standard textbook [62].

## 1.2 Empirical risk minimization

In this section, we will look more closely at the linear ERM problem and its different forms. Let us arrange all the data vectors $\mathbf{x}^j$ as columns of a matrix $\mathbf{X} \in \mathbb{R}^{d \times n}$, i.e., $\mathbf{X}_{:j} = \mathbf{x}^j$. The vector of labels stays the same, specifically $\mathbf{y} \in \mathbb{R}^n$. An example-label pair therefore consists of $(\mathbf{X}_{:j}, y_j)$. Let us start with the basic formulation, which we already introduced earlier.

### 1.2.1 Primal problem

Let $\phi_j : \mathbb{R} \to \mathbb{R}$ be the loss function, and let $r : \mathbb{R}^d \to \bar{\mathbb{R}}$ be the regularizer. Then the linear ERM function can be written as a function $P : \mathbb{R}^d \to \mathbb{R}$ using the intuition behind (1.2) in the

following way

$$P(\mathbf{w}) := \frac{1}{n}\sum_{j=1}^{n}\phi_j(\langle \mathbf{X}_{:j}, \mathbf{w}\rangle) + \lambda r(\mathbf{w}). \tag{1.13}$$

Our goal is to find

$$\mathbf{w}^* = \arg\min_{\mathbf{w}\in\mathbb{R}^d} P(\mathbf{w}). \tag{1.14}$$

We will refer to this problem as the primal problem. In the next part, we introduce the dual formulation of the linear ERM.

## 1.2.2 Dual problem

We derive the dual problem using the primal formulation (1.13) from above. First, we reformulate (1.13) to the equivalent

$$\min_{\mathbf{w}\in\mathbb{R}^d, \mathbf{u}\in\mathbb{R}^n} \left[\frac{1}{n}\sum_{j=1}^{n}\phi_j(u_j) + \lambda r(\mathbf{w})\right], \quad \text{s.t.} \quad \langle \mathbf{X}_{:j}, \mathbf{w}\rangle = u_j \tag{1.15}$$

with its corresponding Lagrangian function given by

$$L(\mathbf{w},\mathbf{u},\boldsymbol{\alpha}) := \frac{1}{n}\sum_{j=1}^{n}\phi_j(u_j) + \lambda r(\mathbf{w}) + \frac{1}{n}\sum_{j=1}^{n}\alpha_j\bigl(u_j - \langle \mathbf{X}_{:j}, \mathbf{w}\rangle\bigr). \tag{1.16}$$

It follows that the dual function $D(\boldsymbol{\alpha})$ is then given by

$$\begin{aligned}
D(\boldsymbol{\alpha}) &:= -\inf_{\mathbf{w},\mathbf{u}} L(\mathbf{w},\mathbf{u},\boldsymbol{\alpha}) \\
&= -\inf_{\mathbf{w},\mathbf{u}} \left[\frac{1}{n}\sum_{j=1}^{n}\phi_j(u_j) + \lambda r(\mathbf{w}) + \frac{1}{n}\sum_{j=1}^{n}\alpha_j\bigl(u_j - \langle \mathbf{X}_{:j}, \mathbf{w}\rangle\bigr)\right] \\
&= -\inf_{\mathbf{w}}\left[\lambda r(\mathbf{w}) - \frac{1}{n}\sum_{j=1}^{n}\alpha_j\langle \mathbf{X}_{:j}, \mathbf{w}\rangle\right] - \inf_{\mathbf{u}}\left[\frac{1}{n}\sum_{j=1}^{n}(\phi_j(u_j) + \alpha_j u_j)\right] \\
&= -\lambda \sup_{\mathbf{w}}\left[\left\langle \frac{1}{\lambda n}\mathbf{X}\boldsymbol{\alpha},\mathbf{w}\right\rangle - r(\mathbf{w})\right] - \sup_{\mathbf{u}}\left[\frac{1}{n}((-\alpha_j)u_i - \phi_j(u_j))\right] \\
&\stackrel{(1.12)}{=} -\lambda r^*\left(\frac{1}{\lambda n}\mathbf{X}\boldsymbol{\alpha}\right) - \sum_{j=1}^{n}\phi^*(-\alpha_j). \tag{1.17}
\end{aligned}$$

The dual ERM problem can be formulated as

$$\boldsymbol{\alpha}^* = \arg\max_{\boldsymbol{\alpha}\in\mathbb{R}^n} D(\boldsymbol{\alpha}), \tag{1.18}$$

where $D(\boldsymbol{\alpha})$ is given by (1.17). We will refer to the vector $\boldsymbol{\alpha}$ as the dual variable. Note, that $\boldsymbol{\alpha}$ has the dimension of the number of examples instead of the dimensionality of the original space.

Often it is assumed, that the function $r$ is 1-strongly convex, so that the dual function $r^*$ is 1-smooth. Another common assumption is that the functions $\phi_j$ are $\frac{1}{\gamma}$-smooth, which would lead to the functions $\phi_j^*$ being $\gamma$-strongly convex. Some of the presented methods will be restricted to these cases.

It is well-known by standard theory [62] that the optimal solutions $\mathbf{w}^*$ and $\boldsymbol{\alpha}^*$ to the primal

(1.14) and dual problems (1.18) respectively, satisfy the optimality conditions given by

$$\begin{aligned} \mathbf{w}^* &= \nabla r^* \left( \frac{1}{\lambda n} \mathbf{X} \boldsymbol{\alpha}^* \right) \\ \alpha_j^* &= -\nabla \phi_j(\langle \mathbf{X}_{:j}, \mathbf{w}^* \rangle), \qquad \forall j \in \{1, \ldots, n\}. \end{aligned}$$

Observe that once we have the solution to the dual problem $\boldsymbol{\alpha}^*$, we can recover the solution to the primal problem $\mathbf{w}^*$ using the optimality conditions. The authors in [72] took advantage of this and showed that in many cases the dual ERM problem is actually easier to solve than the primal ERM problem.

In the following chapters we will take advantage of this as well and we will often aim to solve the dual problem instead of the primal.

### 1.2.3 Loss functions

In this section we introduce all the loss functions used thorough this thesis, with their corresponding dual formulations.

- **Quadratic/squared loss:** By far the most popular loss for regression problems, most notably linear regression. In the definition of the quadratic loss we assume that the labels are real numbers $y_j \in \mathbb{R}$.

$$\phi_j(s) = \frac{1}{2}(s - y_j)^2 \quad \text{and} \quad \phi_j^*(-u) = -u y_i + \frac{u^2}{2} \tag{1.19}$$

  The functions $\phi_j$ and $\phi_j^*$ are both 1-strongly convex and 1-smooth.

- **Logistic/sigmoid loss:** A very popular loss used for classification problems. On top of logistic regression, it is also one of the sensible choices for the activation function in neural networks. For the definition we assume that $y_j \in \{-1, 1\}$, i.e., we have a two-class classification problem.

$$\begin{aligned} \phi_j(s) &= \log(1 + \exp(-y_j s)) \\ \phi_j^*(-u) &= u y_j \log(u y_j) + (1 - u y_j) \log(1 - u y_j) \quad \text{for } u y_j \in [0, 1] \end{aligned} \tag{1.20}$$

  The above function $\phi_j$ is 1/4-smooth and convex, while $\phi_j^*$ is a 4-strongly convex and non-smooth function, defined only on $[0, 1]$.

- **Smoothed hinge loss:** As the name suggests, the smoothed hinge loss is a smoothed version of the hinge loss widely used in support vector machines [6]. The smoothing can be observed in the dual, where we add the term $\gamma u^2 / 2$ to the standard hinge loss dual formula.

$$\phi_j(s) = \begin{cases} 0 & s y_j > 1 \\ 1 - s y_j - \gamma/2 & s y_j < 1 - \gamma \\ \frac{1}{2\gamma}(1 - s y_j)^2 & \text{otherwise} \end{cases} \quad \text{and} \quad \phi_j^*(-u) = -u y_j + \frac{\gamma}{2} u^2 \tag{1.21}$$

  Again, for the definitions to make sense, we assume a two-class classification problem using $y_j \in \{-1, 1\}$. The above function $\phi_j$ is $1/\gamma$-smooth and convex, while $\phi_j^*$ is $\gamma$-strongly convex and defined only on $[0, 1]$.

### 1.2.4 Regularizers

In this part we introduce the most common regularizers used for the ERM problem, also with their dual versions.

- **$L_1$ regularizer:** The $L_1$ regularizer is a simple $L_1$ norm applied to the input vector

$$r(\mathbf{w}) = \|\mathbf{w}\|_1 = \sum_{i=1}^{d} |w_i| = \sum_{i=1}^{d} r_i(w_i). \tag{1.22}$$

One of the popular use cases of the $L_1$ regularizer is its combination with the quadratic loss which leads to an ERM problem commonly referred to as LASSO [78]. It is very popular in signal processing, as it has sparsity-inducing properties. However, it is non-smooth and non-strongly convex (but still convex), therefore it is more difficult to optimize. Its Fenchel conjugate is the $L_\infty$-norm defined below in (1.25).

- **$L_2$ regularizer:** By far the most popular regularizer is the squared $L_2$ norm of the input vector.

$$r(\mathbf{w}) \quad = \quad \frac{1}{2}\|\mathbf{w}\|_2^2 \quad = \quad \sum_{i=1}^{d} \frac{1}{2}w_i^2 \quad = \quad \sum_{i=1}^{d} r_i(w_i) \tag{1.23}$$

Problems containing an $L_2$ regularizers are often referred to as *ridge* problems, e.g., ridge regression in the case of quadratic loss or ridge logistic regression in the case of a logistic loss. The $L_2$ regularizer is 1-strongly convex and 1-smooth and therefore including it usually makes the problem easier to optimize. Due to this behavior it is often the default choice for many models. The $L_2$ regularizer a self-conjugate, i.e., its Fenchel conjugate is the $L_2$ regularizer itself.

- **Box constraints:**

$$r(\mathbf{w}) \quad = \quad \begin{cases} 0 & \max(\mathbf{w}) \leq b \\ \infty & \text{otherwise} \end{cases} \quad = \quad \sum_{i=1}^{d} \begin{cases} 0 & |w_i| \leq b \\ \infty & \text{otherwise} \end{cases} \quad = \quad \sum_{i=1}^{d} r_i(w_i) \tag{1.24}$$

The box constraint is constant on $[-b, b]^d$ and undefined otherwise. The main applications of the box constraints are in cases where the prediction vector has a hard constraint given by an environment.

- **$L_\infty$ regularizer:** As with all the $L_p$-norm regularizers, the $L_\infty$ regularizer is the $L_\infty$-norm applied to the input vector

$$r(\mathbf{w}) \quad = \quad \max_i \{w_i\}. \tag{1.25}$$

It is commonly encountered in the dual problems, as the $L_\infty$ regularizer is the Fenchel conjugate of the $L_1$ regularizer [91].

All of the above regularizers are separable, as it can be seen from their definitions. Also, the above regularizers can be combined to obtain more complex regularizers, e.g., by combining the $L_1$ and the $L_2$ regularizer we get the so-called elastic-net regularizer.

### 1.2.5 Setup for general ERM

To reach the general version of ERM (1.1), we assume that the data vectors are matrices. Specifically, let $\mathbf{X}_1, \ldots, \mathbf{X}_n$ be a sequence of data matrices with $\mathbf{X}_j \in \mathbb{R}^{d \times m}$, with their corresponding labels stored in $\mathbf{y} \in \mathbb{R}^n$. An example-label pair therefore consists of $(\mathbf{X}_j, y_j)$.

Consider the loss functions $\phi_j : \mathbb{R}^m \to \mathbb{R}$ and the regularizer $r : \mathbb{R}^d \to \mathbb{R}$. Then the general primal function $P : \mathbb{R}^d \to \mathbb{R}$ can be written as

$$P(\mathbf{w}) \quad := \quad \frac{1}{n} \sum_{j=1}^{n} \phi_j(\mathbf{X}_j^\top \mathbf{w}) + \lambda r(\mathbf{w}). \tag{1.26}$$

Note that in the special case of $m = d$ and $\mathbf{X}_j = \mathbf{I}_d$, i.e., the $d \times d$ identity matrix, we get the exact formulation (1.1). We will use the above formulation in Chapter 3.

For the loss functions $\phi_j$ and the regularizer $r$ we can consider the generalized dual function $D : \mathbb{R}^{m \times n} \to \mathbb{R}$, which takes the form

$$D(\mathbf{A}) \quad := \quad -\lambda r^* \left( \frac{1}{\lambda n} \sum_{j=1}^{n} \mathbf{X}_j \mathbf{A}_{:j} \right) - \frac{1}{n} \sum_{j=1}^{n} \phi_j^*(-\mathbf{A}_{:j}), \tag{1.27}$$

where the matrix $\mathbf{A} \in \mathbb{R}^{m \times n}$ contains the dual variables.

### 1.2.6 Datasets

For the numerical experiments in this thesis we use several publicly available[1] datasets. A summary of all the datasets is available in Table 1.1.

| Dataset | #samples | #features | sparsity |
|---|---|---|---|
| url | 2,396,130 | 3,231,962 | 0.04 % |
| news | 1,355,191 | 19,996 | 0.03% |
| cov1 | 581,012 | 54 | 22% |
| aloi | 108,000 | 129 | 24.6% |
| w8a | 49,749 | 301 | 4.2% |
| ijcnn1 | 35,000 | 23 | 60.1% |
| rcv1 | 20,242 | 47,237 | 0.2% |
| protein | 17,766 | 358 | 29.1% |
| mushrooms | 8124 | 112 | 18.8% |
| leukemia | 7,129 | 38 | 100.00% |
| dorothea | 800 | 100,000 | 0.9% |

Table 1.1: Summary of real data sets used in this thesis.

The sparsity column reports on the number of nonzeros divided by the number of all entries in the dataset, i.e.,

$$\text{sparsity} = \frac{\#\text{nonzeros}}{\#\text{samples} \times \#\text{features}}.$$

### 1.2.7 Notes on notation

Multi-dimensional objects are always typeset as bold, e.g., we have a vector $\mathbf{x} \in \mathbb{R}^d$ or a matrix $\mathbf{X} \in \mathbb{R}^{d \times n}$. We use slice notation to access elements of vectors and matrices, i.e., $x_i$ is the $i$-th element of vector $\mathbf{x}$, $X_{ij}$ is the element at the $i$-th row and $j$-th column. We index the dimensions from 1 to $d$ by the index $i$ and we index the data examples from 1 to $n$ by the index $j$. The data matrix has always the examples as columns, i.e., $\mathbf{X} \in \mathbb{R}^{d \times n}$ and $\mathbf{X}_{:j}$ is the $j$-th example, while $\mathbf{X}_{i:}$ is the $i$-th feature.

By $\langle \mathbf{u}, \mathbf{v} \rangle$ we denote the inner product between two vectors and by $\|\mathbf{u}\|_2 = \sqrt{\langle \mathbf{u}, \mathbf{u} \rangle}$ we denote the Euclidean norm of a vector. Also, $\|\mathbf{u}\|_1 = \sum_{i=1}^{d} |u_i|$ stands for the $L_1$ norm and $\|\mathbf{u}\|_0 = |\{i : u_i \neq 0\}|$ counts the number of nonzeros in a given vector/matrix. If not specified otherwise, $\|\mathbf{X}\|$ stands for the Frobenius norm for matrices, which is given by $\|\mathbf{X}\| = \sqrt{\sum_{i=1}^{d} \sum_{j=1}^{n} X_{ij}}$.

As a shorthand, we often denote the set of integers $\{1, \ldots, n\}$ as $[n]$. Similarly $[d] := \{1, \ldots, d\}$.

If not specified otherwise, the function $P : \mathbb{R}^d \to \mathbb{R}$ specifies the primal function for linear ERM (1.13) and $D : \mathbb{R}^n \to \mathbb{R}$ refers to its corresponding dual function (1.17). By $F : \mathbb{R}^d \to \bar{\mathbb{R}}$ we denote a general function of the form (1.4), where $f : \mathbb{R}^d \to \mathbb{R}$ and $g : \mathbb{R}^d \to \bar{\mathbb{R}}$ are the smooth and non-smooth parts, respectively.

A sequence of vectors is always indexed in their superscripts, while all the other sequences are indexed in their subscripts. This is due to the confusion between the $i$-th element of vector $\mathbf{x}$, i.e., $x_i$ and the $i$-th term in the sequence of vectors $\mathbf{x}$, i.e., $\mathbf{x}^i$, as well as the ambiguity between a scalar $s$ taken into the power $i$, i.e., $s^i$ and the $i$-th term in a sequence of scalars $s$, i.e., $s_i$.

## 1.3 Baseline methods

In this section we introduce the baseline methods, which we will compare against in the following chapters. Also, we introduce some additional notions connected to these methods, which will be used throughout this work.

---

[1] https://www.csie.ntu.edu.tw/~cjlin/libsvmtools/datasets/

### 1.3.1 Gradient descent

Suppose we want to find the minimizer $\mathbf{x}^*$ of a smooth and convex function $f : \mathbb{R}^d \to \mathbb{R}$. The baseline method for solving this problem dates back to Cauchy in 1847 [5]. The method is known as the gradient descent, and it takes the general form

$$\mathbf{x}^{k+1} \quad = \quad \mathbf{x}^{k+1} \quad = \quad \arg\min_{\mathbf{y} \in \mathbb{R}^d} \{\langle \nabla f(\mathbf{x}^k), \mathbf{y} - \mathbf{x}^k \rangle + \frac{\beta}{2}\|\mathbf{y} - \mathbf{x}^k\|_2^2\} \quad = \quad \mathbf{x}^k - \frac{1}{\beta}\nabla f(\mathbf{x}^k).$$

If $\beta$ is chosen as the smoothness parameter, then we can guarantee that if the number of iterations $K$ satisfies $\mathcal{O}(\beta/\epsilon) < K$, then $f(\mathbf{x}^K) - f(\mathbf{x}^*) < \epsilon$. This behavior is commonly referred to as sublinear convergence rate.

In the case that the function $f$ is also $\mu$-strongly convex, we can guarantee an $\epsilon$-convergence in $\mathcal{O}((\beta/\mu)\log(1/\epsilon))$ iterations, which is referred to as the linear convergence rate. The number $\beta/\mu$ is called the condition number and plays a key role in the convergence rates of many optimization methods. For more details on gradient descent and its variants we recommend [50].

### 1.3.2 Stochastic gradient descent

One of the the main disadvantages of using gradient descent for ERM-type problems is that we need to do a lot of computation for a single iteration. Indeed, we have to compute $d$ partial derivatives of $n$ different loss functions. The steady rise of stochastic methods was based on the idea that one can get a reasonable progress with much less information. An example of such a method was already mentioned in the introduction as the **S**tochastic **G**radient **D**escent [61]. In addition to SGD being commonly applied to purely stochastic problems, it is also used to minimize finite-sum functions, which have the form

$$F(\mathbf{x}) \quad := \quad \frac{1}{n}\sum_{j=1}^{n} f_j(\mathbf{x}).$$

We can observe, that the general ERM problem (1.1) fits into this framework. As an iteration of SGD, we randomly sample a function $f_j$ and we perform the update

$$\mathbf{x}^{k+1} \quad = \quad \mathbf{x}^k - \eta_k \nabla f_j(\mathbf{x}^k). \tag{1.28}$$

There are couple of things to say about this update rule. First, if we pick the function $f_j$ uniformly at random, we are going in the full-gradient direction in expectation, as $\mathbf{E}[\nabla f_j(\mathbf{x})] = F(\mathbf{x})$. Second, the sequence of stepsizes $\eta_k$ has to satisfy the common rule in stochastic approximation theory that $\sum_{k=1}^{\infty} \eta_k \to \infty$, while $\sum_{k=1}^{\infty} \eta_k^2 < \infty$. The most common sequence used for SGD satisfying these condition is $\eta_k = 1/k$. It is also worth noting, that in contrast to all the other methods mentioned in this section, SGD will not stay in the optimum once it reaches it. Indeed, for gradient descent type methods the update is equal to zero once we reach the optimum, but this is not the case for SGD. Due to this behavior, the SGD method is very good at converging to a medium accuracy, while for higher accuracy solutions other methods are prefered.

As for the convergence rates, for $\mu$-strongly convex functions SGD needs $\mathcal{O}(B^2/\mu\epsilon)$ iterations to converge to $\epsilon$-precision, where $B$ is a global bound on the magnitude of the individual gradient. This is a sublinear rate, which is clearly worse than the linear rate of gradient descent. However, for medium accuracy this is competitive, due to the iterations being $n$-times cheaper than the iterations of gradient descent. In the next sections we present alternative stochastic methods for dealing with expensive iterations.

### 1.3.3 Coordinate descent

Randomized coordinate descent plays an important role in the current optimization approaches. It was originally considered in [35] for quadratic objectives and analyzed for general functions in [49]. As the name suggests, it differs from gradient descent by only updating a single coordinate at each iteration, which is a lot cheaper. The main hope of coordinate descent is, that performing $d$ coordinate descent updates is going to make more progress than performing one full-gradient update.

In general, during an iteration $k$ we randomly select a coordinate $i_k$ and we perform the update

$$x_i^{k+1} \quad = \quad \begin{cases} x_i^k - \frac{1}{\beta_i}\nabla_i f(\mathbf{x}^k) & i = i_k \\ x^{k_i} & \text{otherwise} \end{cases},$$

where by $\nabla_i f(\mathbf{x})$ we denote the $i$-th partial derivative of $f$ at $\mathbf{x}$. By standard theory the stepsize parameters $\beta_i$ can be set to be the smoothness parameters of the individual partial derivatives, i.e.,

$$|\nabla_i f(\mathbf{x}) - \nabla_i f(\mathbf{y})| \quad \leq \quad \beta_i \|\mathbf{x} - \mathbf{y}\|_2$$

for each $i$. The number of iterations needed for coordinate descent to converge to $\epsilon$ accuracy is $\mathcal{O}(d \max_i \beta_i/\epsilon)$ in the case of a non-strongly convex objective and $\mathcal{O}(d \max_i \beta_i/\mu \log(1/\epsilon))$ in the case of a $\mu$-strongly convex objective. The above results hold in a case, that the coordinate was chosen uniformly at random. This was shown to be sub-optimal, as by sampling the coordinates according to the probabilities $p_i = \beta_i / \sum_j \beta_j$ one can replace the term $n \max_i \beta_i$ in the complexity by $\sum_i \beta_i$, which is a strictly better result.

Note, that the iteration of coordinate descent is $d$-times cheaper than the iteration of gradient descent – we need to compute only a single partial derivative instead of all the partial derivatives required for full gradient descent. Therefore, although the number of iterations required for coordinate descent to converge is larger than for gradient descent, the total computation time is usually smaller.

### 1.3.4 Proximal methods

The setup described in (1.4) – a function consisting of a smooth and a non-smooth part – is commonly reffered to as the proximal or composite setup. The baseline method for minimizing such a function is the so-called proximal gradient descent method, which takes the form of

$$\mathbf{x}^{k+1} \quad = \quad \arg\min_{\mathbf{y} \in \mathbb{R}^d} \{\langle \nabla f(\mathbf{x}^k), \mathbf{y} - \mathbf{x}^k \rangle + \frac{\beta}{2}\|\mathbf{y} - \mathbf{x}^k\|_2^2 + g(\mathbf{y})\}. \tag{1.29}$$

The above only makes sense if the minimization during each iteration is an easy problem – at the very least, easier than the original problem (1.4). This is indeed the case for many common non-smooth functions encountered in real-world applications.

As an example, suppose the non-smooth function takes the form

$$g(\mathbf{x}) \quad = \quad \begin{cases} 0 & \text{if } \|\mathbf{x}\|_2 \leq R \\ \infty & \text{otherwise} \end{cases}.$$

It is straightforward to show that in this case the proximal iteration (1.29) takes the form

$$\mathbf{x}^{k+1} \quad = \quad R \frac{\mathbf{x}^k - \frac{1}{\beta}\nabla f(\mathbf{x}^k)}{\|\mathbf{x}^k - \frac{1}{\beta}\nabla f(\mathbf{x}^k)\|_2},$$

which is indeed easy to compute.

As for the convergence rates, they remain of the same form as for smooth gradient descent, described in the previous section – for smooth and convex functions we have a sublinear rate of convergence, while for smooth and strongly-convex functions we have a linear convergence rate. Therefore, if we are able to solve (1.29) in a reasonably fast time, adding a non-smooth term to a smooth function does not make the minimization problem any more complex.

The proximal methods can be easily extended to coordinate descent [40]. The only additional requirement on the function is, that the function $g$ is separable (1.11). This is satisfied by many non-smooth functions, e.g. the dual functions $\phi_j^*$ in (1.17) as well as the non-smooth regularizers considered in Section 1.2.4. The separability is required for the update, which takes the form

$$x_i^{k+1} \;=\; \arg\min_{y \in \mathbb{R}} \left\{ \nabla_i f(\mathbf{x}^k)(y - x_i^k) + \frac{\beta}{2}(y - x_i^k)^2 + g_i(y) \right\}, \qquad (1.30)$$

where $i$ is a randomly sampled coordinate. Similarly as for the gradient descent counterpart of proximal methods, the complexity rates match the smooth case, with a single difference that instead of the quantity $\max_i \beta_i$ we have the standard full smoothness parameter $\beta$.

### 1.3.5 Stochastic dual coordinate ascent

In this part we present a final method, which plays a key role in the current theory of solving ERM (1.1), and especially linear ERM (1.2). As the title of the section suggests, we are going to apply the coordinate descent method to the dual of the ERM (1.17). However, in general the functions $\phi_j^*$ are non-smooth and therefore our update is of a proximal form. Directly plugging in the formulation (1.17) in to the proximal coordinate update (1.30) we get the update rule

$$\alpha_j^{k+1} \;=\; \arg\max_{y \in \mathbb{R}} \left\{ -\frac{1}{n}\mathbf{X}_{:j}^\top \nabla r^* \left(\tfrac{1}{\lambda n}\mathbf{X}\boldsymbol{\alpha}^k\right)(y - \alpha_j^k) - \frac{\|\mathbf{X}_{:j}\|_2^2}{2\lambda n^2} y^2 - \frac{1}{n}\phi_j^*(-(\alpha_j^k + y)) \right\},$$

where we maximize instead of minimization because of the dual problem (1.18) formulation. This update is reminiscent of SGD, as each update uses a single example, instead of a single feature. However, as it is simply a coordinate descent applied to the dual problem, SDCA retains the same rate of convergence as proximal coordinate descent, i.e., a linear rate in the strongly convex case and a sublinear rate in the non-strongly convex case. As we can see, this is superior to the rates of SGD and we can see similar behavior also in practice [72].

## 1.4 Chapter summaries

In this part we summarize the contribution in each of the chapters. Each of this chapters is based on a single paper/preprint, which is summarized in the table below.

| Chapter | 2 | 3 | 4 | 5 | 6 |
|---|---|---|---|---|---|
| Paper | [10] | [7] | [8] | [9] | (not yet online) |

During the time of my PhD studies I have co-authored an additional preprint [26] during my internship at the Machine Learning group at Amazon. However, the work does not fill well into the general topic of this work and therefore it is not part of it.

### 1.4.1 Summary of Chapter 2

In this chapter we propose a new randomized method for solving linear ERM (1.2) based on stochastic dual coordinate ascent, which we call *AdaSDCA*. The main novelty of AdaSDCA is that the sampling of dual variables is being adapted on each iteration. Although similar approaches were already considered before, our adaptive sampling is the first with a proven convergence guarantee superior to all known sampling approaches.

The adaptivity is motivated by the optimality conditions for the dual variables which are

$$\alpha_j^* \;=\; -\phi'(\langle \mathbf{X}_{:j}, \mathbf{w}^* \rangle), \qquad \forall j \in [n].$$

Using thess conditions we can measure how far a given dual variable is from its optimum. We do this by introducing the *dual residue*, which is simply

$$\kappa_j \;:=\; |\alpha_j + \phi'(\langle \mathbf{X}_{:j}, \mathbf{w}^* \rangle)|$$

for each $j \in [n]$. Observe that a given dual variable is optimal only if $\kappa_j = 0$. Therefore, our method takes into account the residues $\kappa_j$ and samples the dual variables with larger $\kappa_j$ more often.

However, the resulting sampling is only theoretical as the dual residues are expensive to compute on each iteration. With this in mind we introduced a practical version of our method which we call *AdaSDCA+*. The practical version does cheap updates to the sampling procedure which are approximating the adaptive sampling of AdaSDCA.

In numerical experiments we show the superiority of AdaSDCA in the convergence measured by the iteration count and also the superiority of AdaSDCA+ in wall-clock time over the sampling approaches previously regarded as the state-of-the-art.

The work which this chapter is based on [10] received an award[2] and the proposed method AdaSDCA+ is currently a part of the Tensorflow core library[3].

### 1.4.2 Summary of Chapter 3

In this part we propose a new variance-reduced primal method *dfSDCA* for solving the general ERM (1.1) with an $L_2$ regularizer. The "df" stands for "dual-free" and together with "SDCA" describes the main characteristics of the method – it is a primal-only version of SDCA.

The method is motivated by the optimality conditions of the dual problem which are

$$\boldsymbol{\alpha}_j^* = -\nabla \phi(\mathbf{X}_j^\top \mathbf{w}^*), \qquad \forall j \in [n]$$

for the general ERM. The method maintains auxilary variables $\boldsymbol{\alpha}_j \in \mathbb{R}^m$ corresponding to dual variables in the dual approaches and updates some of them at random on each iteration. The update to each variable is a convex combination between the previous iteration of the variable $\boldsymbol{\alpha}_j^k$ and the current best estimate for the variable which is $-\nabla \phi(\mathbf{X}_j^\top \mathbf{w}^k)$, resulting in the update rule

$$\boldsymbol{\alpha}_j^{k+1} = (1-\theta)\boldsymbol{\alpha}_j^k + \theta(-\nabla \phi(\mathbf{X}_j^\top \mathbf{w}^k)).$$

A very similar method for a slightly less general setup was already considered in [67], therefore we do not claim any significant novelty in the above described method. The main novelty in contrast to the previous method is that dfSDCA has been extended to work for arbitrary fixed sampling strategies.

A sampling strategy $\hat{S}$ is a set-valued random variable defined over all the subsets of $[n]$. Each subset $S$ corresponds to one way of choosing, which variables to update, e.g., $S = \{1\}$ means that we only update the first variable, while $S = [n]$ means that we update all the variables. An example of a sampling itself can be defined by $P(S = [n]) = 1$, i.e., we select all the variables with probability 1, which corresponds to batch updates. On the other hand, we can define a sampling by $P(S = \{j\}) = 1/n$ for each $j \in [n]$, which corresponds to picking a single variable uniformly at random. Also, we can take advantage of a multi-core machine and define a sampling which will uniformly at random distributed one variable per core, such that they can process the variables in parallel.

One can devise many useful samplings and our method gives a convergence rate for every such fixed sampling. Methods with arbitrary sampling were already considered before, but dfSDCA is the first primal variance-reduced method with such property. Also, our analysis extends to a case, where the individual loss functions can be non-convex, as long as the average function remains convex. As an example, this case can be observed in local convergences of deep learning models. Our method is the first to show convergence for arbitrary sampling strategies for such objectives.

To showcase the power of arbitrary sampling, we devise a new sampling called *chunking*. The main idea of chunking is to offer a load-balancing scheme for parallel computing with a convergence rate guarantee, analyzed through our framework. In numerical experiments we observe empirical speed-up in using chunking against standard parallel processing schemes.

---

[2]Best contribution award - 2$^{\text{nd}}$ place at Optimization and Big Data 2015, Edinburgh, award committee: Arkadi Nemirovski (Georgia Institute of Technology), Rodolphe Jenatton (Amazon Berlin)

[3]`https://github.com/tensorflow/tensorflow/blob/master/tensorflow/core/ops/sdca_ops.cc`

### 1.4.3 Summary of Chapter 4

One of the important developments for sampling strategies in the recent years was the introduction of the so-called importance sampling [49, 57, 88]. The importance sampling is currently the best fixed sampling scheme selecting a single variable at each iteration. It outperforms standard uniform sampling both theoretically and empirically.

The main idea of importance sampling is to precompute a specific quantity $v_j$ based on the data (usually $\|\mathbf{X}_{:j}\|_2^2$) for each of the variables and then sample the variables proportional to $v_j$. In practice this means, that the $j$-th variable will have the probability $v_j/\sum_i v_i$ of being sampled. The quantities $v_j$ also appear in the iteration complexity. For uniform sampling, the guarantee is proportional to $\max_j v_j$, while for importance sampling, the guarantee is proportional to the mean $\frac{1}{n}\sum_{j=1}^n v_j$. It follows that the importance sampling has always a lower iteration complexity than the uniform sampling.

Although importance sampling is already known for many years, it was always restricted to the case of serial sampling, i.e., sampling of a single variable on each iteration. In this chapter we proposed and analyzed the first ever importance sampling for minibatches. The main advantage of this approach over the standard serial importance sampling is that it can be used for parallel samplings with multi-core machines. By introducing this sampling, we fill in an important hole in the theory, as both importance sampling and parallel samplings were already known for some time, but they were never successfully combined.

The sampling we propose is not connected to any specific method. It can be combined with any method which allows for arbitrary sampling strategies, even the method dfSDCA we introduced in the previous chapter.

At the end of the chapter we perform numerical experiments to showcase that our approach indeed outperforms standard uniform parallel strategies. We observe speedup corresponding to the difference between $\max_j v_j$ and $\frac{1}{n}\sum_{j=1}^n v_j$ both in serial and parallel setting. This verifies that we managed to successfully generalize the importance sampling to parallel setting.

### 1.4.4 Summary of Chapter 5

In Chapter 5 we compare the performance of primal coordinate descent and dual coordinate ascent on the task of linear ERM. We are not the first ones who decided to do this comparison, but we are the first ones to do an analysis which takes into account the sparsity pattern of the data matrix.

The main difference in this comparison for dense and sparse data comes from the fact that the iteration cost can differ based on the variable being updated. Each variable corresponds to a single row or column of the data matrix, depending on whether we consider the primal or the dual approach. Each nonzero entry of this row/column is used for the computation of the update and therefore the iteration cost is proportional to the number of nonzeros in the given row/column of the data matrix. Previously, this was not taken into account and each iteration was considered to have the same iteration cost.

This small difference has a potentially large impact. For example, the importance sampling, which is considered to be the optimal fixed sampling for serial methods, might be selecting the rows/columns with a large amount of nonzeros much more often than the rows/columns with a small number of nonzeros. This can potentially lead to importance sampling to be sub-optimal for sparse matrices. For a proper comparison of the two approaches, we need to make sure that both of them are using the optimal samplings.

The first important result of this chapter is to show, surprisingly, that the importance sampling is optimal even if we have taken the sparsity pattern into account. Once we have the optimal sampling, we delve deeply into the comparison to draw some conclusions.

The first conclusion is that in the case of a dense data matrix, the comparison boils down to the comparison of the number of examples $n$ vs. the number of features $d$. If $n > d$, then the dual approach should be used and vice-versa, if $d > n$, then the primal approach should be used. This confirms the community belief about the comparison of these two methods.

As for the sparse case, the situation gets more involved. We show that the comparison boils to the difference between the data-dependent quantities

$$\sum_{i=1}^{d} \|\mathbf{X}_{i:}\|_2^2 \cdot \|\mathbf{X}_{i:}\|_0 \quad \text{vs} \quad \sum_{j=1}^{n} \|\mathbf{X}_{:j}\|_2^2 \cdot \|\mathbf{X}_{:j}\|_0.$$

We spend the rest of the chapter on analyzing this difference in many different scenarios and we draw some insightful conclusions. As an example we show that even for binary matrices, the primal methods can outperform dual methods when $n = \mathcal{O}(d^2)$, as far as the structure of the nonzero entries has a specific form – and vice-versa.

Lastly, we include numerical experiments to showcase that our theory meaningfully predicts the behaviour of these two methods in practice. We also include experiments on synthetic and real data in which $d \gg n$ and the dual method vastly outperforms the primal method – and also the other way around.

### 1.4.5 Summary of Chapter 6

Coordinate descent methods are well studied for convex and strongly convex objectives, but their behavior is not yet understood well for non-convex objectives. A recent work [29] revived an old idea [52] to show that for a certain class of non-convex objectives, we can achieve the same rate of convergence as for strongly-convex objectives. This realisation uses the fact that in the analysis of coordinate descent one only uses a certain consequence of strong convexity, in particular that there exists $\mu > 0$ such that

$$\|\nabla f(\mathbf{x})\|_2^2 \geq \mu(f(\mathbf{x}) - f(\mathbf{x}^*)), \qquad \forall \mathbf{x} \in \mathbb{R}^n,$$

where $\mathbf{x}^*$ is the minimizer of $f$. The above consequence is referred to as the *Polyak-Łojasiewicz* (PL) condtion. By assuming that the objective function only satisfies the PL condition and it is not strongly convex in general, we get the same convergence guarantee as for strongly convex functions, i.e., linear convergence of the form $\mathcal{O}(\log(1/\epsilon))$. One can also observe that the PL condition does not imply convexity and there are non-convex functions satisfying the above condition, which leads to a global convergence result in non-convex optimization.

While this is a remarkable result, it is still very restrictive, as it only applies to a very narrow class of functions. In the final chapter we aim to identify additional classes of non-convex functions, for which we can give global convergence guarantees. As a natural extension of the PL condition we introduce the so-called *Weak Polyak-Łojasiewicz* (WPL) condition which is a counterpart of the standard PL condition for (weakly) convex functions. Similarly as for the PL condition, we identify that the standard convergence analysis of coordinate descent for convex objectives only uses a consequence of convexity and not the convexity itself. This consequence is that there exists $\mu > 0$ such that

$$\|\nabla f(\mathbf{x})\| \cdot \|\mathbf{x} - \mathbf{x}^*\| \geq \sqrt{\mu}(f(\mathbf{x}) - f(\mathbf{x}^*)), \qquad \forall \mathbf{x} \in \mathbb{R}^n$$

and we refer to it as the WPL condition. Any smooth function satisfying the WPL condition will share the convergence rate with convex functions, i.e., sublinear convergence $\mathcal{O}(1/\epsilon)$.

Additionally, we provide convergence guarantees for an even more general class of non-convex functions, for which we only assume smoothness. However, for this class of functions we can no longer prove a global convergence rate guarantee, as they can potentially have many local extrema. Nevertheless, we show that the iterative process generated by a coordinate descent method will either convergence to the global optimum or to a point with a bounded gradient $\|\nabla f(\mathbf{x})\| \leq \epsilon$ in $\mathcal{O}(\frac{1}{\epsilon} \log(\frac{1}{\epsilon}))$ iterations.

Further, we introduce a new notion of a *proportion function*. For an **M**-smooth function (1.5) it is defined as

$$\theta(S, \mathbf{x}) \quad := \quad \frac{(\nabla_S f(\mathbf{x}))^\top \mathbf{M}_S^{-1} \nabla_S f(\mathbf{x})}{\|\nabla f(\mathbf{x})\|^2},$$

where $S \subset [n]$ is a given set of coordinates. The proportion function plays a key role in the analysis of various block-selection rules and procedures. Specifically, let $\{S_k\}_{k=0}^{\infty}$ be a sequence of sets, where $S_k$ is the set of coordinates updated on the $k$-th iteration. Let $\{\mathbf{x}^k\}_{k=0}^{\infty}$ be the corresponding sequence of iterates, i.e., we have $\mathbf{x}_{S_k}^{k+1} = \mathbf{x}_{S_k}^{k+1} - (\mathbf{M}_{S_k}^{-1})\nabla_{S_k} f(\mathbf{x}^{k+1})$ and $\mathbf{x}_{[n]\setminus S_k}^{k+1} = \mathbf{x}_{[n]\setminus S_k}^{k}$ otherwise. Additionally, let $\{c_k\}_{k=0}^{\infty}$ be a sequence of scalars such that $\theta(S_k, \mathbf{x}^k) \geq c_k$ or $\mathbf{E}[\theta(S_k, \mathbf{x}^k) \mid \mathbf{x}^k] \geq c_k$ for each $k$. In this case, we show that convergence is guaranteed once $\sum_{k=0}^{\infty} c_k \to \infty$. As we can see, this is condition is very weak and it is easily satisfied, e.g., it is sufficient to bound $c_k \geq c > 0$ for every $k$. As an example consider randomized coordinate descent with uniform sampling. We have

$$\mathbf{E}[\theta(S_k, \mathbf{x}^k) \mid \mathbf{x}^k] \quad = \quad \mathbf{E}_i[\theta(\{i\}, \mathbf{x}^k) \mid \mathbf{x}^k] \quad = \quad \frac{1}{n} \sum_{i=1}^{n} \frac{(\nabla_i f(\mathbf{x}^k))^2}{M_{ii} \|\nabla f(\mathbf{x}^k)\|^2} \quad \geq \quad \frac{1}{n \cdot \max_i \{M_{ii}\}},$$

where the last expression is a constant independent of the iterates and $k$ and therefore guarantees convergence.

Using the proportion function, we recover the state-of-the-art results of coordinate descent methods for previously known and also novel block selection rules. These rules include a class of randomized methods with uniform, non-uniform or adaptive sampling rules, as well as a class of methods with deterministic selection procedures as batch, greedy or cyclic.

Finally, we consider all the above theory for the general proximal setting (1.3), therefore it can be applied to both the primal (1.14) and the dual (1.18) problems in their full generality. The flexibility of the above approach offers a lot of potential for future research, as a new block selection procedure will have a convergence guarantee for all objectives considered in our framework, while a new objective analyzed under our approach will have a whole fleet of block selection rules with convergence guarantees readily available. We showcase this by introducing a new block-selection rule called *greedy minibatches*, which is a minibatch version of the greedy coordinate descent. By lower bounding its proportion function by a constant we give it competitive convergence guarantees in strongly-convex, convex, strongly PL, weakly PL, and general non-convex setups, which can be both smooth and proximal. The other way around, we introduced the weakly PL setup, which can be paired by any method with a bounded proportion function, e.g., gradient descent, coordinate descent with uniform or importance sampling, greedy coordinate descent, minibatch coordinate descent and many more. We believe that the flexibility of our approach opens up many future research opportunities.

# Chapter 2

# Stochastic Dual Coordinate Ascent with Adaptive Probabilities

## 2.1 Introduction

In this chapter we consider the regularized empirical risk minimization problem:

$$\min_{\mathbf{w}\in\mathbb{R}^d} \left[ P(\mathbf{w}) \quad := \quad \frac{1}{n}\sum_{j=1}^{n} \phi_j(\langle \mathbf{X}_{:j}, \mathbf{w}\rangle) + \lambda g(\mathbf{w}) \right]. \tag{2.1}$$

In the context of supervised learning, $\mathbf{w}$ is a linear predictor, the columns of the matrix $\mathbf{X}$, i.e., $\mathbf{X}_{:1}, \ldots, \mathbf{X}_{:n} \in \mathbb{R}^d$ are samples, $\phi_1, \ldots, \phi_n : \mathbb{R}^d \to \mathbb{R}$ are loss functions, $g : \mathbb{R}^d \to \mathbb{R}$ is a regularizer and $\lambda > 0$ a regularization parameter. Hence, we are seeking to identify the predictor which minimizes the average (empirical) loss $P(\mathbf{w})$.

We assume throughout that the loss functions $\phi_1, \ldots, \phi_n$ are $1/\gamma$-smooth for some $\gamma > 0$. That is, we assume they every $\phi_j$ is differentiable and have Lipschitz derivative with Lipschitz constant $1/\gamma$:

$$|\phi'_j(a) - \phi'_j(b)| \quad \leq \quad \frac{1}{\gamma}|a-b|$$

for all $a, b \in \mathbb{R}$. Moreover, we assume that $g$ is 1-strongly convex with respect to the L2 norm:

$$g(\mathbf{w}) \quad \leq \quad \beta g(\mathbf{w}_1) + (1-\beta)g(\mathbf{w}_2) - \frac{\beta(1-\beta)}{2}\|\mathbf{w}_1 - \mathbf{w}_2\|^2$$

for all $\mathbf{w}_1, \mathbf{w}_2 \in \text{dom}\, g$, $0 \leq \beta \leq 1$ and $\mathbf{w} = \beta\mathbf{w}_1 + (1-\beta)\mathbf{w}_2$.

The ERM problem (2.1) has received considerable attention in recent years due to its widespread usage in supervised statistical learning [72]. Often, the number of samples $n$ is very large and it is important to design algorithms that would be efficient in this regime.

**Modern stochastic algorithms for ERM.** Several highly efficient methods for solving the ERM problem were proposed and analyzed recently. These include primal methods such as SAG [64], SVRG [27], S2GD [32], SAGA [12], mS2GD [30] and MISO [42]. Importance sampling was considered in ProxSVRG [84] and S2CD [31].

**Stochastic Dual Coordinate Ascent.** One of the most successful methods in this category is *stochastic dual coordinate ascent (SDCA)*, which operates on the dual of the ERM problem (2.1):

$$\max_{\boldsymbol{\alpha}\in\mathbb{R}^n} \left[ D(\boldsymbol{\alpha}) \quad := \quad -\lambda g^*\left(\frac{1}{\lambda n}\sum_{j=1}^{n}\mathbf{X}_{:j}\alpha_j\right) - \frac{1}{n}\sum_{j=1}^{n}\phi_j^*(-\alpha_j) \right], \tag{2.2}$$

where $g^*$ and $\phi_j^*$ are the convex conjugates[1] of $g$ and $\phi_j$, respectively. Note that in dual problem, there are as many variables as there are samples in the primal: $\boldsymbol{\alpha} \in \mathbb{R}^n$.

SDCA in each iteration randomly selects a dual variable $\alpha_j$, and performs its update, usually via closed-form formula – this strategy is know as randomized coordinate descent. Methods based on updating randomly selected dual variables enjoy, in our setting, a linear convergence rate [72, 70, 75, 71, 88, 55]. These methods have attracted considerable attention in the past few years, and include SCD [69], RCDM [49], UCDC [59], ICD [77], PCDM [60], SPCDM [14], SPDC [86], APCG [37], RCD [44], APPROX [15], QUARTZ [55] and ALPHA [53]. Recent advances on mini-batch and distributed variants can be found in [38], [90], [58], [16], [79], [25], [43] and [41]. Other related work includes [47, 13, 1, 89, 17, 76]. We also point to [83] for a review on coordinate descent algorithms.

**Selection Probabilities.** Naturally, both the theoretical convergence rate and practical performance of randomized coordinate descent methods depends on the probability distribution governing the choice of individual coordinates. While most existing work assumes uniform distribution, it was shown by Richtárik and Takáč [59], Necoara et al. [45], Zhao and Zhang [88] that coordinate descent works for an arbitrary fixed probability distribution over individual coordinates and even subsets of coordinates [57, 55, 53, 54]. In all of these works the theory allows the computation of a fixed probability distribution, known as *importance sampling*, which optimizes the complexity bounds. However, such a distribution often depends on unknown quantities, such as the distances of the individual variables from their optimal values [59, 53]. In some cases, such as for smooth strongly convex functions or in the primal-dual setup we consider here, the probabilities forming an importance sampling can be explicitly computed [57, 88, 55, 53, 54]. Typically, the theoretical influence of using the importance sampling is in the replacement of the maximum of certain data-dependent quantities in the complexity bound by the average.

**Adaptivity.** Despite the striking developments in the field, there is virtually no literature on methods using an *adaptive* choice of the probabilities. We are aware of a few pieces of work; but all resort to heuristics unsupported by theory [21, 63, 2, 39], which unfortunately also means that the methods are sometimes effective, and sometimes not. *We observe that in the primal-dual framework we consider, each dual variable can be equipped with a natural measure of progress which we call "dual residue". We propose that the selection probabilities be constructed based on these quantities.*

*Outline:* In Section 2.2 we summarize the contributions of our work. In Section 2.3 we describe our first, theoretical methods (Algorithm 1) and describe the intuition behind it. In Section 2.4 we provide convergence analysis. In Section 2.5 we introduce Algorithm 2: an variant of Algorithm 1 containing heuristic elements which make it efficiently implementable. We conclude with numerical experiments in Section 2.6. Technical proofs and additional numerical experiments can be found in the appendix.

## 2.2 Contributions

We now briefly highlight the main contributions of this work.

**Two algorithms with adaptive probabilities.** We propose two new stochastic dual ascent algorithms: AdaSDCA (Algorithm 1) and AdaSDCA+ (Algorithm 2) for solving (2.1) and its dual problem (2.2). The novelty of our algorithms is in adaptive choice of the probability distribution over the dual coordinates.

**Complexity analysis.** We provide a convergence rate analysis for the first method, showing that *AdaSDCA enjoys better rate than the best known rate for SDCA with a fixed sampling* [88, 55]. The probabilities are proportional to a certain measure of dual suboptimality associated with each variable.

**Practical method.** AdaSDCA requires the same computational effort per iteration as the batch gradient algorithm. To solve this issue, we propose AdaSDCA+ (Algorithm 2): an efficient heuristic variant of the AdaSDCA. The computational effort of the heuristic method in a single iteration is low, which makes it very competitive with methods based on impor-tance sampling, such as IProx-SDCA [88]. We support this with computational experiments in Section 2.6.

---
[1] By the convex conjugate we mean the fenchel conjugate defined in Chapter 1

*Outline:* In Section 2.2 we summarize the contributions of our work. In Section 2.3 we describe our first, theoretical methods (AdaSDCA) and describe the intuition behind it. In Section 2.4 we provide convergence analysis. In Section 2.5 we introduce AdaSDCA+: a variant of AdaSDCA containing heuristic elements which make it efficiently implementable. We conclude with numerical experiments in Section 2.6. Technical proofs and additional numerical experiments can be found in the appendix.

## 2.3 The algorithm: AdaSDCA

It is well known that the optimal primal-dual pair $(\mathbf{w}^*, \boldsymbol{\alpha}^*) \in \mathbb{R}^d \times \mathbb{R}^n$ satisfies the following *optimality conditions*:

$$\mathbf{w}^* = \nabla g^* \left( \frac{1}{\lambda n} \mathbf{X} \boldsymbol{\alpha}^* \right) \tag{2.3}$$

$$\alpha_j^* = -\nabla \phi_j(\langle \mathbf{X}_{:j}, \mathbf{w}^* \rangle), \quad \forall j \in [n] := \{1, \ldots, n\}, \tag{2.4}$$

where $\mathbf{X}$ is the $d$-by-$n$ matrix with columns $\mathbf{X}_{:1}, \ldots, \mathbf{X}_{:n}$.

**Definition 2.1** (Dual residue). The *dual residue*, $\boldsymbol{\kappa}^t = (\kappa_1^t, \ldots, \kappa_n^t) \in \mathbb{R}^n$, associated with $(\mathbf{w}^t, \boldsymbol{\alpha}^t)$ is given by:

$$\kappa_j^t := \alpha_j^t + \nabla \phi_j(\langle \mathbf{X}_{:j}^\top \mathbf{w}^t \rangle). \tag{2.5}$$

Note, that $\kappa_i^t = 0$ if and only if $\alpha_j$ satisfies (2.4). This motivates the design of AdaSDCA (Algorithm 1) as follows: whenever $|\kappa_j^t|$ is large, the $j$th dual coordinate $\alpha_j$ is suboptimal and hence should be updated more often.

**Definition 2.2** (Coherence). We say that the probability vector $\mathbf{p}^t \in \mathbb{R}^n$ is *coherent* with the dual residue $\boldsymbol{\kappa}^t$ if for all $j \in [n]$ we have

$$\kappa_j^t \neq 0 \quad \Rightarrow \quad p_j^t > 0.$$

Alternatively, $\mathbf{p}^t$ is coherent with $\boldsymbol{\kappa}^t$ if for

$$J_t := \{ j \in [n] : \kappa_j^t \neq 0 \} \subseteq [n].$$

we have $\min_{j \in J_t} p_j^t > 0$.

---

**Algorithm 1** AdaSDCA

**Init:** $v_j := \|\mathbf{X}_{:j}\|_2^2$ for $j \in [n]$; $\boldsymbol{\alpha}^0 \in \mathbb{R}^n$; $\bar{\boldsymbol{\alpha}}^0 = \frac{1}{\lambda n} \mathbf{X} \boldsymbol{\alpha}^0$
**for** $t \geq 0$ **do**
  Primal update: $\mathbf{w}^t = \nabla g^*(\bar{\boldsymbol{\alpha}}^t)$
  Set: $\boldsymbol{\alpha}^{t+1} = \boldsymbol{\alpha}^t$
  Compute residue $\boldsymbol{\kappa}^t$: $\kappa_i^t = \alpha_j^t + \nabla \phi_j(\langle \mathbf{X}_{:j}, \mathbf{w}^t \rangle), \forall j \in [n]$
  Compute probability distribution $\mathbf{p}^t$ coherent with $\boldsymbol{\kappa}^t$
  Generate random $j_t \in [n]$ according to $\mathbf{p}^t$
  Compute: $\quad \Delta \alpha_{j_t}^t = \arg\max_{\Delta \in \mathbb{R}} \left\{ -\phi_{j_t}^*(-(\alpha_{j_t}^t + \Delta)) - \langle \mathbf{X}_{:j_t}, \mathbf{w}^t \rangle \Delta - \frac{v_{j_t}}{2\lambda n} |\Delta|^2 \right\}$
  Dual update: $\alpha_{j_t}^{t+1} = \alpha_{j_t}^t + \Delta \alpha_{j_t}^t$
  Average update: $\bar{\boldsymbol{\alpha}}^t = \bar{\boldsymbol{\alpha}}^t + \frac{\Delta \alpha_{j_t}}{\lambda n} \mathbf{X}_{:j_t}$
**end for**
**Output:** $\mathbf{w}^t, \boldsymbol{\alpha}^t$

---

AdaSDCA is a stochastic dual coordinate ascent method, with an adaptive probability vector $\mathbf{p}^t$, which could potentially change at every iteration $t$. The primal and dual update rules are exactly the same as in standard SDCA [72], which instead uses uniform sampling probability at every iteration and does not require the computation of the dual residue $\boldsymbol{\kappa}$.

Our first result highlights a key technical tool which ultimately leads to the development of good adaptive sampling distributions $\mathbf{p}^t$ in AdaSDCA. For simplicity we denote by $\mathbb{E}_t$ the expectation with respect to the random index $j_t \in [n]$ generated at iteration $t$.

**Lemma 2.3.** *Consider the AdaSDCA algorithm during iteration $t \geq 0$ and assume that $\mathbf{p}^t$ is coherent with $\boldsymbol{\kappa}^t$. Then*

$$\mathbb{E}_t\left[D(\boldsymbol{\alpha}^{t+1}) - D(\boldsymbol{\alpha}^t)\right] - \theta\left(P(\mathbf{w}^t) - D(\boldsymbol{\alpha}^t)\right) \geq -\frac{\theta}{2\lambda n^2} \sum_{j \in J_t} \left(\frac{\theta(v_j + n\lambda\gamma)}{p_j^t} - n\lambda\gamma\right)|\kappa_j^t|^2, \quad (2.6)$$

*for arbitrary*

$$0 \leq \theta \leq \min_{j \in J_t} p_j^t. \quad (2.7)$$

*Proof.* Lemma 2.3 is proved similarly to Lemma 2 in [88], but in a slightly more general setting. For completeness, we provide the proof in the appendix. □

Lemma 2.3 plays a key role in the analysis of stochastic dual coordinate methods [72, 88, 71]. Indeed, if the right-hand side of (2.6) is positive, then the primal dual error $P(\mathbf{w}^t) - D(\boldsymbol{\alpha}^t)$ can be bounded by the expected dual ascent $\mathbb{E}_t[D(\boldsymbol{\alpha}^{t+1}) - D(\boldsymbol{\alpha}^t)]$ times $1/\theta$, which yields the contraction of the dual error at the rate of $1 - \theta$ (see Theorem 2.7). In order to make the right-hand side of (2.6) positive we can take any $\theta$ smaller than $\theta(\boldsymbol{\kappa}^t, \mathbf{p}^t)$ where the function $\theta(\cdot, \cdot) : \mathbb{R}_+^n \times \mathbb{R}_+^n \to \mathbb{R}$ is defined by:

$$\theta(\boldsymbol{\kappa}, \mathbf{p}) \equiv \frac{n\lambda\gamma \sum_{j:\kappa_j \neq 0} |\kappa_j|^2}{\sum_{j:\kappa_j \neq 0} p_j^{-1} |\kappa_j|^2 (v_j + n\lambda\gamma)}. \quad (2.8)$$

We also need to make sure that $0 \leq \theta \leq \min_{j \in J_t} p_j^t$ in order to apply Lemma 2.3. A "good" adaptive probability $\mathbf{p}^t$ should then be the solution of the following optimization problem:

$$\max_{\mathbf{p} \in \mathbb{R}_+^n} \theta(\boldsymbol{\kappa}^t, \mathbf{p}) \quad (2.9)$$

$$\text{s.t.} \quad \sum_{j=1}^n p_i = j$$

$$\theta(\boldsymbol{\kappa}^t, \mathbf{p}) \leq \min_{j:\kappa_j^t \neq 0} p_j$$

A feasible solution to (2.9) is the *importance sampling* (also known as optimal serial sampling) $\mathbf{p}^{\text{imp}}$ defined by:

$$p_j^{\text{imp}} := \frac{v_j + n\lambda\gamma}{\sum_{k=1}^n (v_k + n\lambda\gamma)}, \quad \forall j \in [n], \quad (2.10)$$

which was proposed in [88] to obtain proximal stochastic dual coordinate ascent method with importance sampling (IProx-SDCA). The same optimal probability vector was also deduced, via different means and in a more general setting in [55]. Note that in this special case, since $\mathbf{p}^t$ is independent of the residue $\boldsymbol{\kappa}^t$, the computation of $\boldsymbol{\kappa}^t$ is unnecessary and hence the complexity of each iteration does not scale up with $n$.

It seems difficult to identify other feasible solutions to program (2.9) apart from $\mathbf{p}^{\text{imp}}$, not to mention solve it exactly. However, by relaxing the constraint $\theta(\boldsymbol{\kappa}^t, \mathbf{p}) \leq \min_{j:\kappa_j^t \neq 0} p_j$, we obtain an explicit optimal solution.

**Lemma 2.4.** *The optimal solution* $\mathbf{p}^*(\boldsymbol{\kappa}^t)$ *of*

$$\max_{\mathbf{p} \in \mathbb{R}^n_+} \quad \theta(\boldsymbol{\kappa}^t, \mathbf{p}) \tag{2.11}$$

$$\text{s.t.} \quad \sum_{j=1}^n p_j = 1$$

*is:*

$$(\mathbf{p}^*(\boldsymbol{\kappa}^t))_i = \frac{|\kappa_j^t|\sqrt{v_j + n\lambda\gamma}}{\sum_{k=1}^n |\kappa_k^t|\sqrt{v_k + n\lambda\gamma}}, \quad \forall k \in [n]. \tag{2.12}$$

*Proof.* The proof is deferred to the appendix. □

The suggestion made by (2.12) is clear: we should update more often those dual coordinates $\alpha_j$ which have large absolute dual residue $|\kappa_j^t|$ and/or large Lipschitz constant $v_j$.

If we let $\mathbf{p}^t = \mathbf{p}^*(\boldsymbol{\kappa}^t)$ and $\theta = \theta(\boldsymbol{\kappa}^t, \mathbf{p}^t)$, the constraint (2.7) may not be sastified, in which case (2.6) does not necessarily hold. However, as shown by the next lemma, the constraint (2.7) is not required for obtaining (2.6) when all the functions $\{\phi_j\}_j$ are quadratic.

**Lemma 2.5.** *Suppose that all* $\{\phi_j\}_j$ *are quadratic. Let* $t \geq 0$. *If* $\min_{j \in J_t} p_j^t > 0$, *then* (2.6) *holds for any* $\theta \in [0, +\infty)$.

The proof is deferred to the Appendix.

## 2.4 Convergence results

In this section we present our theoretical complexity results for AdaSDCA. The main results are formulated in Theorem 2.7, covering the general case, and in Theorem 2.11 in the special case when $\{\phi_j\}_{j=1}^n$ are all quadratic.

### 2.4.1 General loss functions

We derive the convergence result from Lemma 2.3.

**Proposition 2.6.** *Let* $t \geq 0$. *If* $\min_{j \in J_t} p_j^t > 0$ *and* $\theta(\boldsymbol{\kappa}^t, \mathbf{p}^t) \leq \min_{j \in J_t} p_j^t$, *then*

$$\mathbb{E}_t\left[D(\boldsymbol{\alpha}^{t+1}) - D(\boldsymbol{\alpha}^t)\right] \geq \theta(\boldsymbol{\kappa}^t, \mathbf{p}^t)\left(P(\mathbf{w}^t) - D(\boldsymbol{\alpha}^t)\right).$$

*Proof.* This follows directly from Lemma 2.3 and the fact that the right-hand side of (2.6) equals 0 when $\theta = \theta(\boldsymbol{\kappa}^t, \mathbf{p}^t)$. □

**Theorem 2.7.** *Consider AdaSDCA. If at each iteration* $t \geq 0$, $\min_{j \in J_t} p_j^t > 0$ *and* $\theta(\boldsymbol{\kappa}^t, \mathbf{p}^t) \leq \min_{j \in J_t} p_j^t$, *then*

$$\mathbb{E}[P(\mathbf{w}^t) - D(\boldsymbol{\alpha}^t)] \leq \frac{1}{\tilde{\theta}_t} \prod_{k=0}^t (1 - \tilde{\theta}_k)\left(D(\boldsymbol{\alpha}^*) - D(\boldsymbol{\alpha}^0)\right), \tag{2.13}$$

*for all* $t \geq 0$ *where*

$$\tilde{\theta}_t := \frac{\mathbb{E}[\theta(\boldsymbol{\kappa}^t, \mathbf{p}^t)(P(\mathbf{w}^t) - D(\boldsymbol{\alpha}^t))]}{\mathbb{E}[P(\mathbf{w}^t) - D(\boldsymbol{\alpha}^t)]}. \tag{2.14}$$

*Proof.* By Proposition 2.6, we know that

$$\begin{aligned}
\mathbb{E}[D(\boldsymbol{\alpha}^{t+1}) - D(\boldsymbol{\alpha}^t)] &\geq \mathbb{E}[\theta(\boldsymbol{\kappa}^t, \mathbf{p}^t)(P(\mathbf{w}^t) - D(\boldsymbol{\alpha}^t))] \\
&\stackrel{(2.14)}{=} \tilde{\theta}_t \, \mathbb{E}[P(\mathbf{w}^t) - D(\boldsymbol{\alpha}^t)] \\
&\geq \tilde{\theta}_t \, \mathbb{E}[D(\boldsymbol{\alpha}^*) - D(\boldsymbol{\alpha}^t)],
\end{aligned} \tag{2.15}$$

whence
$$\mathbb{E}[D(\boldsymbol{\alpha}^*) - D(\boldsymbol{\alpha}^{t+1})] \leq (1 - \tilde{\theta}_t)\mathbb{E}[D(\boldsymbol{\alpha}^*) - D(\boldsymbol{\alpha}^t)].$$

Therefore,
$$\mathbb{E}[D(\boldsymbol{\alpha}^*) - D(\boldsymbol{\alpha}^t)] \leq \prod_{k=0}^{t}(1 - \tilde{\theta}_k)\left(D(\boldsymbol{\alpha}^*) - D(\boldsymbol{\alpha}^0)\right).$$

By plugging the last bound into (2.15) we get the bound on the primal dual error:
$$\begin{aligned}
\mathbb{E}[P(\mathbf{w}^t) - D(\boldsymbol{\alpha}^t)] &\leq \frac{1}{\tilde{\theta}_t}\mathbb{E}[D(\boldsymbol{\alpha}^{t+1}) - D(\boldsymbol{\alpha}^t)] \\
&\leq \frac{1}{\tilde{\theta}_t}\mathbb{E}[D(\boldsymbol{\alpha}^*) - D(\boldsymbol{\alpha}^t)] \\
&\leq \frac{1}{\tilde{\theta}_t}\prod_{k=0}^{t}(1 - \tilde{\theta}_k)\left(D(\boldsymbol{\alpha}^*) - D(\boldsymbol{\alpha}^0)\right). \quad \square
\end{aligned}$$

As mentioned in Section 2.3, by letting every sampling probability $\mathbf{p}^t$ be the importance sampling (optimal serial sampling) $\mathbf{p}^{\text{imp}}$ defined in (2.10), AdaSDCA reduces to IProx-SDCA proposed in [88]. The convergence theory established for IProx-SDCA in [88], which can also be derived as a direct corollary of our Theorem 2.7, is stated as follows.

**Theorem 2.8** ([88]). *Consider AdaSDCA with $\mathbf{p}^t = \mathbf{p}^{imp}$ defined in (2.10) for all $t \geq 0$. Then*
$$\mathbb{E}[P(\mathbf{w}^t) - D(\boldsymbol{\alpha}^t)] \leq \frac{1}{\theta^{imp}}(1 - \theta^{imp})^t\left(D(\boldsymbol{\alpha}^*) - D(\boldsymbol{\alpha}^0)\right),$$

where
$$\theta^{imp} = \frac{n\lambda\gamma}{\sum_{i=1}^{n}(v_j + \lambda\gamma n)}.$$

The next corollary suggests that a *better convergence rate than IProx-SDCA can be achieved by using properly chosen adaptive sampling probability*.

**Corollary 2.9.** *Consider AdaSDCA. If at each iteration $t \geq 0$, $\mathbf{p}_t$ is the optimal solution of (2.9), then (2.13) holds and $\tilde{\theta}_t \geq \theta^{imp}$ for all $t \geq 0$.*

However, solving (2.9) requires large computational effort, because of the dimension $n$ and the non-convex structure of the program. We show in the next section that when all the loss functions $\{\phi_j\}_j$ are quadratic, then we can get better convergence rate in theory than IProx-SDCA by using the optimal solution of (2.11).

### 2.4.2 Quadratic loss functions

The main difficulty of solving (2.9) comes from the inequality constraint, which originates from (2.7). In this section we mainly show that the constraint (2.7) can be relaxed if all $\{\phi_j\}_j$ are quadratic.

**Proposition 2.10.** *Suppose that all $\{\phi_j\}_j$ are quadratic. Let $t \geq 0$. If $\min_{j \in J_t} p_j^t > 0$, then*
$$\mathbb{E}_t\left[D(\boldsymbol{\alpha}^{t+1}) - D(\boldsymbol{\alpha}^t)\right] \geq \theta(\boldsymbol{\kappa}^t, \mathbf{p}^t)\left(P(\mathbf{w}^t) - D(\boldsymbol{\alpha}^t)\right).$$

*Proof.* This is a direct consequence of Lemma 2.5 and the fact that the right-hand side of (2.6) equals 0 when $\theta = \theta(\boldsymbol{\kappa}^t, \mathbf{p}^t)$. $\square$

**Theorem 2.11.** *Suppose that all $\{\phi_j\}_j$ are quadratic. Consider AdaSDCA. If at each iteration $t \geq 0$, $\min_{j \in J_t} p_j^t > 0$, then (2.13) holds for all $t \geq 0$.*

*Proof.* We only need to apply Proposition 2.10. The rest of the proof is the same as in Theorem 2.7. $\square$

**Corollary 2.12.** *Suppose that all $\{\phi_j\}_j$ are quadratic. Consider AdaSDCA. If at each iteration $t \geq 0$, $\mathbf{p}_t$ is the optimal solution of (2.11), which has a closed form (2.12), then (2.13) holds and $\tilde{\theta}_t \geq \theta_*$ for all $t \geq 0$.*

## 2.5 Efficient heuristic variant

Corollary 2.9 and 2.12 suggest how to choose adaptive sampling probability in AdaSDCA which yields a theoretical convergence rate at least as good as IProx-SDCA [88]. However, there are two main implementation issues of AdaSDCA:

1. The update of the dual residue $\boldsymbol{\kappa}^t$ at each iteration costs $\mathcal{O}(\text{nnz}(\mathbf{X}))$ where $\text{nnz}(\mathbf{X})$ is the number of nonzero elements of the matrix $\mathbf{X}$;

2. We do not know how to compute the optimal solution of (2.9).

In this section, we propose a heuristic variant of AdaSDCA, which avoids the above two issues while staying close to the 'good' adaptive sampling distribution.

### 2.5.1 Description of algorithm

---
**Algorithm 2** AdaSDCA+

  **Parameter** a number $m > 1$
  **Initialization** Choose $\boldsymbol{\alpha}^0 \in \mathbb{R}^n$, set $\bar{\boldsymbol{\alpha}}^0 = \frac{1}{\lambda n}\mathbf{X}\boldsymbol{\alpha}^0$
  **for** $t \geq 0$ **do**
    Primal update: $\mathbf{w}^t = \nabla g^*(\bar{\boldsymbol{\alpha}}^t)$
    Set: $\boldsymbol{\alpha}^{t+1} = \boldsymbol{\alpha}^t$
    **if** $\mod(t, n) == 0$ **then**
      **Option I:** Adaptive probability
        Compute: $\kappa_j^t = \alpha_j^t + \nabla\phi_j(\langle\mathbf{X}_{:j}, \mathbf{w}^t\rangle), \ \forall j \in [n]$
        Set: $p_j^t \sim |\kappa_i^t|\sqrt{v_j + n\lambda\gamma}, \ \forall j \in [n]$
      **Option II:** Optimal Importance probability
        Set: $p_j^t \sim (v_j + n\lambda\gamma), \ \forall j \in [n]$
    **end if**
    Generate random $j_t \in [n]$ according to $\mathbf{p}^t$
    Compute: $\quad \Delta\alpha_{j_t}^t = \arg\max\limits_{\Delta \in \mathbb{R}} \left\{-\phi_{j_t}^*(-(\alpha_{j_t}^t + \Delta)) - \langle\mathbf{X}_{:j_t}, \mathbf{w}^t\rangle\Delta - \frac{v_{j_t}}{2\lambda n}|\Delta|^2\right\}$
    Dual update: $\alpha_{j_t}^{t+1} = \alpha_{j_t}^t + \Delta\alpha_{j_t}^t$
    Average update: $\bar{\boldsymbol{\alpha}}^t = \bar{\boldsymbol{\alpha}}^t + \frac{\Delta\alpha_{j_t}}{\lambda n}\mathbf{X}_{:j_t}$
    Probability update: $p_{j_t}^{t+1} \sim p_{j_t}^t/m, \ \ p_k^{t+1} \sim p_k^t, \ \forall k \neq j_t$
  **end for**
  **Output:** $\mathbf{w}^t, \boldsymbol{\alpha}^t$

---

AdaSDCA+ has the same structure as AdaSDCA with a few important differences.

**Epochs** AdaSDCA+ is divided into epochs of length $n$. At the beginning of every epoch, sampling probabilities are computed according to one of two options. During each epoch the probabilities are cheaply updated at the end of every iteration to approximate the adaptive model. The intuition behind is as follows. After $j$ is sampled and the dual coordinate $\alpha_j$ is updated, the residue $\kappa_i$ naturally decreases. We then decrease also the probability that $j$ is chosen in the next iteration, by setting $\mathbf{p}^{t+1}$ to be proportional to $(p_1^t, \ldots p_{j-1}^t, p_j^t/m, p_{j+1}^t, \ldots, p_n^t)$. By doing this we avoid the computation of $\boldsymbol{\kappa}$ at each iteration (issue 1) which costs as much as the full gradient algorithm, while following closely the changes of the dual residue $\boldsymbol{\kappa}$. We reset the adaptive sampling probability after every epoch of length $n$.

**Parameter $m$** The setting of parameter $m$ in AdaSDCA+ directly affects the performance of the algorithm. If $m$ is too large, the probability of sampling the same coordinate twice during an epoch will be very small. This will result in a random permutation through all coordinates every epoch. On the other hand, for $m$ too small the coordinates having larger probabilities at the beginning of an epoch could be sampled more often than it should, even after their corresponding dual residues become sufficiently small. We don't have a definitive rule on the choice of $m$ and we leave this to future work. Experiments with different choices of $m$ can be found in Section 2.6.

Table 2.1: One epoch computational cost of different algorithms

| Algorithm | Cost of an epoch |
|---|---|
| SDCA & Quartz(uniform) | $\mathcal{O}(\text{nnz})$ |
| IProx-SDCA | $\mathcal{O}(\text{nnz} + n \log(n))$ |
| AdaSDCA | $\mathcal{O}(n \cdot \text{nnz})$ |
| AdaSDCA+ | $\mathcal{O}(\text{nnz} + n \log(n))$ |

**Option I & Option II** At the beginning of each epoch, one can choose between two options for resetting the sampling probability. Option I corresponds to the optimal solution of (2.11), given by the closed form (2.12). Option II is the optimal serial sampling probability (2.10), the same as the one used in IProx-SDCA [88]. However, AdaSDCA+ differs significantly with IProx-SDCA since we also update iteratively the sampling probability, which as we show through numerical experiments yields a faster convergence than IProx-SDCA.

### 2.5.2 Computational cost

**Sampling and probability update** During the algorithm we sample $j \in [n]$ from non-uniform probability distribution $\mathbf{p}^t$, which changes at each iteration. This process can be done efficiently using the Random Counters algorithm introduced in Section 6.2 of [49], which takes $\mathcal{O}(n \log(n))$ operations to create the probability tree and $\mathcal{O}(\log(n))$ operations to sample from the distribution or change one of the probabilities.

**Total computational cost** We can compute the computational cost of one epoch. At the beginning of an epoch, we need $\mathcal{O}(\text{nnz}(\mathbf{X}))$ operations to calculate the dual residue $\boldsymbol{\kappa}$. Then we create a probability tree using $\mathcal{O}(n \log(n))$ operations. At each iteration we need $\mathcal{O}(\log(n))$ operations to sample a coordinate, $\mathcal{O}(\text{nnz}(\mathbf{X})/n)$ operations to calculate the update to $\boldsymbol{\alpha}$ and a further $\mathcal{O}(\log(n))$ operations to update the probability tree. As a result an epoch needs $\mathcal{O}(\text{nnz} + n \log(n))$ operations. For comparison purpose we list in Table 2.1 the one epoch computational cost of comparable algorithms.

## 2.6 Numerical experiments

In this section we present results of numerical experiments.

### 2.6.1 Loss functions

We test AdaSDCA and AdaSDCA+, SDCA, and IProx-SDCA for two different types of loss functions $\{\phi_j\}_{j=1}^n$: quadratic loss and smoothed Hinge loss. Let $y \in \mathbb{R}^n$ be the vector of labels. The quadratic loss is given by

$$\phi_j(x) = \frac{1}{2\gamma}(x - y_j)^2$$

and the smoothed Hinge loss is:

$$\phi_j(x) = \begin{cases} 0 & y_j x \geq 1 \\ 1 - y_j x - \gamma/2 & y_j x \leq 1 - \gamma \\ \frac{(1 - y_j x)^2}{2\gamma} & \text{otherwise,} \end{cases}$$

Table 2.2: Dimensions and nonzeros of the datasets

| Dataset | $d$ | $n$ | nnz($\mathbf{X}$)/($nd$) |
|---|---|---|---|
| w8a | 300 | 49,749 | 3.9% |
| dorothea | 100,000 | 800 | 0.9% |
| mushrooms | 112 | 8,124 | 18.8% |
| cov1 | 54 | 581,012 | 22% |
| ijcnn1 | 22 | 49,990 | 41% |

In both cases we use $L_2$-regularizer, i.e.,

$$g(\mathbf{w}) \quad = \quad \frac{1}{2}\|\mathbf{w}\|^2.$$

Quadratic loss functions appear usually in regression problems, and smoothed Hinge loss can be found in linear support vector machine (SVM) problems [71].

### 2.6.2 Numerical results

We used 5 different datasets: w8a, dorothea, mushrooms, cov1 and ijcnn1 (see Table 2.2).

In all our experiments we used $\gamma = 1$ and $\lambda = 1/n$.

**AdaSDCA** The results of the theory developed in Section 2.4 can be observed through Figure 2.1a to Figure 2.1d. AdaSDCA needs the least amount of iterations to converge, confirming the theoretical result.

**AdaSDCA+ vs. others** We can observe through in Figure 2.1 and Figure 2.2, that both options of AdaSDCA+ outperforms SDCA and IProx-SDCA (labelled "Optimal" in the experiments), in terms of number of iterations, for quadratic loss functions and for smoothed Hinge loss functions. One can observe similar results in terms of time in Figure 2.3 and Figure 2.4.

**Option I vs. Option II** Despite the fact that Option I is not theoretically supported for smoothed hinge loss, it still converges faster than Option II on every dataset and for every loss function. The biggest difference can be observed on Figure 2.4d, where Option I converges to the machine precision in just 15 seconds. The default value of $m$ used for Option II is $m = 10$.

**Different choices of $m$** To show the impact of different choices of $m$ on the performance of AdaSDCA+, in Figures 2.5a to 2.6d we compare the results of the two options of AdaSDCA+ using different $m$ equal to 2, 10 and 50. It is hard to draw a clear conclusion here because clearly the optimal $m$ shall depend on the dataset and the problem type.

## 2.A  Proofs

For the sake of the analysis let us define

$$f(\boldsymbol{\alpha}) \;=\; -\lambda g^*\left(\frac{1}{\lambda n}\sum_{j=1}^{n}\mathbf{X}_{:j}\alpha_j\right). \tag{2.16}$$

We shall need the following inequality.

**Lemma 2.13.** *Function $f:\mathbb{R}^n \to \mathbb{R}$ defined in (2.16) satisfies the following inequality:*

$$f(\boldsymbol{\alpha}+\mathbf{h}) \;\leq\; f(\boldsymbol{\alpha}) + \langle \nabla f(\boldsymbol{\alpha}), \mathbf{h}\rangle + \frac{1}{2\lambda n^2}\mathbf{h}^\top \mathbf{X}^\top \mathbf{X}\mathbf{h}, \tag{2.17}$$

*holds for $\forall \boldsymbol{\alpha}, \mathbf{h} \in \mathbb{R}^n$.*

*Proof.* Since $g$ is 1-strongly convex, $g^*$ is 1-smooth. Pick $\boldsymbol{\alpha}, \mathbf{h} \in \mathbb{R}^n$. Since, $f(\boldsymbol{\alpha}) = \lambda g^*(\frac{1}{\lambda n}\mathbf{X}\boldsymbol{\alpha})$, we have

$$\begin{aligned}
f(\boldsymbol{\alpha}+\mathbf{h}) &= \lambda g^*\left(\frac{1}{\lambda n}\mathbf{X}\boldsymbol{\alpha} + \frac{1}{\lambda n}\mathbf{X}\mathbf{h}\right) \\
&\leq \lambda\left[g^*\left(\frac{1}{\lambda n}\mathbf{X}\boldsymbol{\alpha}\right) + \left\langle \nabla g^*\left(\frac{1}{\lambda n}\mathbf{X}\boldsymbol{\alpha}\right), \frac{1}{\lambda n}\mathbf{X}\mathbf{h}\right\rangle + \frac{1}{2}\left\|\frac{1}{\lambda n}\mathbf{X}\mathbf{h}\right\|^2\right] \\
&= f(\boldsymbol{\alpha}) + \langle \nabla f(\boldsymbol{\alpha}), \mathbf{h}\rangle + \frac{1}{2\lambda n^2}\mathbf{h}^\top \mathbf{X}^\top \mathbf{X}\mathbf{h}.
\end{aligned}$$

$\square$

*Proof of Lemma 2.3.* It can be easily checked that the following relations hold

$$\nabla_j f(\boldsymbol{\alpha}^t) \;=\; \frac{1}{n}\langle \mathbf{X}_{:j}, \mathbf{w}^t\rangle, \;\; \forall t \geq 0, \; j \in [n], \tag{2.18}$$

$$g(\mathbf{w}^t) + g^*(\bar{\boldsymbol{\alpha}}^t) \;=\; \langle \mathbf{w}^t, \bar{\boldsymbol{\alpha}}^t\rangle, \;\; \forall t \geq 0, \tag{2.19}$$

where $\{\mathbf{w}^t, \boldsymbol{\alpha}^t, \bar{\boldsymbol{\alpha}}^t\}_{t\geq 0}$ is the output sequence of Algorithm 1. Let $t \geq 0$ and $\theta \in [0, \min_j p_j^t]$. For each $j \in [n]$, since $\phi_j$ is $1/\gamma$-smooth, $\phi_j^*$ is $\gamma$-strongly convex and thus for arbitrary $s_j \in [0,1]$,

$$\begin{aligned}
\phi_j^*(-\alpha_j^t + s_j \kappa_j^t) &= \phi_j^*\big((1-s_j)(-\alpha_j^t) + s_j \nabla \phi_j(\langle \mathbf{X}_{:j}, \mathbf{w}^t\rangle)\big) \\
&\leq (1-s_j)\phi_j^*(-\alpha_j^t) + s_j \phi_j^*(\nabla \phi_j(\langle \mathbf{X}_{:j}, \mathbf{w}^t\rangle)) - \frac{\gamma s_j(1-s_j)|\kappa_j^t|^2}{2}. \quad (2.20)
\end{aligned}$$

We have:

$$\begin{aligned}
f(\boldsymbol{\alpha}^{t+1}) - f(\boldsymbol{\alpha}^t) &\overset{(2.17)}{\leq} \langle \nabla f(\boldsymbol{\alpha}^t), \boldsymbol{\alpha}^{t+1} - \boldsymbol{\alpha}^t\rangle + \frac{1}{2\lambda n^2}\langle \boldsymbol{\alpha}^{t+1} - \boldsymbol{\alpha}^t, \mathbf{X}^\top \mathbf{X}(\boldsymbol{\alpha}^{t+1} - \boldsymbol{\alpha}^t)\rangle \\
&= \nabla_j f(\boldsymbol{\alpha}^t) \Delta \alpha_{j_t}^t + \frac{v_j}{2\lambda n^2}|\Delta \alpha_{j_t}^t|^2 \quad (2.21) \\
&\overset{(2.18)}{=} \frac{1}{n}\langle \mathbf{X}_{:j_t}, \mathbf{w}^t\rangle \Delta \alpha_{j_t}^t + \frac{v_j}{2\lambda n^2}|\Delta \alpha_{j_t}^t|^2 \quad (2.22)
\end{aligned}$$

Thus,

$$
\begin{aligned}
D(\boldsymbol{\alpha}^{t+1}) - D(\boldsymbol{\alpha}^t) &\overset{(2.22)}{\geq} -\frac{1}{n}\langle \mathbf{X}_{:j_t}, \mathbf{w}^t\rangle \Delta\alpha_{j_t}^t - \frac{v_{j_t}}{2\lambda n^2}|\Delta\alpha_{j_t}^t|^2 + \frac{1}{n}\sum_{j=1}^n \phi_j^*(-\alpha_j^t) - \frac{1}{n}\sum_{j=1}^n \phi_j^*(-\alpha_j^{t+1}) \\
&= -\frac{1}{n}\langle \mathbf{X}_{:j_t}, \mathbf{w}^t\rangle \Delta\alpha_{j_t}^t - \frac{v_{j_t}}{2\lambda n^2}|\Delta\alpha_{j_t}^t|^2 + \frac{1}{n}\phi_{j_t}^*(-\alpha_{j_t}^t) - \frac{1}{n}\phi_{j_t}^*\bigl(-(\alpha_{j_t}^t + \Delta\alpha_{j_t}^t)\bigr) \\
&= \max_{\Delta\in\mathbb{R}}\left\{-\frac{1}{n}\langle \mathbf{X}_{:j_t}, \mathbf{w}^t\rangle \Delta - \frac{v_{j_t}}{2\lambda n^2}|\Delta|^2 + \frac{1}{n}\phi_{j_t}^*(-\alpha_{j_t}^t) - \frac{1}{n}\phi_{j_t}^*\bigl(-(\alpha_{j_t}^t + \Delta)\bigr)\right\},
\end{aligned}
$$

where the last equality follows from the definition of $\Delta\alpha_{j_t}^t$ in Algorithm 1. Then by letting $\Delta = -s_j \kappa_{j_t}^t$ for some arbitrary $s_j \in [0,1]$ we get:

$$
\begin{aligned}
D(\boldsymbol{\alpha}^{t+1}) - D(\boldsymbol{\alpha}^t) &\geq \frac{s_j\langle \mathbf{X}_{:j_t}, \mathbf{w}^t\rangle \kappa_{j_t}^t}{n} - \frac{s_j^2 v_{j_t}|\kappa_{j_t}^t|^2}{2\lambda n^2} + \frac{1}{n}\phi_{j_t}^*(-\alpha_{j_t}^t) - \frac{1}{n}\phi_{j_t}^*(-\alpha_{j_t}^t + s_j \kappa_{j_t}^t) \\
&\overset{(2.20)}{\geq} \frac{s_j}{n}\left(\phi_{j_t}^*(-\alpha_{j_t}^t) - \phi_{j_t}^*(\nabla\phi_{j_t}(\langle \mathbf{X}_{:j_t}, \mathbf{w}^t\rangle)) + \langle \mathbf{X}_{:j_t}, \mathbf{w}^t\rangle \kappa_{j_t}^t\right) \\
&\quad - \frac{s_j^2 v_{j_t}|\kappa_{j_t}^t|^2}{2\lambda n^2} + \frac{\gamma s_j(1-s_j)|\kappa_{j_t}^t|^2}{2n}.
\end{aligned}
$$

By taking expectation with respect to $j_t$ we get:

$$
\begin{aligned}
\mathbb{E}_t\left[D(\boldsymbol{\alpha}^{t+1}) - D(\boldsymbol{\alpha}^t)\right] &\geq \sum_{j=1}^n \frac{p_j^t s_j}{n}\left[\phi_j^*(-\alpha_j^t) - \phi_j^*(\nabla\phi_j(\langle \mathbf{X}_{:j}, \mathbf{w}^t\rangle)) + \langle \mathbf{X}_{:j}, \mathbf{w}^t\rangle \kappa_j^t\right] \\
&\quad - \sum_{j=1}^n \frac{p_j^t s_j^2 |\kappa_j^t|^2 (v_j + \lambda\gamma n)}{2\lambda n^2} + \sum_{j=1}^n \frac{p_j^t \gamma s_j |\kappa_j^t|^2}{2n}. \quad (2.23)
\end{aligned}
$$

Set

$$
s_j = \begin{cases} 0, & j \notin J_t \\ \theta/p_j^t, & j \in J_t \end{cases} \quad (2.24)
$$

Then $s_j \in [0,1]$ for each $j \in [n]$ and by plugging it into (2.23) we get:

$$
\begin{aligned}
\mathbb{E}_t\left[D(\boldsymbol{\alpha}^{t+1}) - D(\boldsymbol{\alpha}^t)\right] &\geq \frac{\theta}{n}\sum_{j\in J_t}\left[\phi_j^*(-\alpha_j^t) - \phi_j^*(\nabla\phi_j(\langle \mathbf{X}_{:j}, \mathbf{w}^t\rangle)) + \langle \mathbf{X}_{:j}, \mathbf{w}^t\rangle \kappa_j^t\right] \\
&\quad - \frac{\theta}{2\lambda n^2}\sum_{j\in J_t}\left(\frac{\theta(v_j + n\lambda\gamma)}{p_j^t} - n\lambda\gamma\right)|\kappa_j^t|^2
\end{aligned}
$$

Finally note that:

$$
\begin{aligned}
P(\mathbf{w}^t) - D(\boldsymbol{\alpha}^t) &= \frac{1}{n}\sum_{j=1}^{n}\left[\phi_j(\langle \mathbf{X}_{:j}, \mathbf{w}^t\rangle) + \phi_j^*(-\alpha_j^t)\right] + \lambda\left(g(\mathbf{w}^t) + g^*(\bar{\boldsymbol{\alpha}}^t)\right) \\
&\stackrel{(2.19)}{=} \frac{1}{n}\sum_{j=1}^{n}\left[\phi_j^*(-\alpha_i^t) + \phi_j(\langle \mathbf{X}_{:j}, \mathbf{w}^t\rangle)\right] + \frac{1}{n}\langle \mathbf{w}^t, \mathbf{X}\boldsymbol{\alpha}^t\rangle \\
&= \frac{1}{n}\sum_{j=1}^{n}\big[\phi_j^*(-\alpha_j^t) + \langle \mathbf{X}_{:j}, \mathbf{w}^t\rangle \nabla\phi_j(\langle \mathbf{X}_{:j}, \mathbf{w}^t\rangle) \\
&\qquad\qquad -\phi_j^*(\nabla\phi_j(\langle \mathbf{X}_{:j}, \mathbf{w}^t\rangle)) + \langle \mathbf{X}_{:j}, \mathbf{w}^t\rangle \alpha_j^t\big] \\
&= \frac{1}{n}\sum_{j=1}^{n}\left[\phi_j^*(-\alpha_j^t) - \phi_j^*(\nabla\phi_j(\langle \mathbf{X}_{:j}, \mathbf{w}^t\rangle)) + \langle \mathbf{X}_{:j}, \mathbf{w}^t\rangle \kappa_j^t\right] \\
&= \frac{1}{n}\sum_{j\in J_t}\left[\phi_j^*(-\alpha_j^t) - \phi_j^*(\nabla\phi_j(\langle \mathbf{X}_{:j}, \mathbf{w}^t\rangle)) + \langle \mathbf{X}_{:j}, \mathbf{w}^t\rangle \kappa_j^t\right]
\end{aligned}
$$

$\square$

*Proof of Lemma 2.4.* Note that (2.11) is a standard constrained maximization problem, where everything independent of $\mathbf{p}$ can be treated as a constant. We define the Lagrangian

$$
L(\mathbf{p}, \eta) = \theta(\boldsymbol{\kappa}, \mathbf{p}) - \eta\left(\sum_{j=1}^{n}p_j - 1\right)
$$

and get the following optimality conditions:

$$
\begin{aligned}
\frac{|\kappa_j^t|^2(v_j + n\lambda\gamma)}{p_j^2} &= \frac{|\kappa_k^t|^2(v_k + n\lambda\gamma)}{p_k^2}, \quad \forall j, k \in [n] \\
\sum_{j=1}^{n}p_j &= 1 \\
p_j &\geq 0, \quad \forall j \in [n],
\end{aligned}
$$

the solution of which is (2.12). $\square$

*Proof of Lemma 2.5.* Note that in the proof of Lemma 2.3, the condition $\theta \in [0, \min_{j\in J_t} p_j^t]$ is only needed to ensure that $s_j$ defined by (2.24) is in $[0, 1]$ so that (2.20) holds. If $\phi_j$ is quadratic function, then (2.20) holds for arbitrary $s_j \in \mathbb{R}$. Therefore in this case we only need $\theta$ to be positive and the same reasoning holds. $\square$

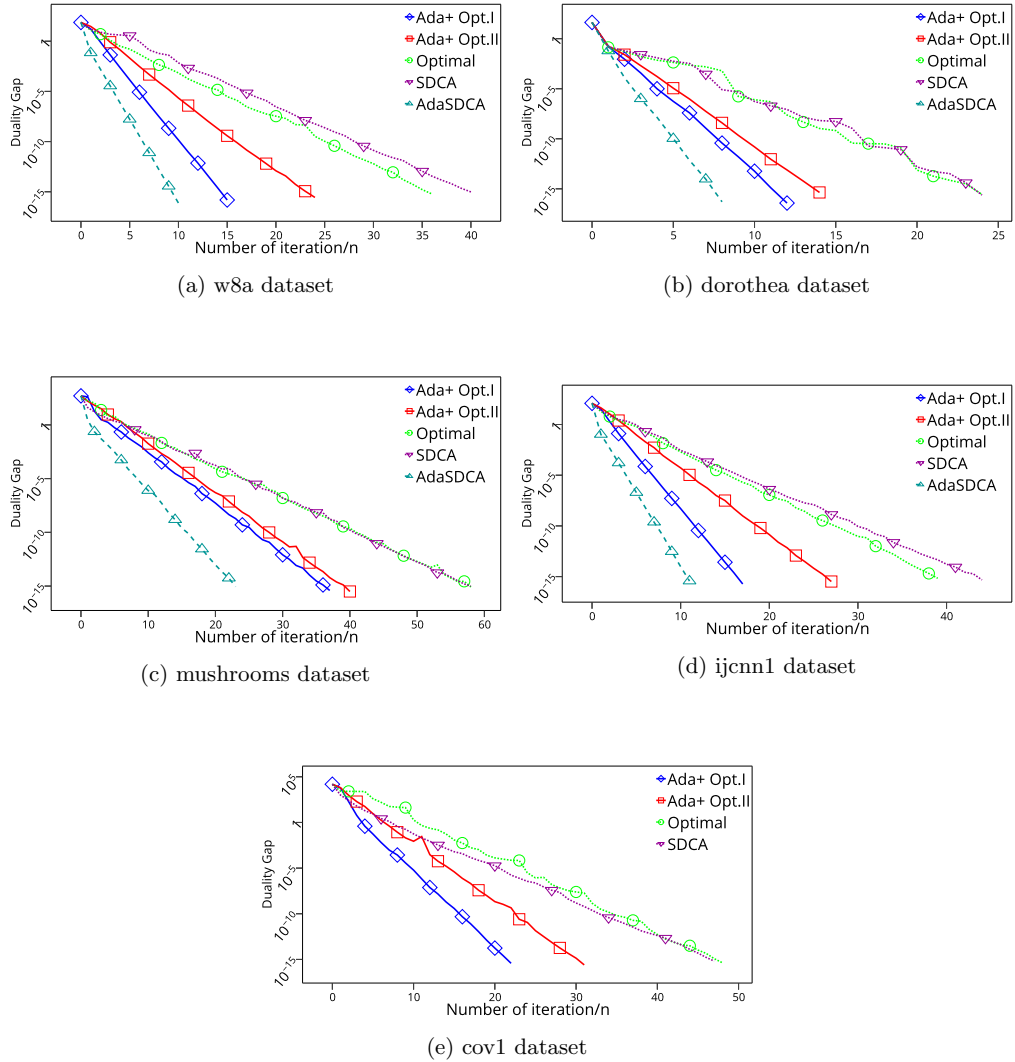

Figure 2.1: Quadratic loss with $L_2$ regularizer, comparing number of iterations with known algorithms.

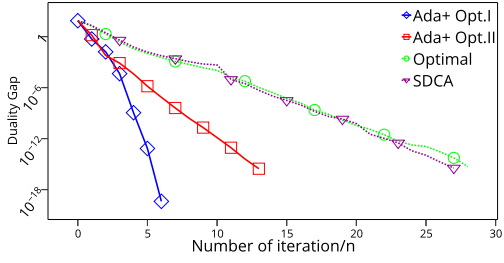
(a) w8a dataset

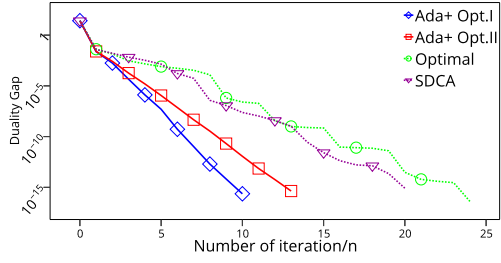
(b) dorothea dataset

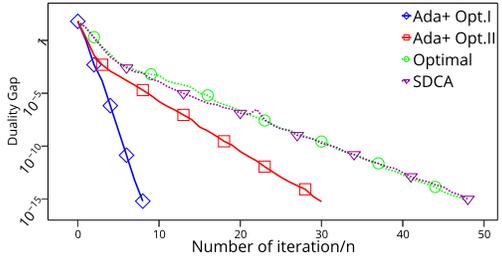
(c) mushrooms dataset

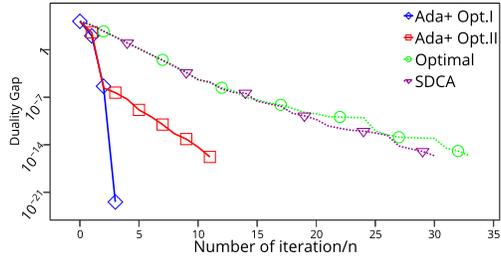
(d) cov1 dataset

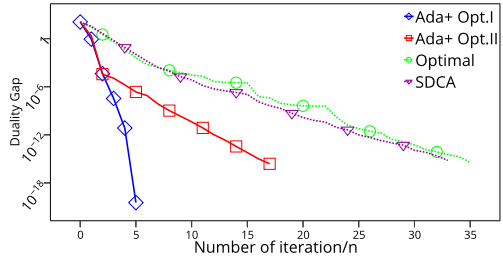
(e) ijcnn1 dataset

Figure 2.2: Smooth hinge loss with $L_2$ regularizer, comparing number of iterations with known algorithms.

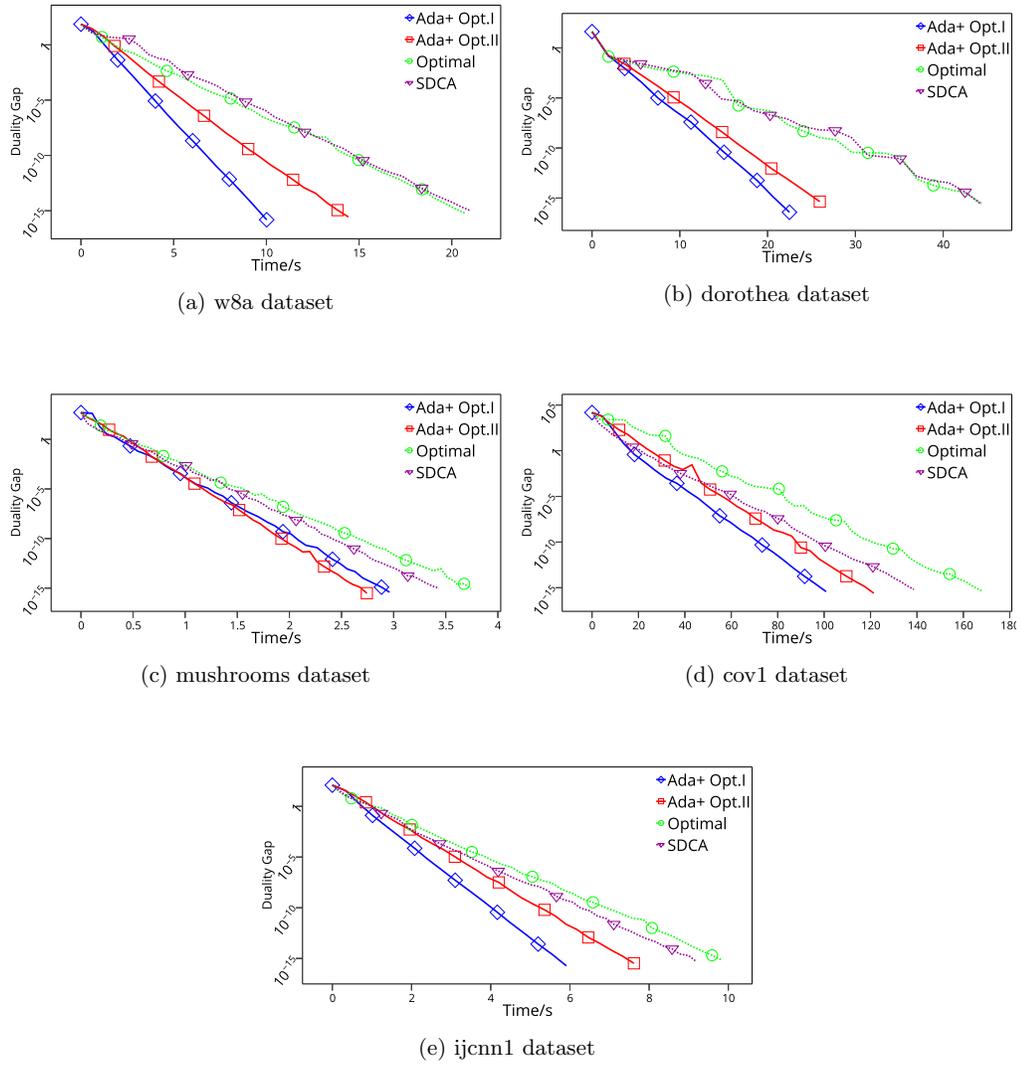

Figure 2.3: Quadratic loss with $L_2$ regularizer, comparing real time with known algorithms.

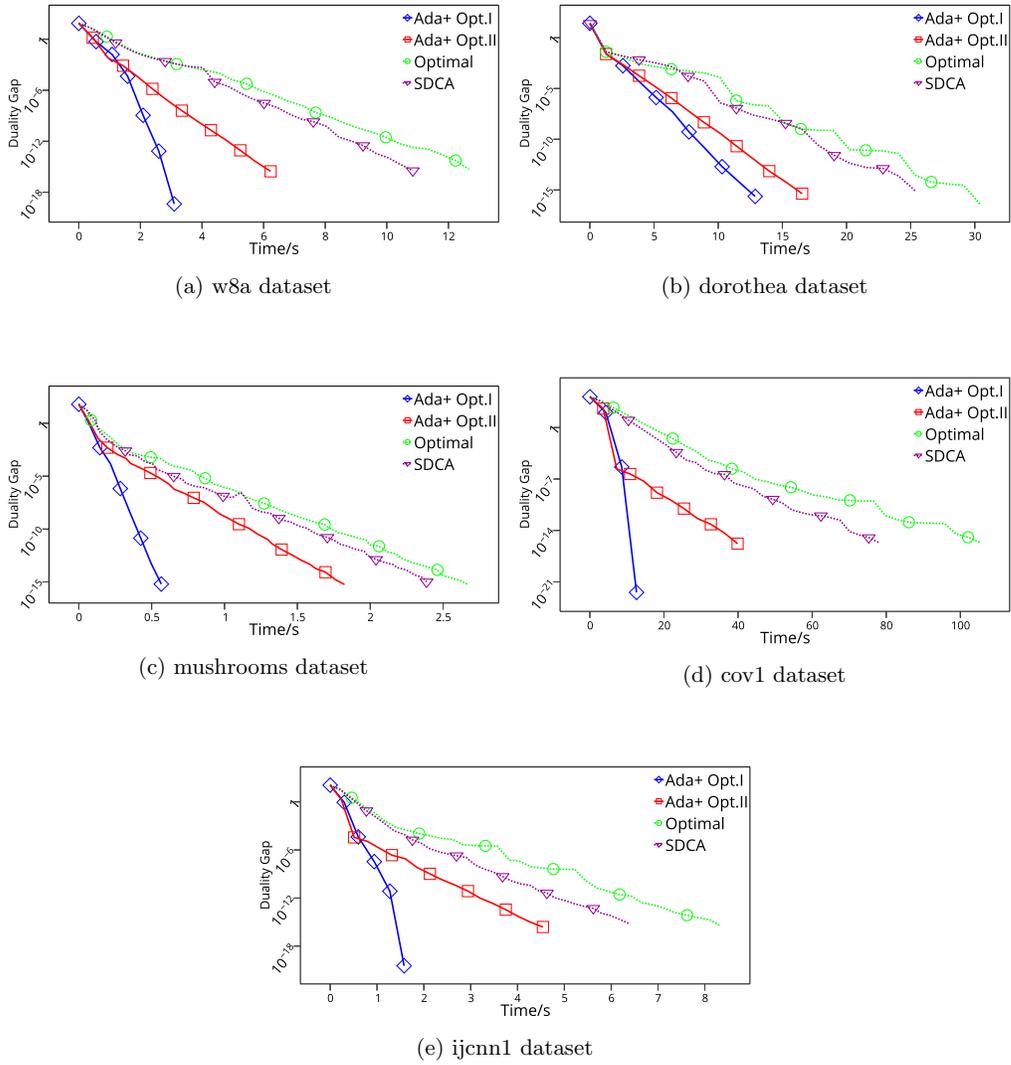

Figure 2.4: Smooth hinge loss with $L_2$ regularizer, comparing real time with known algorithms.

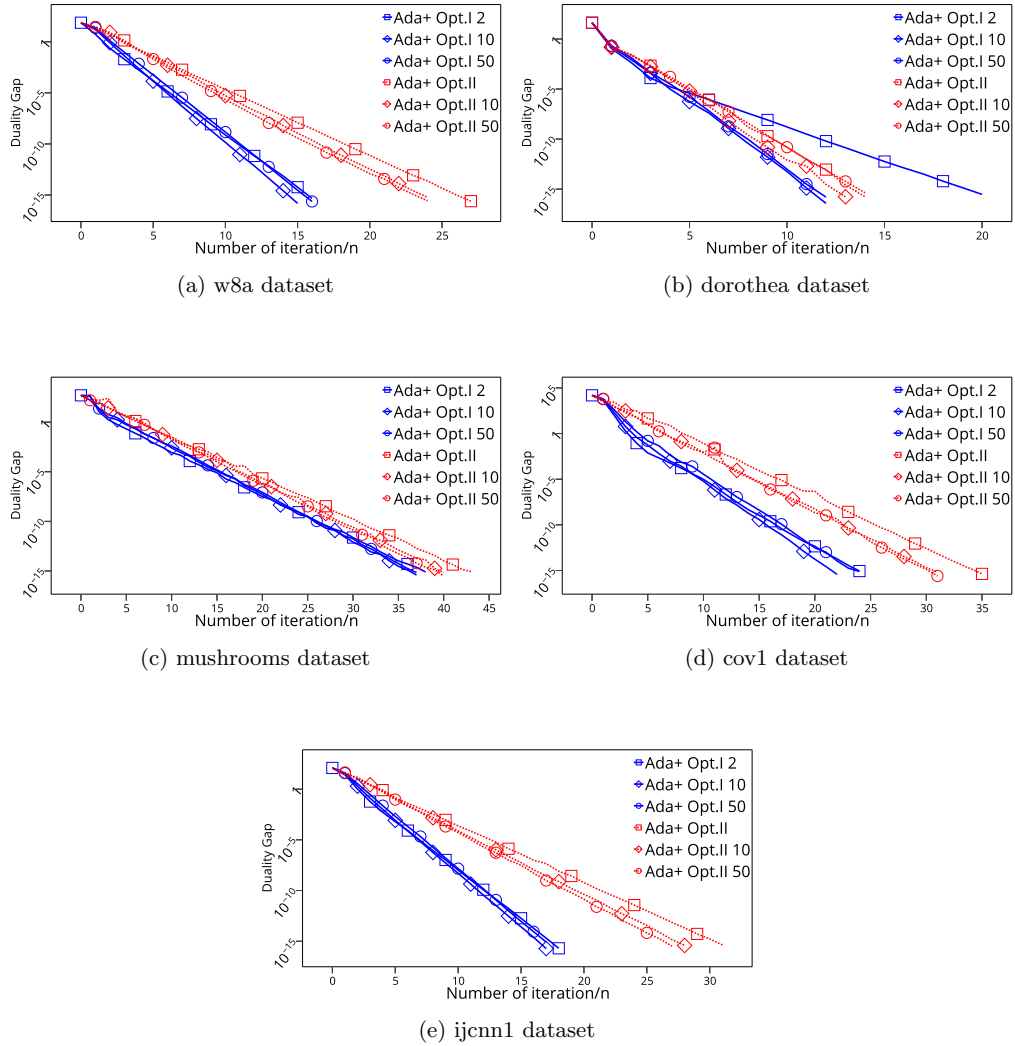

Figure 2.5: Quadratic loss with $L_2$ regularizer, comparison of different choices of the constant $m$

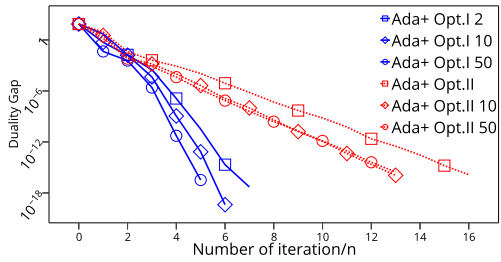
(a) w8a dataset

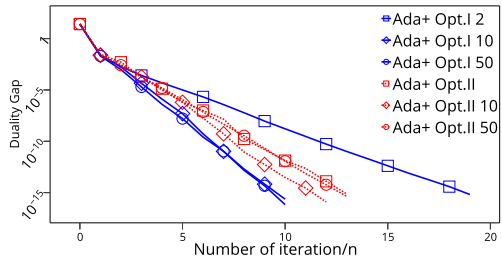
(b) dorothea dataset

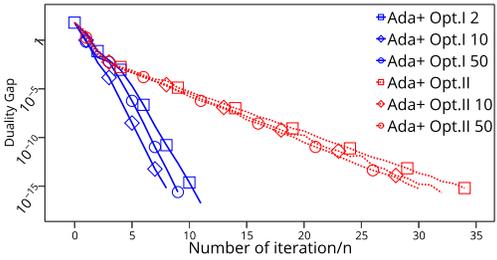
(c) mushrooms dataset

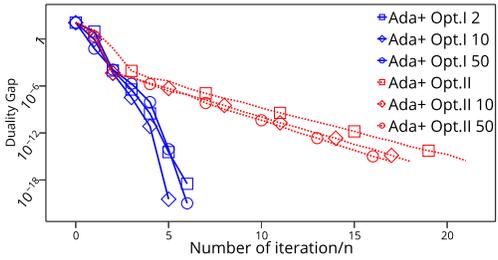
(d) ijcnn1 dataset

Figure 2.6: Smooth hinge loss with $L_2$ regularizer, comparison of different choices of the constant $m$.

# Chapter 3

# Primal Method for ERM with Flexible Mini-batching Schemes and Non-convex Losses

## 3.1 Introduction

Empirical risk minimization (ERM) is a very successful and immensely popular paradigm in machine learning, used to train a variety of prediction and classification models. Given examples $\mathbf{X}_1, \ldots, \mathbf{X}_n \in \mathbb{R}^{d \times m}$, loss functions $\phi_1, \ldots, \phi_n : \mathbb{R}^m \to \mathbb{R}$ and a regularization parameter $\lambda > 0$, the L2-regularized ERM problem is an optimization problem of the form

$$\min_{\mathbf{w} \in \mathbb{R}^d} \left[ P(\mathbf{w}) \quad := \quad \frac{1}{n} \sum_{j=1}^{n} \phi_j(\mathbf{X}_j^\top \mathbf{w}) + \frac{\lambda}{2} \|\mathbf{w}\|^2 \right] \tag{3.1}$$

Note that the examples are considered matrices instead of standard vectors, as this is a more general setup. Throughout the chapter we shall assume that for each $j$, the loss function $\phi_j$ is $l_i$-smooth with $l_j > 0$. That is, for all $\mathbf{x}, \mathbf{y} \in \mathbb{R}^m$ and all $j \in [n] := \{1, 2, \ldots, n\}$, we have

$$\|\nabla \phi_j(\mathbf{x}) - \nabla \phi_j(\mathbf{y})\| \quad \leq \quad l_j \|\mathbf{x} - \mathbf{y}\|. \tag{3.2}$$

In the last few years, a lot of research effort was put into designing new efficient algorithms for solving this problem (and some of its modifications). The frenzy of activity was motivated by the realization that SGD [61], not so long ago considered the state-of-the-art method for ERM, was far from being optimal, and that new ideas can lead to algorithms which are far superior to SGD in both theory and practice. The methods that belong to this category include SAG [64], SDCA [72], SVRG [27], S2GD [32], mS2GD [30], SAGA [12], S2CD [31], QUARTZ [55], ASDCA [71], prox-SDCA [70], IPROX-SDCA [88], A-PROX-SDCA [73], AdaSDCA [10], SDNA [56]. Methods analyzed for arbitrary mini-batching schemes include NSync [57], ALPHA [53] and QUARTZ [55].

In order to find an $\epsilon$-solution in expectation, state of the art (and non-accelerated) methods for solving (3.1) only need $\mathcal{O}((n+\kappa)\log(1/\epsilon))$ steps, where each step involves the computation of the gradient $\nabla \phi_j(\mathbf{X}_j^\top \mathbf{w})$ for some randomly selected example $j$. The quantity $\kappa$ is the condition number. Typically one has $\kappa = \frac{\max_j l_j \|\mathbf{X}_j\|^2}{\lambda}$ for methods picking $j$ uniformly at random, and $\kappa = \frac{\sum_j l_j \|\mathbf{X}_j\|^2}{n\lambda}$ for methods picking $j$ using a carefully designed data-dependent importance sampling. Computation of such a gradient typically involves work which is equivalent to reading the example $\mathbf{X}_j$, that is, $\mathcal{O}(\text{nnz}(\mathbf{X}_j)) \leq dm$ arithmetic operations.

**Contributions.** In this work we develop a new algorithm for the L2-regularized ERM problem (3.1). Our method extends a technique recently introduced by Shalev-Shwartz [67], which enables a dual-free analysis of SDCA, to *arbitrary mini-batching schemes*. That is, our method works at each iteration with a random subset of examples, chosen in an i.i.d. fashion from an arbitrary distribution. Such flexible schemes are useful for various reasons, including i) the development of distributed or robust variants of the method, ii) design of importance sampling for improving the complexity rate, iii) design of a sampling which is aimed at obtaining efficiencies elsewhere, such us utilizing NUMA (non-uniform memory access) architectures, and iv) streamlining and speeding up the processing of each mini-batch by means of assigning to each processor approximately even workload so as to reduce idle time (we do experiments with the latter setup). In comparison with [67], our method is able to better *utilize the information in the data examples* $\mathbf{X}_1, \ldots, \mathbf{X}_n$, leading to a better *data-dependent* bound. For convex loss functions, our complexity results match those of QUARTZ [55] in terms of the rate (the logarithmic factors differ). QUARTZ is a primal-dual method also allowing for arbitrary mini-batching schemes. However, while [55] only characterize the decay of expected risk, we also give bounds for the *sequence of iterates*. In particular, we show that for convex loss functions, our method enjoys the rate $\mathcal{O}(\max_j\{1/p_j + l_j v_j/(\lambda p_j n)\})$, where $p_j$ and $v_j$ are some parameters of the method associated with the sampling and data. For instance, in the special case picking a single example at a time uniformly at random, we have $p_j = 1/n$ and $v_j = \|\mathbf{X}_j\|_{\text{op}}^2$, whereby we obtain one of the $\mathcal{O}(n + \kappa)\log(1/\epsilon)$ rates mentioned above. The other rate can be recovered using importance sampling. The advantage of a dual-free analysis comes from the fact that it guarantees convergence even for *non-convex* loss functions, as long as the average loss is convex. This is a step toward understanding non-convex models. We illustrate through experiments the utility of being able to design arbitrary mini-batching schemes.

## 3.2 Algorithm

We shall now describe the method (Algorithm 3).

---

**Algorithm 3** dfSDCA: Dual-Free SDCA with Arbitrary Sampling

---

**Parameters:** Sampling $\hat{S}$, stepsize $\theta$
**Initialization** $\mathbf{A}^0 \in \mathbb{R}^{m \times n}$, set $\mathbf{w}^0 = \frac{1}{\lambda n}\sum_{j=1}^n \mathbf{X}_j \mathbf{A}_{:j}^0$, $p_j = \mathbf{Prob}(j \in S_t)$
**for** $t \geq 1$ **do**
  Sample a set $S_t$ according to $\hat{S}$
  **for** $j \in S_t$ in parallel **do**
    $\mathbf{A}_{:j}^t = \mathbf{A}_{:j}^{t-1} - \theta p_j^{-1}(\nabla \phi_j(\mathbf{X}_j^\top \mathbf{w}^{t-1}) + \mathbf{A}_{:j}^{t-1})$
  **end for**
  $\mathbf{w}^t = \mathbf{w}^{t-1} - \sum_{j \in S_t} \theta(n\lambda p_j)^{-1} \mathbf{X}_j(\nabla \phi_j(\mathbf{X}_j^\top \mathbf{w}^{t-1}) + \mathbf{A}_{:j}^{t-1})$
**end for**

---

The method encodes a family of algorithms, depending on the choice of the sampling $\hat{S}$, which encodes a particular *mini-batching scheme*. Formally, a sampling $\hat{S}$ is a set-valued random variable with values being the subsets of $[n]$, i.e., subsets of examples. In this chapter, we use the terms "mini-batching scheme" and "sampling" interchangeably. A sampling is defined by the collection of probabilities $\mathbf{Prob}(S)$ assigned to every subset $S \subseteq [n]$ of the examples.

The method maintains $n$ vectors $\mathbf{A}_{:j} \in \mathbb{R}^m$ as columns of the matrix $\mathbf{A} \in \mathbb{R}^{m \times n}$ and a vector $\mathbf{w} \in \mathbb{R}^d$. At the beginning of step $t$, we have $\mathbf{A}_{:j}^{t-1}$ for all $j$ and $\mathbf{w}^{t-1}$ computed and stored in memory. We then pick a random subset $S_t$ of the examples, according to the mini-batching scheme, and update variables $\mathbf{A}_{:j}$ for $j \in S_t$, based on the computation of the gradients $\nabla \phi_j(\mathbf{X}_j^\top \mathbf{w}^{t-1})$ for $j \in S_t$. This is followed by an update of the vector $\mathbf{w}$, which is performed so as to maintain the relation

$$\mathbf{w}^t = \frac{1}{\lambda n} \sum_{j=1}^n \mathbf{X}_j \mathbf{A}_{:j}^t. \tag{3.3}$$

This relation is maintained for the following reason. If $\mathbf{w}^*$ is the optimal solution to (3.1), then

$$0 = \nabla P(\mathbf{w}^*) = \frac{1}{n}\sum_{j=1}^n \mathbf{X}_j \nabla \phi_j(\mathbf{X}_j^\top \mathbf{w}^*) + \lambda \mathbf{w}^*, \tag{3.4}$$

and hence $\mathbf{w}^* = \frac{1}{\lambda n}\sum_{j=1}^n \mathbf{X}_j \mathbf{A}^*_{:j}$, where $\mathbf{A}^*_{:j} := -\nabla\phi_j(\mathbf{X}_j^\top \mathbf{w}^*)$. So, if we believe that the variables $\mathbf{A}_{:j}$ converge to $-\nabla\phi_j(\mathbf{X}_j^\top \mathbf{w}^*)$, it indeed does make sense to maintain (3.3). Why should we believe this? This is where the specific update of the "dual variables" $\mathbf{A}_{:j}$ comes from: $\mathbf{A}_{:j}$ is set a convex combination of its previous value and our best estimate so far of $-\nabla\phi_j(\mathbf{X}_j^\top \mathbf{w}^*)$, namely, $-\nabla\phi_j(\mathbf{X}_j^\top \mathbf{w}^{t-1})$. Indeed, the update can be written as $\mathbf{A}_{:j}^t = (1-\theta n p_j^{-1})\mathbf{A}_{:j}^{t-1} + \theta n p_j^{-1}(-\nabla\phi_j(\mathbf{X}_j^\top \mathbf{w}^{t-1}))$. Why does *this* make sense? Because we believe that $\mathbf{w}^{t-1}$ converges to $\mathbf{w}^*$. Admittedly, this reasoning is somewhat "circular". However, a better word to describe this reasoning would be: "iterative".

## 3.3 Main results

Let $p_j := \mathbb{P}(j \in S_t)$. We say that a sampling $\hat{S}$ is *proper* if $p_j > 0$ for all $j \in [n]$. In the rest of the chapter we assume that all samplings are proper.

We assume the knowledge of parameters $v_1, \ldots, v_n > 0$ for which

$$\mathbf{E}\left[\Big\|\sum_{j\in S_t} \mathbf{X}_j \mathbf{A}_{:j}\Big\|^2\right] \leq \sum_{j=1}^n p_j v_j \|\mathbf{X}_{:j}\|^2 \qquad (3.5)$$

for all values of $\mathbf{A} \in \mathbb{R}^{m\times n}$. Tight and easily computable formulas for such parameters can be found in [54]. For instance, whenever $\mathbf{Prob}(|S_t| \leq \tau) = 1$, inequality (3.5) holds with $v_j = \tau\|\mathbf{X}_j\|^2$.

Further, let $L_1, \ldots, L_n$ be constants for which the inequality

$$\|\nabla\phi_j(\mathbf{X}_j^\top \mathbf{w}) - \nabla\phi_j(\mathbf{X}_j^\top \mathbf{z})\| \leq L_j \|\mathbf{w} - \mathbf{z}\| \qquad (3.6)$$

holds for all $\mathbf{w}, \mathbf{z} \in \mathbb{R}^d$ and all $i$ and let $L := \max_j L_j$. Note that we can always bound $L_j \leq l_j\|\mathbf{X}_j\|$. However, $L_j$ can be better (smaller) than $l_j\|\mathbf{X}_j\|$, which is why we will use $L_j$ in our results whenever possible.

To simplify the exposure, we will write

$$W_t := \|\mathbf{w}^t - \mathbf{w}^*\|^2, \qquad a_j^t := \|\mathbf{A}_{:j}^t - \mathbf{A}_{:j}^*\|^2, \quad j = 1, 2, \ldots, n. \qquad (3.7)$$

### 3.3.1 Non-convex loss functions

Our result will be expressed in terms of the decay of the potential $D_t := \frac{\lambda}{2} W_t + \frac{\lambda}{2n}\sum_{j=1}^n \frac{1}{L_j^2} a_j^t$, where $W_t$ and $a_j^t$ are defined in (3.7).

**Theorem 3.1.** *Assume the average loss function $\frac{1}{n}\sum_{j=1}^n \phi_j$ is convex. If (3.6) holds and we let*

$$\theta \leq \min_j \frac{p_j n \lambda^2}{L_j^2 v_j + n\lambda^2}, \qquad (3.8)$$

*then the for $t \geq 0$ the potential $D_t$ decays exponentially to zero as*

$$\mathbf{E}[D_t] \leq e^{-\theta t} D_0. \qquad (3.9)$$

*Moreover, if we set $\theta$ equal to the upper bound in (3.8), then*

$$T \geq \max_j \left(\frac{1}{p_j} + \frac{L_j^2 v_j}{\lambda^2 p_j n}\right) \log\left(\frac{(L+\lambda)D_0}{\lambda\epsilon}\right) \quad \Rightarrow \quad \mathbf{E}[P(\mathbf{w}^T) - P(\mathbf{w}^*)] \leq \epsilon.$$

### 3.3.2 Convex loss functions

Our result will be expressed in terms of the decay of the potential $E_t := \frac{\lambda}{2} W_t + \frac{1}{2n}\sum_{j=1}^n \frac{1}{l_j} a_j^t$, where $W_t$ and $a_j^t$ are defined in (3.7).

**Theorem 3.2.** *Assume that all loss functions $\{\phi_j\}$ are convex and satisfy (3.2). If we run Algorithm 3 with parameter $\theta$ satisfying the inequality*

$$\theta \leq \min_j \frac{p_j n \lambda}{l_j v_j + n \lambda}, \tag{3.10}$$

*then the for $t \geq 0$ the potential $E_t$ decays exponentially to zero as*

$$\mathbf{E}\left[E_t\right] \leq e^{-\theta t} E_0. \tag{3.11}$$

*Moreover, if we set $\theta$ equal to the upper bound in (3.10), then*

$$T \geq \max_j \left(\frac{1}{p_j} + \frac{l_j v_j}{\lambda p_j n}\right) \log\left(\frac{(L+\lambda)E_0}{\lambda \epsilon}\right) \quad \Rightarrow \quad \mathbf{E}[P(\mathbf{w}^T) - P(w^*)] \leq \epsilon$$

The rate, $\theta$, precisely matches that of the QUARTZ algorithm [55]. Quartz is the only other method for ERM which has been analyzed for an arbitrary mini-batching scheme. Our algorithm is dual-free, and as we have seen above, allows for an analysis covering the case of non-convex loss functions.

## 3.4 Chunking

In this section we illustrate one usage of the ability of our method to work with an arbitrary mini-batching scheme. Further examples include the ability to design distributed variants of the method [58], or the use of importance/adaptive sampling to lower the number of iterations [59][88][55][10].

One marked disadvantage of standard mini-batching ("choose a subset of examples, uniformly at random") used in the context of parallel processing on multicore processors is the fact that in a synchronous implementation there is a loss of efficiency due to the fact that the computation time of $\nabla \phi(\mathbf{X}_j^\top \mathbf{w})$ may differ through $j$. This is caused by the data examples having varying degree of sparsity. We hence introduce a new sampling which mitigates this issue.

**Chunks:** Let nnz $\mathbf{X}_j$ be the numebr of nonzero entries in $\mathbf{X}_j$. Choose sets $G_1, \ldots, G_k \subset [n]$, such that $\cup_{i=1}^k G_i = [n]$ and $G_i \cap G_j = \emptyset \; \forall i, j$ and $\psi(i) := \sum_{j \in G_i} \text{nnz}(\mathbf{X}_j)$ is similar for every $i$, i.e. $\psi(1) \approx \cdots \approx \psi(k)$. Instead of sampling $\tau$ coordinates we propose a new sampling, which on each iteration $t$ samples $\tau$ sets $G_1^t, \ldots, G_\tau^t$ out of $G_1, \ldots, G_k$ and uses coordinates $j \in \cup_{i=1}^\tau G_i^t$ as the sampled set. We assign each core one of the sets $G_i^t$ for parallel computation. The advantage of this sampling lies in the fact, that the load of computing $\nabla \phi(\mathbf{X}_j^\top \mathbf{w})$ for all $j \in G_i$ is similar for all $i \in [k]$. Hence, using this sampling we minimize the waiting time of processors.

**How to choose $G_1, \ldots, G_k$:** We introduce the following algorithm:

---
**Algorithm 4** Naive Chunks
---
**Parameters:** vector of nnz $\mathbf{u}$
**Initialization** $n = \text{length}(\mathbf{u})$; Empty vector $\mathbf{g}$ and $\mathbf{s}$ of length $n$; $m = \max(\mathbf{u})$
$\mathbf{g}[1] = 1, \quad \mathbf{s}[1] = \mathbf{u}[1], \quad i = 1$
**for** $t = 2 : n$ **do**
  **if** $\mathbf{g}[i] + \mathbf{u}[t] \leq m$ **then**
    $\mathbf{g}[i] = \mathbf{g}[i] + 1, \mathbf{s}[i] = \mathbf{s}[i] + \mathbf{u}[t]$
  **else**
    $i = i + 1, \mathbf{g}[i] = 1, \mathbf{s}[i] = \mathbf{u}[t]$
  **end if**
**end for**

---

The algorithms returns the partition of $[n]$ into $G_1, \ldots, G_k$ in a sense, that the first $\mathbf{g}[1]$ coordinates belong to $G_1$, next $\mathbf{g}[2]$ coordinates belong to $G_2$ and so on. The main advantage of this approach is, that it makes a preprocessing step on the dataset which takes just one pass through the data. On Figure 3.1a through Figure 3.1f we show the impact of Algorithm 4 on the probability of the waiting time of a single core, which we measure by the difference $\max_{j \in S_t}\{\mathrm{nnz}(\mathbf{X}_j)\} - \frac{1}{\tau}\sum_{j \in S_t} \mathrm{nnz}(\mathbf{X}_j)$ and $\max_{j \in [\tau]}\{\mathrm{nnz}(G_j^t)\} - \frac{1}{\tau}\sum_{j=1}^{\tau} \mathrm{nnz}(G_j^t)$ for the initial and preprocessed dataset respectively. We can observe, that the waiting time is smaller using the preprocessing.

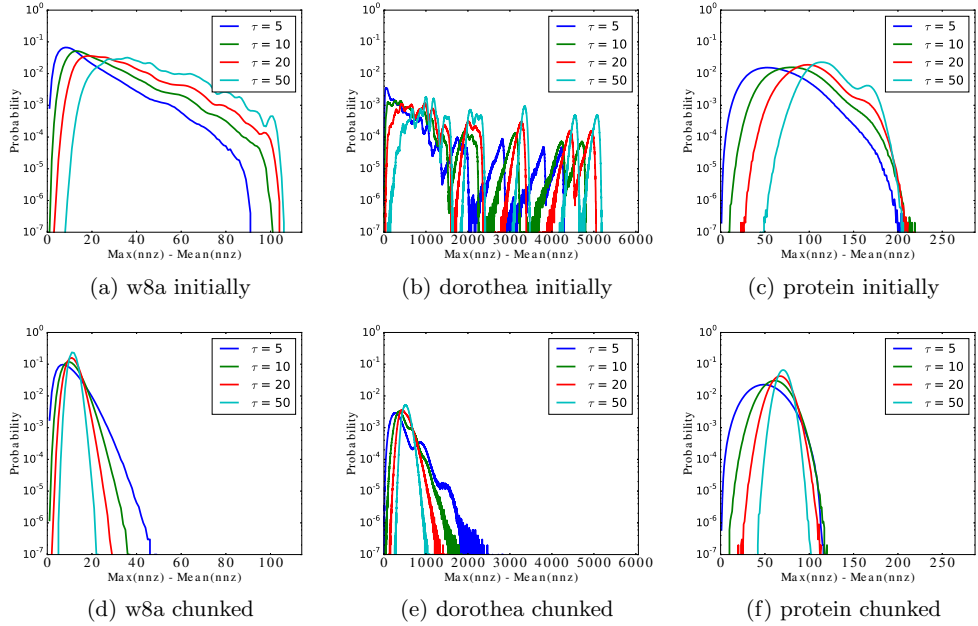

Figure 3.1: Distribution of the difference between the maximum number of nonzeros processed by a single core and the mean of all nonzeros processed by each core. This difference shows us, how much time is wasted per core waiting on the slowest core to finish its task, therefore smaller numbers are better. The first row corresponds to the initial distribution while the second row shows the distribution after using Algorithm 4.

## 3.5 Experiments

In all our experiments we used ridge logistic regression. We normalized the datasets so that $\max_j \|\mathbf{X}_j\| = 1$, and fixed $\lambda = 1/n$. The datasets used for experiments are summarized in Table 3.1.

| Dataset | #samples | #features | sparsity |
| --- | --- | --- | --- |
| w8a | 49,749 | 300 | 3.8% |
| dorothea | 800 | 100,000 | 0.9% |
| protein | 17,766 | 358 | 29% |
| rcv1 | 20,242 | 47,237 | 0.2% |
| cov1 | 581,012 | 54 | 22% |

Table 3.1: Datasets used in the experiments.

**Experiment 1.** In Figure 3.2a we compared the performance of Algorithm 3 with uniform serial sampling against state of the art algorithms such as SGD [61], SAG[64] and S2GD [32] in number of epochs. The real running time of the algorithms was 0.46s for S2GD, 0.79s for SAG, 0.47s for SDCA and 0.58s for SGD. In Figure 3.2b we show the convergence rate for different regularization parameters $\lambda$. In Figure 3.2c we show convergence rates for different serial samplings: uniform, importance [88] and also 4 different randomly generated serial samplings. These samplings were generated in a controlled manner, such that *random c* has $(\max_j\ p_j)/(\min_j\ p_j) < c$. All of these samplings have linear convergence as shown in the theory.

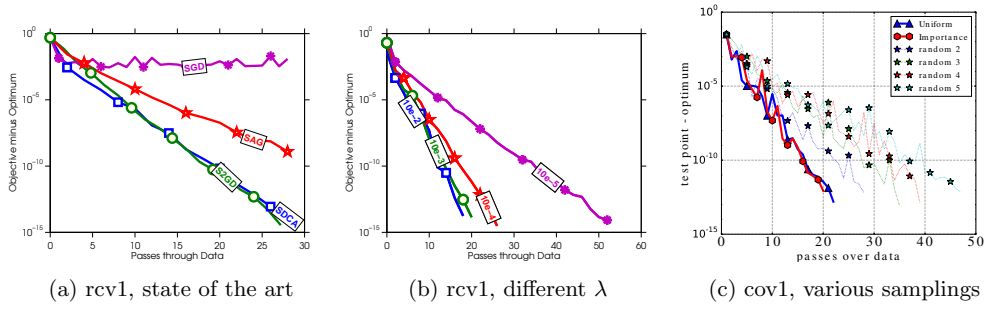

(a) rcv1, state of the art  (b) rcv1, different $\lambda$  (c) cov1, various samplings

Figure 3.2: LEFT: Comparison of SDCA with other state of the art methods. MIDDLE: SDCA for various values of $\lambda$. RIGHT: SDCA run with various samplings $\hat{S}$.

**Experiment 2: New sampling vs. old sampling.** In Figure 3.3a through Figure 3.3l we compare the performance of a standard parallel sampling against sampling of blocks $G_1, \ldots, G_k$ output by Algorithm 4. To simulate the synchronization of parallel processes, we measure the wall time by looking at the thread which has to process the most amount of nonzero entries, as the other threads will have to wait for it to finish. Therefore, in each iteration we approximate the actual wall time by $\max_{j \in S_t}\{\text{nnz}(\mathbf{X}_j)\}$ and $\max_{j \in [\tau]}\{\text{nnz}(G_j)\}$ for the standard and new sampling respectively. The x-axis on the plots stands for the amount of data passed by a single thread as a proxy for time.

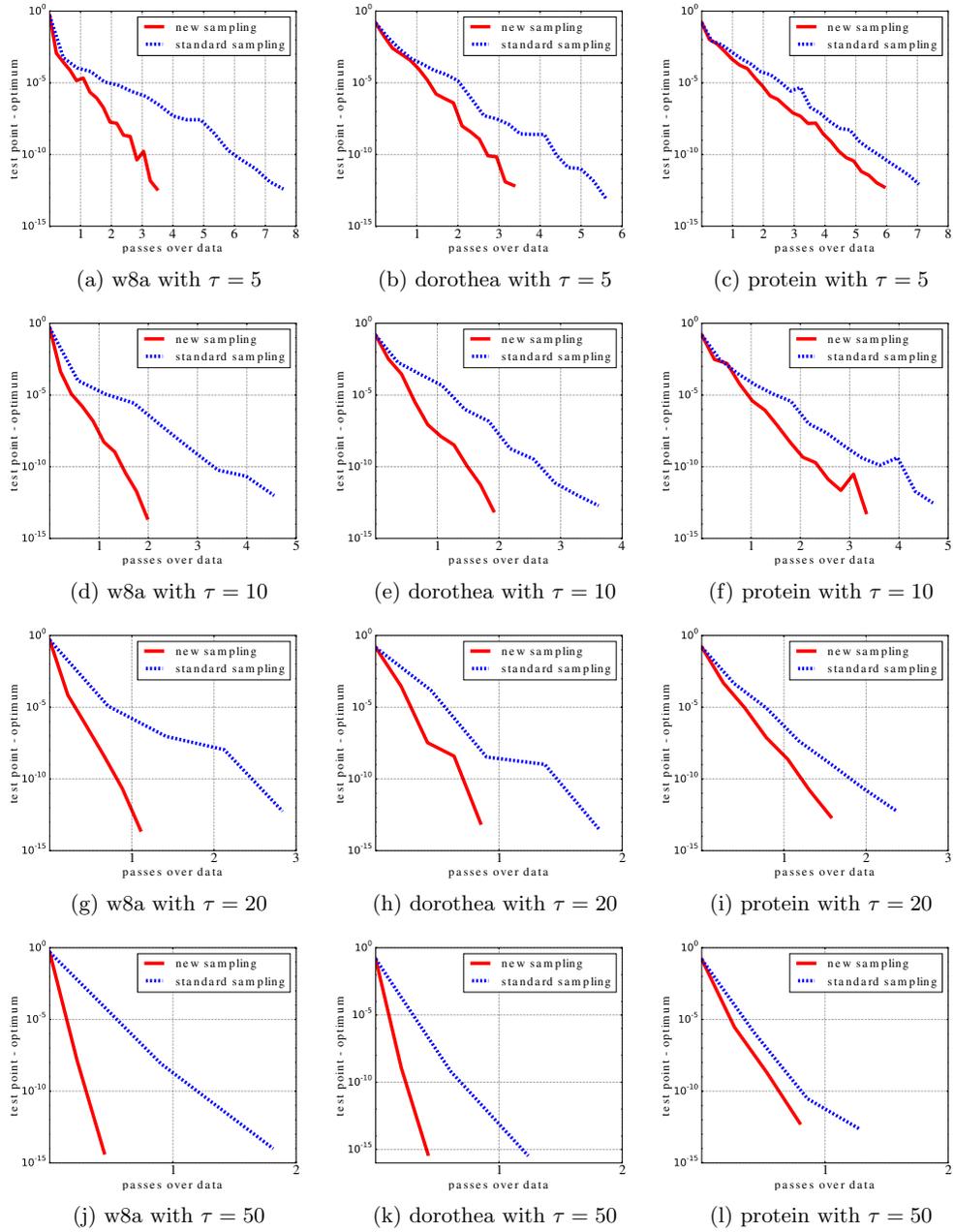

Figure 3.3: Logistic regression with $\lambda = 1/n$. Comparison between new and standard sampling with fine-tuned stepsizes for different values of $\tau$.

## 3.A  Proofs

As a first approximation, our proof is an extension of the proof of Shalev-Shwartz [67] to accommodate an arbitrary sampling [57][53][55][56]. For all $j$ and $t$ we let $\mathbf{U}_{:j} = -\nabla \phi_j(\mathbf{X}_j^\top \mathbf{w}^{t-1})$ and $\mathbf{Z}_{:j} = \mathbf{A}_{:j}^{t-1} - \mathbf{U}_{:j}$. We will use the following lemma.

**Lemma 3.3** (Evolution of $a_j^t$ and $W^t$). *For a fixed iteration $t$ the following holds:*

$$\mathbf{E}_{\hat{S}}\left[a_j^{t-1} - a_j^t\right] = \theta \left[\|\mathbf{A}_{:j}^{t-1} - \mathbf{A}_{:j}^*\|^2 - \|\mathbf{U}_{:j} - \mathbf{A}_{:j}^*\|^2 + (1 - \theta p_j^{-1})\|\mathbf{Z}_{:j}\|^2\right] \qquad (3.12)$$

$$\mathbf{E}_{\hat{S}}\left[W_{t-1} - W_t\right] \geq \frac{2\theta}{\lambda}\langle \mathbf{w}^{t-1} - \mathbf{w}^*, \nabla P(\mathbf{w}^{t-1})\rangle - \frac{\theta^2}{n^2 \lambda^2} \sum_{j=1}^n \frac{v_j}{p_j}\|\mathbf{Z}_{:j}\|^2. \qquad (3.13)$$

*Proof.* It follows that for $j \in S_t$ using the definition (3.7) we have

$$\begin{aligned}
a_j^{t-1} - a_j^t &= \|\mathbf{A}_{:j}^{t-1} - \mathbf{A}_{:j}^*\|^2 - \|\mathbf{A}_{:j}^t - \mathbf{A}_{:j}^*\|^2 \\
&= \|\mathbf{A}_{:j}^{t-1} - \mathbf{A}_{:j}^*\|^2 - \|(1 - \theta p_i^{-1})(\mathbf{A}_{:j}^{t-1} - \mathbf{A}_{:j}^*) + \theta p_i^{-1}(\mathbf{U}_{:j} - \mathbf{A}_{:j}^*)\|^2 \\
&= \|\mathbf{A}_{:j}^{t-1} - \mathbf{A}_{:j}^*\|^2 - (1 - \theta p_i^{-1})\|\mathbf{A}_{:j}^{t-1} - \mathbf{A}_{:j}^*\|^2 - \theta p_i^{-1}\|\mathbf{U}_{:j} - \mathbf{A}_{:j}^*\|^2 \\
&\quad + \theta p_i^{-1}(1 - \theta p_i^{-1})\|\mathbf{A}_{:j}^{t-1} - \mathbf{U}_{:j}\|^2 \\
&= \theta p_i^{-1}\left[\|\mathbf{A}_{:j}^{t-1} - \mathbf{A}_{:j}^*\|^2 - \|\mathbf{U}_{:j} - \mathbf{A}_{:j}^*\|^2 + (1 - \theta p_i^{-1})\|\mathbf{Z}_{:j}\|^2\right]
\end{aligned}$$

and for $j \notin S_t$ we have $a_j^{t-1} - a_j^t = 0$. Taking the expectation over $S_t$ we get the result.

For the second potential we get

$$\begin{aligned}
W_{t-1} - W_t &= \|\mathbf{w}^{t-1} - \mathbf{w}^*\|^2 - \|\mathbf{w}^t - \mathbf{w}^*\|^2 \\
&= \frac{2\theta}{n\lambda} \sum_{j \in S_t} p_j^{-1}\langle \mathbf{w}^{t-1} - \mathbf{w}^*, \mathbf{X}_j \mathbf{Z}_{:j}\rangle - \frac{\theta^2}{n^2 \lambda^2}\|\sum_{j \in S_t} p_j^{-1} \mathbf{X}_j \mathbf{Z}_{:j}\|^2
\end{aligned}$$

Taking the expectation over $S_t$, using inequality (3.5), and noting that

$$\frac{1}{n}\sum_{j=1}^n \mathbf{X}_j \mathbf{Z}_{:j} = \frac{1}{n}\sum_{j=1}^n \mathbf{X}_j \nabla \phi(\mathbf{X}_j^\top \mathbf{w}^{t-1}) + \lambda \mathbf{w}^{t-1} = \nabla P(\mathbf{w}^{t-1}), \qquad (3.14)$$

we get

$$\begin{aligned}
\mathbf{E}\left[W_{t-1} - W_t\right] &= \frac{2\theta}{n\lambda}\sum_{j=1}^n \langle \mathbf{w}^{t-1} - \mathbf{w}^*, \mathbf{X}_j \mathbf{Z}_{:j}\rangle - \frac{\theta^2}{n^2\lambda^2}\mathbf{E}\left[\|\sum_{j \in S_t} \mathbf{X}_j(p_j^{-1}\mathbf{Z}_{:j})\|^2\right] \\
&\stackrel{(3.5)}{\geq} \frac{2\theta}{n\lambda}\sum_{j=1}^n \langle \mathbf{w}^{t-1} - \mathbf{w}^*, \mathbf{X}_j \mathbf{Z}_{:j}\rangle - \frac{\theta^2}{n^2\lambda^2}\sum_{j=1}^n p_j v_j \|p_j^{-1}\mathbf{Z}_{:j}\|^2 \\
&\stackrel{(3.14)}{=} \frac{2\theta}{\lambda}\langle \mathbf{w}^{t-1} - \mathbf{w}^*, \nabla P(\mathbf{w}^{t-1})\rangle - \frac{\theta^2}{n^2\lambda^2}\sum_{j=1}^n \frac{v_j}{p_j}\|\mathbf{Z}_{:j}\|^2 \qquad \square
\end{aligned}$$

### 3.A.1 Proof of Theorem 3.1 (nonconvex case)

Combining (3.12) and (3.13), we obtain

$$\mathbf{E}[D_{t-1} - D_t] \geq \frac{\theta\lambda}{2n} \sum_{j=1}^{n} \frac{1}{L_j^2} \left[ \|\mathbf{A}_{:j}^{t-1} - \mathbf{A}_{:j}^*\|^2 - \|\mathbf{U}_{:j} - \mathbf{A}_{:j}^*\|^2 + (1 - \theta p_j^{-1})\|\mathbf{Z}_{:j}\|^2 \right]$$

$$+ \frac{\lambda}{2} \left[ \frac{2\theta}{\lambda} \langle \mathbf{w}^{t-1} - \mathbf{w}^*, \nabla P(\mathbf{w}^{t-1}) \rangle - \frac{\theta^2}{n^2\lambda^2} \sum_{j=1}^{n} \frac{v_j}{p_j} \|\mathbf{Z}_{:j}\|^2 \right]$$

$$= \frac{\theta}{2n} \sum_{j=1}^{n} \left[ \frac{\lambda}{L_j^2} \left( a_j^{t-1} - \|\mathbf{U}_{:j} - \mathbf{A}_{:j}^*\|^2 \right) + \left( \frac{\lambda(1 - \theta p_j^{-1})}{L_j^2} - \frac{\theta v_j}{n\lambda p_j} \right) \|\mathbf{Z}_{:j}\|^2 \right]$$

$$+ \theta \langle \mathbf{w}^{t-1} - \mathbf{w}^*, \nabla P(\mathbf{w}^{t-1}) \rangle$$

$$\overset{(3.8)}{\geq} \frac{\theta}{2n} \sum_{j=1}^{n} \frac{\lambda}{L_j^2} \left( a_j^{t-1} - \|\mathbf{U}_{:j} - \mathbf{A}_{:j}^*\|^2 \right) + \theta \langle \mathbf{w}^{t-1} - \mathbf{w}^*, \nabla P(\mathbf{w}^{t-1}) \rangle.$$

Using (3.6) we have $\|\mathbf{U}_{:j} - \mathbf{A}_{:j}^*\|^2 = \|\nabla\phi_j(\mathbf{X}_j^\top \mathbf{w}^{t-1}) - \nabla\phi_j(\mathbf{X}_j^\top \mathbf{w}^*)\|^2 \leq L_j^2 \|\mathbf{w}^{t-1} - \mathbf{w}^*\|^2$. Using the strong convexity of $P$ we have $\langle \mathbf{w}^{t-1} - \mathbf{w}^*, \nabla P(\mathbf{w}^{t-1}) \rangle \geq P(\mathbf{w}^{t-1}) - P(\mathbf{w}^*) + \frac{\lambda}{2}\|\mathbf{w}^{t-1} - \mathbf{w}^*\|^2$ and $P(\mathbf{w}^{t-1}) - P(\mathbf{w}^*) \geq \frac{\lambda}{2}\|\mathbf{w}^{t-1} - \mathbf{w}^*\|^2$, which together yields $\langle \mathbf{w}^{t-1} - \mathbf{w}^*, \nabla P(\mathbf{w}^{t-1}) \rangle \geq \lambda \|\mathbf{w}^{t-1} - \mathbf{w}^*\|^2$. Therefore,

$$\mathbf{E}[D_{t-1} - D_t] \geq \theta \left[ \frac{1}{n} \sum_{j=1}^{n} \frac{\lambda}{2L_j^2} a_j^{t-1} + \left( -\frac{\lambda}{2} + \lambda \right) W_{t-1} \right] = \theta D_{t-1}.$$

It follows that $\mathbf{E}[D_t] \leq (1-\theta)D_{t-1}$, and repeating this recursively we end up with $\mathbf{E}[D_{t-1}] \leq (1-\theta)^t D_0 \leq e^{-\theta t} D_{(0)}$. This concludes the proof of the first part of Theorem 3.1. The second part of the proof follows by observing that $P$ is $(L+\lambda)$-smooth, which gives $P(\mathbf{w}) - P(\mathbf{w}^*) \leq \frac{L+\lambda}{2}\|\mathbf{w} - \mathbf{w}^*\|^2$.

### 3.A.2 Convex case

For the next theorem we need an additional lemma:

**Lemma 3.4.** *Assume that $\phi_j$ are $L_j$-smooth and convex. Then, for every $\mathbf{w}$,*

$$\frac{1}{n}\sum_{j=1}^{n} \frac{1}{L_j} \|\nabla\phi_j(\mathbf{X}_j^\top \mathbf{w}) - \nabla\phi_j(\mathbf{X}_j^\top \mathbf{w}^*)\|^2 \leq 2\left(P(\mathbf{w}) - P(\mathbf{w}^*) - \frac{\lambda}{2}\|\mathbf{w} - \mathbf{w}^*\|^2\right) \quad (3.15)$$

*Proof.* Let

$$g_j(\mathbf{x}) = \phi_j(\mathbf{x}) - \phi_j(\mathbf{X}_j^\top \mathbf{w}^*) - \langle \nabla\phi_j(\mathbf{X}_j^\top \mathbf{w}^*), \mathbf{x} - \mathbf{X}_j^\top \mathbf{w}^* \rangle.$$

Clearly, $g_j$ is also $l_j$-smooth. By convexity of $\phi_j$ we have $g_j(\mathbf{x}) \geq 0$ for all $\mathbf{x}$. It follows that $g_j$ satisfies $\|\nabla g_j(\mathbf{x})\|^2 \leq 2l_j g_j(\mathbf{x})$. Using the definition of $g_j$, we obtain

$$\|\nabla\phi_j(\mathbf{X}_j^\top \mathbf{w}) - \nabla\phi_j(\mathbf{X}_j^\top \mathbf{w}^*)\|^2 = \|\nabla g_j(\mathbf{X}_j^\top \mathbf{w})\|^2$$
$$\leq 2l_j[\phi_j(\mathbf{X}_j^\top \mathbf{w}) - \phi_j(\mathbf{X}_j^\top \mathbf{w}^*) - \langle \nabla\phi_j(\mathbf{X}_j^\top \mathbf{w}^*), \mathbf{X}_j^\top \mathbf{w} - \mathbf{X}_j^\top \mathbf{w}^* \rangle]. \quad (3.16)$$

Summing these terms up weighted by $1/l_j$ and using (3.4) we get

$$\frac{1}{n}\sum_{j=1}^n \frac{1}{l_j}\|\nabla\phi_j(\mathbf{X}_j^\top \mathbf{w}) - \nabla\phi_j(\mathbf{X}_j^\top \mathbf{w}^*)\|^2$$

$$\stackrel{(3.16)}{\leq} \sum_{j=1}^n \frac{2}{n}[\phi_j(\mathbf{X}_j^\top \mathbf{w}) - \phi_j(\mathbf{X}_j^\top \mathbf{w}^*) - \langle \mathbf{X}_j \nabla\phi_j(\mathbf{X}_j^\top \mathbf{w}^*), \mathbf{w} - \mathbf{w}^*\rangle]$$

$$\stackrel{(3.4)}{=} 2\left[P(\mathbf{w}) - \frac{\lambda}{2}\|\mathbf{w}\|^2 - P(\mathbf{w}^*) + \frac{\lambda}{2}\|\mathbf{w}^*\|^2 + \lambda\langle\mathbf{w}^*, \mathbf{w}-\mathbf{w}^*\rangle\right]$$

$$= 2\left[P(\mathbf{w}) - P(\mathbf{w}^*) - \frac{\lambda}{2}\|\mathbf{w}-\mathbf{w}^*\|^2\right]. \quad \square$$

### 3.A.3   Proof of Theorem 3.2

Combining (3.12) and (3.13), we obtain

$$\mathbf{E}[E_{t-1} - E_t] \geq \frac{\theta}{n}\sum_{j=1}^n \frac{1}{2l_j}\left[\|\mathbf{A}_{:j}^{t-1} - \mathbf{A}_{:j}^*\|^2 - \|\mathbf{U}_{:j} - \mathbf{A}_{:j}^*\|^2 + (1-\theta p_j^{-1})\|\mathbf{Z}_{:j}\|^2\right]$$

$$+ \frac{\lambda}{2}\left[\frac{2\theta}{\lambda}\langle\mathbf{w}^{t-1}-\mathbf{w}^*, \nabla P(\mathbf{w}^{t-1})\rangle - \frac{\theta^2}{n^2\lambda^2}\sum_{j=1}^n \frac{v_j}{p_j}\|\mathbf{Z}_{:j}\|^2\right]$$

$$= \frac{\theta}{n}\sum_{j=1}^n \left[\frac{1}{2l_j}(a_j^{t-1} - \|\mathbf{U}_{:j} - \mathbf{A}_{:j}^*\|^2) + \left(\frac{(1-\theta p_j^{-1})}{2l_j} - \frac{\theta v_j}{2p_j\lambda n}\right)\right]$$

$$+ \theta\langle\mathbf{w}^{t-1}-\mathbf{w}^*, \nabla P(\mathbf{w}^{t-1})\rangle$$

$$\stackrel{(3.10)}{\geq} \frac{\theta}{n}\sum_{j=1}^n\left[\frac{1}{2l_j}(a_j^{t-1} - \|\mathbf{U}_{:j} - \mathbf{A}_{:j}^*\|^2)\right] + \theta\langle\mathbf{w}^{t-1}-\mathbf{w}^*, \nabla P(\mathbf{w}^{t-1})\rangle$$

Using the convexity of $P$ we have $P(\mathbf{w}^{t-1}) - P(\mathbf{w}^*) \leq \langle\mathbf{w}^{t-1} - \mathbf{w}^*, \nabla P(\mathbf{w}^{t-1})\rangle$ and using Lemma 3.4, we have

$$\mathbf{E}[E_{t-1} - E^t] \stackrel{(3.15)}{\geq} \frac{\theta}{n}\sum_{j=1}^n \frac{1}{2l_j}a_j^{t-1} - \theta\left(P(\mathbf{w}^{t-1}) - P(\mathbf{w}^*) - \frac{\lambda}{2}\|\mathbf{w}^{t-1}-\mathbf{w}^*\|^2\right)$$

$$+\theta\langle\mathbf{w}^{t-1}-\mathbf{w}^*, \nabla P(\mathbf{w}^{t-1})\rangle$$

$$\geq \theta\left[\frac{1}{n}\sum_{j=1}^n \frac{1}{2l_j}a_j^{t-1} + \frac{\lambda}{2}W_{t-1}\right] = \theta E_{t-1}.$$

This gives $\mathbf{E}[E_t] \leq (1-\theta)E_{t-1}$, which concludes the first part of the Theorem 3.2. The second part follows by observing, that $P$ is $(L+\lambda)$-smooth, which gives $P(\mathbf{w}) - P(\mathbf{w}^*) \leq \frac{L+\lambda}{2}\|\mathbf{w}-\mathbf{w}^*\|^2$.

# Chapter 4

# Importance Sampling for Minibatches

## 4.1 Introduction

Supervised learning is a widely adopted learning paradigm with important applications such as regression, classification and prediction. The most popular approach to training supervised learning models is via empirical risk minimization (ERM). In ERM, the practitioner collects data composed of example-label pairs, and seeks to identify the best predictor by minimizing the empirical risk, i.e., the average risk associated with the predictor over the training data.

With ever increasing demand for accuracy of the predictors, largely due to successful industrial applications, and with ever more sophisticated models that need to be trained, such as deep neural networks [23, 33], or multiclass classification [24], increasing volumes of data are used in the training phase. This leads to huge and hence extremely computationally intensive ERM problems.

Batch algorithms—methods that need to look at all the data before taking a single step to update the predictor—have long been known to be prohibitively impractical to use. Typical examples of batch methods are gradient descent and classical quasi-Newton methods. One of the most popular algorithms for overcoming the deluge-of-data issue is stochastic gradient descent (SGD), which can be traced back to a seminal work of Robbins and Monro [61]. In SGD, a single random example is selected in each iteration, and the predictor is updated using the information obtained by computing the gradient of the loss function associated with this example. This leads to a much more fine-grained iterative process, but at the same time introduces considerable stochastic noise, which eventually—typically after one or a few passes over the data—-effectively halts the progress of the method, rendering it unable to push the training error (empirical risk) to the realm of small values.

### 4.1.1 Strategies for dealing with stochastic noise

Several approaches have been proposed to deal with the issue of stochastic noise in the finite-data regime. The most important of these are i) decreasing stepsizes, ii) minibatching, iii) importance sampling and iv) variance reduction via "shift", listed here from historically first to the most modern.

The first strategy, *decreasing stepsizes*, takes care of the noise issue by a gradual and direct scale-down process, which ensures that SGD converges to the ERM optimum [85]. However, an unwelcome side effect of this is a considerable slowdown of the iterative process [3]. For instance, the convergence rate is sublinear even if the function to be minimized is strongly convex.

The second strategy, *minibatching*, deals with the noise by utilizing a random set of examples in the estimate of the gradient, which effectively decreases the variance of the estimate [74]. However, this has the unwelcome side-effect of requiring more computation. On the other hand, if a parallel processing machine is available, the computation can be done concurrently, which ultimately leads to speedup. This strategy does not result in an improvement of the convergence rate (unless progressively larger minibatch sizes are used, at the cost of further computational burden [18]), but can lead to massive improvement of the leading constant, which ultimately means acceleration (almost linear speedup for sparse data) [75].

The third strategy, *importance sampling*, operates by a careful data-driven design of the probabilities of selecting examples in the iterative process, leading to a reduction of the variance of the stochastic gradient thus selected. Typically, the overhead associated with computing the sampling probabilities and with sampling from the resulting distribution is negligible, and hence the net effect is speedup. In terms of theory, for standard SGD this improves a non-dominant term in the complexity. On the other hand, when SGD is combined with variance reduction, then this strategy leads to the improvement of the leading constant in the complexity estimate, typically via replacing the maximum of certain data-dependent quantities by their average [57, 31, 88, 55, 46, 7, 10].

Finally, and most recently, there has been a considerable amount of research activity due to the ground-breaking realization that one can gain the benefits of SGD (cheap iterations) without having to pay through the side effects mentioned above (e.g., halt in convergence due to decreasing stepsizes or increase of workload due to the use of minibatches). The result, in theory, is that for strongly convex losses (for example), one does not have to suffer sublinear convergence any more, but instead a fast linear rate "kicks in". In practice, these methods dramatically surpass all previous existing approaches.

The main algorithmic idea is to change the search direction itself, via a properly designed and cheaply maintainable *"variance-reducing shift"* (control variate). Methods in this category are of two types: those operating in the primal space (i.e., directly on ERM) and those operating in a dual space (i.e., with the dual of the ERM problem). Methods of the primal variety include SAG [64], SVRG [27], S2GD [32], proxSVRG [84], SAGA [12], mS2GD [30] and MISO [42]. Methods of the dual variety work by updating randomly selected dual variables, which correspond to examples. These methods include SCD [69], RCDM [49, 59], SDCA [72], Hydra [58, 16], mSDCA [75], APCG [37], AsySPDC [38], RCD [44], APPROX [15], SPDC [86], ProxSDCA [70], ASDCA [71], IProx-SDCA [88], and QUARTZ [55].

### 4.1.2 Combining strategies

We wish to stress that the key strategies, mini-batching, importance sampling and variance-reducing shift, should be seen as orthogonal tricks, and as such they can be combined, achieving an amplification effect. For instance, the first primal variance-reduced method allowing for mini-batching was [30]; while dual-based methods in this category include [71, 55, 7]. Variance-reduced methods with importance sampling include [49, 59, 57, 53] for general convex minimization problems, and [88, 55, 46, 7] for ERM.

## 4.2 Contributions

Despite considerable effort of the machine learning and optimization research communities, no importance sampling for minibatches was previously proposed, nor analyzed. The reason for this lies in the underlying theoretical and computational difficulties associated with the design and successful implementation of such a sampling. One needs to come up with a way to focus on a reasonable set of subsets (minibatches) of the examples to be used in each iteration (issue: there are many subsets; which ones to choose?), assign meaningful data-dependent non-uniform probabilities to them (issue: how?), and then be able to sample these subsets according to the chosen distribution (issue: this could be computationally expensive).

The tools that would enable one to consider these questions did not exist until recently. However, due to a recent line of work on analyzing variance-reduced methods utilizing what is known as *arbitrary sampling* [57, 55, 53, 54, 7], we are able to ask these questions and provide answers. In this work we design a novel family of samplings—*bucket samplings*—and a particular member of this family—*importance sampling for minibatches*. We illustrate the power of this sampling in combination with the reduced-variance dfSDCA method for ERM. This method is a primal variant of SDCA, first analyzed by Shalev-Shwartz [67], and extended by Csiba and Richtárik [7] to the arbitrary sampling setting. However, our sampling can be combined with any stochastic method for ERM, such as SGD or S2GD, and extends beyond the realm of ERM, to convex optimization problems in general. However, for simplicity, we do not discuss these extensions in this work.

We analyze the performance of the new sampling theoretically, and by inspecting the results we are able to comment on when can one expect to be able to benefit from it. We illustrate on synthetic datasets with varying distributions of example sizes that our approach can lead to *dramatic speedups* when compared against standard (uniform) minibatching, of *one or more degrees of magnitude*. We then test our method on real datasets and confirm that the use of importance minibatching leads to up to an order of magnitude speedup. Based on our experiments and theory, we predict that for real data with particular shapes and distributions of example sizes, importance sampling for minibatches will operate in a favourable regime, and can lead to speedup higher than one order of magnitude.

### 4.2.1 Related work

The idea of using non-uniform sampling in the parallel regime is by no means new. In the following we highlight several recent approaches in a chronological order and we describe their main differences to our method.

The first attempt for a potential speed-up using a non-uniform parallel sampling was proposed in [57]. However, to compute the optimal probability vector one has to solve a linear programming problem, which can easily be more complex than the original problem. The authors do not propose a practical version, which would overcome this issue.

The approach described in [87] uses the idea of a stratified sampling, which is a well-known strategy in statistics. The authors use clustering to group the examples into several partitions and sample an example from each of the partitions uniformly. This approach is similar to ours, with two main differences: i) we do not need clustering for our approach (it can be computationally very expensive) ii) we allow non-uniform sampling inside each of the partitions, which leads to the main speed-up in our work.

Instead of directly improving the convergence rate of the methods, the authors in [7] propose a strategy to improve the synchronized parallel implementation of a method by a load-balancing scheme. The method divides the examples into groups, which have similar sum of the amount of nonzero entries. When each core processes a single group, it should take the same time to finish as all the other groups, which leads to shorter waiting time in synchronization. Although this is a non-uniform parallel sampling, this approach takes a completely different direction than our method.

Lastly, in [22] the authors actually propose a scheme for importance sampling with minibatches. In the chapter they assume that they can sample a minibatch with a fixed size (without repetition), such that the probabilities of sampling individual examples will be proportional to some given values. However, this is easier said than done – until our work there was no sampling scheme which would allow for such minibatches. Therefore, the authors theoretically described an idea, which can be used in practice using our scheme.

## 4.3 The problem

Let $\mathbf{X} \in \mathbb{R}^{d \times n}$ be a data matrix in which features are represented in rows and examples in columns, and let $\mathbf{y} \in \mathbb{R}^n$ be a vector of labels corresponding to the examples. Our goal is to find a linear predictor $\mathbf{w} \in \mathbb{R}^d$ such that $\langle \mathbf{X}_{:j}, \mathbf{w} \rangle \sim y_j$, where the pair $(\mathbf{X}_{:j}, y_j) \in \mathbb{R}^d \times \mathbb{R}$ is sampled from the underlying distribution over data-label pairs. In the L2-regularized Empirical Risk Minimization problem, we find $\mathbf{w}$ by solving the optimization problem

$$\min_{\mathbf{w} \in \mathbb{R}^d} \left[ P(\mathbf{w}) \quad := \quad \frac{1}{n} \sum_{j=1}^n \phi_j(\langle \mathbf{X}_{:j}, \mathbf{w} \rangle) + \frac{\lambda}{2} \|\mathbf{w}\|_2^2 \right], \tag{4.1}$$

where $\phi_j : \mathbb{R} \to \mathbb{R}$ is a loss function associated with example-label pair $(\mathbf{X}_{:j}, y_j)$, and $\lambda > 0$. For instance, the square loss function is given by $\phi_j(t) = 0.5(t - y_j)^2$. Our results are not limited to L2-regularized problems though: an arbitrary strongly convex regularizer can be used instead [55]. We shall assume throughout that the loss functions are convex and $1/\gamma$-smooth, where $\gamma > 0$. The latter means that for all $s, t \in \mathbb{R}$ and all $j \in [n] := \{1, 2, \ldots, n\}$, we have

$$|\phi'_j(s) - \phi'_j(t)| \leq \frac{1}{\gamma}|s - t|.$$

This setup includes ridge and logistic regression, smoothed hinge loss, and many other problems as special cases [72]. Again, our sampling can be adapted to settings with non-smooth losses, such as the hinge loss.

## 4.4 The algorithm

In this chapter we illustrate the power of our new sampling in tandem with Algorithm 5 (dfS-DCA) for solving (4.1).

---

**Algorithm 5** dfSDCA [7]

---

**Parameters:** Sampling $\hat{S}$, stepsize $\theta > 0$
**Initialization:** Choose $\boldsymbol{\alpha}^0 \in \mathbb{R}^n$,
set $\mathbf{w}^0 = \frac{1}{\lambda n}\sum_{j=1}^n \mathbf{X}_{:j}\alpha_j^0$, $p_j = \mathbf{Prob}(j \in \hat{S})$
**for** $t \geq 1$ **do**
  Sample a fresh random set $S_t$ according to $\hat{S}$
  **for** $j \in S_t$ **do**
    $\Delta_j = \phi'_j(\langle \mathbf{X}_{:j}, \mathbf{w}^{t-1}\rangle) + \alpha_j^{t-1}$
    $\alpha_j^t = \alpha_j^{t-1} - \theta p_j^{-1}\Delta_j$
  **end for**
  $\mathbf{w}^t = \mathbf{w}^{t-1} - \sum_{j \in S_t} \theta(n\lambda p_j)^{-1}\Delta_j \mathbf{X}_{:j}$
**end for**

---

The method has two parameters. A "sampling" $\hat{S}$, which is a random set-valued mapping [60] with values being subsets of $[n]$, the set of examples. No assumptions are made on the distribution of $\hat{S}$ apart from requiring that $p_j$ is positive for each $j$, which simply means that each example has to have a chance of being picked. The second parameter is a stepsize $\theta$, which should be as large as possible, but not larger than a certain theoretically allowable maximum depending on $P$ and $\hat{S}$, beyond which the method could diverge.

Algorithm 5 maintains $n$ "dual" variables, $\alpha_1^t, \ldots, \alpha_n^t \in \mathbb{R}$, which act as variance-reduction shifts. This is most easily seen in the case when we assume that $S_t = \{j\}$ (no minibatching). Indeed, in that case we have

$$\mathbf{w}^t = \mathbf{w}^{t-1} - \frac{\theta}{n\lambda p_j}(g_j^{t-1} + \mathbf{X}_{:j}\alpha_j^{t-1}),$$

where $g_j^{t-1} := \mathbf{X}_{:j}\Delta_j$ is the stochastic gradient. If $\theta$ is set to a proper value, as we shall see next, then it turns out that for all $j \in [n]$, $\alpha_j$ is converging $\alpha_j^* := -\phi'_j(\langle \mathbf{X}_{:j}, \mathbf{w}^*\rangle)$, where $\mathbf{w}^*$ is the solution to (4.1), which means that the shifted stochastic gradient converges to zero. This means that its variance is progressively vanishing, and hence no additional strategies, such as decreasing stepsizes or minibatching are necessary to reduce the variance and stabilize the process. In general, dfSDCA in each step picks a random subset of the examples, denoted as $S_t$, updates variables $\alpha_j^t$ for $j \in S_t$, and then uses these to update the predictor $\mathbf{w}$.

### 4.4.1 Complexity of dfSDCA

In order to state the theoretical properties of the method, we define

$$E_t := \frac{\lambda}{2}\|\mathbf{w}^t - \mathbf{w}^*\|_2^2 + \frac{\gamma}{2n}\|\boldsymbol{\alpha}^t - \boldsymbol{\alpha}^*\|_2^2.$$

Most crucially to this chapter, we assume the knowledge of parameters $v_1, \ldots, v_n > 0$ for which the following ESO[1] inequality holds

$$\mathbf{E}\left[\|\sum_{j \in S_t} h_j \mathbf{X}_{:j}\|^2\right] \leq \sum_{j=1}^n p_j v_j h_j^2 \qquad (4.2)$$

holds for all $\mathbf{h} \in \mathbb{R}^n$. Tight and easily computable formulas for such parameters can be found in [54]. For instance, whenever $\mathbf{Prob}(|S_t| \leq \tau) = 1$, inequality (4.2) holds with $v_j = \tau \|\mathbf{X}_{:j}\|^2$. However, this is a conservative choice of the parameters, as better choices for $v_j$ exist in some specific cases [54]. Convergence of dfSDCA is described in the next theorem.

**Theorem 4.1** ([7]). *Assume that all loss functions $\{\phi_j\}$ are convex and $1/\gamma$ smooth. If we run Algorithm 5 with parameter $\theta$ satisfying the inequality*

$$\theta \leq \min_j \frac{p_j n \lambda \gamma}{v_j + n \lambda \gamma}, \qquad (4.3)$$

*where $\{v_j\}$ satisfy (4.2), then the potential $E_t$ decays exponentially to zero as*

$$\mathbf{E}[E_t] \leq e^{-\theta t} E_0.$$

*Moreover, if we set $\theta$ equal to the upper bound in (4.3) so that*

$$\frac{1}{\theta} = \max_j \left(\frac{1}{p_j} + \frac{v_j}{p_j n \lambda \gamma}\right) \qquad (4.4)$$

*then*

$$t \geq \frac{1}{\theta} \log\left(\frac{(1+\lambda\gamma)E_0}{\lambda\gamma\epsilon}\right) \quad \Rightarrow \quad \mathbf{E}[P(\mathbf{w}^t) - P(\mathbf{w}^*)] \leq \epsilon.$$

## 4.5 Bucket sampling

We shall first explain the concept of "standard" importance sampling.

### 4.5.1 Standard importance sampling

Assume that $\hat{S}$ always picks a single example only. In this case, (4.2) holds for $v_j = \|\mathbf{X}_{:j}\|^2$, independently of $\mathbf{p} := (p_1, \ldots, p_n)$ [54]. This allows us to choose the sampling probabilities as $p_j \sim v_j + n\lambda\gamma$, which ensures that (4.4) is minimized. This is *importance sampling*. The number of iterations of dfSDCA is in this case proportional to

$$\frac{1}{\theta^{(\text{imp})}} := n + \frac{\sum_{j=1}^n v_j}{n\lambda\gamma}.$$

If uniform probabilities are used, the average in the above formula gets replaced by the maximum:

$$\frac{1}{\theta^{(\text{unif})}} := n + \frac{\max_j v_j}{\lambda\gamma}.$$

Hence, one should expect the following *speedup* when comparing the importance and uniform samplings:

$$\sigma := \frac{\max_j \|\mathbf{X}_{:j}\|^2}{\frac{1}{n}\sum_{j=1}^n \|\mathbf{X}_{:j}\|^2}. \qquad (4.5)$$

If $\sigma = 10$ for instance, then dfSDCA with importance sampling is $10\times$ faster than dfSDCA with uniform sampling.

---

[1]ESO = Expected Separable Overapproximation [60, 54].

### 4.5.2 Uniform minibatch sampling

In machine learning, the term "minibatch" is virtually synonymous with a special sampling, which we shall here refer to by the name $\tau$-nice sampling [60]. Sampling $\hat{S}$ is $\tau$-nice if it picks uniformly at random from the collection of all subsets of $[n]$ of cardinality $\tau$. Clearly, $p_j = \tau/n$ and, moreover, it was shown by Qu and Richtárik [54] that (4.2) holds with $\{v_j\}$ defined by

$$v_j^{(\tau\text{-nice})} \quad = \quad \sum_{i=1}^{d}\left(1 + \frac{(\|\mathbf{X}_{i:}\|_0 - 1)(\tau - 1)}{n-1}\right) X_{ij}^2. \tag{4.6}$$

In the case of $\tau$-nice sampling we have the stepsize and complexity given by

$$\theta^{(\tau\text{-nice})} \quad = \quad \min_j \frac{\tau\lambda\gamma}{v_j^{(\tau\text{-nice})} + n\lambda\gamma}, \tag{4.7}$$

$$\frac{1}{\theta^{(\tau\text{-nice})}} \quad = \quad \frac{n}{\tau} + \frac{\max_j v_j^{(\tau\text{-nice})}}{\tau\lambda\gamma}. \tag{4.8}$$

Learning from the difference between the uniform and importance sampling of single example (Section 4.5.1), one would ideally wish the importance minibatch sampling, which we are yet to define, to lead to complexity of the type (4.8), where the maximum is replaced by an average.

### 4.5.3 Bucket sampling: definition

We now propose a family of samplings, which we call *bucket samplings*. Let $B_1, \ldots, B_\tau$ be a partition of $[n] = \{1, 2, \ldots, n\}$ into $\tau$ nonempty sets ("buckets").

**Definition 4.2** (Bucket sampling). We say that $\hat{S}$ is a bucket sampling if for all $k \in [\tau]$, $|\hat{S} \cap B_k| = 1$ with probability 1.

Informally, a bucket sampling picks one example from each of the $\tau$ buckets, forming a minibatch. Hence, $|\hat{S}| = \tau$ and $\sum_{j \in B_l} p_j = 1$ for each $l = 1, 2\ldots, \tau$, where, as before, $p_j := \mathbf{Prob}(j \in \hat{S})$. Notice that given the partition, the vector $\mathbf{p} = (p_1, \ldots, p_n)$ *uniquely determines* a bucket sampling. Hence, we have a family of samplings indexed by a single $n$-dimensional vector. Let $\mathcal{P}_B$ be the set of all vectors $\mathbf{p} \in \mathbb{R}^n$ describing bucket samplings associated with partition $B = \{B_1, \ldots, B_\tau\}$. Clearly,

$$\mathcal{P}_B \quad = \quad \left\{\mathbf{p} \in \mathbb{R}^n : \sum_{j \in B_l} p_j = 1 \text{ for all } l \ \& \ p_j \geq 0 \text{ for all } j\right\}.$$

Note, that the sampling inside each bucket $B_j$ can be performed in $\mathcal{O}(\log |B_j|)$ time using a binary tree, with an initial overhead and memory of $\mathcal{O}(|B_j| \log |B_j|)$, as explained in [49]

### 4.5.4 Optimal bucket sampling

The optimal bucket sampling is that for which (4.4) is minimized, which leads to a complicated optimization problem:

$$\min_{\mathbf{p} \in \mathcal{P}_B} \max_j \frac{1}{p_j} + \frac{v_j}{p_j n \lambda \gamma} \quad \text{subject to } \{v_j\} \text{ satisfy (4.2).}$$

A particular difficulty here is the fact that the parameters $\{v_j\}$ depend on the vector $\mathbf{p}$ in a complicated way. In order to resolve this issue, we prove the following result.

**Theorem 4.3.** *Let $\hat{S}$ be a bucket sampling described by partition $B = \{B_1, \ldots, B_\tau\}$ and vector $\mathbf{p}$. Then the ESO inequality (4.2) holds for parameters $\{v_j\}$ set to*

$$v_j \quad = \quad \sum_{i=1}^{d}\left(1 + \left(1 - \tfrac{1}{\omega_i'}\right)\delta_i\right) X_{ij}^2, \tag{4.9}$$

where $J_i := \{j \in [n] \; : \; X_{ij} \neq 0\}$, $\delta_i := \sum_{j \in J_i} p_j$ and $\omega_i' := |\{l \; : \; J_i \cap B_l \neq \emptyset\}|$.

Observe that $J_i$ is the set of examples which express feature $i$, and $\omega_i'$ is the number of buckets intersecting with $J_i$. Clearly, that $1 \leq \omega_i' \leq \tau$ (if $\omega_i' = 0$, we simply discard this feature from our data as it is not needed). Note that the effect of the quantities $\{\omega_i'\}$ on the value of $v_j$ is small. Indeed, unless we are in the extreme situation when $\omega_i' = 1$, which has the effect of neutralizing $\delta_i$, the quantity $1 - 1/\omega_i'$ is between $1 - 1/2$ and $1 - 1/\tau$. Hence, for simplicity, we could instead use the slightly more conservative parameters:

$$v_j = \sum_{i=1}^{d} \left(1 + \left(1 - \frac{1}{\tau}\right)\delta_i\right)X_{ij}^2.$$

### 4.5.5 Uniform bucket sampling

Assume all buckets are of the same size: $|B_l| = n/\tau$ for all $l$. Further, assume that $p_j = 1/|B_l| = \tau/n$ for all $i$. Then $\delta_i = \tau|J_i|/n$, and hence Theorem 4.3 says that

$$v_j^{(\text{unif})} = \sum_{i=1}^{d} \left(1 + \left(1 - \frac{1}{\omega_i'}\right)\frac{\tau|J_i|}{n}\right)\mathbf{X}_{ij}^2, \tag{4.10}$$

and in view of (4.4), the complexity of dfSDCA with this sampling becomes

$$\frac{1}{\theta^{(\text{unif})}} = \frac{n}{\tau} + \frac{\max_j v_j^{(\text{unif})}}{\tau \lambda \gamma}. \tag{4.11}$$

Formula (4.6) is very similar to the one for $\tau$-nice sampling (4.10), despite the fact that the sets/minibatches generated by the uniform bucket sampling have a special structure with respect to the buckets. Indeed, it is easily seen that the difference between between $1 + \frac{\tau|J_i|}{n}$ and $1 + \frac{(\tau-1)(|J_i|-1)}{(n-1)}$ is negligible. Moreover, if either $\tau = 1$ or $|J_i| = 1$ for all $i$, then $\omega_i' = 1$ for all $i$ and hence $v_j = \|\mathbf{X}_{:j}\|^2$. This is also what we get for the $\tau$-nice sampling.

## 4.6 Importance minibatch sampling

In the light of Theorem 4.3, we can formulate the problem of searching for the optimal bucket sampling as

$$\min_{\mathbf{p} \in \mathcal{P}_B} \max_j \frac{1}{p_j} + \frac{v_j}{p_j n \lambda \gamma} \quad \text{subject to } \{v_j\} \text{ satisfy (4.9)}. \tag{4.12}$$

Still, this is not an easy problem. *Importance minibatch sampling* arises as an approximate solution of (4.12). Note that the uniform minibatch sampling is a feasible solution of the above problem, and hence we should be able to improve upon its performance.

### 4.6.1 Approach 1: alternating optimization

Given a probability distribution $\mathbf{p} \in \mathcal{P}_B$, we can easily find $\mathbf{v}$ using Theorem 4.3. On the other hand, for any fixed $\mathbf{v}$, we can minimize (4.12) over $\mathbf{p} \in \mathcal{P}_B$ by choosing the probabilities in each group $B_l$ and for each $j \in B_l$ via

$$p_j = \frac{n\lambda\gamma + v_j}{\sum_{k \in B_l} n\lambda\gamma + v_k}. \tag{4.13}$$

This leads to a natural alternating optimization strategy. An example of the standard convergence behaviour of this scheme is showed in Table 4.1. Empirically, this strategy converges to a pair $(\mathbf{p}^*, \mathbf{v}^*)$ for which (4.13) holds. Therefore, the resulting complexity will be

$$\frac{1}{\theta^{(\tau\text{-imp})}} = \frac{n}{\tau} + \max_{l \in [\tau]} \frac{\frac{\tau}{n}\sum_{j \in B_l} v_j^*}{\tau \lambda \gamma}. \tag{4.14}$$

| quantity \ iteration | 1 | 2 | 3 | 4 | 5 | 6 |
|---|---|---|---|---|---|---|
| $\max_j(|p_j^{\text{new}} - p_j^{\text{old}}|)$ | $7 \cdot 10^{-5}$ | $7 \cdot 10^{-6}$ | $7 \cdot 10^{-7}$ | $8 \cdot 10^{-8}$ | $8 \cdot 10^{-9}$ | $9 \cdot 10^{-10}$ |
| $\|\mathbf{p}^{\text{new}} - \mathbf{p}^{\text{old}}\|_2$ | $1 \cdot 10^{-3}$ | $2 \cdot 10^{-4}$ | $2 \cdot 10^{-5}$ | $2 \cdot 10^{-6}$ | $2 \cdot 10^{-7}$ | $2 \cdot 10^{-8}$ |

Table 4.1: Example of the convergence speed of the alternating optimization scheme for $w8a$ dataset (see Table 4.5) with $\tau = 8$. The table demonstrates the difference in probabilities for two successive iterations ($\mathbf{p}^{\text{old}}$ and $\mathbf{p}^{\text{new}}$). We observed a similar behaviour for all datasets and all choices of $\tau$.

We can compare this result against the complexity of $\tau$-nice in (4.8). We can observe that the terms are very similar, up to two differences. First, the importance minibatch sampling has a maximum over group averages instead of a maximum over everything, which leads to speedup, other things equal. On the other hand, $\mathbf{v}^{(\tau\text{-nice})}$ and $\mathbf{v}^*$ are different quantities. The alternating optimization procedure for computation of $(\mathbf{v}^*, \mathbf{p}^*)$ is costly, as one iteration takes a pass over all data. Therefore, in the next subsection we propose a closed form formula which, as we found empirically, offers nearly optimal convergence rate.

### 4.6.2 Approach 2: practical formula

For each group $B_l$, let us choose for all $j \in B_l$ the probabilities as follows:

$$p_j^* = \frac{n\lambda\gamma + v_j^{(\text{unif})}}{\sum_{k \in B_l} n\lambda\gamma + v_k^{(\text{unif})}} \qquad (4.15)$$

where $v_j^{(\text{unif})}$ is given by (4.10). Note that computing all $v_j^{(\text{unif})}$ can be done directly in one pass over the data, and computing all $p_j^*$ is a simple re-weighting. This is the same computational cost as for standard serial importance sampling. Also, this process can be straightforwardly parallelized, which leads to additional savings in time. The overhead of using this sampling approach is therefore at most one pass over the data, which is negligible in most scenarios considered.

After doing some simplifications, the associated complexity result is

$$\frac{1}{\theta^{(\tau\text{-imp})}} = \max_l \left\{ \left( \frac{n}{\tau} + \frac{\frac{\tau}{n} \sum_{j \in B_l} v_j^{(\text{unif})}}{\tau\lambda\gamma} \right) \beta_l \right\}, \qquad (4.16)$$

where

$$\beta_l := \max_{j \in B_l} \frac{n\lambda\gamma + s_j}{n\lambda\gamma + v_j^{(\text{unif})}}$$

and

$$s_j := \sum_{i=1}^d \left( 1 + \left( 1 - \frac{1}{\omega_i'} \right) \sum_{k \in J_i} p_k^* \right) X_{ij}^2.$$

We would ideally want to have $\beta_l = 1$ for all $l$ (this is what we get for importance sampling without minibatches). If $\beta_l \approx 1$ for all $l$, then the complexity $1/\theta^{(\tau\text{-imp})}$ is an improvement on the complexity of the uniform minibatch sampling since the maximum of group averages is always better than the maximum of all elements $v_j^{(\text{uni})}$:

$$\frac{n}{\tau} + \frac{\max_l \left( \frac{\tau}{n} \sum_{j \in B_l} v_j^{(\text{unif})} \right)}{\tau\lambda\gamma} \leq \frac{n}{\tau} + \frac{\max_j v_j^{(\text{unif})}}{\tau\lambda\gamma}.$$

Indeed, the difference can be very large.

Finally, we would like to comment on the choice of the partitions $B_1, \ldots, B_\tau$, as they clearly affect the convergence rate. The optimal choice of the partitions is given by minimizing in $B_1, \ldots, B_\tau$ the maximum over group sums in (4.16), which is a complicated optimization problem. Instead, we used random partitions of the same size in our experiments, which we believe is a good solution for the partitioning problem. The logic is simple: the minimum of the maximum over the group sums will be achieved, when all the group sums have similar values. If we set the partitions to the same size and we distribute the examples randomly, there is a good chance that the group sums will have similar values (especially for large amounts of data).

## 4.7 Experiments

We now comment on the results of our numerical experiments, with both synthetic and real datasets. We plot the optimality gap $P(\mathbf{w}^t) - P(\mathbf{w}^*)$ and in the case of real data also the test error (vertical axis) against the computational effort (horizontal axis). We measure computational effort by the number of effective passes through the data divided by $\tau$. We divide by $\tau$ as a normalization factor; since we shall compare methods with a range of values of $\tau$. This is reasonable as it simply indicates that the $\tau$ updates are performed in parallel. Hence, what we plot is an implementation-independent model for time.

We compared two algorithms:

1) $\tau$-**nice**: dfSDCA using the $\tau$-nice sampling with stepsizes given by (4.7) and (4.6),

2) $\tau$-**imp**: dfSDCA using $\tau$-importance sampling (i.e., importance minibatch sampling) defined in Subsection 4.6.2.

As the methods are randomized, we always plot the average over 5 runs. For each dataset we provide two plots. In the left figure we plot the convergence of $\tau$-nice for different values of $\tau$, and in the right figure we do the same for $\tau$-importance. The horizontal axis has the same range in both plots, so they are easily comparable. The values of $\tau$ we used to plot are $\tau \in \{1, 2, 4, 8, 16, 32\}$. In all experiments we used the logistic loss: $\phi_j(z) = \log(1 + e^{-y_j z})$ and set the regularizer to $\lambda = \max_j \|\mathbf{X}_{:j}\|/n$. We will observe the theoretical and empirical ratio $\theta^{(\tau\text{-imp})}/\theta^{(\tau\text{-nice})}$. The theoretical ratio is computed from the corresponding theory. The empirical ratio is the ratio between the horizontal axis values at the moments when the algorithms reached the precision $10^{-10}$.

### 4.7.1 Artificial data

We start with experiments using artificial data, where we can control the sparsity pattern of $\mathbf{X}$ and the distribution of $\{\|\mathbf{X}_{:j}\|^2\}$. We fix $n = 50,000$ and choose $d = 10,000$ and $d = 1,000$. For each feature we sampled a random sparsity coefficient $\omega'_i \in [0, 1]$ to have the average sparsity $\omega' := \frac{1}{d}\sum_j^d \omega'_i$ under control. We used two different regimes of sparsity: $\omega' = 0.1$ (10% nonzeros) and $\omega' = 0.8$ (80% nonzeros). After deciding on the sparsity pattern, we rescaled the examples to match a specific distribution of norms $L_j = \|\mathbf{X}_{:j}\|^2$; see Table 4.2. The code column shows the corresponding code in Julia to create the vector of norms $L$. The distributions can be also observed as histograms in Figure 4.1.

| label | code | $\sigma$ |
|---|---|---|
| extreme | L = ones(n);L[1] = 1000 | 980.4 |
| chisq1 | L = rand(chisq(1),n) | 17.1 |
| chisq10 | L = rand(chisq(10),n) | 3.9 |
| chisq100 | L = rand(chisq(100),n) | 1.7 |
| uniform | L = 2*rand(n) | 2.0 |

Table 4.2: Distributions of $\|\mathbf{X}_{:j}\|^2$ used in artificial experiments.

The corresponding experiments can be found in Figure 4.4 and Figure 4.5. The theoretical and empirical speedup are also summarized in Tables 4.3 and 4.4.

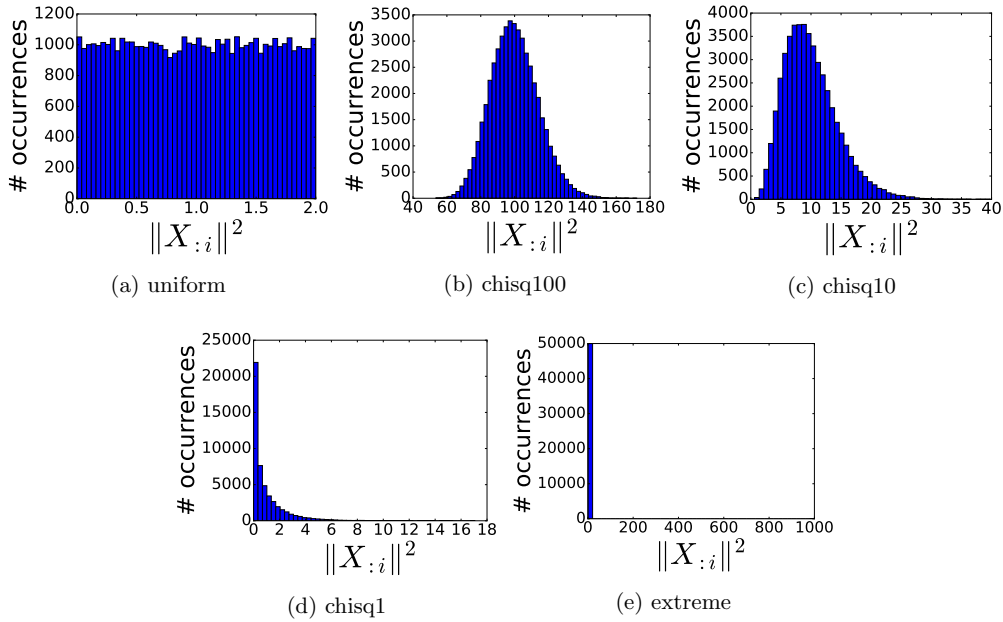

Figure 4.1: The distribution of $\|\mathbf{X}_{:i}\|^2$ for synthetic data

| Data | $\tau = 1$ | $\tau = 2$ | $\tau = 4$ | $\tau = 8$ | $\tau = 16$ | $\tau = 32$ |
|---|---|---|---|---|---|---|
| uniform | 1.2 : 1.0 | 1.2 : 1.1 | 1.2 : 1.1 | 1.2 : 1.1 | 1.3 : 1.1 | 1.4 : 1.1 |
| chisq100 | 1.5 : 1.3 | 1.5 : 1.3 | 1.5 : 1.4 | 1.6 : 1.4 | 1.6 : 1.4 | 1.6 : 1.4 |
| chisq10 | 1.9 : 1.4 | 1.9 : 1.5 | 2.0 : 1.4 | 2.2 : 1.5 | 2.5 : 1.6 | 2.8 : 1.7 |
| chisq1 | 1.9 : 1.4 | 2.0 : 1.4 | 2.2 : 1.5 | 2.5 : 1.6 | 3.1 : 1.6 | 4.2 : 1.7 |
| extreme | 8.8 : 4.8 | 9.6 : 6.6 | 11 : 6.4 | 14 : 6.4 | 20 : 6.9 | 32 : 6.1 |

Table 4.3: The **theoretical** : **empirical** ratios $\theta^{(\tau\text{-imp})}/\theta^{(\tau\text{-nice})}$ for sparse artificial data ($\omega' = 0.1$)

| Data | $\tau = 1$ | $\tau = 2$ | $\tau = 4$ | $\tau = 8$ | $\tau = 16$ | $\tau = 32$ |
|---|---|---|---|---|---|---|
| uniform | 1.2 : 1.1 | 1.2 : 1.1 | 1.4 : 1.2 | 1.5 : 1.2 | 1.7 : 1.3 | 1.8 : 1.3 |
| chisq100 | 1.5 : 1.3 | 1.6 : 1.4 | 1.6 : 1.5 | 1.7 : 1.5 | 1.7 : 1.6 | 1.7 : 1.6 |
| chisq10 | 1.9 : 1.3 | 2.2 : 1.6 | 2.7 : 2.1 | 3.1 : 2.3 | 3.5 : 2.5 | 3.6 : 2.7 |
| chisq1 | 1.9 : 1.3 | 2.6 : 1.8 | 3.7 : 2.3 | 5.6 : 2.9 | 7.9 : 3.2 | 10 : 3.9 |
| extreme | 8.8 : 5.0 | 15 : 7.8 | 27 : 12 | 50 : 16 | 91 : 21 | 154 : 28 |

Table 4.4: The **theoretical** : **empirical** ratios $\theta^{(\tau\text{-imp})}/\theta^{(\tau\text{-nice})}$. Artificial data with $\omega' = 0.8$ (dense)

### 4.7.2 Real data

We used several publicly available datasets[2], summarized in Table 4.5, which we randomly split into a train (80%) and a test (20%) part. The test error is measured by the empirical risk (4.1) on the test data without a regularizer. The resulting test error was compared against the best achievable test error, which we computed by minimizing the corresponding risk. Experimental results are in Figure 4.2 and Figure 4.3. The theoretical and empirical speedup table for these datasets can be found in Table 4.6.

---

[2]https://www.csie.ntu.edu.tw/ cjlin/libsvmtools/datasets/

| Dataset | #samples | #features | sparsity | $\sigma$ |
|---------|----------|-----------|----------|------|
| ijcnn1  | 35,000   | 23        | 60.1%    | 2.04 |
| protein | 17,766   | 358       | 29.1%    | 1.82 |
| w8a     | 49,749   | 301       | 4.2%     | 9.09 |
| url     | 2,396,130 | 3,231,962 | 0.04 % | 4.83 |
| aloi    | 108,000  | 129       | 24.6%    | 26.01 |

Table 4.5: Summary of real data sets ($\sigma$ = predicted speedup).

| Data | $\tau = 1$ | $\tau = 2$ | $\tau = 4$ | $\tau = 8$ | $\tau = 16$ | $\tau = 32$ |
|------|-----------|-----------|-----------|-----------|------------|------------|
| ijcnn1  | 1.2 : 1.1 | 1.4 : 1.1 | 1.6 : 1.3 | 1.9 : 1.6 | 2.2 : 1.6 | 2.3 : 1.8 |
| protein | 1.3 : 1.2 | 1.4 : 1.2 | 1.5 : 1.4 | 1.7 : 1.4 | 1.8 : 1.5 | 1.9 : 1.5 |
| w8a     | 2.8 : 2.0 | 2.9 : 1.9 | 2.9 : 1.9 | 3.0 : 1.9 | 3.0 : 1.8 | 3.0 : 1.8 |
| url     | 3.0 : 2.3 | 2.6 : 2.1 | 2.0 : 1.8 | 1.7 : 1.6 | 1.8 : 1.6 | 1.8 : 1.7 |
| aloi    | 13 : 7.8  | 12 : 8.0  | 11 : 7.7  | 9.9 : 7.4 | 9.3 : 7.0 | 8.8 : 6.7 |

Table 4.6: The **theoretical** : **empirical** ratios $\theta^{(\tau\text{-imp})}/\theta^{(\tau\text{-nice})}$.

### 4.7.3 Conclusion

In all experiments, $\tau$-importance sampling performs significantly better than $\tau$-nice sampling. The theoretical speedup factor computed by $\theta^{(\tau\text{-imp})}/\theta^{(\tau\text{-nice})}$ provides an excellent estimate of the actual speedup. We can observe that on denser data the speedup is higher than on sparse data. This matches the theoretical intuition for $v_j$ for both samplings. Similar behaviour can be also observed for the test error, which is pleasing. As we observed for artificial data, for extreme datasets the speedup can be arbitrary large, even several orders of magnitude. *A rule of thumb: if one has data with large $\sigma$, practical speedup from using importance minibatch sampling will likely be dramatic.*

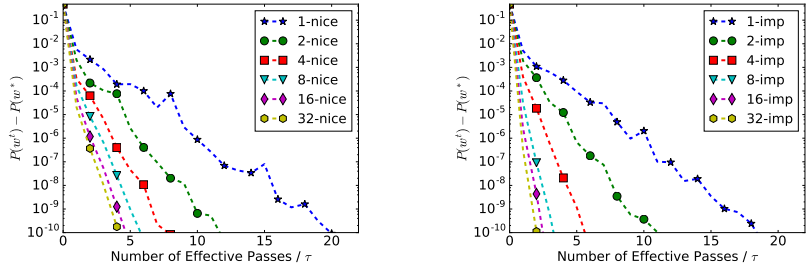

(a) ijcnn1, $\tau$-nice (left), $\tau$-importance (right)

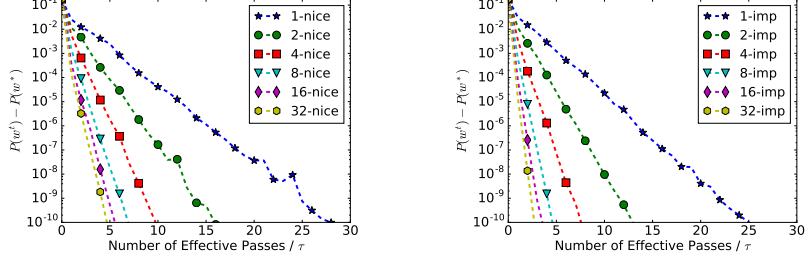

(b) protein, $\tau$-nice (left), $\tau$-importance (right)

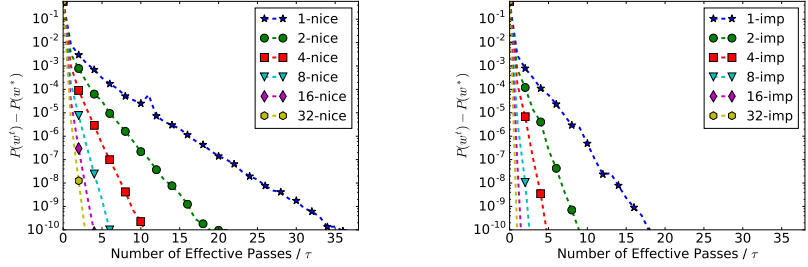

(c) w8a, $\tau$-nice (left), $\tau$-importance (right)

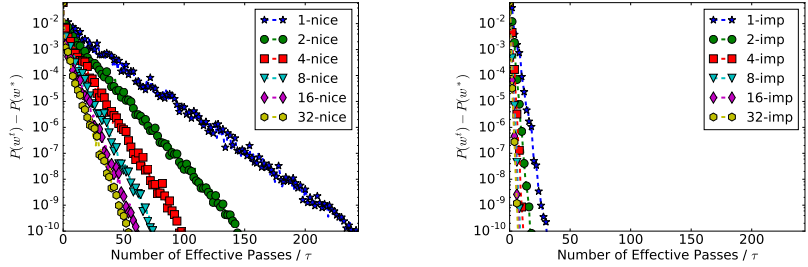

(d) aloi, $\tau$-nice (left), $\tau$-importance (right)

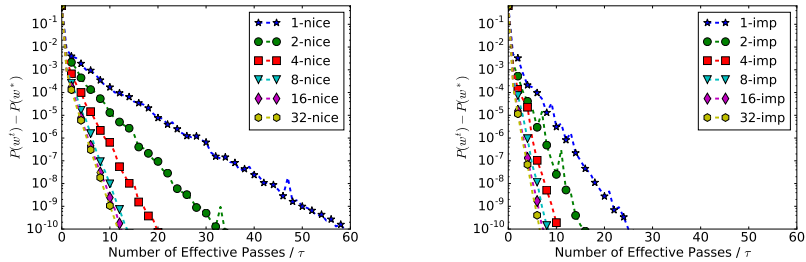

(e) url, $\tau$-nice (left), $\tau$-importance (right)

Figure 4.2: Train error for datasets from Table 4.5

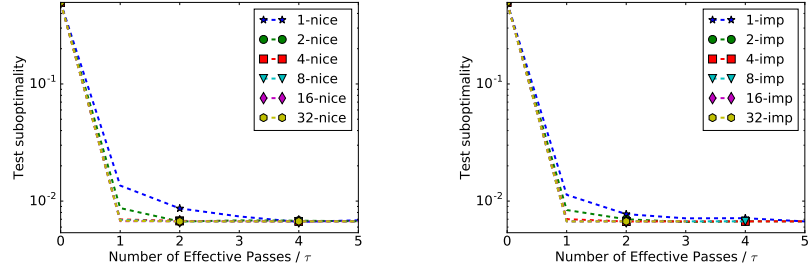

(a) ijcnn1, $\tau$-nice (left), $\tau$-importance (right)

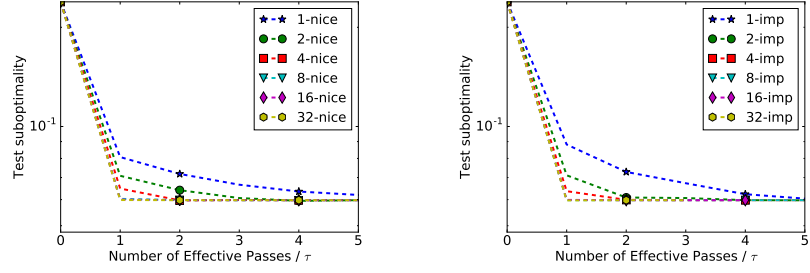

(b) protein, $\tau$-nice (left), $\tau$-importance (right)

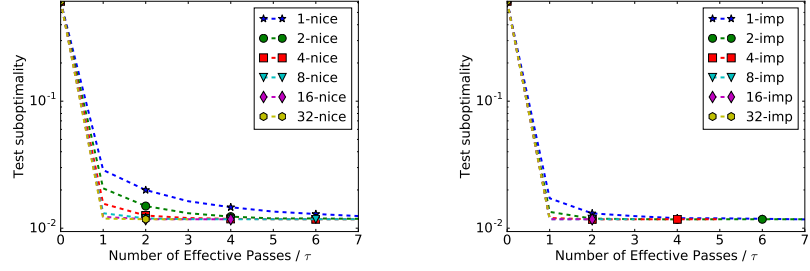

(c) w8a, $\tau$-nice (left), $\tau$-importance (right)

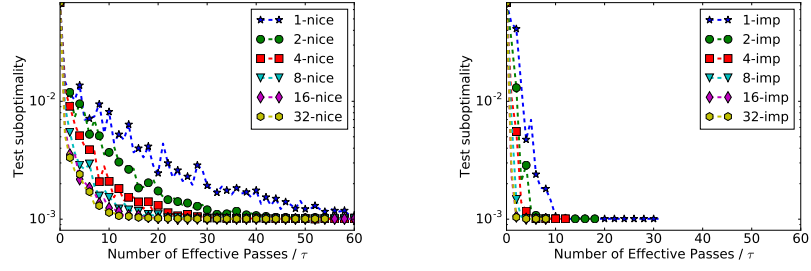

(d) aloi, $\tau$-nice (left), $\tau$-importance (right)

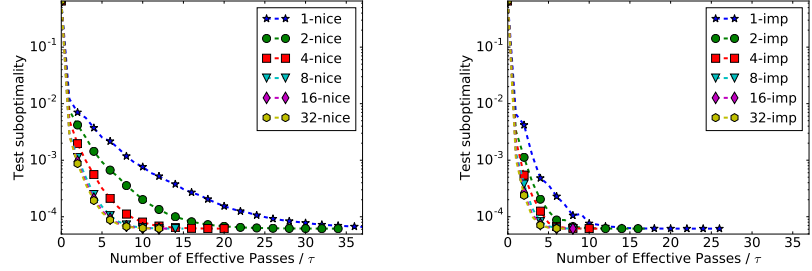

(e) url, $\tau$-nice (left), $\tau$-importance (right)

Figure 4.3: Test error for datasets from Table 4.5

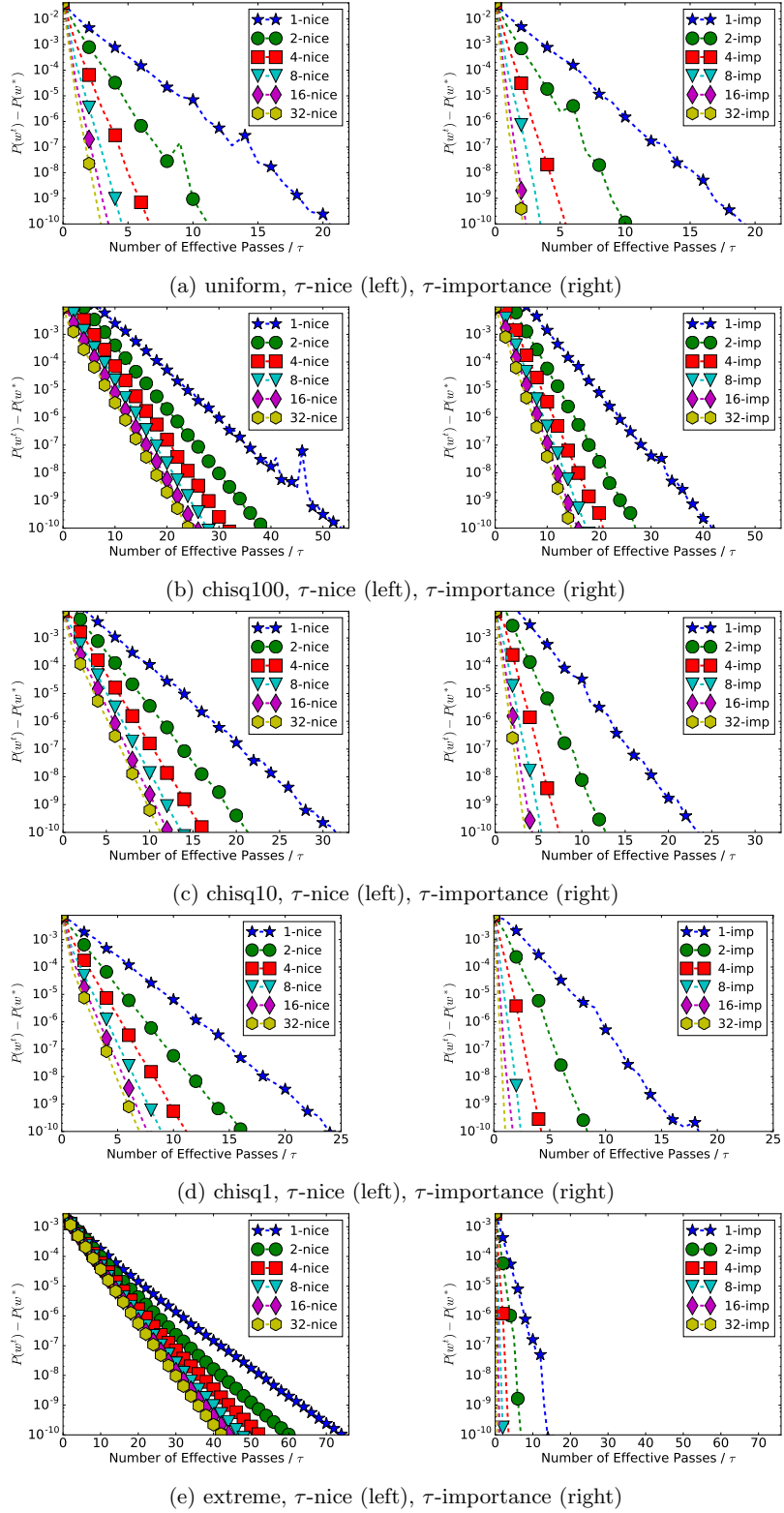

(a) uniform, $\tau$-nice (left), $\tau$-importance (right)

(b) chisq100, $\tau$-nice (left), $\tau$-importance (right)

(c) chisq10, $\tau$-nice (left), $\tau$-importance (right)

(d) chisq1, $\tau$-nice (left), $\tau$-importance (right)

(e) extreme, $\tau$-nice (left), $\tau$-importance (right)

Figure 4.4: Artificial datasets from Table 4.2 with $\omega = 0.8$

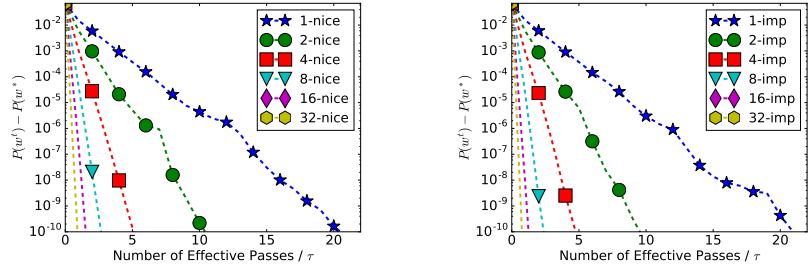
(a) uniform, $\tau$-nice (left), $\tau$-importance (right)

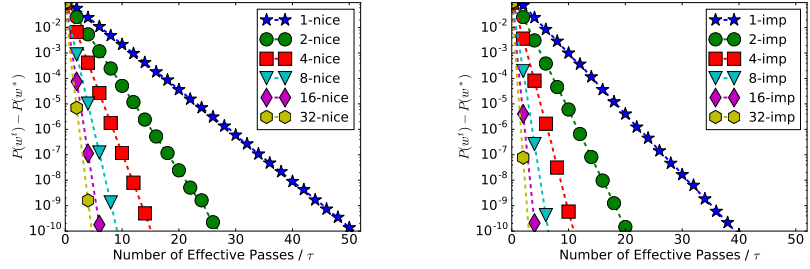
(b) chisq100, $\tau$-nice (left), $\tau$-importance (right)

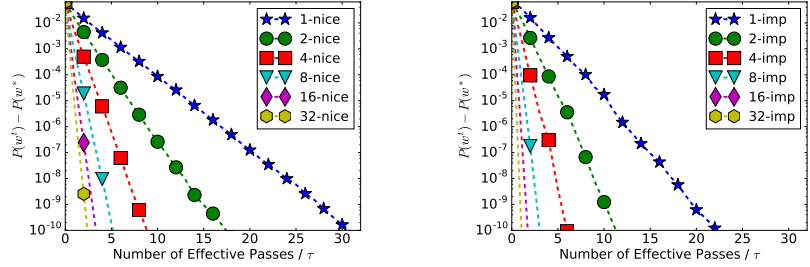
(c) chisq10, $\tau$-nice (left), $\tau$-importance (right)

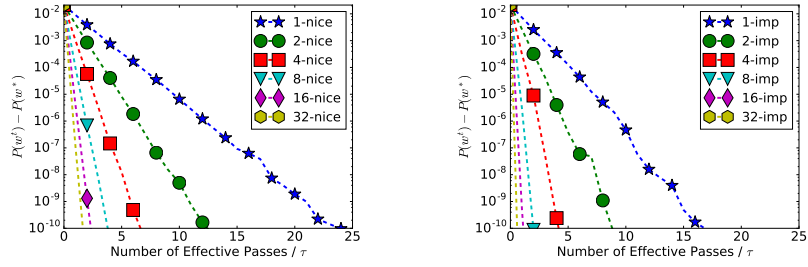
(d) chisq1, $\tau$-nice (left), $\tau$-importance (right)

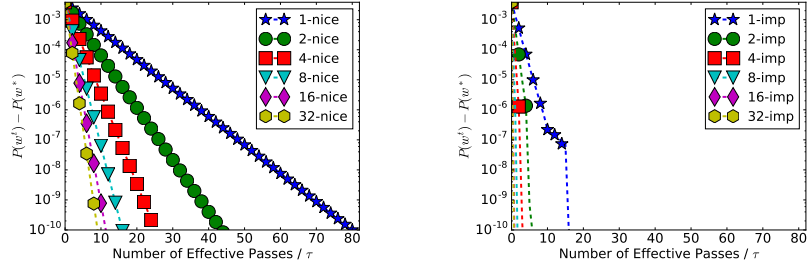
(e) extreme, $\tau$-nice (left), $\tau$-importance (right)

Figure 4.5: Artificial datasets from Table 4.2 with $\omega = 0.1$

## 4.A Proof of Theorem 4.3

### 4.A.1 Three lemmas

We first establish three lemmas, and then proceed with the proof of the main theorem. With each sampling $\hat{S}$ we associate an $n \times n$ "probability matrix" defined as follows: $P_{ij}(\hat{S}) = \mathbf{Prob}(i \in \hat{S}, j \in \hat{S})$. Our first lemma characterizes the probability matrix of the bucket sampling.

**Lemma 4.4.** *If $\hat{S}$ is a bucket sampling, then*

$$\mathbf{P}(\hat{S}) \quad = \quad \mathbf{p}\mathbf{p}^\top \circ (\mathbf{E} - \mathbf{B}) + \mathrm{Diag}(\mathbf{p}), \tag{4.17}$$

*where $\mathbf{E} \in \mathbb{R}^{n \times n}$ is the matrix of all ones,*

$$\mathbf{B} \quad := \quad \sum_{l=1}^{\tau} \mathbf{P}(B_l), \tag{4.18}$$

*and $\circ$ denotes the Hadamard (elementwise) product of matrices. Note that $\mathbf{B}$ is the 0-1 matrix given by $B_{ij} = 1$ if and only if $i, j$ belong to the same bucket $B_l$ for some $l$.*

*Proof.* Let $\mathbf{P} = \mathbf{P}(\hat{S})$. By definition

$$P_{ij} = \begin{cases} p_i & i = j \\ p_i p_j & i \in B_l, \ j \in B_k, \ l \neq k \\ 0 & \text{otherwise.} \end{cases}$$

It only remains to compare this to (4.17). □

**Lemma 4.5.** *Let $J$ be a nonempty subset of $[n]$, let $\mathbf{B}$ be as in Lemma 4.4 and put $\omega'_J := |\{l : J \cap B_l \neq \emptyset\}|$. Then*

$$\mathbf{P}(J) \circ \mathbf{B} \quad \succeq \quad \frac{1}{\omega'_J} \mathbf{P}(J). \tag{4.19}$$

*Proof.* For any $\mathbf{h} \in \mathbb{R}^n$, we have

$$\mathbf{h}^\top \mathbf{P}(J) \mathbf{h} \quad = \quad \left(\sum_{j \in J} h_j\right)^2 \quad = \quad \left(\sum_{l=1}^{\tau} \sum_{j \in J \cap B_l} h_j\right)^2$$

$$\leq \quad \omega'_J \sum_{l=1}^{\tau} \left(\sum_{j \in J \cap B_l} h_j\right)^2 \quad = \quad \omega'_J \sum_{l=1}^{\tau} \mathbf{h}^\top \mathbf{P}(J \cap B_l) \mathbf{h},$$

where we used the Cauchy-Schwarz inequality. Using this, we obtain

$$\mathbf{P}(J) \circ \mathbf{B} \quad \stackrel{(4.18)}{=} \quad \mathbf{P}(J) \circ \sum_{l=1}^{\tau} \mathbf{P}(B_l) = \sum_{l=1}^{\tau} \mathbf{P}(J) \circ \mathbf{P}(B_l)$$

$$= \quad \sum_{l=1}^{\tau} \mathbf{P}(J \cap B_l) \stackrel{(4.A.1)}{\succeq} \frac{1}{\omega'} \mathbf{P}(J).$$

□

**Lemma 4.6.** *Let $J$ be any nonempty subset of $[n]$ and $\hat{S}$ be a bucket sampling. Then*

$$\mathbf{P}(J) \circ \mathbf{p}\mathbf{p}^\top \quad \preceq \quad \left(\sum_{j \in J} p_j\right) \mathrm{Diag}(\mathbf{P}(J \cap \hat{S})). \tag{4.20}$$

*Proof.* Choose any $\mathbf{h} \in \mathbb{R}^n$ and note that

$$\mathbf{h}^\top (\mathbf{P}(J) \circ \mathbf{p}\mathbf{p}^\top)\mathbf{h} = \left(\sum_{j \in J} p_j h_i\right)^2 = \left(\sum_{j \in J} s_j t_j\right)^2,$$

where $s_j = \sqrt{p_j} h_j$ and $t_j = \sqrt{p_j}$. It remains to apply the Cauchy-Schwarz inequality:

$$\sum_{j \in J} s_j t_j \leq \sum_{j \in J} s_j^2 \sum_{j \in J} t_j^2$$

and notice that the $j$-th element on the diagonal of $\mathbf{P}(J \cap \hat{S})$ is $p_j$ for $j \in J$ and 0 for $j \notin J$ □

### 4.A.2 Proof of Theorem 4.3

By Theorem 5.2 in [54], we know that inequality (4.2) holds for parameters $\{v_j\}$ set to

$$v_j = \sum_{i=1}^{d} \lambda'(\mathbf{P}(J_i \cap \hat{S})) X_{ij}^2,$$

where $\lambda'(\mathbf{M})$ is the largest normalized eigenvalue of symmetric matrix $\mathbf{M}$ defined as

$$\lambda'(\mathbf{M}) := \max_{\mathbf{h}} \left\{ \mathbf{h}^\top \mathbf{M} \mathbf{h} \ : \ \mathbf{h}^\top \mathrm{Diag}(\mathbf{M})\mathbf{h} \leq 1 \right\}.$$

Furthermore,

$$\begin{aligned}
\mathbf{P}(J_i \cap \hat{S}) &= \mathbf{P}(J_i) \circ \mathbf{P}(\hat{S}) \\
&\stackrel{(4.17)}{=} \mathbf{P}(J_i) \circ \mathbf{p}\mathbf{p}^\top - \mathbf{P}(J_i) \circ \mathbf{p}\mathbf{p}^\top \circ \mathbf{B} + \mathbf{P}(J_i) \circ \mathrm{Diag}(\mathbf{p}) \\
&\stackrel{(4.19)}{\preceq} \left(1 - \frac{1}{\omega'_J}\right) \mathbf{P}(J_i) \circ \mathbf{p}\mathbf{p}^\top + \mathbf{P}(J_i) \circ \mathrm{Diag}(\mathbf{p}) \\
&\stackrel{(4.20)}{\preceq} \left(1 - \frac{1}{\omega'_J}\right) \delta_i \, \mathrm{Diag}(\mathbf{P}(J_i \cap \hat{S})) + \mathrm{Diag}(\mathbf{P}(J_i \cap \hat{S})),
\end{aligned}$$

whence $\lambda'(\mathbf{P}(J_i \cap \hat{S})) \leq 1 + (1 - 1/\omega'_J) \delta_i$, which concludes the proof.

# Chapter 5

# Coordinate Descent Faceoff: Primal or Dual?

## 5.1 Introduction

In the last 5 years or so, randomized coordinate descent (RCD) methods [69, 49, 59, 60] have become immensely popular in a variety of machine learning tasks, with supervised learning being a prime example. The main reasons behind the rise of RCD-type methods is that they can be easily implemented, have intuitive appeal, and enjoy superior theoretical and practical behaviour when compared to classical methods such as SGD [61], especially in high dimensions, and in situations when solutions of medium to high accuracy are needed. One of the most important success stories of RCD is in the domain of training linear predictors via regularized empirical risk minimization (ERM).

The highly popular SDCA algorithm [72] arises as the application of RCD [59] to the *dual problem* associated with the (primal) ERM problem[1]. In practice, SDCA is most effective in situations where the number of examples ($n$) exceeds the number of features ($d$). Since the dual of ERM is an $n$ dimensional problem, it makes intuitive sense to apply RCD to the dual. Indeed, RCD can be seen as a randomized decomposition strategy, reducing the $n$ dimensional problem to a sequence of (randomly generated) one-dimensional problems.

However, if the number of features exceeds the number of examples, and especially when the difference is very large, RCD methods [60] have been found very attractive for solving the *primal problem* (i.e., the ERM problem) directly. For instance, distributed variants of RCD, such as Hydra [58] and its accelerated cousin Hydra$^2$ [16] have been successfully applied to solving problems with billions of features.

Recently, a variety of novel primal methods for ERM have been designed, including SAG [64], SVRG [27], S2GD [32], proxSVRG [84], mS2GD [30], SAGA [12], MISO [42] and S2CD [31]. As SDCA, all these methods improve dramatically on SGD [61] as a benchmark, which they achieve by employing one of a number of variance-reduction strategies. These methods enjoy essentially identical theoretical complexity bounds as SDCA. In this sense, conclusions based on our study complexity of primal RCD vs dual RCD are valid also when comparing primal RCD with appropriate variants of any of the above mentioned methods (e.g., SVRG). For simplicity, we do not explore this further in this chapter, and instead focus on comparing primal versus dual RCD.

### 5.1.1 Contributions

In this chapter we provide the first joint study of these two approaches—applying RCD to the primal vs dual problems—and we do so in the context of L2-regularized linear ERM. First, we show through a rigorous theoretical analysis that for dense data, the intuition that the primal approach is better than the dual approach when $n \leq d$, and vice versa, is precisely correct.

---

[1]Indeed, the analysis of SDCA in [72] proceeds by applying the complexity result from [59] to the *dual problem*, and then arguing that the same rate applies to the primal suboptimality as well.

However, we show that for sparse data, this does not need to be the case: primal RCD can significantly outperform dual RCD even if $d \ll n$, and vice versa, dual RCD can be much faster than primal RCD even if $n \ll d$. In particular, we identify that the face-off between primal and dual RCD boils down to the comparison of as single quantity associated with the data matrix and its transpose. Moreover, we show that, surprisingly, a single sampling strategy minimizes both the (bound on the) number of iterations and the overall expected complexity of RCD. Note that the latter complexity measure takes into account also the average cost of the iterations, which depends on the structure and sparsity of the data, and on the sampling strategy employed. We confirm our theoretical findings using extensive experiments with both synthetic and real data sets.

## 5.2 Primal and dual formulations of ERM

Let $\mathbf{X} \in \mathbb{R}^{d \times n}$ be a data matrix, with $n$ referring to the number of examples and $d$ to the number of features. With each example $\mathbf{X}_{:j} \in \mathbb{R}^d$ we associate a loss function $\phi_j : \mathbb{R} \to \mathbb{R}$, and pick a regularization constant $\lambda > 0$. The key problem of this chapter is the L2-regularized ERM problem for linear models

$$\min_{\mathbf{w} \in \mathbb{R}^d} \left[ P(\mathbf{w}) \quad := \quad \frac{1}{n} \sum_{j=1}^n \phi_j(\langle \mathbf{X}_{:j}, \mathbf{w} \rangle) + \frac{\lambda}{2} \|\mathbf{w}\|_2^2 \right], \tag{5.1}$$

where $\langle \cdot, \cdot \rangle$ denotes the standard Euclidean inner product and $\|\mathbf{w}\|_2 := \sqrt{\langle \mathbf{w}, \mathbf{w} \rangle}$. We refer to (5.1) as the *primal problem*. We assume throughout that the functions $\{\phi_j\}$ are convex and $1/\gamma$-smooth, which is given by the bounds

$$\phi_j(s) + \phi_j'(s)t \quad \leq \quad \phi_j(s+t) \quad \leq \quad \phi_j(s) + \phi_j'(s)t + \frac{1}{2\gamma}t^2, \tag{5.2}$$

for all $s, t \in \mathbb{R}$. The *dual problem* of (5.1) is

$$\max_{\boldsymbol{\alpha} \in \mathbb{R}^n} \left[ D(\boldsymbol{\alpha}) \quad := \quad -\frac{1}{2\lambda n^2} \|\mathbf{X}\boldsymbol{\alpha}\|_2^2 - \frac{1}{n} \sum_{j=1}^n \phi_j^*(-\alpha_j) \right], \tag{5.3}$$

where $\phi_j^* : \mathbb{R} \to \mathbb{R}$ is the convex conjugate of $\phi_j$, defined by $\phi_j^*(s) := \sup\{st - \phi_j(t) \,:\, t \in \mathbb{R}\}$. It is well known that that $P(\mathbf{w}) \geq D(\boldsymbol{\alpha})$ for every pair $(\mathbf{w}, \boldsymbol{\alpha}) \in \mathbb{R}^d \times \mathbb{R}^n$ and $P(\mathbf{w}^*) = D(\boldsymbol{\alpha}^*)$ [72, 55]. Moreover, the primal and dual optimal solutions, $\mathbf{w}^*$ and $\boldsymbol{\alpha}^*$, respectively, are unique, and satisfy the relations $\mathbf{w}^* = \frac{1}{\lambda n}\mathbf{X}\boldsymbol{\alpha}^*$ and $\alpha_j^* = \phi_j'(\langle \mathbf{X}_{:j}, \mathbf{w}^* \rangle)$ for all $j \in [n] := \{1, \ldots, n\}$, which also uniquely characterize them.

### 5.2.1 Note on the setup

We choose the above setup, because linear ERM offers a good balance between the level of developed theory and practical interest. The coordinate descent methods for primal/dual linear ERM have been around for couple of years and there is no doubt that they are well suited for this task. Their convergence rates are well estabilshed and therefore we can confidently build upon them. We could also consider quadratic problems, where the bounds are known to be tighter, but the setup is less general and therefore of smaller importance to the machine learning community. For this reasons we believe that linear ERM is the most appropriate setup for the direct comparison between primal and dual approaches.

As we constrain our analysis to this setup, we do not claim any general conclusions about the advantage of one approach over the other. The results in this setup are only meant to offer us new insight into the comparison, which we believe is enlightening.

## 5.3 Primal and dual RCD

In its general "arbitrary sampling" form [57], RCD applied to the primal problem (5.1) has the form

$$w_i^{k+1} \leftarrow \begin{cases} w_i^k - \frac{1}{u_i'} \nabla_i P(\mathbf{w}^k) & \text{for } i \in S_k, \\ w_i^k & \text{for } i \notin S_k, \end{cases} \quad (5.4)$$

where $u_1', \ldots, u_d' > 0$ are parameters of the method and $\nabla_i P(\mathbf{w}) = \frac{1}{n}\sum_{j=1}^n \phi_j'(\langle \mathbf{X}_{:j}, \mathbf{w}\rangle)\mathbf{X}_{ij} + \lambda w_i$ is the $i$th partial derivative of $P$ at $\mathbf{w}$. This update is performed for a random subset of the coordinates $i \in S_k \subseteq [d]$ chosen in an i.i.d. fashion according to some sampling $\hat{S}_P$. The parameters $u_i'$ are usually computed ahead of the iterative process and need to be selected carefully in order for the method to work [57, 54]. A standard result is that one can set $u_i' := \frac{1}{\gamma n} u_i + \lambda$, where $\mathbf{u} = (u_1, \ldots, u_d)$ is chosen so as to satisfy the Expected Separable Overapproximation (ESO) inequality

$$\mathbf{P} \circ \mathbf{X}\mathbf{X}^\top \preceq \text{Diag}(\mathbf{p} \circ \mathbf{u}), \quad (5.5)$$

where $\mathbf{P}$ is the $d \times d$ matrix with entries $\mathbf{P}_{ij} = \mathbf{Prob}(i \in \hat{S}_P, j \in \hat{S}_P)$, $\mathbf{p} = \text{Diag}(\mathbf{P}) \in \mathbb{R}^d$ and $\circ$ denotes the Hadamard (element-wise) product of matrices. The resulting method is formally described as Algorithm 6. Note, that there are ways to run the method without precomputing $u_i'$ (e.g. [49]), but we will focus on the scenario where we compute them upfront, as this is the more standard way of developing the theory.

We will focus on applying serial coordinate descent to (5.1) and (5.3). For the case of generality we include them in the arbitrary sampling form.

---

**Algorithm 6** Primal RCD: NSync [57]

    **Input:** initial iterate $\mathbf{w}^0 \in \mathbb{R}^d$; sampling $\hat{S}_P$;
    ESO parameters $u_1, \ldots, u_d > 0$
    **Initialize:** $\mathbf{z}^0 = \mathbf{X}^\top \mathbf{w}^0$
    **for** $k = 0, 1, \ldots$ **do**
      Sample $S_k \subseteq [d]$ according to $\hat{S}_P$
      **for** $i \in S_k$ **do**
        $\Delta_i^k = -\frac{\gamma n}{u_i + n\lambda\gamma}\left(\frac{1}{n}\sum_{j=1}^n \phi_j'(z_j^k)X_{ij} + \lambda w_i^k\right)$
        Update $w_i^{k+1} = w_i^k + \Delta_i^k$
      **end for**
      **for** $i \notin S_k$ **do**
        $w_i^{k+1} = w_i^k$
      **end for**
      Update $\mathbf{z}^{k+1} = \mathbf{z}^k + \sum_{i \in S_k} \Delta_i^k \mathbf{X}_{i:}^\top$
    **end for**

---

When applying RCD to the dual problem (5.3), we can't proceed as above since the functions $\phi_j^*$ are not necessarily smooth, and hence we can't compute the partial derivatives of the dual objective. The standard approach here is to use a proximal variant of RCD [60]. In particular, Algorithm 7 has been analyzed in [55]. Like Algorithm 6, Algorithm 7 is also capable to work with an arbitrary sampling, which in this case is a random subset of $[n]$. The ESO parameters $\mathbf{v} = (v_1, \ldots, v_j)$ must in this case satisfy the ESO inequality

$$\mathbf{Q} \circ \mathbf{X}^\top \mathbf{X} \preceq \text{Diag}(\mathbf{q} \circ \mathbf{v}), \quad (5.6)$$

where $\mathbf{Q}$ is an $n \times n$ matrix with entries $Q_{ij} = \mathbf{Prob}(i \in \hat{S}_D, j \in \hat{S}_D)$ and $\mathbf{q} = \text{Diag}(\mathbf{Q}) \in \mathbb{R}^n$.

If we assume that $|\hat{S}_P| = 1$ (resp. $|\hat{S}_D| = 1$) with probability 1 (i.e., of the samplings are "serial"), then it is trivial to observe that (5.5) (resp. (5.6)) holds with

$$\mathbf{u} = \text{Diag}(\mathbf{X}\mathbf{X}^\top) \quad \text{and} \quad \mathbf{v} = \text{Diag}(\mathbf{X}^\top \mathbf{X}). \quad (5.7)$$

The proof of the above and other easily computable expressions for $\mathbf{u}$ and $\mathbf{v}$ for more complicated

**Algorithm 7** Dual RCD: Quartz [55]

**Input:** initial dual variables $\boldsymbol{\alpha}^0 \in \mathbb{R}^n$, sampling $\hat{S}_D$; ESO parameters $v_1, \ldots, v_n > 0$
**Initialize:** set $\mathbf{w}^0 = \frac{1}{\lambda n}\mathbf{X}\boldsymbol{\alpha}^0$
**for** $k = 0, 1, \ldots$ **do**
    Sample $S_k \subseteq [n]$ according to $\hat{S}_D$
    **for** $j \in S_k$ **do**
        $\Delta_j^k = \arg\max_{h \in \mathbb{R}}\{-\phi_j^*(-(\alpha_j + h)) - h\langle \mathbf{X}_{:j}, \mathbf{w}^k\rangle - \frac{v_j h^2}{2\lambda n}\}$
        Update $\alpha_j^{k+1} = \alpha_j^k + \Delta_j^k$
    **end for**
    **for** $j \notin S_k$ **do**
        $\alpha_j^{k+1} = \alpha_j^k$
    **end for**
    Update $\mathbf{w}^{k+1} = \mathbf{w}^k + \frac{1}{\lambda n}\sum_{j \in S_k} \Delta_j^k \mathbf{X}_{:j}$
**end for**

samplings can be found in [54].

### 5.3.1 Note on the methods

To understand the key differences in convergence properties of these two approaches, we analyse their behaviour in their most basic formulations. In practice, both methods can be extended in many different ways, including possibilities as: line-search, adaptive probabilities [10], local smoothness [80], and more. All these extensions offer empirical speed-up, but the theoretical speed-up cannot be quantified. At the same time all of them use additional computations, which in combination with the last point renders them uncomparable with standard version of RCD.

## 5.4 Iteration complexity and total arithmetic complexity

In this section we give expressions for the total expected arithmetic complexity of the two algorithms.

### 5.4.1 Number of iterations

Iteration complexity of Algorithms 6 and 7 is described in the following theorem. We do not claim novelty here, the results follow by applying theorems in [57] and [55] to the problems (5.1) and (5.3), respectively. We include a proof sketch in the appendix.

**Theorem 5.1.** *(Complexity: Primal vs Dual RCD) Let $\{\phi_j\}$ be convex and $1/\gamma$-smooth.*
*(i) Let*
$$K_P(\hat{S}_P, \epsilon) := \max_{i \in [d]} \left(\frac{u_i + n\lambda\gamma}{p_i n\lambda\gamma}\right) \log\left(\frac{c_P}{\epsilon}\right), \tag{5.8}$$

*where $c_P$ is a constant depending on $\mathbf{w}^0$ and $\mathbf{w}^*$. If $\hat{S}_P$ is proper (i.e., $p_i > 0$ for all i), and $\mathbf{u}$ satisfies (5.5), then iterates of primal RCD satisfy*
$$k \geq K_P(\hat{S}_P, \epsilon) \quad \Rightarrow \quad \mathbf{E}[P(\mathbf{w}^k) - P(\mathbf{w}^*)] \leq \epsilon.$$

*(ii) Let*
$$K_D(\hat{S}_D, \epsilon) := \max_{j \in [n]} \left(\frac{v_j + n\lambda\gamma}{q_j n\lambda\gamma}\right) \log\left(\frac{c_D}{\epsilon}\right), \tag{5.9}$$

*where $c_D$ is a constant depending on $\mathbf{w}^0$ and $\mathbf{w}^*$. If $\hat{S}_D$ is proper (i.e., $q_i > 0$ for all i), and $\mathbf{v}$* satisfies (5.6), then iterates of dual RCD satisfy

$$k \geq K_D(\hat{S}_D, \epsilon) \quad \Rightarrow \quad \mathbf{E}[P(\mathbf{w}^k) - P(\mathbf{w}^*)] \leq \epsilon.$$

The above results are the standard state-of-the-art bounds for primal and dual coordinate descent. From now on we will use the shorthand $K_P := K_P(\hat{S}_P, \epsilon)$ and $K_D := K_D(\hat{S}_D, \epsilon)$, when the quantity $\epsilon$ and the samplings $\hat{S}_P$ and $\hat{S}_D$ are clear from the context.

### 5.4.2 Average cost of a single iteration

Let $\|\cdot\|_0$ be the number of nonzeros in a matrix/vector. We can observe, that the computational cost associated with one iteration of Algorithm 6 is $\mathcal{O}(\|\mathbf{X}_{i:}\|_0)$ assuming that we picked the dimension $i$. As the dimension was picked randomly, we have to take the expectation over all the possible dimensions to get the average cost of an iteration. This leads us to the cost

$$\begin{aligned} W_P(\mathbf{X}, \hat{S}_P) &:= \mathcal{O}\left(\mathbf{E}\left[\sum_{i \in \hat{S}_P} \|\mathbf{X}_{i:}\|_0\right]\right) \\ &= \mathcal{O}\left(\sum_{i=1}^d p_i \|\mathbf{X}_{i:}\|_0\right), \end{aligned} \quad (5.10)$$

for Algorithm 6 and similarly for Algorithm 7 the average cost is

$$\begin{aligned} W_D(\mathbf{X}, \hat{S}_D) &:= \mathcal{O}\left(\mathbf{E}\left[\sum_{j \in \hat{S}_D} \|\mathbf{X}_{:j}\|_0\right]\right) \\ &= \mathcal{O}\left(\sum_{j=1}^n q_i \|\mathbf{X}_{:j}\|_0\right). \end{aligned} \quad (5.11)$$

From now on we will use the shorthand $W_P := W_P(\mathbf{X}, \hat{S}_P)$ and $W_D := W_D(\mathbf{X}, \hat{S}_D)$, when the matrix $\mathbf{X}$ and the samplings $\hat{S}_P$ and $\hat{S}_D$ are clear from the context.

We remark that the constant hidden in $\mathcal{O}$ may be larger for Algorithm 6 than for Algorithm 7. The reason for this is that for Algorithm 6 we compute the one-dimensional derivative $\phi'_j$ for every nonzero term in the sum, while for Algorithm 7 we do this only once. Depending on the loss $\phi_j$, this may lead to slower iterations. There is no difference if we use the squared loss as $\phi_j$. On the other hand, if $\phi_j$ is the logistic loss and we compute $\phi'_j$ directly, experimentation shows that the constant can be around 50. However, in practice this constant can be often completely diminished, for example by using a look-up table.

### 5.4.3 Total complexity

By combining the bounds on the number of iterations provided by Theorem 5.1 with the formulas (5.10) and (5.11) for the cost of a single iteration we obtain the following expressions for the *total complexity* of the two algorithms, where we ignore the logarithmic terms and drop the $\tilde{\mathcal{O}}$ symbol:

$$\begin{aligned} T_P(\mathbf{X}, \hat{S}_P) &:= K_P W_P \\ &\stackrel{(5.8)+(5.10)}{=} \left(\max_{i \in [d]} \frac{u_i + n\lambda\gamma}{p_i n\lambda\gamma}\right)\left(\sum_{i=1}^d p_i \|\mathbf{X}_{i:}\|_0\right), \end{aligned} \quad (5.12)$$

$$\begin{aligned} T_D(\mathbf{X}, \hat{S}_D) &:= K_D W_D \\ &\stackrel{(5.9)+(5.11)}{=} \left(\max_{j \in [n]} \frac{v_j + n\lambda\gamma}{q_j n\lambda\gamma}\right)\left(\sum_{j=1}^n q_j \|\mathbf{X}_{:j}\|_0\right). \end{aligned} \quad (5.13)$$

Again, from now on we will use the shorthand $T_P := WT_P(\mathbf{X}, \hat{S}_P)$ and $T_D := T_D(\mathbf{X}, \hat{S}_D)$, when the matrix $\mathbf{X}$ and the samplings $\hat{S}_P$ and $\hat{S}_D$ are clear from the context.

## 5.5 Choosing a sampling that minimizes the total Complexity

In this section we identify the *optimal sampling* in terms of the *total complexity*. This is different from previous results on *importance sampling*, which neglect to take into account the cost of the iterations [57, 55, 88, 46]. For simplicity, we shall only consider *serial* samplings, i.e., samplings which only pick a single coordinate at a time. The situation is much more complicated with non-serial samplings where first importance sampling results have only been derived recently [8].

### 5.5.1 Uniform sampling

The simplest serial sampling is the *uniform sampling*: it selects every coordinate with the same probability, i.e. $p_i = 1/d$, $\forall i \in [d]$ and $q_j = 1/n$, $\forall j \in [n]$. In view of (5.12), (5.13) and (5.7), we get

$$T_P = \|\mathbf{X}\|_0 \left(1 + \frac{1}{n\lambda\gamma} \max_{i \in [d]} \|\mathbf{X}_{i:}\|_2^2\right)$$

and

$$T_D = \|\mathbf{X}\|_0 \left(1 + \frac{1}{n\lambda\gamma} \max_{j \in [n]} \|\mathbf{X}_{:j}\|_2^2\right).$$

We can now clearly see that whether $T_P \leq T_D$ or $T_P \geq T_D$ does not simply depend on $d$ vs $n$, but instead depends on the relative value of the quantities $\max_{i \in [d]} \|\mathbf{X}_{i:}\|_2^2$ and $\max_{j \in [n]} \|\mathbf{X}_{:j}\|_2^2$. Having said that, we shall not study these quantities in this chapter. The reason for this is that for the sake of brevity, we shall instead focus on comparing the primal and dual RCD methods for optimal sampling which minimizes the total complexity, in which case we will obtain different quantities.

### 5.5.2 Importance sampling

By *importance sampling* we mean the serial sampling $\hat{S}_P$ (resp. $\hat{S}_D$) which minimizes the bounds $K_P$ in 5.8 (resp. $K_D$ in (5.9)). It can easily be seen (see also [57], [55], [88]), that importance sampling probabilities are given by

$$p_i^* = \frac{u_i + n\lambda\gamma}{\sum_l (u_l + n\lambda\gamma)} \quad \text{and} \quad q_j^* = \frac{v_j + n\lambda\gamma}{\sum_l (v_l + n\lambda\gamma)}. \tag{5.14}$$

On the other hand, one can observe that the average iteration cost of importance sampling may be larger than the average iteration cost of uniform serial sampling. Therefore, it is a natural question to ask, whether it is necessarily better. In view of (5.12), (5.13) and (5.14), the total complexities for importance sampling are

$$T_P = \|\mathbf{X}\|_0 + \frac{1}{n\lambda\gamma} \sum_{i=1}^{d} \|\mathbf{X}_{i:}\|_0 \|\mathbf{X}_{i:}\|_2^2 \tag{5.15}$$

$$T_D = \|\mathbf{X}\|_0 + \frac{1}{n\lambda\gamma} \sum_{j=1}^{n} \|\mathbf{X}_{:j}\|_0 \|\mathbf{X}_{:j}\|_2^2. \tag{5.16}$$

Since a weighted average is smaller than the maximum, the total complexity of both methods with importance sampling is always better than with uniform sampling. However, this does not mean that importance sampling is the sampling that minimizes total complexity.

### 5.5.3 Optimal sampling

The next theorem states that, in fact, importance sampling *does* minimize the total complexity.

**Theorem 5.2.** *The optimal serial sampling (i.e., the serial sampling minimizing the total expected complexity $T_P$ (resp, $T_D$)) is the importance sampling* (5.14).

## 5.6 The face-off

In this section we investigate the two quantities in (5.15) and (5.16), $T_P$ and $T_D$, measuring the total complexity of the two methods as functions of the data $\mathbf{X}$. Clearly, it is enough to focus on the quantities

$$C_P(\mathbf{X}) \quad := \quad \sum_{i=1}^{d} \|\mathbf{X}_{i:}\|_0 \|\mathbf{X}_{i:}\|^2 \tag{5.17}$$

$$C_D(\mathbf{X}) \quad := \quad \sum_{j=1}^{n} \|\mathbf{X}_{:j}\|_0 \|\mathbf{X}_{:j}\|^2. \tag{5.18}$$

We shall ask questions such as: when is $C_P(\mathbf{X})$ larger/smaller than $C_D(\mathbf{X})$, and by how much. In this regard, it is useful to note that $C_P(\mathbf{X}) = C_D(\mathbf{X}^\top)$. Our first result gives tight lower and upper bounds on their ratio.

**Theorem 5.3.** *For any $\mathbf{X} \in \mathbb{R}^{d \times n}$ with no zero rows or columns, we have the bounds $\|\mathbf{X}\|_F^2 \leq C_P(\mathbf{X}) \leq n\|\mathbf{X}\|_F^2$ and $\|\mathbf{X}\|_F^2 \leq C_D(\mathbf{X}) \leq d\|\mathbf{X}\|_F^2$. It follows that $1/d \leq C_P(\mathbf{X})/C_D(\mathbf{X}) \leq n$. Moreover, all these bounds are tight.*

Since $C_P(\mathbf{X})$ (resp. $C_D(\mathbf{X})$) can dominate the expression (5.12) (resp. (5.13)) for total complexity, it follows that, depending on the data matrix $\mathbf{X}$, *the primal method can be up to $d$ times faster than the dual method, and up to $n$ times slower than the dual method.*
Note, that for the above result to hold, we need to have the magnitudes of the individual entries in $\mathbf{X}$ potentially unbounded. However, this is not the case in practice. In the following sections we study more restricted classes of matrices, for which we are still able to claim some theoretical results.

### 5.6.1 Dense data

If $\mathbf{X}$ is a dense deterministic matrix ($X_{ij} \neq 0$ for all $i, j$), then $C_P(\mathbf{X}) = n\|\mathbf{X}\|_F^2$ and $C_D(\mathbf{X}) = d\|\mathbf{X}\|_F^2$, and we reach the conclusion that everything boils down to $d$ vs $n$.

### 5.6.2 Binary data

In Theorem 5.3 we showed, that without further constraints on the data we cannot say directly from $d$ and $n$, which of the approaches will perform better. The main argument in the proof of Theorem 5.3 is based on the possibility of arbitrary magnitudes of the individual data entries. In this part we go to the other extreme – we assume that all the magnitudes of the non-zero entries are the same. Surprisingly, using only the structure of the nonzero entries of the data we show, that one can have $C_P(\mathbf{X}) \leq C_D(\mathbf{X})$ even if $d \ll n$ (specifically $n \leq \frac{d^2}{4} - \frac{3}{2}d - 1$). On the other hand, we show that if the difference between $d$ and $n$ is too big (specifically $n > d^2 + 3d$), then $C_P(\mathbf{X}) > C_D(\mathbf{X})$. Naturally, we get the symmetric results as well. The class of binary matrices is only an example which we use to illustrate the phenomenon, one can easily come up with similar bounds for other classes of matrices as well.

Let $\mathbb{B}^{d \times n}$ denote the set of $d \times n$ matrices $\mathbf{X}$ with (signed) binary elements, i.e., with $X_{ij} \in \{-1, 0, 1\}$ for all $i, j$. Note, that the following results trivially hold also for entries in $\{-a, 0, a\}$, for any $a \neq 0$. For $\mathbf{X} \in \mathbb{B}^{d \times n}$, the expressions in (5.17) and (5.18) can be also written in the form $C_P(\mathbf{X}) = \sum_{i=1}^{d} \|\mathbf{X}_{i:}\|_0^2$ and $C_D(\mathbf{X}) = \sum_{j=1}^{n} \|\mathbf{X}_{:j}\|_0^2$. By $\mathbb{B}_{\neq 0}^{d \times n}$ we denote the set of all matrices in $\mathbb{B}^{d \times n}$ with nonzero columns and rows.

For positive integers $a, b$ we write $\bar{a}_b := b \lfloor \frac{a}{b} \rfloor$ (i.e., $a$ rounded down to the closest multiple of $b$). Further, we write
$$R(\alpha, d, n) \quad := \quad U(\alpha, d, n)/L(\alpha, n),$$
where
$$L(\alpha, n) \quad := \quad \frac{1}{n}(\bar{\alpha}_n^2 + (\alpha - \bar{\alpha}_n)(2\bar{\alpha}_n + n))$$
and
$$\begin{aligned} U(\alpha, d, n) \quad := \quad & (n+1)\overline{(\alpha - d)}_{n-1} + d - 1 \\ & + [\alpha - d + 1 - \overline{(\alpha - d)}_{n-1}]^2. \end{aligned}$$

The following is a refinement of Theorem 5.3 for binary matrices of fixed cardinality $\alpha$.

**Theorem 5.4.** *For all $\mathbf{X} \in \mathbb{B}_{\neq 0}^{d \times n}$ with $\alpha = \|\mathbf{X}\|_0$ we have the bounds*
$$\frac{1}{R(\alpha, n, d)} \quad \leq \quad \frac{C_P(\mathbf{X})}{C_D(\mathbf{X})} \quad \leq \quad R(\alpha, d, n).$$
*Moreover, these bounds are tight.*

The above theorem follows from Lemma 5.10, which we formulate and prove in the Appendix. This lemma establishes formulas for the minimum and maximum of $C_D$ and $C_P$, subject to the constraint $\|\mathbf{X}\|_0 = \alpha$, in terms of the functions $L$ and $U$. Using the same Lemma 5.10 we can further show the following theorem:

**Theorem 5.5.** *Let $d \leq n \leq \frac{d^2}{4} - \frac{3}{2}d - 1$. Then there exists a matrix $\mathbf{X} \in \mathbb{B}_{\neq 0}^{d \times n}$ such that $C_P(\mathbf{X}) < C_D(\mathbf{X})$. Symetrically, if $n \leq d \leq \frac{n^2}{4} - \frac{3}{2}n - 1$ then there exists a matrix $\mathbf{X} \in \mathbb{B}_{\neq 0}^{d \times n}$ such that $C_D(\mathbf{X}) < C_P(\mathbf{X})$.*

The above theorem shows, that even if $n = \mathcal{O}(d^2)$, the primal method can be better than the dual method – and vice-versa.

Further, as we show in Lemma 5.11 in the Appendix, if $d \geq n$ and $\alpha \geq n^2 + 3n$, then $R(\alpha, d, n) \leq 1$. Likewise, if $n \geq d$ and $\alpha \geq d^2 + 3d$, then $R(\alpha, n, d) \leq 1$. Combined with Theorem 5.4, this has an interesting consequence, spelled out in the next theorem and its corollary.

**Theorem 5.6.** *Let $\mathbf{X} \in \mathbb{B}_{\neq 0}^{d \times n}$. If $d \geq n$ and $\|\mathbf{X}\|_0 \geq n^2 + 3n$, then $C_P(\mathbf{X}) \leq C_D(\mathbf{X})$. By symmetry, if $n \geq d$ and $\|\mathbf{X}\|_0 \geq d^2 + 3d$, then $C_D(\mathbf{X}) \leq C_P(\mathbf{X})$.*

This result says that for binary data, and $d \geq n$, the primal method is better than the dual method even for non-dense data, as long as the the data is "dense enough". Observe that as long as $d \geq n^2 + 3n$, all matrices $\mathbf{X} \in \mathbb{B}_{\neq 0}^{d \times n}$ satisfy $\|\mathbf{X}\|_0 \geq d \geq n^2 + 3n \geq n$. This leads to the following corollary.

**Corollary 5.7.** *If $d \geq n^2 + 3n$, then for all $\mathbf{X} \in \mathbb{B}_{\neq 0}^{d \times n}$ we have $C_P(\mathbf{X}) \leq C_D(\mathbf{X})$. By symmetry, if $n \geq d^2 + 3d$, then for all $\mathbf{X} \in \mathbb{B}_{\neq 0}^{d \times n}$ we have $C_D(\mathbf{X}) \leq C_P(\mathbf{X})$.*

In words, the corollary states that for binary data where the number of features ($d$) is large enough in comparison with the number of examples ($n$), the primal method will be always better. On the other hand, if $n$ is large enough, the dual method will be always better. This behavior can be observed in Figure 5.1. For large enough $d$, all the values $R(\alpha, d, n)$ are below 1, and therefore the primal method is always better than the dual in this regime.

We note, that the gap in the results of Theorem 5.5 and Theorem 5.6 arises from using the bounds $a - b \leq \bar{a}_b := b \lfloor \frac{a}{b} \rfloor \leq a$, which can get very loose. A tighter bounds can be achieved, but as they have mostly theoretical purpose, we do not feel they are needed. The bounds well convey the main message of the chapter as they are.

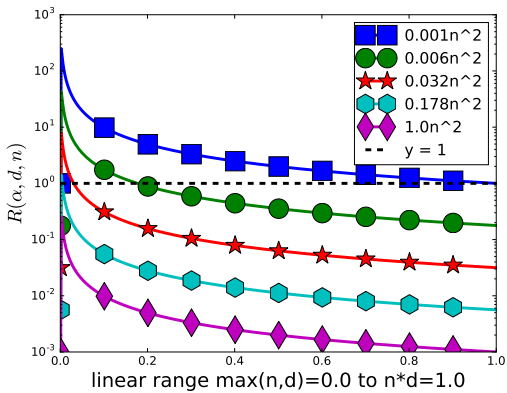

Figure 5.1: The value $R(\alpha, d, n)$ plotted for $n = 10^3$, $n \leq d \leq n^2$ and $\max\{d, n\} \leq \alpha \leq nd$.

Table 5.1: Details on the datasets used in the experiments

| dataset | $d$ | $n$ | density | $\|\mathbf{X}\|_0$ | $C_P$ | $C_D$ | $T_P/T_D$ |
|---|---:|---:|---:|---:|---:|---:|---:|
| news | 1,355,191 | 19,996 | 0.03% | 9,097,916 | $3 \times 10^7$ | $9 \times 10^6$ | 2.0 |
| leukemia | 7,129 | 38 | 100.00% | 270,902 | $1 \times 10^7$ | $2 \times 10^9$ | 0.5 |

## 5.7 Experiments

We conducted experiments on both real and synthetic data. The problem we were interested in is a standard logistic regression with an L2-regularizer, i.e.,

$$P(\mathbf{w}) = \frac{1}{n} \sum_{j=1}^{n} \log(1 + \exp(-y_j \langle \mathbf{X}_{:j}, \mathbf{w} \rangle)) + \frac{\lambda}{2} \|\mathbf{w}\|_2^2.$$

In all our experiments we used $\lambda = 1/n$ and we normalized all the entries of $\mathbf{X}$ by the average column norm. Note that for logistic loss there is no closed form solution for $\Delta_j^k$ in Algorithm 7. Therefore we use a variant of Algorithm 7 where $\Delta_j^k = \eta(\phi_j'(\langle \mathbf{X}_{:j}, \mathbf{w} \rangle) + \alpha_j^k)$ with the step size $\eta$ defined as $\eta = \min_{j \in [n]}(q_j n \lambda \gamma)/(v_j + n \lambda \gamma)$. This variant has the same convergence rate guarantees as Algorithm 7 and does not require exact minimization. Details can be found in [55].

We plot the training error against the number of passes through the data. The number of passes is calculated according to the number of visited nonzero entries in the matrix $\mathbf{X}$. One pass means that we look at $\|\mathbf{X}\|_0$ nonzero entries of $\mathbf{X}$, but not necessarily all of them, as we can visit some of them multiple times.

We showcase the conclusions from the theory on two real datasets and multiple synthetic datasets. We constructed all the synthetic experiments in a way, that according to the theory the primal approach should be better . We note, that the same plots could be generated symmetrically for the dual approach.

### 5.7.1 General data

We look at the matrices which give the worst-case bounds for general matrices (Theorem 5.3) and their empirical properties for different choices of $d$ and $n$. These matrices have highly non-uniform distribution of the nonzeros and moreover require the entries to have their magnitudes differ by many orders (see the proof of Theorem 5.3). We performed 2 experiments, where we showed the potential empirical speedup for the primal method for $d = n$ and also for $d \ll n$ (which is highly unfavourable for the primal method). The corresponding figures are Figure 5.2a and 5.2b. For a square dataset, we can clearly observe a large speed-up. For $d \ll n$ we can observe, that the theory holds and the primal method is still faster, but because of numerical issues (as mentioned, the magnitudes of the entries differ by many orders) and the fact that the optimal value is very close to an "initial guess" of the algorithm, the difference in speed is more difficult to observe.

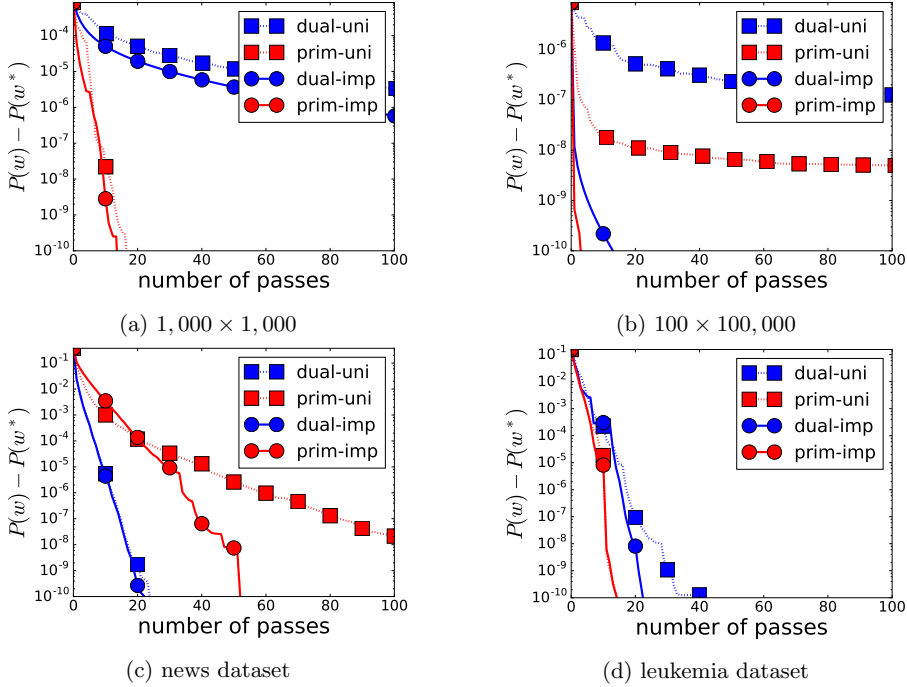

Figure 5.2: Testing the worst case for general matrices and real datasets

### 5.7.2 Synthetic binary data

We looked at matrices with all entries in $\{a, -a, 0\}$ for some $a \neq 0$. We fixed the number of features to be $d = 100$ and we varied the number of examples $n$ and the sparsity level $\alpha = \|\mathbf{X}\|_0$. For each triplet $[d, n, \alpha]$ we produced the worst-case matrix for dual RCD according to the developed theory formalized in Theorem 5.4 (see the proof for more details on the matrix structure). The results are in Figure 5.3.

Each row corresponds to one fixed value of $n$, while each column corresponds to one sparsity level $\alpha$ given by the proportion of nonzero entries, e.g., nnz $\sim 1\%$ stands for $\alpha \sim 0.01 \cdot nd$.

In the experiments we can observe the behaviour described in Theorem 5.6. While $n$ is comparable to $d$, the primal method outperforms the dual method. When the sparsity level $\alpha$ reaches values of $\sim d^2$, the dual method outperforms the primal although the matrix structure is much better suited for the primal method. Also note, that the right column corresponds to dense matrices, where larger $n$ is the only dominant factor.

### 5.7.3 Real data

We used two real datasets to showcase our theory: news and leukemia[2]. The news dataset in Figure 5.2c is a nice example of our theory in practice. As shown in Table 5.1 we have $d \gg n$, but the dual method is empirically faster than the primal one. The reason is simple: the news dataset uses a bag of words representation of news articles. If we look at the distribution of features (words), there are many words which appear just very rarely and there are words commonly used in many articles. The features have therefore a very skewed distribution of their nonzero entries. On the other hand, the examples have close to a uniform distribution, as the number of distinct words in an article usually does not take on extreme values. As shown in the theory (proof of Theorem 5.4), this distribution of nonzero entries highly favors the dual approach.

---
[2]both datasets are available from https://www.csie.ntu.edu.tw/ cjlin/libsvmtools/datasets/

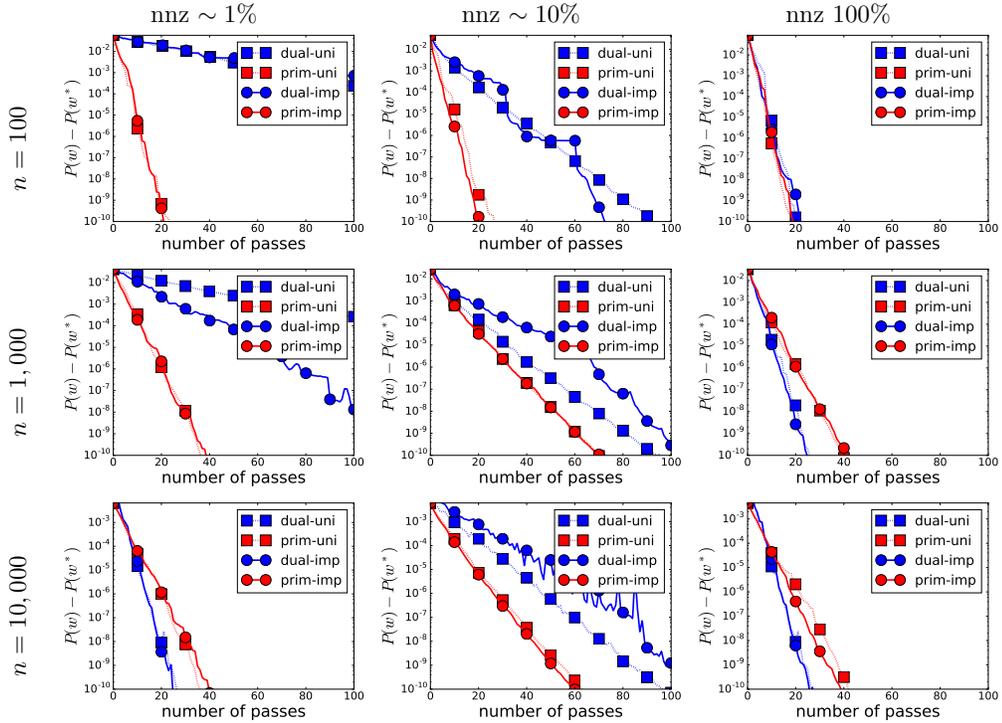

Figure 5.3: Worst-case experiments with various dimensions and sparsity levels for $d = 100$

The leukemia dataset in Figure 5.2d is a fully dense dataset and $d \gg n$. Therefore, as our theoretical analysis shows, the primal approach should be better. The ratio between the runtimes is not very large, as the constant $\|\mathbf{X}\|_0$ is of similar order as the additional term in the computation of the true runtime (recall (5.12), (5.13)). The empirical speedup in Figures 5.2c and Figures 5.2d matches the theoretical predictions in the last column of Table 5.1.

## 5.8 Conclusions and extensions

We have shown that the question whether RCD should be applied to the primal or the dual problem depends on the structure of the training dataset. For dense data, this simply boils down to whether we have more data or parameters, which is intuitively appealing. We have shown, both theoretically, and through experiments with synthetic and real datasets, that contrary to what seems to be a popular belief, primal RCD can outperform dual RCD even if $n \gg d$ and vice-versa. If a user is willing to invest one pass over the data, we recommend to compare the quantities $T_P$ and $T_D$ (or $C_P$ and $C_D$) to figure out which approach has faster convergence according to the theory.

In order to focus on the main message, we have chosen to present our results for simple (as opposed to "accelerated") variants of RCD. However, our results can be naturally extended to accelerated variants of RCD, such as APPROX [15], ASDCA [71], APCG [37], ALPHA [53] and SPDC [86].

## 5.A Proof of Theorem 5.1

We say that $P \in \mathcal{C}^1(\mathbf{M})$, if

$$P(\mathbf{w}+\mathbf{h}) \leq P(\mathbf{w}) + \langle \nabla P(\mathbf{w}), \mathbf{h}\rangle + \frac{1}{2}\mathbf{h}^\top \mathbf{M}\mathbf{h}, \quad \forall \mathbf{w}, \mathbf{h} \in \mathbb{R}^d.$$

For three vectors $\mathbf{a}, \mathbf{b}, \mathbf{c} \in \mathbb{R}^n$ we define $\langle \mathbf{a}, \mathbf{b}\rangle_\mathbf{c} := \sum_{i=1}^d a_i b_i c_i$ and $\|\mathbf{a}\|_\mathbf{c}^2 := \langle \mathbf{a}, \mathbf{a}\rangle_\mathbf{c} = \sum_{i=1}^d c_i a_i^2$. Also, let for $\emptyset \neq S \subseteq [d]$ and $\mathbf{h} \in \mathbb{R}^d$, we write $h_S := \sum_{i \in S} h_i \mathbf{e}_i$, where $\mathbf{e}_i$ is the $i$-th coordinate vector (i.e., standard basis vector) in $\mathbb{R}^d$.

We will need the following two lemmas.

**Lemma 5.8.** *The primal objective $P$ satisfies $P \in \mathcal{C}^1(\mathbf{M})$, where $\mathbf{M} = \lambda \mathbf{I} + \frac{1}{n\gamma}\mathbf{X}\mathbf{X}^\top$.*

*Proof.*

$$\begin{aligned}
P(\mathbf{w}+\mathbf{h}) &\stackrel{(5.1)}{=} \frac{1}{n}\sum_{j=1}^n \phi_j(\langle \mathbf{X}_{:j}, \mathbf{w}\rangle + \langle \mathbf{X}_{:j}, \mathbf{h}\rangle) + \frac{\lambda}{2}\|\mathbf{w}+h\|^2 \\
&\stackrel{(5.2)}{\leq} \frac{1}{n}\sum_{j=1}^n \left[\phi_j(\langle \mathbf{X}_{:j}, \mathbf{w}\rangle) + \phi_j'(\langle \mathbf{X}_{:j}, \mathbf{w}\rangle) \cdot \langle \mathbf{X}_{:j}, \mathbf{h}\rangle + \frac{1}{2\gamma}\langle \mathbf{X}_{:j}, \mathbf{h}\rangle^2\right] + \frac{\lambda}{2}\|\mathbf{w}\|^2 \\
&\quad + \lambda \langle \mathbf{w}, \mathbf{h}\rangle + \frac{\lambda}{2}\|\mathbf{h}\|^2 \\
&= \frac{1}{n}\sum_{j=1}^n \phi_j(\langle \mathbf{X}_{:j}, \mathbf{w}\rangle) + \frac{\lambda}{2}\|\mathbf{w}\|^2 + \left\langle \frac{1}{n}\sum_{j=1}^n \phi_j'(\langle \mathbf{X}_{:j}, \mathbf{w}\rangle)\mathbf{X}_{:j} + \lambda \mathbf{w}, \mathbf{h}\right\rangle \\
&\quad + \frac{1}{2}\mathbf{h}^\top \left(\frac{1}{n\gamma}\sum_{j=1}^n \mathbf{X}_{:j}(\mathbf{X}_{:j})^\top + \lambda \mathbf{I}\right)\mathbf{h} \\
&= P(\mathbf{w}) + \langle \nabla P(\mathbf{w}), \mathbf{h}\rangle + \frac{1}{2}\mathbf{h}^\top \mathbf{M}\mathbf{h}.
\end{aligned}$$

□

**Lemma 5.9.** *If $P \in \mathcal{C}^1(\mathbf{M})$ and $\mathbf{u}' \in \mathbb{R}^d$ is such that $\mathbf{P} \circ \mathbf{M} \preceq \mathrm{Diag}(\mathbf{p} \circ \mathbf{u}')$, then*

$$\mathbf{E}[P(\mathbf{w}+\mathbf{h}_{[\hat{S}_P]})] \leq P(\mathbf{w}) + \langle \nabla P(w), \mathbf{h}\rangle_\mathbf{p} + \frac{1}{2}\|\mathbf{h}\|_{\mathbf{p}\circ\mathbf{u}'}^2.$$

*Proof.* See [54], Section 3. □

We can now proceed to the proof of Theorem 5.1.
First, note that

$$\mathbf{P}\circ\mathbf{M} = \lambda\mathrm{Diag}(\mathbf{p}) + \frac{1}{n\gamma}(\mathbf{P}\circ\mathbf{X}\mathbf{X}^\top) \preceq \lambda\mathrm{Diag}(\mathbf{p}) + \frac{1}{n\gamma}\mathrm{Diag}(\mathbf{p}\circ\mathbf{u})$$

with $\mathbf{u}$ defined as in (5.5). We now separately establish the two complexity results; (i) for primal RCD and (ii) for dual RCD.

(i) The proof is a consequence of the proof of the main theorem of [57]. Assumption 1 from [57] holds with $w_i := \lambda + \frac{1}{n\gamma}u_i$ (Lemma 5.8 & Lemma 5.9) and Assumption 2 from [57] holds with standard Euclidean norm and $\gamma := \lambda$. We follow the proof all the way to the bound

$$\mathbf{E}[P(\mathbf{w}^k) - P(\mathbf{w}^*)] \leq (1-\mu)^k(P(\mathbf{w}^0) - P(\mathbf{w}^*))$$

which holds for $\mu$ defined by

$$\mu := \frac{\lambda}{\max_i \frac{u_i + n\lambda\gamma}{p_i n\gamma}}$$

by direct substitution of the quantities. The result follows by standard arguments. Note that $C_P = P(\mathbf{w}^0) - P(\mathbf{w}^*)$.

(ii) The proof is a direct consequence of the proof of the main theorem of [55], using the fact that $P(\mathbf{w}^k) - P(\mathbf{w}^*) \leq P(\mathbf{w}^k) - D(\boldsymbol{\alpha}^k)$, as the weak duality holds. Note that $C_D = P(\mathbf{w}^0) - D(\boldsymbol{\alpha}^0)$.

## 5.B  Proof of Theorem 5.2

The proofs for Algorithm 6 and Algorithm 7 are analogous, and hence we will establish the result for Algorithm 6 only. For brevity, denote $s_i = u_i + n\lambda\gamma$. We aim to solve the optimization problem:

$$\mathbf{p}^* \leftarrow \underset{\mathbf{p} \in \mathbb{R}_+^d \,:\, \sum_i p_i = 1}{\arg\min} \quad T_P \stackrel{(5.12)}{=} \left( \max_{i \in [d]} \frac{s_i}{p_i n\lambda\gamma} \right) \cdot \sum_{i=1}^d p_i \|\mathbf{X}_{i:}\|_0. \tag{5.19}$$

First observe, that the problem is homogeneous in $\mathbf{p}$, i.e., if $\mathbf{p}$ is optimal, also $c\mathbf{p}$ will be optimal for $c > 0$, as the solution will be the same. Using this argument, we can remove the constraint $\sum_i p_i = 1$. Also, we can remove the multiplicative factor $1/(n\lambda\gamma)$ from the denominator as it does not change the arg min. Hence we get the simpler problem

$$\mathbf{p}^* \leftarrow \underset{\mathbf{p} \in \mathbb{R}_+^d}{\arg\min} \quad \left[ \left( \max_{i \in [d]} \frac{s_i}{p_i} \right) \cdot \sum_{i=1}^d p_i \|\mathbf{X}_{i:}\|_0 \right]. \tag{5.20}$$

Now choose optimal $\mathbf{p}$ and assume that there exist $j, k$ such that $s_j/p_j < s_k/p_k$. By a small decrease in $p_j$, we will still have $s_j/p_j \leq s_k/p_k$, and hence the term $\max_i s_i/p_i$ stays unchanged. However, the term $\sum_i p_i \|\mathbf{X}_{i:}\|_0$ decreased. This means that the optimal sampling must satisfy $s_i/p_i = const$ for all $i$. However, this is precisely the importance sampling.

## 5.C  Proof of Theorem 5.3

By assumption, all rows and columns of $\mathbf{X}$ are nonzero. Therefore, $1 \leq \|\mathbf{X}_{i:}\|_0 \leq n$ and $1 \leq \|\mathbf{X}_{:j}\|_0 \leq d$, and the bounds on $C_P$ and $C_D$ follow by applying this to (5.17) and (5.18) respectively. The bounds for the ratio follow immediately by combining the previous bounds. It remains to establish tightness. For $a, b, c \in \mathbb{R}$, let $\mathbf{X}(a, b, c) \in \mathbb{R}^{d \times n}$ be the matrix defined as follows:

$$X_{ij}(a, b, c) = \begin{cases} a & i \neq 1 \wedge j = 1 \\ b & i = 1 \wedge j \neq 1 \\ c & i = 1 \wedge j = 1 \\ 0 & \text{otherwise.} \end{cases}$$

Notice that $\mathbf{X}(a, b, c)$ does not have any zero rows nor columns as long as $a, b, c$ are nonzero. Since $C_P(\mathbf{X}(a, b, c)) = (d-1)a^2 + n(n-1)b^2 + nc^2$ and $C_D(\mathbf{X}(a, b, c)) = d(d-1)a^2 + (n-1)b^2 + dc^2$, one readily sees that

$$\lim_{\substack{b \to 0 \\ c \to 0}} \frac{C_P(\mathbf{X}(a,b,c))}{C_D(\mathbf{X}(a,b,c))} = \frac{1}{d} \quad \text{and} \quad \lim_{\substack{a \to 0 \\ c \to 0}} \frac{C_P(\mathbf{X}(a,b,c))}{C_D(\mathbf{X}(a,b,c))} = n.$$

## 5.D  Proof of Theorem 5.4

We first need a lemma.

**Lemma 5.10.** *Let $\alpha$ be an integer satisfying $\max\{d, n\} \leq \alpha \leq dn$ and let $L$ and $U$ be the functions defined in Section 5.6.2. We have the following identities:*

$$L(\alpha, n) = \min_{\mathbf{X} \in \mathbb{B}_{\neq 0}^{d \times n}} \{C_D(\mathbf{X}) \;:\; \|\mathbf{X}\|_0 = \alpha\} \tag{5.21}$$

$$L(\alpha, d) = \min_{\mathbf{X} \in \mathbb{B}_{\neq 0}^{d \times n}} \{C_P(\mathbf{X}) \;:\; \|\mathbf{X}\|_0 = \alpha\} \tag{5.22}$$

$$U(\alpha, n, d) = \max_{\mathbf{X} \in \mathbb{B}_{\neq 0}^{d \times n}} \{C_D(\mathbf{X}) \;:\; \|\mathbf{X}\|_0 = \alpha\} \tag{5.23}$$

$$U(\alpha, d, n) = \max_{\mathbf{X} \in \mathbb{B}_{\neq 0}^{d \times n}} \{C_P(\mathbf{X}) \;:\; \|\mathbf{X}\|_0 = \alpha\}. \tag{5.24}$$

*Proof.* Let $\mathbf{X} \in \mathbb{B}_{\neq 0}^{d \times n}$ be an arbitrary matrix and let $\boldsymbol{\omega} = (\omega_1, \ldots, \omega_n)$, where $\omega_j := \|\mathbf{X}_{:j}\|_0$. Let $\alpha = \|\mathbf{X}\|_0 = \sum_j \omega_j$. Observe that $C_D(\mathbf{X}) = \sum_{j=1}^n \|\mathbf{X}_{:j}\|_0^2 = \|\boldsymbol{\omega}\|_2^2$.

(i) We shall first establish (5.21). Assume that there exist two columns $j, k$ of $\mathbf{X}$, such that $\omega_j + 2 \leq \omega_k$, i.e., their difference in the number of nonzeros is at least 2. Because $\omega_k > \omega_j$, there has to exist a row which has a nonzero entry in the $k$-th column and a zero entry in the $j$-th column. Let $\mathbf{X}'$ be the matrix obtained from $\mathbf{X}$ by switching these two entries. Note that $C_P(\mathbf{X}) = C_P(\mathbf{X}')$. However, we have

$$C_D(\mathbf{X}) - C_D(\mathbf{X}') \;=\; \omega_j^2 + \omega_k^2 - (\omega_j + 1)^2 - (\omega_k - 1)^2 \;=\; 2\omega_k - 2\omega_j - 2 > 0.$$

It follows that while there exist two such columns, the minimum is not achieved. So, we only need to consider matrices $\mathbf{X}$ for which there exists integer $a$ such that $\omega_j = a$ or $\omega_j = a + 1$ for every $j$. Let $b = |\{j \;:\; \omega_j = a\}|$.

We can now without loss of generality assume that $0 \leq b \leq n - 1$. Indeed, we can do this is because the choices $b = 0$ and $b = n$ lead to the same matrices, and hence by focusing on $b = 0$ we have not removed any matrices from consideration. With simple calculations we get

$$\alpha \;=\; ba + (n - b)(a + 1) \;=\; n(a + 1) - b.$$

Note that $\alpha + b$ is a multiple of $n$. It follows that $b = n - \alpha + \bar{\alpha}_n$ and $a = \bar{\alpha}_n/n$. Up to the ordering of the columns (which does not affect $C_D(\mathbf{X})$) we have just one candidate $\mathbf{X}$, therefore it has to be the minimizer of $C_D$. Finally, we can easily calculate the minimum as

$$\sum_{j=1}^n \omega_j^2 \;=\; ba^2 + (n - b)(a + 1)^2 \;=\; (n - \alpha + \bar{\alpha}_n)\left(\frac{\bar{\alpha}_n}{n}\right)^2 + (\alpha - \bar{\alpha}_n)\left(\frac{\bar{\alpha}_n}{n} + 1\right)^2$$

$$=\; \frac{1}{n}\left(\bar{\alpha}_n^2 + (\alpha - \bar{\alpha}_n)(2\bar{\alpha}_n + n)\right) \;=\; L(\alpha, n).$$

(ii) Claim (5.22) follows from part (5.21) via symmetry: $C_P(\mathbf{X}) = C_D(\mathbf{X}^\top)$ and $\|\mathbf{X}\|_0 = \|\mathbf{X}^\top\|_0$.

(iii) We now establish claim (5.23). Assume that there exist a pair of columns $j, k$ such that $1 < \omega_j \leq \omega_k < d$. Let $\mathbf{X}'$ be the matrix obtained from $\mathbf{X}$ by zeroing out an entry in the $j$-th column and putting a nonzero inside the $k$-th column. Then

$$C_D(\mathbf{X}') - C_D(\mathbf{X}) \;=\; (\omega_j - 1)^2 + (\omega_k + 1)^2 - \omega_j^2 - \omega_k^2 \;=\; 2\omega_k - 2\omega_j + 2 > 0.$$

It follows that while there exist such a pair of columns, the maximum is not achieved. This condition leaves us with matrices $\mathbf{X}$ where at most one column $j$ has $\omega_j$ *not* equal to 1 or $d$.

Formally, let $a = |\{j \;:\; \omega_j = d\}|$. Then we have $n - a - 1$ columns with 1 nonzero and 1 column with $b$ nonzeros, where $1 \leq b < d$. This is correct, as $b = d$ is the same as $b = 1$ and $a$ being one more. We can compute $a$ and $b$ from the equation

$$(n - a - 1) \cdot 1 + 1 \cdot b + a \cdot d = \alpha$$
$$b + a(d - 1) = \alpha - n + 1$$

as the only solution to the division with remainder of $\alpha - n + 1$ by $d - 1$, with the difference that $b \in \{1, \ldots, d - 1\}$ instead of the standard $\{0, \ldots, d - 2\}$. We get

$$a = \left\lfloor \frac{a - n}{d - 1} \right\rfloor \quad \text{and} \quad b = \alpha - n + 1 - \overline{(\alpha - n)}_{d-1}.$$

The maximum can now be easily computed as follows:

$$\sum_{j=1}^{n} \omega_j^2 = (n - a - 1) + b^2 + ad^2$$
$$= n - \left\lfloor \frac{a - n}{d - 1} \right\rfloor - 1 + \left(\alpha - n + 1 - \overline{(\alpha - n)}_{d-1}\right)^2 + \left\lfloor \frac{a - n}{d - 1} \right\rfloor d^2$$
$$= U(\alpha, n, d).$$

(iv) Again, claim (5.24) follows from (5.23) via symmetry.

$\square$

We can now proceed to the proof of the theorem.

The quantity is the ratio between the maximal value of $C_P$ and the minimal value of $C_D$, we have to show that there exists a matrix $\mathbf{X}$ such that this is achieved. Assume we have a matrix $\mathbf{X}$ which has the maximal $C_P$. In the proof of Lemma 5.10 we showed, that by switching entries in $\mathbf{X}$ we can get the minimal value of $C_D$ without changing $C_P$. Therefore we can achieve maximal $C_P$ and minimal $C_D$ at the same time. Analogically for the other case.

## 5.E  Proof of Theorem 5.5

Let us establish the proof for the first case – the second will follow by symetricity. The main idea of the proof is to choose $\alpha$ such that the best possible case for the primal method is going to be better than the worst possible case for the dual method, i.e., $L(\alpha, d) < U(\alpha, n, d)$. First, we get some bounds on $L$ and $U$ based on the trivial bounds $a - b \leq \bar{a}_b := b\lfloor \frac{a}{b} \rfloor \leq a$ for positive integers $a, b$. First, the bound on $L$:

$$
\begin{aligned}
L(\alpha, d) &= \frac{1}{d}\left[\bar{\alpha}_d^2 + (\alpha - \bar{\alpha}_d)(2\bar{\alpha}_d + d)\right] \\
&\leq \frac{1}{d}\left[\alpha^2 + (\alpha - (\alpha - d))(2\alpha + d)\right] \\
&= \frac{1}{d}\left[\alpha^2 + 2\alpha d + d^2\right]
\end{aligned}
$$

Second, the bound on $U$:

$$
\begin{aligned}
U(\alpha, n, d) &= (d+1)\overline{(\alpha - n)}_{d-1} + n - 1 + \left[\alpha - n + 1 - \overline{(\alpha - n)}_{d-1}\right]^2 \\
&\geq (d+1)(\alpha - n - d + 1) + n - 1 + [\alpha - n + 1 - \alpha + n]^2 \\
&= d\alpha - dn - d^2 + \alpha + 1.
\end{aligned}
$$

Therefore, it is sufficient to show

$$\alpha^2 + 2\alpha d + d^2 \quad < \quad d^2\alpha - d^2 n - d^3 + d\alpha + d,$$

We will show that this holds for $\alpha = \frac{d(d-1)}{2}$. Observe, that $d(d-1)$ is even and therefore $\alpha$ is an integer. From the definition of $\mathbb{B}_{\neq 0}^{d \times n}$ we know, that $n + d - 1 \leq \alpha \leq nd$ has to hold. To verify this, we use the conditions from the assumptions:

$$n + d - 1 \leq \left(\frac{d^2}{4} - \frac{3}{2}d - 1\right) + d - 1 = \frac{d(d-1)}{2} - \frac{d^2}{4} - 2 < \frac{d(d-1)}{2} = \alpha$$

and the other inequality

$$\alpha = \frac{d(d-1)}{2} < d^2 \leq nd.$$

Using the above defined $\alpha$ we can proceed to prove the claim:

$$
\begin{aligned}
\alpha^2 + 2\alpha d + d^2 &= \frac{d^4}{4} + \frac{d^3}{2} + \frac{d^2}{4} \\
&< \frac{d^4}{4} + \frac{d^3}{2} + d^2 + d \\
&= \frac{d^4}{2} - d^3 + d - d^2\left(\frac{d^2}{2} - \frac{3}{2}d - 1\right) \\
&\leq \frac{d^4}{2} - d^3 + d - d^2 n \\
&= d^2\alpha - d^2 n - d^3 + d\alpha + d,
\end{aligned}
$$

which finishes the proof.

## 5.F  Proof of Theorem 5.6

As shown in the main text, the theorem follows from the following lemma. Hence, we only need to prove the lemma.

**Lemma 5.11.** *If $d \geq n$ and $\alpha \geq n^2 + 3n$, then $R(\alpha, d, n) \leq 1$. If $n \geq d$ and $\alpha \geq d^2 + 3d$, then $R(\alpha, n, d) \leq 1$.*

*Proof.* We focus on the first part, the second follows in an analogous way. Using the two assumptions, we have $\alpha(n^2 + 3n) + n^3 \leq \alpha^2 + dn^2$. By adding $n^2 + n$ to the right hand side and after reshuffling, we obtain the inequality

$$n\left[(n+1)(\alpha - d) + d - 1 + n^2\right] \quad \leq \quad (\alpha - n)^2.$$

For positive integers $a, b$, we have the trivial estimates $a - b \leq \bar{a}_b := b\lfloor \frac{a}{b} \rfloor \leq a$. We use them to bound four expressions:

$$\begin{aligned}
(\alpha - d) &\geq \overline{(\alpha - d)}_{n-1} \\
n^2 &\geq (\alpha - d + 1 - \overline{(\alpha - d)}_{n-1})^2 \\
\bar{\alpha}_n^2 &\geq (\alpha - n)^2 \\
(\alpha - \bar{\alpha}_n)(2\bar{\alpha}_n + n) &\geq 0
\end{aligned}$$

Using these bounds one-by-one we get the result

$$\begin{aligned}
n\left[(n+1)(\alpha - d) + d - 1 + n^2\right] &\leq (\alpha - n)^2 \\
n\left[(n+1)\overline{(\alpha - d)}_{n-1} + d - 1 + n^2\right] &\leq (\alpha - n)^2 \\
n\left[(n+1)\overline{(\alpha - d)}_{n-1} + d - 1 + (\alpha - d + 1 - \overline{(\alpha - d)}_{n-1})^2\right] &\leq (\alpha - n)^2 \\
n\left[(n+1)\overline{(\alpha - d)}_{n-1} + d - 1 + (\alpha - d + 1 - \overline{(\alpha - d)}_{n-1})^2\right] &\leq \bar{\alpha}_n^2 \\
n\left[(n+1)\overline{(\alpha - d)}_{n-1} + d - 1 + (\alpha - d + 1 - \overline{(\alpha - d)}_{n-1})^2\right] &\leq \bar{\alpha}_n^2 + (\alpha - \bar{\alpha}_n)(2\bar{\alpha}_n + n) \\
R(\alpha, d, n) &\leq 1
\end{aligned}$$

□

# Chapter 6

# Global Convergence of Arbitrary-Block Gradient Methods for Generalized Polyak-Łojasiewicz Functions

## 6.1 Introduction

During the last decade, gradient-type methods have become the methods of choice for solving optimization problems of very large sizes arising in fields such as machine learning, data science, engineering, and visual computing.

Consider the optimization problem

$$\min_{\mathbf{x}\in\mathbb{R}^n} f(\mathbf{x}),$$

where $f : \mathbb{R}^n \to \mathbb{R}$ is a differentiable function. Assume that this problem has a nonempty set of global minimizers $\mathcal{X}^*$ (clearly, $\nabla f(x^*) = 0$ for all $x^* \in \mathcal{X}^*$). It is well known [4] that if $f$ is $L$-smooth and $\mu$-strongly convex, where $L \geq \mu > 0$, then the gradient descent method $\mathbf{x}^{k+1} = \mathbf{x}^k - \frac{1}{L}\nabla f(\mathbf{x}^k)$ for all $k \geq 0$ satisfies $\xi(\mathbf{x}^{k+1}) \leq \left(1 - \frac{\mu}{L}\right)\xi(\mathbf{x}^k)$, where $\mathbf{x}^* \in \mathcal{X}^*$ and

$$\xi(\mathbf{x}) \quad := \quad f(\mathbf{x}) - f(\mathbf{x}^*), \tag{6.1}$$

is the *optimality gap function*. Motivated by the rise of nonconvex models in fields such as image and signal processing and deep learning, there is interest in studying the performance of gradient-type methods for nonconvex functions.

As observed by Polyak in 1963 [52], and recently popularized and further studied by Karimi, Nutini and Schmidt [29] in the context of proximal methods, proofs of linear convergence rely on a certain *consequence* of strong convexity known as the *Polyak-Łojasiewicz* (PL) inequality. Since functions satisfying the PL inequality need not be convex, linear convergence of gradient methods to the global optimum extends beyond the realm of convex functions.

The (strong) PL inequality can be written in the form

$$\tfrac{1}{2}\|\nabla f(\mathbf{x})\|^2 \quad \geq \quad \mu \cdot \xi(\mathbf{x}), \qquad \mathbf{x} \in \mathbb{R}^n. \tag{6.2}$$

We write $f \in \mathcal{S}_{PL}(\mu)$ if $f$ satisfies (6.2). The PL inequality and methods based on it have been an inspiration for many researchers in recent years [11, 34, 20]. It is known that in order to guarantee $\xi(\mathbf{x}^k) \leq \epsilon$, it suffices to take $k = \mathcal{O}((L/\mu)\log(1/\epsilon))$.

The starting point of this paper is the realization that *while the PL inequality serves as a generalization of strong convexity, there is no equivalent generalization of (weak) convexity*. One of the key contributions of this paper is to remedy this situation by introducing the *weak*

*PL inequality*:
$$\|\nabla f(\mathbf{x})\| \cdot \|\mathbf{x} - \mathbf{x}^*\| \geq \sqrt{\mu} \cdot \xi(\mathbf{x}), \qquad \mathbf{x} \in \mathbb{R}^n. \tag{6.3}$$

We write $f \in \mathcal{W}_{PL}(\mu)$ if $f$ satisfies (6.3). If $f$ is convex, then $f \in \mathcal{W}_{PL}(1)$. Indeed, by convexity and Cauchy-Schwartz inequality, we have

$$\xi(\mathbf{x}) = f(\mathbf{x}) - f(\mathbf{x}^*) \leq \langle \nabla f(\mathbf{x}), \mathbf{x} - \mathbf{x}^* \rangle \leq \|\nabla f(\mathbf{x})\| \cdot \|\mathbf{x} - \mathbf{x}^*\|.$$

However, $\mathcal{W}_{PL}(1)$ contains nonconvex functions as well. As an example consider the function $f(x_1, x_2) = x_1^2 x_2^2$, for which it is straightforward to show that $f \in \mathcal{W}_{PL}(1)$ and it is apparently nonconvex. If we allow $0 < \mu < 1$, the inequality (6.3) becomes weaker, and holds for a larger family of functions still.

In this paper we prove that $\xi(\mathbf{x}^k) \leq \epsilon$ if $k \geq 2LR^2/\epsilon$, where $R$ is a uniform upper bound on $\|\mathbf{x}^k - \mathbf{x}^0\|$. Since gradient descent is a monotonic methods, such a bound exists if, for instance, the level set $\{\mathbf{x} : f(\mathbf{x}) \leq f(\mathbf{x}^0)\}$ is bounded. This result extends standard convergence result for gradient descent for convex functions to weak PL functions.

Further contributions of this work are:

(i) We consider a large family of gradient type methods. The methods include block coordinate descent with arbitrary block selection rules, such as cyclic, greedy, randomized, adaptive and so on. Gradient descent arises as a special case when the active block at each iteration consists of all coordinates. Also, we introduce a novel method called greedy minibatch descent which we analyze using our developed theory to prove convergence for it in various setups mentioned in the next point.

(ii) We extend all results (weak PL inequality, algorithms and complexity results) to the proximal setup. That is, we consider composite optimization problems of the form $\min F(x) := f(x) + g(x)$, where $f$ is differentiable, and $g$ is a simple nonsmooth function. For instance, the weak PL inequality (6.3) arises as a special case of the new *proximal weak PL inequality* when $g = 0$. The complexity results are the same: $\mathcal{O}(\log(1/\epsilon))$ for PL functions, $\mathcal{O}(1/\epsilon)$ for WPL functions and $\mathcal{O}(1/\epsilon \log(1/\epsilon))$ for general nonconvex functions. The specific rates can be found in Tables 6.1 and 6.2. Note that $\rho(\mathbf{x}^0)$ is a problem dependent constant standing for the distance between the initial iterate and the optimum. To the best of our knowledge, all the rates in the tables are novel, except $\mathcal{S}_{PL}^g$ for gradient descent and uniform and greedy coordinate descents for both smooth and non-smooth case. These were already shown in [29].

| block selection rule | $\mathcal{S}_{PL}(\mu)$ | $\mathcal{W}_{PL}(\rho)$ | general $g = 0$ |
|---|---|---|---|
| gradient descent | $\frac{\lambda_{\max}(\mathbf{M})}{\mu} \log(\frac{\xi(\mathbf{x}^0)}{\epsilon})$ | $\frac{\lambda_{\max}(\mathbf{M})}{\rho(\mathbf{x}^0)\epsilon}$ | $\frac{\lambda_{\max}(\mathbf{M})\xi(\mathbf{x}^0)}{\epsilon} \log(\frac{\xi(\mathbf{x}^0)}{\epsilon})$ |
| uniform coordinate | $\frac{n \max_i \{M_{ii}\}}{\mu} \log(\frac{\xi(\mathbf{x}^0)}{\epsilon})$ | $\frac{n \max_i \{M_{ii}\}}{\rho(\mathbf{x}^0)\epsilon}$ | $\frac{n \max_i \{M_{ii}\}\xi(\mathbf{x}^0)}{\epsilon} \log(\frac{\xi(\mathbf{x}^0)}{\epsilon})$ |
| importance coord. | $\frac{\sum_{i=1}^n M_{ii}}{\mu} \log(\frac{\xi(\mathbf{x}^0)}{\epsilon})$ | $\frac{\sum_{i=1}^n M_{ii}}{\rho(\mathbf{x}^0)\epsilon}$ | $\frac{\sum_{i=1}^n M_{ii}\xi(\mathbf{x}^0)}{\epsilon} \log(\frac{\xi(\mathbf{x}^0)}{\epsilon})$ |
| greedy coordinate | $\frac{\sum_{i=1}^n M_{ii}}{\mu} \log(\frac{\xi(\mathbf{x}^0)}{\epsilon})$ | $\frac{\sum_{i=1}^n M_{ii}}{\rho(\mathbf{x}^0)\epsilon}$ | $\frac{\sum_{i=1}^n M_{ii}\xi(\mathbf{x}^0)}{\epsilon} \log(\frac{\xi(\mathbf{x}^0)}{\epsilon})$ |
| uniform minibatch | $\frac{1}{\mu \lambda_{\min}(\mathbf{E}[\mathbf{M}_{[S]}^{-1}])} \log(\frac{\xi(\mathbf{x}^0)}{\epsilon})$ | $\frac{1}{\rho(\mathbf{x}^0)\lambda_{\min}(\mathbf{E}[\mathbf{M}_{[S]}^{-1}])\epsilon}$ | $\frac{\xi(\mathbf{x}^0)}{\lambda_{\min}(\mathbf{E}[\mathbf{M}_{[S]}^{-1}])\epsilon} \log(\frac{\xi(\mathbf{x}^0)}{\epsilon})$ |
| greedy minibatch | $\frac{1}{\mu \lambda_{\min}(\mathbf{E}[\mathbf{M}_{[S]}^{-1}])} \log(\frac{\xi(\mathbf{x}^0)}{\epsilon})$ | $\frac{1}{\rho(\mathbf{x}^0)\lambda_{\min}(\mathbf{E}[\mathbf{M}_{[S]}^{-1}])\epsilon}$ | $\frac{\xi(\mathbf{x}^0)}{\lambda_{\min}(\mathbf{E}[\mathbf{M}_{[S]}^{-1}])\epsilon} \log(\frac{\xi(\mathbf{x}^0)}{\epsilon})$ |

Table 6.1: Iteration complexity guarantees for $\mathbf{E}[\xi(\mathbf{x}^K)] \leq \epsilon$ in the smooth cases.

### 6.1.1 Outline

We first perform the our general analysis specified for smooth gradient descent in Section 6.2. In Section 6.3 we introduce the general setup considered in the rest of the work. In Section 6.4 we introduce the proportion function, which is a tool for the general analysis of block selection rules. In Sections 6.5.2 to 6.7 we prove the main theory for strongly PL, weakly PL, and general non-convex functions. Finally, in Section 6.8 we perform numerical experiments confirming our theoretical findings.

| block selection rule | $\mathcal{S}^g_{PL}(\mu)$ | $\mathcal{W}^g_{PL}(\rho)$ | general $g \neq 0$ |
|---|---|---|---|
| gradient descent | $\frac{\lambda_{\max}(\mathbf{M})}{\mu} \log(\frac{\xi(\mathbf{x}^0)}{\epsilon})$ | $\frac{\lambda_{\max}(\mathbf{M})}{\rho(\mathbf{x}^0)\epsilon}$ | $\frac{\lambda_{\max}(\mathbf{M})\xi(\mathbf{x}^0)}{\epsilon} \log(\frac{\xi(\mathbf{x}^0)}{\epsilon})$ |
| uniform coordinate | $\frac{n\max_i\{M_{ii}\}}{\mu} \log(\frac{\xi(\mathbf{x}^0)}{\epsilon})$ | $\frac{n\max_i\{M_{ii}\}}{\rho(\mathbf{x}^0)\epsilon}$ | $\frac{n\max_i\{M_{ii}\}\xi(\mathbf{x}^0)}{\epsilon} \log(\frac{\xi(\mathbf{x}^0)}{\epsilon})$ |
| greedy coordinate | $\frac{n\max_i\{M_{ii}\}}{\mu} \log(\frac{\xi(\mathbf{x}^0)}{\epsilon})$ | $\frac{n\max_i\{M_{ii}\}}{\rho(\mathbf{x}^0)\epsilon}$ | $\frac{n\max_i\{M_{ii}\}\xi(\mathbf{x}^0)}{\epsilon} \log(\frac{\xi(\mathbf{x}^0)}{\epsilon})$ |
| uniform minibatch | $\frac{nL_\tau}{\tau\mu} \log(\frac{\xi(\mathbf{x}^0)}{\epsilon})$ | $\frac{nL_\tau}{\rho(\mathbf{x}^0)\tau\epsilon}$ | $\frac{\xi(\mathbf{x}^0)nL_\tau}{\tau\epsilon} \log(\frac{\xi(\mathbf{x}^0)}{\epsilon})$ |
| greedy minibatch | $\frac{nL_\tau}{\tau\mu} \log(\frac{\xi(\mathbf{x}^0)}{\epsilon})$ | $\frac{nL_\tau}{\rho(\mathbf{x}^0)\tau\epsilon}$ | $\frac{\xi(\mathbf{x}^0)nL_\tau}{\tau\epsilon} \log(\frac{\xi(\mathbf{x}^0)}{\epsilon})$ |

Table 6.2: Iteration complexity guarantees for $\mathbf{E}[\xi(\mathbf{x}^K)] \leq \epsilon$ in the non-smooth cases.

### 6.1.2 Notation

We use boldface to denote a multi-dimensional object. As an example, we have a vector $\mathbf{x}$, a matrix $\mathbf{X}$, while a scalar entry of a vector $x_i$ has a normal typeset. By $[n]$ we denote the set $\{1, \ldots, n\}$. $\|\mathbf{x}\| = \langle \mathbf{x}, \mathbf{x}\rangle^{1/2}$ is the L2 norm, where $\langle \mathbf{x}, \mathbf{y}\rangle = \sum_i x_i y_i$ is the standard inner product. refer to Table 6.3 in the appendix for a summary of frequently used notation.

## 6.2 Gradient Descent

We assume throughout this section that $f$ is $L$-smooth for some $L > 0$:

$$f(\mathbf{x} + \mathbf{h}) \leq f(\mathbf{x}) + \langle \nabla f(\mathbf{x}), \mathbf{h}\rangle + \tfrac{L}{2}\|\mathbf{h}\|^2, \qquad \mathbf{x}, \mathbf{h} \in \mathbb{R}^n. \tag{6.4}$$

We shall write $f \in C^1(L)$. In addition to this assumption, in our analysis we consider several classes of *nonconvex* objectives: PL functions, weak PL functions, and gradient dominated functions.

In this section we perform a novel analysis of the gradient descent[1] method for minimizing $f$:

$$\mathbf{x}^{k+1} = \mathbf{x}^k - \tfrac{1}{L}\nabla f(\mathbf{x}^k). \tag{6.5}$$

for the above classes of nonconvex functions. As we shall show, for these classes of objectives gradient descent converges to the global minimizer. By focusing on the notoriously known gradient descent method first, we illuminate some of the key insights of this paper without distractions from additional complications caused by the proximal setup and particularities of other algorithms, making the more general treatment in further sections more easily digestable.

A key role in the analysis is played by the *forcing function* associated with $f$, defined as

$$\mu(\mathbf{x}) := \frac{\|\nabla f(\mathbf{x})\|^2}{2\xi(\mathbf{x})}, \qquad x \in \mathbb{R}^n/\mathcal{X}^*. \tag{6.6}$$

For any fixed value of this function, the smaller the gradient $\|\nabla f(\mathbf{x})\|^2$ is, the smaller the optimality gap $\xi(\mathbf{x})$. In other words, small gradients force the optimality gap to become small. The importance of this function is clear from the following simple lemma, which says that the larger $\mu(\mathbf{x}^k)$ is, the more reduction we get at iteration $k$ in the optimality gap.

**Lemma 6.1.** *Let $f \in C^1(L)$ and let $\{\mathbf{x}^k\}_{k\geq 0}$ be the sequence of iterates produced by the gradient descent method (6.5). As long as $\mathbf{x}^k \notin \mathcal{X}^*$, we have*

$$\xi(\mathbf{x}^{k+1}) \leq \left(1 - \frac{\mu(\mathbf{x}^k)}{L}\right)\xi(\mathbf{x}^k).$$

*Moreover, $\|\nabla f(\mathbf{x})\|^2 \leq 2L\xi(\mathbf{x})$ for all $\mathbf{x}$.*

---

[1] For simplicity, we consider gradient descent with fixed stepsize inversely proportional to the Lipschitz constant: $1/L$. While one can extend our results to other stepsize strategies using standard techniques, we avoid doing so as to present our results in a simple setting.

*Proof.* Let $\mathbf{h}^k = -\frac{1}{L}\nabla f(\mathbf{x}^k)$. Then

$$
\begin{aligned}
\xi(\mathbf{x}^{k+1}) &\stackrel{(6.5)}{=} \xi(\mathbf{x}^k + \mathbf{h}^k) \\
&\stackrel{(6.1)}{=} f(\mathbf{x}^k + \mathbf{h}^k) - f(\mathbf{x}^*) \\
&\stackrel{(6.4)}{\leq} f(\mathbf{x}^k) + \langle \nabla f(\mathbf{x}^k), \mathbf{h}^k \rangle + \tfrac{L}{2}\|\mathbf{h}^k\|^2 - f(\mathbf{x}^*) \\
&\stackrel{(6.1)}{=} \xi(\mathbf{x}^k) + \langle \nabla f(\mathbf{x}^k), \mathbf{h}^k \rangle + \tfrac{L}{2}\|\mathbf{h}^k\|^2 \\
&= \xi(\mathbf{x}^k) - \tfrac{1}{2L}\|\nabla f(\mathbf{x}^k)\|^2 \\
&\stackrel{(6.6)}{=} (1 - \tfrac{\mu(\mathbf{x}^k)}{L})\xi(\mathbf{x}^k).
\end{aligned}
$$

Since $\xi(\mathbf{x}^{k+1}) \geq 0$, it must be the case that $\mu(\mathbf{x}) \leq L$ for all $\mathbf{x} \notin \mathcal{X}^*$. $\square$

It is well known that gradient descent is monotonic: $\xi(\mathbf{x}^{k+1}) \leq \xi(\mathbf{x}^k)$ for all $k$. Note that this property follows from the second-to-last identity in the proof, and relies on the assumption of $L$-smoothness only. If $f$ is $\mu$-strongly convex or, more generally, if $f \in \mathcal{S}_{PL}(\mu)$, then $\mu(\mathbf{x}) \geq \mu$ for all $\mathbf{x} \notin \mathcal{X}^*$, and Lemma 6.1 implies the linear rate $\xi(\mathbf{x}^{k+1}) \leq \left(1 - \frac{\mu}{L}\right)\xi(\mathbf{x}^k)$. This result was shown already by Polyak [52].

### 6.2.1 Weakly Polyak-Łojasiewicz functions

Consider now functions satisfying a weak version of the PL inequality. To the best of our knowledge, this is the first work where such functions are considered.

**Definition 6.2** (Weak Polyak-Łojasiewicz functions). We say that $f : \mathbb{R}^n \to \mathbb{R}$ is a weak Polyak-Łojasiewicz (WPL) function with parameter $\mu \geq 0$ if there exists $\mathbf{x}^* \in \mathcal{X}^*$ such that

$$\sqrt{\mu} \cdot \xi(\mathbf{x}) \leq \|\nabla f(\mathbf{x})\| \cdot \|\mathbf{x} - \mathbf{x}^*\|, \qquad \mathbf{x} \in \mathbb{R}^n. \tag{6.7}$$

For simplicity, we write $f \in \mathcal{W}_{PL}(\mu)$.

Consider the Huber loss given by

$$H(z) := \begin{cases} z^2 & |z| < 1 \\ 2|z| - 1 & \text{otherwise,} \end{cases}$$

and the derived function $f$ given by $f(x_1, x_2) = H(x_1)H(x_2)$. It is straightforward to show from the definitions, that $f$ is an example of a smooth function, such that $f \in \mathcal{W}_{PL}(\mu)$ for some $\mu > 0$ and $f \notin \mathcal{S}_{PL}(\mu)$ for all $\mu > 0$.

Note that all[2] functions belong to $\mathcal{W}_{PL}(0)$. As the next result shows, WPL functions admit a lower bound on $\mu(\mathbf{x})$ which is proportional to $\xi(\mathbf{x})$ and inversely proportional to $\|\mathbf{x} - \mathbf{x}^*\|^2$.

**Lemma 6.3.** *If $f \in \mathcal{W}_{PL}(\mu)$, then*

$$\mu(\mathbf{x}) \geq \frac{\mu \xi(\mathbf{x})}{2\|\mathbf{x} - \mathbf{x}^*\|^2}, \qquad \mathbf{x} \in \mathbb{R}^n / \mathcal{X}^*.$$

*Proof.* We have

$$\mu(\mathbf{x}) \stackrel{(6.6)}{=} \frac{\|\nabla f(\mathbf{x})\|^2}{2\xi(\mathbf{x})} \stackrel{(6.7)}{\geq} \frac{\mu \xi^2(\mathbf{x})/\|\mathbf{x} - \mathbf{x}^*\|^2}{2\xi(\mathbf{x})} = \frac{\mu \xi(\mathbf{x})}{2\|\mathbf{x} - \mathbf{x}^*\|^2}.$$

$\square$

Several basic properties of WPL functions are summarized in Appendix 6.A. Combining Lemma 6.1 and Lemma 6.3, we get the recursion

---
[2] By "all" we implicitly mean all functions for which the definition make sense. That is, differentiable and bounded from below.

$$\xi(\mathbf{x}^{k+1}) \leq \left(1 - \frac{\mu \xi(\mathbf{x}^k)}{2L\|\mathbf{x}^k - \mathbf{x}^*\|^2}\right) \xi(\mathbf{x}^k). \tag{6.8}$$

The next lemma will be useful in the analysis of this recursion.

**Lemma 6.4.** *Let $\{\alpha^t\}_{t=0}^k$ and $\{\beta^t\}_{t=0}^k$ be two sequences of positive numbers satisfying the recursion*

$$\alpha^{t+1} \leq \left(1 - \alpha^t \beta^t\right) \alpha^t. \tag{6.9}$$

*Then for all $k \geq 0$ we have the bound*

$$\alpha^k \leq \frac{\alpha^0}{1 + \alpha^0 \sum_{t=0}^{k-1} \beta^t}.$$

*Proof.* As $\alpha^t$ and $\beta^t$ are positive numbers, we have $\alpha^{t+1} \leq \alpha^t$ for all $t$ using (6.9). Observe that

$$\frac{1}{\alpha^{t+1}} - \frac{1}{\alpha^t} = \frac{\alpha^t - \alpha^{t+1}}{\alpha^{t+1} \alpha^t} \geq \frac{\alpha^t - \alpha^{t+1}}{(\alpha^t)^2} \overset{(6.9)}{\geq} \beta^t,$$

which we can recursively used to show

$$\frac{1}{\alpha^k} \geq \frac{1}{\alpha^{k-1}} + \beta^{k-1} \geq \ldots \geq \frac{1}{\alpha^0} + \sum_{t=0}^{k-1} \beta^t.$$

We get the result by inverting the last equation. $\square$

By applying the above lemma to recursion (6.8), we get a global convergence result for gradient descent applied to a WPL function.

**Theorem 6.5.** *Let $f \in \mathcal{C}^1(L)$, and let $\{\mathbf{x}^k\}_{k \geq 0}$ be the sequence of iterates produced by the gradient descent method (6.5). Assume $f \in \mathcal{W}_{PL}(\mu)$ for $\mu > 0$. Then for all $k \geq 1$ we have*

$$\xi(\mathbf{x}^k) \leq \frac{\xi(\mathbf{x}^0)}{1 + \xi(\mathbf{x}^0) \frac{\mu}{2L} \sum_{t=0}^{k-1} \frac{1}{\|\mathbf{x}^t - \mathbf{x}^*\|^2}} \leq \frac{2L}{\mu} \cdot \frac{1}{\sum_{t=0}^{k-1} \frac{1}{\|\mathbf{x}^t - \mathbf{x}^*\|^2}}. \tag{6.10}$$

The second inequality is obtained from the first by neglecting the additive constant 1 in the denominator.

By monotonicity, all iterates of gradient descent stay in the level set $\mathcal{L}_0 := \{\mathbf{x} \in \mathbb{R}^n : f(\mathbf{x}) \leq f(\mathbf{x}^0)\}$. If this set is bounded, then $R := \max_{\mathbf{x} \in \mathcal{L}_0} \|\mathbf{x} - \mathbf{x}^*\| < +\infty$, and we have $\|\mathbf{x}^k - \mathbf{x}^*\| \leq R$ for all $k$. In this case, the bound (6.10) implies

$$\xi(\mathbf{x}^k) \leq \frac{2LR^2}{\mu k}.$$

### 6.2.2 Gradient dominated functions

We now consider a new class of (not necessarily convex) functions. To the best of our knowledge, this class was not considered in optimization before.

**Definition 6.6.** *We say that function $f : \mathbb{R}^n \to \mathbb{R}$ is $\varphi$-gradient dominated if there exists an increasing function $\varphi : \mathbb{R}_+ \to \mathbb{R}_+$ such that $\varphi(0) = 0$, $\lim_{t \to 0} \varphi(t) = 0$ and*

$$\xi(\mathbf{x}) \leq \varphi(\|\nabla f(\mathbf{x})\|^2), \qquad \mathbf{x} \in \mathbb{R}^n. \tag{6.11}$$

The above definition essentially says that for any sequence $\{\mathbf{x}^k\}$ (not necessarily related to iterates of gradient descent) such that $\|\nabla f(\mathbf{x}^k)\|^2 \to 0$, we must have $f(\mathbf{x}^k) \to f(\mathbf{x}^*)$. In particular, if $f$ has multiple minimizers, all must have the same function value.

As an example of the function $\varphi$, we might consider any function of the form $\varphi(t) = c \cdot |t|^p$, where $c > 0$ and $p > 0$.

**Theorem 6.7.** *Assume that $f \in \mathcal{C}^1(L)$ is $\varphi$-gradient dominated. Pick $\epsilon > 0$ and let*

$$k \geq \frac{2L\xi(\mathbf{x}^0)}{\varphi^{-1}(\epsilon)} \log\left(\frac{\xi(\mathbf{x}^0)}{\varphi(\epsilon)}\right). \tag{6.12}$$

*Then* $\min\{\xi(\mathbf{x}^t) \,:\, t = 0, 1, \ldots, k\} \leq \varphi(\epsilon)$.

*Proof.* If $\varphi(\|\nabla f(\mathbf{x}^t)\|^2) \leq \epsilon$ for some $k = 0, 1, \ldots, k - 1$, we are done by applying (6.11). Otherwise, $\varphi(\|\nabla f(\mathbf{x}^t)\|^2) > \epsilon$ for all $t = 0, 1, \ldots, k - 1$. Then in view of (6.6), $\mu(\mathbf{x}^t) > \varphi^{-1}(\epsilon)/(2\xi(\mathbf{x}^t))$ for all such $t$. By monotonicity, $\xi(\mathbf{x}^{t+1}) \leq \xi(\mathbf{x}^t)$ for all $t$, whence $\mu(\mathbf{x}^t) > \varphi^{-1}(\epsilon)/(2\xi(\mathbf{x}^0))$. By applying Lemma 6.1, we get

$$\xi(\mathbf{x}^{t+1}) \leq \left(1 - \frac{\varphi^{-1}(\epsilon)}{2L\xi(\mathbf{x}^0)}\right) \xi(\mathbf{x}^t), \qquad t = 0, 1, \ldots, k - 1.$$

By unrolling the recurrence, and using the bound $(1 - s)^{1/s} \leq \exp(-1)$ (which holds for $0 < s \leq 1$), we get

$$\xi(\mathbf{x}^k) \leq \left(1 - \frac{\varphi^{-1}(\epsilon)}{2L\xi(\mathbf{x}^0)}\right)^k \xi(\mathbf{x}^0) \leq \exp\left(-\frac{k\varphi^{-1}(\epsilon)}{2L\xi(\mathbf{x}^0)}\right) \xi(\mathbf{x}^0) \overset{(6.12)}{\leq} \varphi(\epsilon).$$

$\square$

### 6.2.3 Literature Review

The original gradient descent method was developed by Cauchy [5] and it has seen a lot of development ever since. This development is documented in detail in [50, 4]. In the recent years, a version of gradient descent called coordinate descent was developed. The first developments of coordinate descent are due to [35] and it was first analyzed for general convex objectives by Nesterov in [49].

An important part of the coordinate descent is its ability to work with arbitrary block selection strategies. In the seminal work of Nesterov [49], there were three strategies introduced, which are known as coordinate descent with uniform probabilities, coordinate descent with importance sampling, and greedy coordinate descent. The first two strategies fall into the family of randomized strategies. These were further developed in [59, 57, 60, 54, 10, 56, 8]. The third selection rule is a deterministic strategy similar in nature to batch gradient descent and it was shown to be superior over randomized methods in terms of iteration complexity in [51].

The Polyak-Łojasiewicz condition was first introduced by Boris Polyak in [52]. It was revived recently in [29] and applied to modern optimization approaches. Since then, multiple papers used the condition to develop new approaches [11, 34, 20].

Lastly, we note that our framework considers the non-accelerated version of coordinate descent methods, although they play a key role in modern theory. If needed, the acceleration can be achieved non-directly by using approaches as proposed in [36] or [19]. We leave the accelerated counterpart of this framework to future work.

## 6.3 General Setup

In this section we move beyond the simplified setup considered in the previous section and introduce the setting considered in this paper in its full generality. Our general treatment differs from that in Section 6.2 in several ways.

First, we consider the *composite* optimization problem

$$\min_{x \in \mathbb{R}^n} \{F(\mathbf{x}) := f(\mathbf{x}) + g(\mathbf{x})\}, \tag{6.13}$$

where $f$ is assumed to be smooth, and $g$ is a simple (possibly nonconvex and nonsmooth) separable function. Second, $f$ is assumed to be smooth in a slightly more general sense than $L$-smoothness of Section 6.2. Third, we go beyond gradient descent and consider a large family of first order methods which include randomized, cyclic, adaptive and greedy coordinate descent, in serial and adaptive settings.

In the following, we will refer to the problem (6.13) with $g = 0$ as the smooth case and otherwise as the non-smooth case.

### 6.3.1 Smoothness and separability

We assume, that $f$ is **M**-smooth, which is formalized by the following assumption:

**Assumption 6.8.** *We say that a function $f : \mathbb{R}^n \to \mathbb{R}$ is **M**-smooth, if there exists a positive definite matrix $\mathbf{M} \in \mathbb{R}^{n \times n}$ such that*

$$f(\mathbf{x} + \mathbf{h}) \leq f(\mathbf{x}) + \langle \nabla f(\mathbf{x}), \mathbf{h} \rangle + \tfrac{1}{2} \langle \mathbf{M}\mathbf{h}, \mathbf{h} \rangle, \qquad \mathbf{x}, \mathbf{h} \in \mathbb{R}^n. \tag{6.14}$$

In the non-smooth case ($g \neq 0$), we will without loss of generality assume **M** to be a multiple of the identity matrix; specifically $\mathbf{M} := L\mathbf{I}$, where **I** is the identity matrix. If (6.8) holds for some **M**, we can always replace **M** by $L\mathbf{I}$, where $L = \lambda_{\max}(\mathbf{M})$. It is easy to verify, that if a function is **M**-smooth, it is also $L\mathbf{I}$-smooth. Note that in some cases we might choose $L$ to be smaller. This will be in detail explained in Section 6.4.2.

For simplicity, we will use the notation with **M** also in the non-smooth case, but we will always treat it as the diagonal matrix $L\mathbf{I}$.

In addition to smoothness of $f$, we assume that the function $g$ is separable, which is defined as follows:

**Assumption 6.9.** *We say that a function $g : \mathbb{R}^n \to \bar{\mathbb{R}} := \mathbb{R} \cup \{+\infty\}$ is separable, if there exist $n$ scalar functions $g_1, \ldots, g_n : \mathbb{R}^n \to \bar{\mathbb{R}}$, such that*

$$g(\mathbf{x}) = \sum_{i=1}^{n} g_i(x_i). \tag{6.15}$$

Note that the function $g$ is treated as the non-smooth part of the problem (e.g., L1 norm, box constraints, and so on). When we refer to a smooth problem, we assume the setup with $g = 0$, while all the other setups are referred as non-smooth problems.

The problem described in (6.13) is encountered in many areas, ranging from machine learning and signal processing to biology and beyond. We believe it does not need to be motivated further, as it was already considered in a lot of previous works.

### 6.3.2 Masking vectors and matrices

Let $\mathbf{x} \in \mathbb{R}^n$ be an arbitrary vector and let $\mathbf{X} \in \mathbb{R}^{n \times n}$ be an arbitrary matrix. We will need to index vectors and matrices by subsets of coordinates $\emptyset \neq S \subseteq [n]$. The indexing has two distinct forms. By $\mathbf{x}_S$ we denote the $|S|$-dimensional vector constructed by taking the entries of $\mathbf{x}$ with indices in $S$, while the notation $\mathbf{x}_{[S]}$ is used to zero out every entry of $\mathbf{x}$ not appearing in $S$ without changing its length. We have similar notation for matrices, where $\mathbf{X}_S$ denotes the $|S| \times |S|$ matrix of the entries with both column and row indices in $S$, while $\mathbf{X}_{[S]}$ is used for the matrix $\mathbf{X}$ with entries zeroes out outside of columns and rows with indices in $S$. As a rule of thumb, the subscript $S$ changes the dimensions of the object, while the subscript $[S]$ maintains its dimensions.

To illustrate the notation, consider the following example.

**Example 1.** *Let $\mathbf{x} \in \mathbb{R}^3$, $\mathbf{X} \in \mathbb{R}^{3 \times 3}$ and $S = \{1, 3\}$. Then*

$$\mathbf{x} = \begin{bmatrix} 1 \\ 2 \\ 3 \end{bmatrix} \quad \Rightarrow \quad \mathbf{x}_S = \begin{bmatrix} 1 \\ 3 \end{bmatrix}, \quad \mathbf{x}_{[S]} = \begin{bmatrix} 1 \\ 0 \\ 3 \end{bmatrix}$$

*and*

$$\mathbf{X} = \begin{bmatrix} 1 & 4 & 7 \\ 2 & 5 & 8 \\ 3 & 6 & 9 \end{bmatrix} \quad \Rightarrow \quad \mathbf{X}_S = \begin{bmatrix} 1 & 7 \\ 3 & 9 \end{bmatrix}, \quad \mathbf{X}_{[S]} = \begin{bmatrix} 1 & 0 & 7 \\ 0 & 0 & 0 \\ 3 & 0 & 9 \end{bmatrix}.$$

### 6.3.3 Algorithm

We now propose and analyze a wide class of block descent algorithms for solving (6.13). For a non-empty block of coordinates $S \subseteq [n]$, we define

$$U_S(\mathbf{x}, \mathbf{u}) := \langle (\nabla f(\mathbf{x}))_{[S]}, \mathbf{u} \rangle + \frac{1}{2}\mathbf{u}^\top \mathbf{M}_{[S]} \mathbf{u} + \sum_{i \in S} \left[g_i(x_i + u_i) - g_i(x_i)\right]. \tag{6.16}$$

We assume that finding a minimizer of $U_S(\mathbf{x}, \mathbf{h})$ in $\mathbf{h}$ is cheap (e.g., there exists a closed form solution, or an efficient algorithm). Given an iterate $\mathbf{x}^k$, in iteration $k$ of our method we select a block $S_k \subseteq [n]$ of active coordinates, according to an *arbitrary block selection procedure* $\mathcal{P}_k$, and subsequently minimize $U_{S_k}(\mathbf{x}^k, \mathbf{h})$ in $\mathbf{h}$. The result is denoted $\mathbf{u}^k$. Due to the structure of problem (6.16), the minimizer of $U_{S_k}$ does not depend on $u_i$ for $i \notin S_k$. Hence, only the active coordinates $\mathbf{u}_i^k$ for $i \in S_k$ are relevant, and we set $\mathbf{x}^{k+1} = \mathbf{x}^k + \mathbf{u}_{[S_k]}^k$. Equivalently, using other alternative notation, we can write this as

$$\mathbf{x}_i^{k+1} = \begin{cases} \mathbf{x}_i^k + \mathbf{u}_i^k, & i \in S_k, \\ \mathbf{x}_i^k, & i \notin S_k, \end{cases}$$

or

$$\mathbf{x}_{S_k}^{k+1} = \mathbf{x}_{S_k}^k + \mathbf{u}_{S_k}^k, \qquad \mathbf{x}_{[n]/S_k}^{k+1} = \mathbf{x}_{[n]/S_k}^k.$$

This is Algorithm 8.

---
**Algorithm 8** Proximal Arbitrary-Block Descent Method

**Input:** Initial iterate $\mathbf{x}^0 \in \text{dom}(g)$, arbitrary block selection procedures $\{\mathcal{P}_k\}_{k=0}^{K-1}$
**for** $k = 0, \ldots, K-1$ **do**
  Pick a non-empty subset of coordinates $S_k \subset [n]$ using the procedure $\mathcal{P}_k$
  Compute $\mathbf{u}^k \in \arg\min_{\mathbf{u}}\{U_{S_k}(\mathbf{x}^k, \mathbf{u})\}$
  Update $\mathbf{x}^{k+1} = \mathbf{x}^k + \mathbf{u}_{[S_k]}^k$
**end for**
Output $\mathbf{x}^K$

---

Observe, that in the smooth case we have

$$\mathbf{u}_{S_k}^k = -\mathbf{M}_{S_k}^{-1} \nabla_{S_k} f(\mathbf{x}^k), \tag{6.17}$$

where $\nabla_S f(\mathbf{x}) := (\nabla f(\mathbf{x}))_S$. The iteration can be computed as a solution of a linear system of a size $|S_k| \times |S_k|$, which is very cheap for small $|S_k|$. In the non-smooth case, the iterate does not have a closed-form solution in general but can be solved fast for a lot of forms of $g$, e.g., box constraints or L1 norm. Since we assume $\mathbf{M} = L\mathbf{I}$, we can write

$$\mathbf{u}_{S_k}^k = \arg\min_{\mathbf{u} \in \mathbb{R}^{|S_k|}} \left\{ \langle \nabla_S f(\mathbf{x}^k), \mathbf{u} \rangle + \frac{L}{2}\|\mathbf{u}\|^2 + \sum_{i \in S_k} g_i(x_i + u_i) \right\}. \tag{6.18}$$

For simplicity, assume $\mathbf{M} = L\mathbf{I}$ for some constant $L > 0$. If $g = 0$, and we always pick $S_k = [n]$, then Algorithm 8 reduces to gradient descent, considered in Section 6.2. If $g \neq 0$, and we always pick $S_k = [n]$, then Algorithm 8 reduces to proximal gradient descent. On the other hand, if we always pick $|S_k| = 1$, then we obtain coordinate descent ($g = 0$) or proximal coordinate descent ($g \neq 0$). The selection procedure $\mathcal{P}_k$ may be set to choose the coordinates in a cyclic manner, greedily, randomly according to any (fixed or evolving) probability law, and even adaptively to the entire history of the iterative process. There are many other possibilities between the two extremes of always selecting $|S_k| = 1$ and $|S_k| = n$. Such methods can be considered block coordinate descent methods, subspace descent methods, or parallel coordinate descent methods (as the updates to individual coordinates can be performed in parallel). We stress that unlike all other methods considered in the literature, in our method we allow for the block selection procedure $\mathcal{P}_k$ to be *arbitrary*, without any restrictions whatsover.

By removing these restrictions, we allow for several new possibilities in the sampling procedure. These are: 1) the block selection procedure might change from iteration to iteration, allowing for adaptive strategies as [10], 2) the procedure might depend on previous iterations, which opens up the possibility of cyclic and other similar selections, and 3) the procedure does not have to be randomized, which allows for greedy selection procedures [51].

In all our convergence results we shall enforce several key common assumptions, together with some additional assumptions. In order to avoid repeating the common core, we shall summarize them here.

**Assumption 6.10** (Common Assumptions). *Let $f : \mathbb{R}^n \to \mathbb{R}$ be an **M**-smooth (6.14) function, let $g : \mathbb{R}^n \to \bar{\mathbb{R}}$ be separable (6.15), and let $F : \mathbb{R}^n \to \bar{\mathbb{R}}$ be function defined using $f$ and $g$ as in (6.13). Assume $F$ has a global minimizer $\mathbf{x}^*$, such that $F(\mathbf{x}^*) > -\infty$. Let $\mathbf{x}^0 \in \text{dom}(g)$ be an initial point and let the sequence $\{\mathbf{x}^k\}_{k=1}^K$ be generated using Algorithm 8, where $\{S_k\}_{k=0}^{K-1}$ is an arbitrary sequence of non-empty (possibly random) subsets of $[n]$.*

### 6.3.4 Forcing function

As before, let us define the *optimality gap function*

$$\xi(\mathbf{x}) \;:=\; F(\mathbf{x}) - F(\mathbf{x}^*), \tag{6.19}$$

where $\mathbf{x}^*$ is a minimizer of $F$. Observe that $\xi(\mathbf{x}) \geq 0$ with equality only if $\mathbf{x} \in \mathcal{X}^* := \arg\min_\mathbf{x} F(\mathbf{x})$.

We now extend the definition of the forcing function to the proximal setting.

**Definition 6.11** (Forcing function: proximal version). Let

$$\lambda(\mathbf{x}) \;:=\; -L \cdot \min_{\mathbf{y} \in \mathbb{R}^n} \left\{ \langle \nabla f(\mathbf{x}), \mathbf{y} \rangle + \frac{L}{2}\|\mathbf{y}\|^2 + g(\mathbf{x} + \mathbf{y}) - g(\mathbf{x}) \right\}, \tag{6.20}$$

where $L$ is chosen such that (6.14) holds with $\mathbf{M} = L\mathbf{I}$. Specific choices of $L$ are explained in detail in Section 6.4.2. The non-negative function

$$\mu(\mathbf{x}) \;:=\; \frac{\lambda(\mathbf{x})}{\xi(\mathbf{x})} \tag{6.21}$$

is the *proximal forcing function*. The domain of $mu$ is $\text{dom}(g)/\mathcal{X}^*$.

Pick any $\mathbf{x} \in \text{dom}(g)/\mathcal{X}^*$. The minimum in (6.20) is non-positive since setting $\mathbf{y} = \mathbf{0}$ gives a zero value. Hence, $\lambda(\mathbf{x}) \geq 0$. Since the denominator in (6.21) is positive, $\mu(\mathbf{x})$ is always nonnegative.

In the smooth case ($g = 0$) we have

$$\mu(\mathbf{x}) \;=\; \frac{\|\nabla f(\mathbf{x})\|^2}{2\xi(\mathbf{x})}, \tag{6.22}$$

which is simply the forcing function (6.6).

### 6.3.5 Proportion function

We introduce one more notion, which we call the *proportion function*. This function plays an important role in our theory.

**Definition 6.12.** Let

$$\mathcal{X} \;:=\; \{\mathbf{x} \in \mathbb{R}^n \;:\; \lambda(\mathbf{x}) \neq 0\}. \tag{6.23}$$

The *proportion function* is defined by

$$\theta(S, \mathbf{x}) \;:=\; \frac{-\min_{\mathbf{u} \in \mathbb{R}^n} \mathbf{U}_S(\mathbf{x}, \mathbf{u})}{\lambda(\mathbf{x})} \tag{6.24}$$

for all $S \in \mathcal{S}$ and $\mathbf{x} \in \mathcal{X}$. For $\mathbf{x} \in \mathcal{X}^*$, we set $\theta(S, \mathbf{x}) = 0$.

Similarly as in the definition of the forcing function, both of the minimizations in (6.24) are non-positive, as $\mathbf{u} = 0$ and $\mathbf{y} = 0$ gives zero value in the numerator and denominator, respectively. Also observe that in the smooth case ($g = 0$) we have

$$\theta(S, \mathbf{x}) \quad = \quad \frac{(\nabla_S f(\mathbf{x}))^\top (\mathbf{M}_S)^{-1} \nabla_S f(\mathbf{x})}{\|\nabla f(\mathbf{x})\|^2}. \tag{6.25}$$

Note that the matrix $(\mathbf{M}_S)^{-1}$ exists since all principal submatrices of a positive definite matrix are also positive definite. A more detailed treatment of the proportion function can be found in Section 6.4.

Also note that $\mathcal{X}$ (6.23) might differ from $\mathcal{X}^*$ in the case of local minimizers.

### 6.3.6 Generic descent lemma

We now formulate a simple but important descent lemma which bounds the progress gained by a single iteration of Algorithm 8. Our bound applies to arbitrary block selection rules, and will enable us to prove global convergence results for Algorithm 8 for new classes of nonconvex functions.

**Lemma 6.13.** *Let $\mathbf{x}^{k+1}$ be the next iterate of Algorithm 8 generated from $\mathbf{x}^k$ by picking a nonempty set of coordinates $S_k \subseteq [n]$. Then*

$$\xi(\mathbf{x}^{k+1}) \quad \leq \quad \left[1 - \theta(S_k, \mathbf{x}^k) \cdot \mu(\mathbf{x}^k)\right] \cdot \xi(\mathbf{x}^k). \tag{6.26}$$

*Applying this repeatedly, for all $K \geq 1$ we obtain the estimate*

$$\xi(\mathbf{x}^K) \quad \leq \quad \left(\prod_{k=0}^{K-1} \left[1 - \theta(S_k, \mathbf{x}^k) \cdot \mu(\mathbf{x}^k)\right]\right) \xi(\mathbf{x}^0). \tag{6.27}$$

*Proof.* See Section 6.B.1 □

Recursion (6.26) is a direct generalization of Lemma 6.1. Indeed, if $\mathbf{M} = L\mathbf{I}$, $g = 0$, and $S_k = [n]$, then in the view of (6.25), $\theta([n], \mathbf{x}^k) = 1/L$. Note that as $\theta \geq 0$ and $\mu \geq 0$, we can be sure that $\xi(\mathbf{x}^{k+1}) \leq \xi(\mathbf{x}^k)$, which means that we are not getting worse by iterating the Algorithm 8. The proof of Lemma 6.13 is straightforward. The difficulty will lie in bounding the forcing and proportion functions so as to obtain convergence. In Sections 6.5.2, 6.6 and 6.7 we apply this lemma to prove the main results of this paper.

The $K$-step bound (6.27) provides us with a compact and generic bound on the optimality gap at the $K$-th iterate of Algorithm 8 dependent on the iterates $\mathbf{x}^0, \mathbf{x}^1, \ldots, \mathbf{x}^{K-1}$ and the selected sets $S_0, S_1, \ldots, S_{k-1}$. This result does not immediately imply convergence as at this level of generality, it is possible for the product appearing in (6.27) not to converge to zero. Indeed, this corollary also covers the situation where $S_k = 1$ for all $k$, which clearly can't result in convergence. We will need to introduce further restrictions in order to establish convergence.

## 6.4 Proportion Function

In this section we show standard bounds on the proportion function, which are independent of a given iterate. This will be important further, to recover the convergence rates given by standard theory. We note, that for stochastic methods we bound the expectation of the proportion function conditioned on the last iterate, instead of directly bounding the proportion function. This quantity will be important in the theory for the convergence of stochastic methods, as specified in following sections.

The proportion function is the only quantity in the convergence theory, which is dependent on the choice of the set of coordinates $S$. Therefore, to analyze a new sampling strategy for coordinate descent, one only has to show a bound on the proportion function. This opens up a possible venue of novel techniques for coordinate selection.

In the following we tackle all the known cases of samplings, which we break down by their smoothness. Also, we introduce and analyze a new sampling procedure to showcase the generality of our framework.

### 6.4.1 Smooth problems

In the case of $g = 0$, we have the proportion function equal to (6.25), i.e.,

$$\theta(S, \mathbf{x}) \quad = \quad \frac{(\nabla_S f(\mathbf{x}))^\top (\mathbf{M}_S)^{-1} \nabla_S f(\mathbf{x})}{\|\nabla f(\mathbf{x})\|^2}$$

for all $\mathbf{x} \in \mathcal{X}^*$, and all $\emptyset \neq S \subseteq [n]$. Let us break-down the cases according to specific choices of the set $S$.

- **Batch Gradient descent:** In the case when $S = [n]$, we recover the standard gradient descent strategy, which dates back to the work of Cauchy [5]. In this case we can lower bound the proportion function by

$$\theta([n], \mathbf{x}) \quad = \quad \frac{(\nabla f(\mathbf{x}))^\top \mathbf{M}^{-1} \nabla f(\mathbf{x})}{\|\nabla f(\mathbf{x})\|^2} \quad \geq \quad \frac{1}{\lambda_{\max}(\mathbf{M})} \qquad (6.28)$$

  for all $\mathbf{x} \in \mathcal{X}^*$.

- **Serial Coordinate Descent:** Suppose $S = \{i\}$, for some given $i$. In this case we get

$$\theta(\{i\}, \mathbf{x}) \quad = \quad \frac{(\nabla_i f(\mathbf{x}))^2}{M_{ii}\|\nabla f(\mathbf{x})\|^2}. \qquad (6.29)$$

  There are multiple strategies for choosing the coordinate $i$, and we tackle them one by one.

  - **Uniform probabilities:** Suppose we choose the coordinate $i$ with the probability $p_i$ given by $p_i = 1/n$, independently of $\mathbf{x}$. This strategy was originally analyzed in [49]. The expectation of $\theta$ can be lower bounded as

$$\mathbf{E}_i[\theta(\{i\}, \mathbf{x}) \mid \mathbf{x}] \quad \stackrel{(6.29)}{=} \quad \frac{1}{n}\sum_{i=1}^{n} \frac{(\nabla_i f(\mathbf{x}))^2}{M_{ii}\|\nabla f(\mathbf{x})\|^2} \quad \geq \quad \frac{1}{n \cdot \max_i\{M_{ii}\}}. \qquad (6.30)$$

  - **Importance sampling:** Suppose we choose the coordinate $i$ with the probability $p_i$ given by

$$p_i \quad = \quad \frac{M_{ii}}{\sum_{i=1}^{n} M_{ii}}, \qquad (6.31)$$

    independently of $\mathbf{x}$. Again, this strategy was originally analyzed in [49]. The expectation of $\theta$ can be lower bounded as

$$\mathbf{E}_i[\theta(\{i\}, \mathbf{x}) \mid \mathbf{x}] \quad \stackrel{(6.29)}{=} \quad \frac{1}{\sum_{i=1}^{n} M_{ii}} \sum_{i=1}^{n} \frac{M_{ii}(\nabla_i f(\mathbf{x}))^2}{M_{ii}\|\nabla f(\mathbf{x})\|^2} \quad \geq \quad \frac{1}{\sum_{i=1}^{n} M_{ii}} \qquad (6.32)$$

  - **Greedy choice:** Suppose we choose the coordinate $i$ deterministically as

$$i \quad = \quad \arg\max_j \left\{ \frac{(\nabla_i f(\mathbf{x}))^2}{M_{ii}} \right\}. \qquad (6.33)$$

    It is straightforward that this strategy maximizes the proportion function (6.29) for a single iteration, given that we choose only a single coordinate. It was originally proposed in [49] and further improved in [51]. However, in this case we do not have a better bound than (6.32), which would be independent of $\mathbf{x}^3$. We can get this bound

---

[3]In [51] they proved a slightly better bound using $\ell_1$-strong convexity, which can be achieved by replacing the $\ell_2$-norm by $\ell_\infty$-norm in the definition of the proportion and forcing functions.

using that the maximum of some quantity is more than its average weighted by $M_{ii}$

$$\theta(\{i\}, \mathbf{x}) \quad = \quad \frac{\max_j\{(\nabla_j f(\mathbf{x}))^2/M_{jj}\}}{\|\nabla f(\mathbf{x})\|^2} \quad \geq \quad \frac{1}{\sum_{j=1}^n M_{jj}}. \tag{6.34}$$

Observe that in the case of $\|\nabla f(\mathbf{x})\|^2 = (\nabla_i f(\mathbf{x}))^2$ we have that the above lower bound (6.34) could be $\sum_{j=1}^n M_{jj}/M_{ii}$ larger to still hold. Therefore on some iterations, greedy coordinate descent is much better than coordinate descent with importance sampling, although their global bounds are the same. This usally leads to superiority of greedy rules in practice, in the case that they can be implemented cheaply.

- **Minibatch Coordinate Descent:** Suppose $|S| = \tau$, for some given $n > \tau > 1$. There are currently two analyzed strategies in this case and we introduce a third.

  - $\tau$**-nice sampling:** Assume that we want to sample a subset of $\tau > 1$ coordinates at each iteration, uniformly at random from all subsets of cardinality $\tau$. It can be inferred from results established in [56] that the expectation of the proportion function can be lower bounded by the quantity

    $$\mathbf{E}[\theta(S, \mathbf{x}) \mid \mathbf{x}] \quad \stackrel{(6.25)}{=} \quad \frac{(\nabla f(\mathbf{x}))^\top \mathbf{E}[\mathbf{M}_{[S]}^{-1}]\nabla f(\mathbf{x})}{\|\nabla f(\mathbf{x})\|^2} \quad \geq \quad \lambda_{\min}\left(\mathbf{E}\left[\mathbf{M}_{[S]}^{-1}\right]\right), \tag{6.35}$$

    where $\mathbf{M}_{[S]}^{-1}$ is the $n \times n$ matrix constructed by putting the matrix $\mathbf{M}_S^{-1}$ on the columns and rows specified by $S$ and zero out the rest of the entries. Additionally, assuming that we have a factorization $\mathbf{M} = \mathbf{A}^\top \mathbf{A}$, where $\mathbf{A} \in \mathbb{R}^{m \times n}$, using results from [54], this can be further lower bounded as

    $$\lambda_{\min}\left(\mathbf{E}\left[\mathbf{M}_{[S]}^{-1}\right]\right) \quad \stackrel{[56]}{\geq} \quad \frac{1}{n \cdot \max_i\{v_i\}}, \tag{6.36}$$

    where

    $$v_i \quad \stackrel{[54]}{:=} \quad \sum_{j=1}^m \left[1 + \frac{(\|\mathbf{A}_{j:}\|_0 - 1)(\tau - 1)}{n - 1}\right] A_{ji}^2.$$

  - **Importance sampling for minibatches:** A minibatch version of importance sampling was recently proposed in [8]. The main idea is as follows: we randomly partition the coordinates into $\tau$ approximately equally sized "buckets", and subsequently and independently perform standard importance sampling (as described above) for each bucket. The sampling is then generated as the union of all sampled coordinates. For specific bounds, we recommend discussing the original paper [8], as they are not available in a compact form.

  - **Greedy minibatches:** To showcase the power of our framework, we introduce a brand-new selection rule called *greedy minibatches*. This selection rule aims to select the set $S$ such that it minimizes the proportion function on the current iteration, i.e.,

    $$S^g \quad := \quad \arg\max_{S:|S|=\tau} \left\{(\nabla_S f(\mathbf{x}))^\top (\mathbf{M}_S)^{-1} \nabla_S f(\mathbf{x})\right\}. \tag{6.37}$$

    The above selection procedure is a difficult problem in general, but it might be feasible for some specific problems, e.g., for diagonal $\mathbf{M}$ or the function $f$ with a special structure (see [51] for examples).

    To get a lower bound on the proportion function independent of the current iterate $\mathbf{x}$, we use the argument that the maximum of some quantities is always at least equal to any weighted mean of the same quantities. Using this, we get that $\theta(S^g, \mathbf{x}) \geq \mathbf{E}[\theta(S, \mathbf{x}) \mid \mathbf{x}]$, where the sets $S$ were selected according to the $\tau$-nice sampling, importance sampling for minibatches introduced above, or any other sampling.

Therefore, using (6.35) we get that

$$\theta(S^g, \mathbf{x}) \overset{(6.25)+(6.35)}{\geq} \lambda_{\min}\left(\mathbf{E}\left[\mathbf{M}_{[S]}^{-1}\right]\right). \tag{6.38}$$

Observe that both the selection rule and the above bound is a generalization of the greedy coordinate sampling to minibatches. Also observe that the above bound is potentially very loose. As an example, consider a diagonal matrix $\mathbf{M}$ with all the elements equal to 1. The right-hand side of (6.38) is then equal to $\tau/n$, while the left-hand side is equal to 1. Even in this very special case, the bound is disregarding a factor of $n/\tau$, which is potentially huge. For this reason, it is expected that the greedy minibatches will actually perform much better in practice than in theory (6.38).

### 6.4.2 Non-smooth problems

Let us define for all $i$ and all $\mathbf{x} \in \mathbb{R}^n$ the function

$$\lambda_i(\mathbf{x}) := -L \cdot \min_{v \in \mathbb{R}} \left\{ \nabla_i f(\mathbf{x}) v + \frac{L}{2} v^2 + g_i(x_i + v) - g_i(x_i) \right\}. \tag{6.39}$$

Using the definition of proportion function (6.24) with the diagonal $\mathbf{M} := L\mathbf{I}$ we get

$$\theta(S, \mathbf{x}) \overset{(6.24)+(6.39)}{=} \frac{\sum_{i \in S} \lambda_i(\mathbf{x})}{L \sum_{j=1}^{n} \lambda_j(\mathbf{x})}. \tag{6.40}$$

As mentioned in Section 6.3.1, the value $L$ can be chosen as $L = \lambda_{\max}(\mathbf{M})$ to satisfy the inequality (6.14). However, this choice of $L$ might be suboptimal in many cases. As an example, when analyzing coordinate descent methods, the vector $\mathbf{h}$ will always be of the form $h\mathbf{e}^i$, where $\mathbf{e}^i$ is the $i$-th coordinate vector (see Appendix 6.B.1). Therefore, it would be sufficient for $L$ to satisfy (6.14) for the specific choice $\mathbf{h} = h\mathbf{e}^i$ which leads to

$$f(\mathbf{x} + h\mathbf{e}^i) \leq f(\mathbf{x}) + \nabla_i f(\mathbf{x})h + \frac{M_{ii}}{2}h^2 \leq f(\mathbf{x}) + \nabla_i f(\mathbf{x})h + \frac{L}{2}h^2$$

for each $i \in \{1, \ldots, n\}$ and $h \in \mathbb{R}$. We can easily observe from the above that this will be satisfied with $L := \max_i\{M_{ii}\} \leq \lambda_{\max}(\mathbf{M})$. To account for this in general, we need to take some additional measures. Specifically, let $\mathcal{S}$ be the collection of all sets $S \subset [n]$ which can be possibly generated during the iterative process by $\mathcal{P}_k$. As an example, $\mathcal{S} = \{\{1\}, \{2\}, \ldots, \{n\}\}$ corresponds to all such sets for coordinate descent. We define $L$ as the smallest number for which

$$f(\mathbf{x} + \mathbf{h}_{[S]}) \leq f(\mathbf{x}) + \langle \nabla f(\mathbf{x}), \mathbf{h}_{[S]} \rangle + \frac{1}{2} \langle \mathbf{M}\mathbf{h}_{[S]}, \mathbf{h}_{[S]} \rangle$$

is satisfied for every $\mathbf{x}, \mathbf{h} \in \mathbb{R}^n$ and every $S \in \mathcal{S}$. It is straightforward to see that if the size of possible sets $S \in \mathcal{S}$ is upper bounded as $|S| \leq \tau$, then we can safely choose

$$L = L_\tau := \max_{S:|S|=\tau} \{\lambda_{\max}(\mathbf{M}_S)\} \leq \max_{S:|S|=\tau} \left\{ \sum_{i \in S} M_{ii} \right\}, \tag{6.41}$$

where the last inequality is due to the eigenvalues being positive and their sum being the trace. We can easily verify that this generalizes to $L = L_n = \lambda_{\max}(\mathbf{M})$ for gradient descent and $L = L_1 = \max_i\{M_{ii}\}$ for coordinate descent.

Now, we will breakdown the cases depending on the block selection procedure.

- **Proximal Gradient Descent:** The first result in the proximal setting was the Iterative Shrinkage Tresholding Algorithm (ISTA) for the $\ell_1$-norm, which selects all the coordinates on each iteration. The result was generalized to as general non-smooth penalty in [48].

The bound on the proportion function is given by

$$\theta([n], \mathbf{x}) \stackrel{(6.40)}{=} \frac{\sum_{i=1}^{n} \lambda_i(\mathbf{x})}{L_n \sum_{j=1}^{n} \lambda_j(\mathbf{x})} = \frac{1}{L_n} \stackrel{(6.41)}{=} \frac{1}{\lambda_{\max}(\mathbf{M})}. \tag{6.42}$$

- **Serial Proximal Coordinate descent:** Suppose $|S| = 1$ and specifically let $S = \{i\}$. Then we have that

$$\theta(\{i\}, \mathbf{x}) \stackrel{(6.40)}{=} \frac{\lambda_i(\mathbf{x})}{L_1 \sum_{j=1}^{n} \lambda_j(\mathbf{x})}, \tag{6.43}$$

as the size of the selected sets are upper bounded by 1. The procedure leading to the choice of $i$ distinguishes between various serial approaches.

  - **Uniform probabilities:** Assume we choose the coordinate $i$ uniformly at random at each iteration. Then we can bound the expectation of the proportion function as

$$\mathbf{E}[\theta(\{i\}, \mathbf{x}) \mid \mathbf{x}] \stackrel{(6.43)}{=} \mathbf{E}\left[\frac{\lambda_i(\mathbf{x})}{L_1 \sum_{j=1}^{n} \lambda_j(\mathbf{x})}\right] = \frac{\frac{1}{n}\sum_{i=1}^{n} \lambda_i(\mathbf{x})}{L_1 \sum_{j=1}^{n} \lambda_j(\mathbf{x})} = \frac{1}{nL_1}$$

$$\stackrel{(6.41)}{=} \frac{1}{n \max_i \{M_{ii}\}}. \tag{6.44}$$

  - **Greedy choice:** Another approach is to pick the coordinate $i$, which maximizes the proportion function in (6.43), which was analyzed in [51]. The best bound independent on $\mathbf{x}$ coincides with the above bound for uniform sampling, using the fact that a mean of some quantity is less than its maximum

$$\theta(\{i\}, \mathbf{x}) \stackrel{(6.43)}{=} \frac{\max_i \{\lambda_i(\mathbf{x})\}}{L_1 \sum_{j=1}^{n} \lambda_j(\mathbf{x})} \geq \frac{\frac{1}{n}\sum_{i=1}^{n} \lambda_i(\mathbf{x})}{L_1 \sum_{j=1}^{n} \lambda_j(\mathbf{x})} = \frac{1}{nL_1}$$

$$\stackrel{(6.41)}{=} \frac{1}{n \max_i \{M_{ii}\}}. \tag{6.45}$$

  Similarly as in the smooth case, observe that the quantity $\max_i\{\lambda_i(\mathbf{x})\}$ is potentially up to $n$ times larger than $\frac{1}{n}\sum_{i=1}^{n} \lambda_i(\mathbf{x})$ and therefore the bound (6.45) can be potentially $n$ times larger in some cases, which usually results in better empirical results.

- **Minibatch Proximal Coordinate Descent:** Suppose $|S| = \tau$, where $\tau$ is given, which implies that we can use $L = L_\tau$ in the bounds. We introduce two options for the block selection procedure.

  - **$\tau$-nice sampling:** Only one variant of a sampling for this setup was considered before [60], and that is a uniform choice of $\tau$ coordinates without repetition. In expectation, each coordinate has a chance of $\tau/n$ to be picked, which is used in the bound to get

$$\mathbf{E}[\theta(S, \mathbf{x}) \mid \mathbf{x}] \stackrel{(6.40)}{=} \mathbf{E}\left[\frac{\sum_{i \in S} \lambda_i(\mathbf{x})}{L_\tau \sum_{j=1}^{n} \lambda_j(\mathbf{x})}\right] = \frac{\frac{\tau}{n}\sum_{i=1}^{n} \lambda_i(\mathbf{x})}{L_\tau \sum_{j=1}^{n} \lambda_j(\mathbf{x})} = \frac{\tau}{nL_\tau}$$

$$\stackrel{(6.41)}{=} \frac{\tau}{n \max_{S:|S|=\tau}\{\lambda_{\max}(\mathbf{M}_S)\}}. \tag{6.46}$$

  To our best knowledge, this bound is new, as the previous bound considered $L_n$ instead of $L_\tau$. As $L_\tau \leq L_n$, the new bound is better.

  - **Greedy minibatches:** Similarly as in the smooth case, we introduce a new selection rule – greedy minibatches. Specifically, the corresponding set $S^g$ is given by

$$S^g := \arg\max_{S:|S|=\tau} \left\{\sum_{i \in S} \lambda_i(\mathbf{x})\right\}. \tag{6.47}$$

Note that for $\tau = 1$ we recover the greedy coordinate descent. To give global bounds independent of $\mathbf{x}$ for this strategy, we again use the fact that maximum is an upper bound for the mean, to get

$$\theta(S^g, \mathbf{x}) \stackrel{(6.40)}{=} \frac{\sum_{i \in S^g} \lambda_i(\mathbf{x})}{L_\tau \sum_{j=1}^n \lambda_j(\mathbf{x})} \geq \frac{\frac{\tau}{n} \sum_{i=1}^n \lambda_i(\mathbf{x})}{L_\tau \sum_{j=1}^n \lambda_j(\mathbf{x})} = \frac{\tau}{nL_\tau}$$

$$\stackrel{(6.41)}{=} \frac{\tau}{n \max_{S:|S|=\tau}\{\lambda_{\max}(\mathbf{M}_S)\}}. \tag{6.48}$$

Note that the the above bound can be potentially very pessimistic, as

$$1 \geq \frac{\sum_{i=1}^n \lambda_i(\mathbf{x})}{\sum_{i \in S^g} \lambda_i(\mathbf{x})} \geq \frac{\tau}{n}$$

by using $\sum_{i \in S^g} \lambda_i(\mathbf{x}) \geq \frac{\tau}{n} \sum_{i=1}^n \lambda_i(\mathbf{x}) \geq \frac{\tau}{n} \sum_{i \in S^g} \lambda_i(\mathbf{x})$. Therefore, the bound (6.48) can be up to $n/\tau$ times better in certain cases.

## 6.5 Strongly Polyak-Łojasiewicz Functions

In this section, we reinvent the strongly PL functions using the proximal forcing function (6.21), and develop the corresponding convergence rates. Also, we show how to recover the known results in this setting by applying our theory.

### 6.5.1 Strongly PL functions

**Definition 6.14** (Strongly PL functions: composite case). We say that $F$ is a strongly PL function there exists a scalar $\mu > 0$ satisfying

$$\mu(\mathbf{x}) \geq \mu \tag{6.49}$$

for all $\mathbf{x} \in \text{dom}(g)/\mathcal{X}^*$. The collection of all functions $F$ satisfying inequality (6.49) will be denoted $\mathcal{S}_{PL}^g(\mu)$, and we say that $F$ is strongly PL with parameter $\mu$.

Recall that in the smooth case we said that a function $f \in \mathcal{S}_{PL}(\mu)$, if $f$ satisfied the condition (6.2), which can be observed to be equivalent to (6.49) for smooth functions. Therefore we have that $\mathcal{S}_{PL}(\mu) \subset \mathcal{S}_{PL}^g(\mu)$.

The above definition is not new, it was originally introduced in a slightly different form by Karimi et al. [29].

### 6.5.2 Strongly convex functions are strongly PL

Let $\lambda \geq 0$. Function $F : \mathbb{R}^n \to \bar{\mathbb{R}}$ is said to be $\lambda$-*strongly convex*, if for all $\mathbf{x}, \mathbf{y} \in \mathbb{R}^n$ and all $\beta \in [0, 1]$ we have

$$F(\beta \mathbf{x} + (1-\beta)\mathbf{y}) \leq \beta F(\mathbf{x}) + (1-\beta) F(\mathbf{y}) - \frac{\lambda \beta(1-\beta)}{2} \|\mathbf{x} - \mathbf{y}\|^2. \tag{6.50}$$

If $F$ is 0-strongly convex, we refer to it simply as *convex*. Consider a differentiable function $f : \mathbb{R}^n \to \mathbb{R}$. If $f$ is $\lambda$-strongly convex for $\lambda \geq 0$, then

$$f(\mathbf{x} + \mathbf{h}) \geq f(\mathbf{x}) + \langle \nabla f(\mathbf{x}), \mathbf{h} \rangle + \frac{\lambda}{2} \|\mathbf{h}\|^2, \qquad \mathbf{x}, \mathbf{h} \in \mathbb{R}^n. \tag{6.51}$$

We will now show that if $F$ is strongly convex, then $F \in \mathcal{S}_{PL}^g(\rho)$ for some specific $\rho$. This means that the class of strongly PL (composite) functions contains the class of strongly convex (composite) functions.

**Theorem 6.15.** *Assume $F$ is $\lambda_F$-strongly convex with $\lambda_F > 0$, and $f$ is $\lambda_f$-strongly convex with $\lambda_f \geq 0$. Then*

$$\mu(\mathbf{x}) \geq \mu := \min\left\{\frac{L}{2}, \frac{L\lambda_F}{\lambda_F - \lambda_f + L}\right\} \tag{6.52}$$

*for all $\mathbf{x} \in \mathrm{dom}(g)/\mathcal{X}^*$, and hence $F \in \mathcal{S}_{PL}^g(\mu)$.*

*Proof.* See Section 6.B.2. □

In the above theorem we do not enforce separability assumption on $g$.

### 6.5.3 Convergence

We have the following convergence result, establishing convergence to a global minimizer.

**Theorem 6.16.** *Invoke Assumption 6.10. Further, assume $F \in \mathcal{S}_{PL}^g(\mu)$ (i.e., $F$ is strongly PL with parameter $\mu$), and let*

$$\mu_k := \frac{\mu \mathbf{E}[\xi(\mathbf{x}^k) \cdot \theta(S_k, \mathbf{x}^k)]}{\mathbf{E}[\xi(\mathbf{x}^k)]}. \tag{6.53}$$

*Then*

$$\sum_{k=0}^{K-1} \mu_k \geq \log\left(\frac{\xi(\mathbf{x}^0)}{\epsilon}\right) \quad \Rightarrow \quad \mathbf{E}[\xi(\mathbf{x}^K)] \leq \epsilon. \tag{6.54}$$

*Proof.* See Section 6.B.3. □

In order to get concrete complexity results from the above theorem, we need to estimate the speed of growth of $\sum_{k=0}^{K-1} \mu_k$ in $K$. There no universal way to do this, which is why we state the above result the way we do. Instead, in each situation this needs to be estimated separately. Typically, this will be done by lower bounding $\mu_k$ for each $k$ separately. Let us illustrate this using a couple examples.

If the blocks are selected *deterministically*, then the expectations in (6.53) do not play any role, and we have $\mu_k = \mu \theta(S_k, \mathbf{x}^k)$. If we have a global lower bound of the form

$$\theta(S_k, \mathbf{x}^k) \geq c > 0$$

readily available, then $\sum_{k=0}^{K-1} \mu_k \geq \mu c K$, which implies the rate

$$K \geq \frac{1}{c\mu} \log\left(\frac{\xi(\mathbf{x}^0)}{\epsilon}\right) \quad \Rightarrow \quad \xi(\mathbf{x}^K) \leq \epsilon \tag{6.55}$$

or, equivalently, $\xi(\mathbf{x}^K) \leq \xi(\mathbf{x}^0) \cdot e^{-c\mu K}$. More generally, convergence is established whenever we can lower bound $\theta(S_k, \mathbf{x}^k) \geq c_k > 0$, where the constants $c_k$ sum up to infinity.

If the blocks are selected *stochastically*, then the sequence of iterates $\mathbf{x}^k$ is also stochastic. In such cases, it is often possible to come up with a bound for the expectation of $\theta(S_k, \mathbf{x}^k)$ conditioned on $\mathbf{x}^k$:

$$\mathbf{E}[\theta(S_k, \mathbf{x}^k) \mid \mathbf{x}^k] \geq c > 0. \tag{6.56}$$

If this is the case, we claim that $\mu_k$ can be lower bounded by $\mu c$, and hence Theorem 6.16 implies the same rate as before: given by (6.55), only with a bound on the expectation $\mathbf{E}[\xi(\mathbf{x}^k)] \leq \epsilon$ on the right-hand side. More generally, if $\mathbf{E}[\theta(S_k, \mathbf{x}^k) \mid \mathbf{x}^k] \geq c_k \geq 0$, then $\mu_k \geq \mu c_k$, and convergence is guaranteed as long as $\sum_k c_k = \infty$.

Let us now return to the claim. Indeed, the bound $\mu_k \geq \mu c$ follows from (6.56) by applying Lemma 6.25 (see Appendix) with $X := \mathbf{x}^k, Y := \theta_k$ and $f(\mathbf{x}) := \xi(\mathbf{x})$.

### 6.5.4 Applications of Theorem 6.16

We now showcase the use of Theorem 6.16 on selected algorithms which arise as special cases of our generic method (Algorithm 6.18).

**Proximal gradient descent.** If for all $k$ we choose $S_k = [n]$ with probability 1, Algorithm 6.18 reduces to (proximal) gradient descent. In both smooth and non-smooth case we have $\theta([n], \mathbf{x}^k) \geq 1/\lambda_{\max}(\mathbf{M})$ (see (6.28) and (6.42)). Substituting into (6.55), we get the rate

$$K \geq \frac{L}{\mu} \log\left(\frac{\xi(\mathbf{x}^0)}{\epsilon}\right) \quad \Rightarrow \quad \xi(\mathbf{x}^K) \leq \epsilon,$$

in both the smooth and non-smooth case. These results were previously obtained for strongly PL functions in [29].

**Randomized coordinate descent with uniform probabilities.** Randomized coordinate descent, analyzed in [49, 59], arises as special case of Algorithm 6.18 by choosing $S_k = \{i_k\}$, where $i_k$ is an index chosen from $[n]$ uniformly at random, and independently of the history of the method.

Note that a bound on $\mathbf{E}[\theta(S_k, \mathbf{x}^k) \mid \mathbf{x}^k]$ is readily available in Section 6.4, specifically in (6.30) for the smooth case and (6.44) in the non-smooth case, and it takes the form

$$\mathbf{E}[\theta(S_k, \mathbf{x}^k) \mid \mathbf{x}^k] \quad \geq \quad c \quad := \quad \frac{1}{n \max_i\{M_{ii}\}}.$$

As we have seen in the discussion immediately following Theorem 6.16, this implies the bound $\mu_k \geq \mu c$. Applying Theorem 6.16, we conclude that

$$K \geq \frac{n \max\{M_{ii}\}}{\mu} \log\left(\frac{\xi(\mathbf{x}^0)}{\epsilon}\right) \quad \Rightarrow \quad \mathbf{E}[\xi(\mathbf{x}^K)] \leq \epsilon,$$

which is the same bound as given in [29] for strongly PL functions in both smooth and non-smooth case.

**Randomized coordinate descent with importance sampling** We now allow for specific nonuniform probabilities: probability of choosing $S_k = \{i\}$ is proportional to $M_{ii}$ (see (6.31)). In view of (6.32), we get $\mathbf{E}[\theta(S_k, \mathbf{x}^k) \mid \mathbf{x}^k] \geq c := 1/\sum_{i=1}^n M_{ii}$ for smooth functions. This leads to the complexity result

$$K \geq \frac{\sum_{i=1}^n M_{ii}}{\mu} \log\left(\frac{\xi(\mathbf{x}^0)}{\epsilon}\right) \quad \Rightarrow \quad \mathbf{E}[\xi(\mathbf{x}^K)] \leq \epsilon,$$

which is *is a new result* for strongly PL functions. However, in the special case of strongly-convex functions, this result is known [49, 59, 57].

**Minibatch coordinate descent.** Assume $S_k$ is a subset of $[n]$ of cardinality $\tau$, chosen uniformly at random. This leads to (strandard) minibatch coordinate descent. In view of (6.35), we have the bound $\mathbf{E}[\theta(S, \mathbf{x}) \mid \mathbf{x}] \geq \lambda_{\min}\left(\mathbf{E}\left[\mathbf{M}_{[S]}^{-1}\right]\right)$ in the smooth case. Substituting into (6.55), we get the rate

$$K \geq \frac{1}{\mu \lambda_{\min}(\mathbf{E}[\mathbf{M}_S^{-1}])} \log\left(\frac{\xi(\mathbf{x}^0)}{\epsilon}\right) \quad \Rightarrow \quad \mathbf{E}[\xi(\mathbf{x}^K)] \leq \epsilon.$$

This result is novel for the strongly PL case, but it was already established before for strongly convex functions in [56].

In the non-smooth case we can use the bound (6.46) given by $\mathbf{E}[\theta(S, \mathbf{x}) \mid \mathbf{x}] \geq \tau/nL_\tau$, where $L_\tau$ is given by (6.41). It follows that we have the guarantee

$$K \geq \frac{nL_\tau}{\mu\tau} \log\left(\frac{\xi(\mathbf{x}^0)}{\epsilon}\right) \quad \Rightarrow \quad \mathbf{E}[\xi(\mathbf{x}^K)] \leq \epsilon,$$

which is novel for strongly PL and also for strongly convex functions, as $L_\tau$ is less or equal to the standard $L_n$ [60].

**Greedy coordinate descent.** Assume a single coordinate $S = \{i\}$ is chosen using the rule stated in (6.33), i.e., the choice maximizes the proportion function on the given iteration. In the smooth case we have the bound (6.34), i.e., $\theta(S_k, \mathbf{x}^k) \geq c := 1/\sum_{i=1}^{n} M_{ii}$. This leads to the rate

$$K \geq \frac{\sum_{i=1}^{n} M_{ii}}{\mu} \log\left(\frac{\xi(\mathbf{x}^0)}{\epsilon}\right) \quad \Rightarrow \quad \xi(\mathbf{x}^K) \leq \epsilon.$$

Note that this is identical to the rate of randomized coordinate descent with importance sampling, with the exception that we have $\xi(\mathbf{x}^K)$ instead of $\mathbf{E}[\xi(\mathbf{x}^K)]$ due to the deterministic nature of the method. Again, this result was already established in [29]. In the special case of strongly convex functions, this was first established by Nesterov [49].

As for the proximal case, we have a bound (6.45) given by $\theta(S, \mathbf{x}^k) \geq c := 1/n \max_i\{M_{ii}\}$ which leads to a similar rate

$$K \geq \frac{n \max_i\{M_{ii}\}}{\mu} \log\left(\frac{\xi(\mathbf{x}^0)}{\epsilon}\right) \quad \Rightarrow \quad \xi(\mathbf{x}^K) \leq \epsilon.$$

This rate is identical to the rate of randomized coordinate descent, except for the expectation on the right-hand side. This result was also established in [29], while the strongly convex version is due [51].

**Greedy minibatches.** Assume the smooth and and that we choose the set of coordinates $S$ according to the rule described in (6.37), which is

$$S := \arg\max_{S:|S|=\tau} \left\{(\nabla_S f(\mathbf{x}))^\top (\mathbf{M}_S)^{-1} \nabla_S f(\mathbf{x})\right\}.$$

Assuming that $|S| = \tau$, we can see that this rule maximizes the proportion function for the given iteration. The corresponding lower bound takes the form (6.38) which leads to the rate

$$K \geq \frac{1}{\lambda_{\min}(\mathbf{E}[\mathbf{M}_{[S]}^{-1}])\mu} \log\left(\frac{\xi(\mathbf{x}^0)}{\epsilon}\right) \quad \Rightarrow \quad \xi(\mathbf{x}^K) \leq \epsilon.$$

Similarly as in the serial case, this bound is identical to uniform minibatches, with the only exception being the dropped expectation on the optimality gap $\xi(\mathbf{x})$, due to this algorithm being deterministic.

In the non-smooth case we have the bound on the proportion function given by (6.48) which leads to the rate

$$K \geq \frac{nL_\tau}{\mu\tau} \log\left(\frac{\xi(\mathbf{x}^0)}{\epsilon}\right) \quad \Rightarrow \quad \xi(\mathbf{x}^K) \leq \epsilon.$$

Both of the above bounds are novel, as greedy minibatches is a novel sampling approach.

## 6.6 Weakly Polyak-Łojasiewicz Functions

In this section we introduce a generalized definition of Weakly PL functions, which specializes to Definition 6.2 for a specific choice of parameters, and also covers the proximal case. We show that proximal convex functions can be analyzed using this framework and we give a general convergence rate guarantee for this class. Lastly, we specify our theory to several known setups, showcasing the generality of our definition.

### 6.6.1 Weakly PL functions

**Definition 6.17** (Weakly PL functions: general case). We say that $F$ is a weakly PL function, if there exists a scalar function $\rho : \mathbb{R}^n \to \mathbb{R}^+$, such that

$$\mu(\mathbf{x}) \quad \geq \quad \rho(\mathbf{x}^0) \cdot \xi(\mathbf{x}) \quad > \quad 0, \tag{6.57}$$

for all $\mathbf{x}^0, \mathbf{x} \in \mathrm{dom}(g)/\mathcal{X}^*$ such that $\xi(\mathbf{x}) \leq \xi(\mathbf{x}^0)$. The collection of all functions $F$ satisfying inequality (6.57) will be denoted $\mathcal{W}_{PL}^g(\rho)$, and we say that $F$ is weakly PL with parameter $\rho$.

Recall that in the smooth case we said that a function $f \in \mathcal{W}_{PL}(\mu)$, if $f$ satisfied the condition (6.11)

$$\mu(\mathbf{x}) \geq \frac{\mu \xi(\mathbf{x})}{2\|\mathbf{x} - \mathbf{x}^*\|^2}, \qquad \mathbf{x} \in \mathbb{R}^n / \mathcal{X}^*.$$

In the proof of convergence of these methods, we further bounded the right-hand side of the above expression by $\mu \xi(\mathbf{x})/2R^2$, where we defined $R := \max_{\mathbf{x} \in \mathcal{L}_0} \|\mathbf{x} - \mathbf{x}^*\| < +\infty$, with $\mathcal{L}_0 := \{\mathbf{x} \in \mathbb{R}^n : f(\mathbf{x}) \leq f(\mathbf{x}^0)\}$. To see the above smooth case in our general framework defined in (6.57), we can use $\rho(\mathbf{x}^0) = \frac{1}{2R^2}$, which satisfies our assumption, as $R$ does not depend on $\mathbf{x}$.

### 6.6.2 Weakly convex functions are weakly PL

We will now show that if $F$ is convex, then $F \in \mathcal{W}^g_{PL}(\rho)$ for some specific $\rho$. This means that weakly PL functions generalize weakly convex functions also in the composite setting.

**Theorem 6.18.** *Assume $f$ and $g$ are convex, and let $\mathbf{x}^*$ be a global minimizer of $F$. Then*

$$\mu(\mathbf{x}) \geq \xi(\mathbf{x}) \cdot \min\left\{\frac{1}{2\xi(\mathbf{x})}, \frac{1}{2L\|\mathbf{x} - \mathbf{x}^*\|^2}\right\} \tag{6.58}$$

*for all $\mathbf{x} \in \mathrm{dom}(g)/\mathcal{X}^*$. Also, $F$ is a weakly PL function with the parameter $\rho$ given by*

$$\rho(\mathbf{x}^0) = \min\left\{\frac{L}{2\xi(\mathbf{x}^0)}, \frac{1}{2R^2}\right\}, \tag{6.59}$$

*where*

$$R := \max_{\mathbf{x} \in \mathbb{R}^n : f(\mathbf{x}) \leq f(\mathbf{x}^0)} \|\mathbf{x} - \mathbf{x}^*\| < +\infty. \tag{6.60}$$

*Proof.* See Section 6.B.4. □

Note that in the smooth case we have $\xi(\mathbf{x}) \leq \frac{L}{2}\|\mathbf{x} - \mathbf{x}^*\|^2$, therefore

$$\min\left\{\frac{1}{2\xi(\mathbf{x})}, \frac{1}{2L\|\mathbf{x} - \mathbf{x}^*\|^2}\right\} = \frac{1}{2L\|\mathbf{x} - \mathbf{x}^*\|^2},$$

which also leads to $\rho(\mathbf{x}^0) = 1/(2R^2)$ in the smooth case.

### 6.6.3 Convergence

For this class, we can show the following convergence result.

**Theorem 6.19.** *Invoke Assumption 6.10. Further, assume $F \in \mathcal{W}^g_{PL}(\rho)$ (i.e., $F$ is weakly PL with parameter $\rho$). Let*

$$\mu_k := \rho(\mathbf{x}^0) \frac{\mathbf{E}[(\xi(\mathbf{x}^k))^2 \cdot \theta(S_k, \mathbf{x}^k)]}{(\mathbf{E}[\xi(\mathbf{x}^k)])^2}. \tag{6.61}$$

*Then*

$$\sum_{k=0}^{K-1} \mu_k \geq \frac{1}{\epsilon} \quad \Rightarrow \quad \mathbf{E}[\xi(\mathbf{x}^K)] \leq \epsilon. \tag{6.62}$$

*Proof.* See Section 6.B.5. □

Similarly as in the previous section, we need to bound the quantity $\sum_{k=0}^{K} \mu_k$ in $K$ to get a complexity result. The standard theory is developed by bounding each of $\mu_k$ separately, which leads to a simpler looking, but less tight bounds.

In the case that the blocks are selected deterministically, then the expectations in (6.61) do not play any role and we have $\mu_k = \rho(\mathbf{x}^0)\theta(S_k, \mathbf{x}^k)$. Additionally, if we have a global lower bound $\theta(S_k, \mathbf{x}^k) \geq c > 0$, then $\sum_{k=0}^{K-1} \mu_k \geq \rho(\mathbf{x}^0)cK$, which implies the rate

$$K \geq \frac{1}{\rho(\mathbf{x}^0)c\epsilon} \quad \Rightarrow \quad \xi(\mathbf{x}^K) \leq \epsilon$$

and equivalently $\xi(\mathbf{x}^K) \leq \frac{1}{\rho(\mathbf{x}^0)cK}$.

We get a similar result for stochastic block selection. Specifically, if we have a global bound for the quantity $\mathbf{E}[\theta(S_k, \mathbf{x}^k) \mid \mathbf{x}^k] \geq c > 0$, we can lower bound the values of $\mu_k$ in Theorem 6.19 by $\rho(\mathbf{x}^0)c$. This can be done by first lower bounding the denominator of (6.61) by $\mathbf{E}[(\xi(\mathbf{x}^k))^2]$ using the trivial bound $\mathbf{E}[X]^2 \leq \mathbf{E}[X^2]$ and further applying Lemma 6.25 with $X := \mathbf{x}^k, Y := \theta_k, f(\mathbf{x}) = (\xi(\mathbf{x}))^2$.

In general, we can claim convergence if we have a sequence of lower bounds $\{c_k\}_{k=0}^\infty$ for $\theta(S_k, \mathbf{x}^k) \geq c_k$ in the deterministic case or $\mathbf{E}[\theta(S_k, \mathbf{x}^k) \mid \mathbf{x}^k] \geq c_k$ in the stochastic case, and additionally $\sum_{k=0}^\infty c_k = \infty$.

### 6.6.4 Applications of Theorem 6.19

We now showcase the use of Theorem 6.19 on several methods which arise as special cases of our generic method (Algorithm 6.18). All the results for general $\mathcal{W}_{PL}^g$ functions are novel, as the notion itself is novel. In most cases, we recover known theory by specializing the results to convex objectives.

**Proximal gradient descent.** In both smooth and non-smooth case we have $\theta([n], \mathbf{x}^k) \geq c := 1/\lambda_{\max}(\mathbf{M})$ (see (6.28) and (6.42)). Substituting into (6.62), we get the rate

$$K \geq \frac{\lambda_{\max}(\mathbf{M})}{\rho(\mathbf{x}^0)\epsilon} \quad \Rightarrow \quad \xi(\mathbf{x}^K) \leq \epsilon,$$

in both the smooth and non-smooth case. This result is novel for Weakly PL functions, but it was established before for convex functions in [50].

**Randomized coordinate descent with uniform probabilities.** The bound on the quantity $\mathbf{E}[\theta(S_k, \mathbf{x}^k) \mid \mathbf{x}^k)]$ is readily available in Section 6.4, specifically in (6.30) for the smooth case and (6.44) in the non-smooth case, and it takes the form $\mathbf{E}[\theta(S_k, \mathbf{x}^k) \mid \mathbf{x}^k] \geq c := \frac{1}{n \max_i \{M_{ii}\}}$. As discussed right after Theorem 6.19, this implies the bound $\mu_k \geq \mu c$. Applying Theorem 6.19, we conclude that

$$K \geq \frac{n \max\{M_{ii}\}}{\rho(\mathbf{x}^0)\epsilon} \quad \Rightarrow \quad \mathbf{E}[\xi(\mathbf{x}^K)] \leq \epsilon,$$

which is novel for weakly PL functions, but it is well known for convex objectives [49].

**Randomized coordinate descent with importance sampling** In view of (6.32), we get the lower bound $\mathbf{E}[\theta(S_k, \mathbf{x}^k) \mid \mathbf{x}^k)] \geq c := 1/\sum_{i=1}^n M_{ii}$ for smooth functions. This leads to the complexity result

$$K \geq \frac{\sum_{i=1}^n M_{ii}}{\rho(\mathbf{x}^0)\epsilon} \quad \Rightarrow \quad \mathbf{E}[\xi(\mathbf{x}^K)] \leq \epsilon,$$

which is a novel result for weakly PL functions. However, in the special case of convex functions, this result is known [49].

**Minibatch coordinate descent.** In view of (6.35), we have the lower bound $\mathbf{E}[\theta(S, \mathbf{x}) \mid \mathbf{x}] \geq \lambda_{\min}\left(\mathbf{E}\left[\mathbf{M}_{[S]}^{-1}\right]\right)$ in the smooth case. Substituting into (6.62), we get the rate

$$K \geq \frac{1}{\rho(\mathbf{x}^0)\lambda_{\min}(\mathbf{E}[\mathbf{M}_S^{-1}])\epsilon} \quad \Rightarrow \quad \mathbf{E}[\xi(\mathbf{x}^K)] \leq \epsilon.$$

This result is novel for the weakly PL case and to the best of our knowledge also for the convex case, as the results in [56] only consider strongly convex objectives.

In the non-smooth case we can use the bound (6.46) given by $\mathbf{E}[\theta(S, \mathbf{x}) \mid \mathbf{x}] \geq \tau/nL_\tau$, where $L_\tau$ is given by (6.41). It follows that we have the guarantee

$$K \geq \frac{nL_\tau}{\rho(\mathbf{x}^0)\tau\epsilon} \quad \Rightarrow \quad \mathbf{E}[\xi(\mathbf{x}^K)] \leq \epsilon,$$

which is novel for weakly PL and to the best of our knowledge also for weakly convex functions, as $L_\tau$ is less or equal to the standard $L_n$ [60].

**Greedy coordinate descent.** In the smooth case we have the bound (6.34), i.e., $\theta(S_k, \mathbf{x}^k) \geq c := 1/\sum_{i=1}^n M_{ii}$. This leads to the rate

$$K \geq \frac{\sum_{i=1}^n M_{ii}}{\rho(\mathbf{x}^0)\epsilon} \quad \Rightarrow \quad \xi(\mathbf{x}^K) \leq \epsilon.$$

This is identical to the above rate of randomized coordinate descent with importance sampling, with a dropped expectation on $\xi(\mathbf{x}^0)$. This result is new for weakly PL functions and in the special case of weakly convex functions, it was first established by Nesterov [49].

In the proximal case, we have a bound (6.45) given by $\theta(S, \mathbf{x}^k) \geq c := 1/n \max_i\{M_{ii}\}$ which leads to the rate

$$K \geq \frac{n \max_i\{M_{ii}\}}{\rho(\mathbf{x}^0)\epsilon} \quad \Rightarrow \quad \xi(\mathbf{x}^K) \leq \epsilon.$$

This rate is identical to the rate of randomized coordinate descent, except for the expectation on the right-hand side. This result is novel for weakly PL and to the best of our knowledge it is also novel for the special case of convex functions.

**Greedy minibatches.** In the smooth case we have the bound (6.38) which leads to the rate

$$K \geq \frac{1}{\rho(\mathbf{x}^0)\lambda_{\min}(\mathbf{E}[\mathbf{M}_S^{-1}])\epsilon} \quad \Rightarrow \quad \mathbf{E}[\xi(\mathbf{x}^K)] \leq \epsilon.$$

In the non-smooth case we have the bound on the proportion function given by (6.48) which leads to the rate

$$K \geq \frac{nL_\tau}{\rho(\mathbf{x}^0)\tau\epsilon} \quad \Rightarrow \quad \mathbf{E}[\xi(\mathbf{x}^K)] \leq \epsilon.$$

Both of the above results are shared with uniform minibatch coordinate descent, except that we are bounding the optimality gap itself instead of its expectation.

Both of the above rates are novel, as greedy minibatches was introduced as a novel sampling approach.

## 6.7 General Nonconvex Functions

In this section we establish a generic convergence result applicable to general nonconvex functions. This is done at the expense of losing global optimality: we will show that either $\lambda(\mathbf{x}^k)$ gets small, or that $F(\mathbf{x}^k)$ is close to the global minimum $F(\mathbf{x}^*)$. Recall that in the smooth case ($g = 0$) we have $\lambda(\mathbf{x}^k) = \frac{1}{2}\|\nabla f(\mathbf{x}^k)\|^2$.

**Theorem 6.20.** *Invoke Assumption 6.10. Let $\epsilon > 0$ be fixed. Further, let*

$$\mu_k \quad := \quad \frac{\mathbf{E}[\xi(\mathbf{x}^k)\theta(S_k, \mathbf{x}^k)]}{\mathbf{E}[\xi(\mathbf{x}^k)]}. \tag{6.63}$$

*If the inequality*

$$\frac{\xi(\mathbf{x}^0)}{\epsilon} \log\left(\frac{\xi(\mathbf{x}^0)}{\epsilon}\right) \quad \leq \quad \sum_{k=0}^{K-1} \mu_k, \tag{6.64}$$

*holds, then at least one of the following conclusions holds:*

*(i) $\lambda(\mathbf{x}^k) < \epsilon$ for at least one $k \in \{0, \ldots, K-1\}$,*

*(ii) $\mathbf{E}[\xi(\mathbf{x}^K)] \leq \epsilon$.*

*Proof.* See Section 6.B.6. □

In the smooth case ($g = 0$), condition ($i$) reduces to $\frac{1}{2}\|\nabla f(\mathbf{x}^k)\|^2 \leq \epsilon$ for at least one $k \in \{0, \ldots, K-1\}$. Also note that the same strategies for bounding $\mu_k$ as those outlined in 6.5.3 apply here. To sum it up, if we have a global bound $\theta(S_k, \mathbf{x}^k) \geq c > 0$ in the deterministic case, then it follows from Theorem 6.20 that

$$K \geq \frac{\xi(\mathbf{x}^0)}{c\epsilon} \log\left(\frac{\xi(\mathbf{x}^0)}{\epsilon}\right) \quad \Rightarrow \quad \left((\xi(\mathbf{x}^K) \leq \epsilon) \vee (\exists k \in [K] : \lambda(\mathbf{x}^k) \leq \epsilon)\right).$$

Similarly, if we have a bound $\mathbf{E}[\theta(S_k, \mathbf{x}^k) \mid \mathbf{x}^k] \geq c > 0$ in the stochastic case, we get the same result as above with an expectation over the optimality gap $\mathbf{E}[\xi(\mathbf{x}^K)]$, as the optimality gap becomes a random variable.

### 6.7.1 Applications of Theorem 6.20

In this part we apply the results from Theorem 6.20 to several known methods to acquire local convergence guarantees. To the best of our knowledge, all the results in this section are novel.

**Proximal gradient descent.** In both smooth and non-smooth case we have $\theta([n], \mathbf{x}^k) \geq c := 1/\lambda_{\max}(\mathbf{M})$ (see (6.28) and (6.42)). Substituting into (6.64), we get the rate

$$K \geq \frac{\lambda_{\max}(\mathbf{M})\xi(\mathbf{x}^0)}{\epsilon} \log\left(\frac{\xi(\mathbf{x}^0)}{\epsilon}\right) \quad \Rightarrow \quad \left((\xi(\mathbf{x}^K) \leq \epsilon) \vee (\exists k \in [K] : \lambda(\mathbf{x}^k) \leq \epsilon)\right)$$

in both the smooth and non-smooth case.

**Randomized coordinate descent with uniform probabilities.** In both smooth and non-smooth case we have the bound $\mathbf{E}[\theta(S_k, \mathbf{x}^k) \mid \mathbf{x}^k] \geq c := \frac{1}{n \max_i \{M_{ii}\}}$ (see (6.30) and (6.44)). This implies the rate

$$K \geq \frac{\max_i\{M_{ii}\}\xi(\mathbf{x}^0)}{\epsilon} \log\left(\frac{\xi(\mathbf{x}^0)}{\epsilon}\right) \quad \Rightarrow \quad \left((\mathbf{E}[\xi(\mathbf{x}^K)] \leq \epsilon) \vee (\exists k \in [K] : \lambda(\mathbf{x}^k) \leq \epsilon)\right)$$

for both smooth and non-smooth functions.

**Randomized coordinate descent with importance sampling** In view of (6.32), we have that $\mathbf{E}[\theta(S_k, \mathbf{x}^k \mid \mathbf{x}^k)] \geq c := 1/\sum_{i=1}^n M_{ii}$ for smooth functions. This leads to the complexity result

$$K \geq \frac{\xi(\mathbf{x}^0)\sum_{i=1}^n M_{ii}}{\epsilon} \log\left(\frac{\xi(\mathbf{x}^0)}{\epsilon}\right) \quad \Rightarrow \quad \left((\mathbf{E}[\xi(\mathbf{x}^K)] \leq \epsilon) \vee (\exists k \in [K] : \lambda(\mathbf{x}^k) \leq \epsilon)\right).$$

**Minibatch coordinate descent.** We have the bound $\mathbf{E}[\theta(S, \mathbf{x}) \mid \mathbf{x}] \geq \lambda_{\min}\left(\mathbf{E}\left[\mathbf{M}_{[S]}^{-1}\right]\right)$ in the smooth case from (6.35). Substituting into (6.64), we get the rate

$$K \geq \frac{\xi(\mathbf{x}^0)}{\epsilon \lambda_{\min}(\mathbf{E}[\mathbf{M}_{[S]}^{-1}])} \log\left(\frac{\xi(\mathbf{x}^0)}{\epsilon}\right) \quad \Rightarrow \quad \left((\mathbf{E}[\xi(\mathbf{x}^K)] \leq \epsilon) \vee (\exists k \in [K] : \lambda(\mathbf{x}^k) \leq \epsilon)\right)$$

for smooth objectives. In the non-smooth case we can use the bound (6.46) which is given by $\mathbf{E}[\theta(S, \mathbf{x}) \mid \mathbf{x}] \geq \tau/nL_\tau$, where $L_\tau$ is given by (6.41). It follows that we get the rate

$$K \geq \frac{\xi(\mathbf{x}^0)nL_\tau}{\tau\epsilon} \log\left(\frac{\xi(\mathbf{x}^0)}{\epsilon}\right) \quad \Rightarrow \quad \left((\mathbf{E}[\xi(\mathbf{x}^K)] \leq \epsilon) \vee (\exists k \in [K] : \lambda(\mathbf{x}^k) \leq \epsilon)\right)$$

for non-smooth functions.

**Greedy coordinate descent.** In the smooth case we have the bound (6.34), i.e., $\theta(S_k, \mathbf{x}^k) \geq c := 1/\sum_{i=1}^{n} M_{ii}$. This leads to the rate

$$K \geq \frac{\xi(\mathbf{x}^0)\sum_{i=1}^{n} M_{ii}}{\epsilon} \log\left(\frac{\xi(\mathbf{x}^0)}{\epsilon}\right) \quad \Rightarrow \quad \left((\xi(\mathbf{x}^K) \leq \epsilon) \vee (\exists k \in [K] : \lambda(\mathbf{x}^k) \leq \epsilon)\right)$$

for smooth functions. In the proximal case, we have a bound (6.45) given by $\theta(S, \mathbf{x}^k) \geq c := 1/n \max_i\{M_{ii}\}$ which leads to the rate

$$K \geq \frac{\max_i\{M_{ii}\}\xi(\mathbf{x}^0)}{\epsilon} \log\left(\frac{\xi(\mathbf{x}^0)}{\epsilon}\right) \quad \Rightarrow \quad \left((\xi(\mathbf{x}^K) \leq \epsilon) \vee (\exists k \in [K] : \lambda(\mathbf{x}^k) \leq \epsilon)\right)$$

for non-smooth functions.

**Greedy minibatches.** In the smooth case we have the bound (6.38) which leads to the rate

$$K \geq \frac{\xi(\mathbf{x}^0)}{\epsilon \lambda_{\min}(\mathbf{E}[\mathbf{M}_{[S]}^{-1}])} \log\left(\frac{\xi(\mathbf{x}^0)}{\epsilon}\right) \quad \Rightarrow \quad \left((\xi(\mathbf{x}^K) \leq \epsilon) \vee (\exists k \in [K] : \lambda(\mathbf{x}^k) \leq \epsilon)\right)$$

In the non-smooth case we have the bound on the proportion function given by (6.48) which leads to the rate

$$K \geq \frac{\xi(\mathbf{x}^0)nL_\tau}{\tau\epsilon} \log\left(\frac{\xi(\mathbf{x}^0)}{\epsilon}\right) \quad \Rightarrow \quad \left((\xi(\mathbf{x}^K) \leq \epsilon) \vee (\exists k \in [K] : \lambda(\mathbf{x}^k) \leq \epsilon)\right)$$

Uniform minibatch coordinate descent has bounds of the same form, except that the expectation is missing due to the deterministic nature of the method.

Also, both of the above rate are new, as greedy minibatches was introduced as a novel sampling approach.

## 6.8 Experiments

In this section we show results of sample numerical experiments. We will focus on showcasing the theory of the general non-convex optimization methods introduced in Section 6.7.

### 6.8.1 Setup

For our experiments, we consider a function defined as in (6.13) with $f$ defined as

$$f(\mathbf{x}) \quad := \quad \frac{1}{2m}\|\mathbf{A}\mathbf{x} - \mathbf{b}\|^2 + \frac{1}{m}\cos(\langle \mathbf{c}, \mathbf{x}\rangle) \tag{6.65}$$

with $\mathbf{A} \in \mathbb{R}^{m \times n}$ and $\mathbf{b}, \mathbf{c} \in \mathbb{R}^n$, and the function $g$ defined as

$$g(\mathbf{x}) \quad := \quad \lambda\|\mathbf{x}\|_1. \tag{6.66}$$

While it is not motivated by any specific problem, it is clearly non-convex and we can easily control it to observe the behavior of the proposed method.

### 6.8.2 Global convergence of serial coordinate descent

In the first part, we consider the setup defined in the above section using (6.65) and (6.66) with $m = 1000$ and $n = 100$. We generate $\mathbf{A} \in \mathbb{R}^{m \times n}$ as a random matrix with fixed singular values linearly spaced between $\frac{1}{m}$ and 1. The vector $\mathbf{b} \in \mathbb{R}^m$ is set to $\mathbf{b} = \mathbf{A}\mathbf{y}$ for a vector $\mathbf{y} \in \mathbb{R}^n$ randomly generated from a normalized gaussian distribution. Similarly, $\mathbf{c}$ is also randomly generated from a normalized gaussian. We consider two different problems, based on the value of $\lambda$. We have a smooth problem for $\lambda = 0$ and a non-smooth problem for $\lambda = \frac{1}{2m}$.

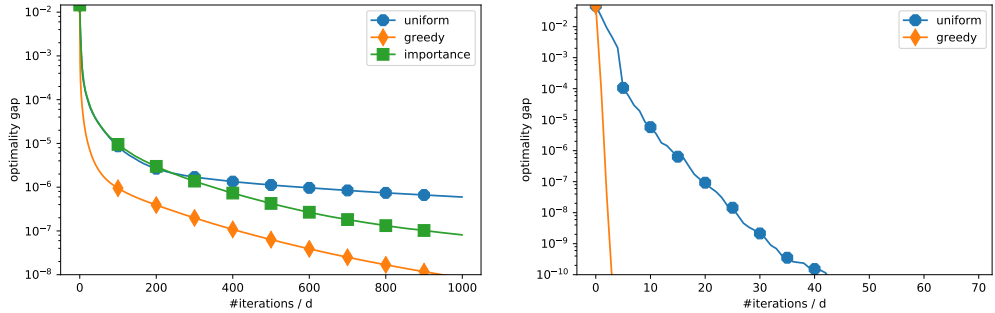

Figure 6.1: Plots of the convergences of various sampling methods for the smooth (left) and non-smooth (right) experiment.

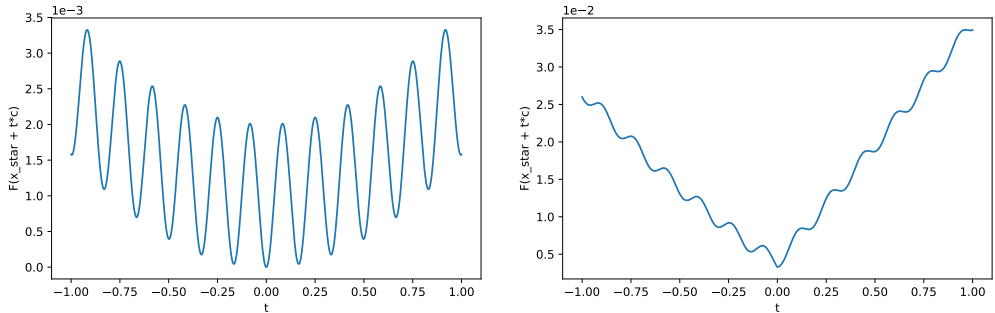

Figure 6.2: Plots of the functions in a 1-dimensional slice around the optimum for the smooth (left) and non-smooth (right) experiment.

We measure the performance of the various coordinate descent approaches described in Section 6.7.1. The convergence behaviors corresponding to the smooth and non-smooth setup can be found on Figure 6.1, on the left and right plots respectively.

To make sure that the functions being optimized are indeed non-convex, we plotted a 1-dimensional slice of the function being optimized around its optimum. These plots can be found on Figure 6.2, for the smooth and non-smooth version respectively.

### 6.8.3 Local convergence of the gradient

The theory of the non-convex case in Theorem 6.20 does not always guarantee convergence to the global optimum, but it at least guarantees a convergence of the magnitude of the gradient. To showcase this scenario, we focused on a 1-dimensional instance of the smooth problem defined in (6.65). We set $\mathbf{A} = [1]$, $\mathbf{b} = \frac{\pi}{c}$ and we have chosen $c$ such that the function $f$ has a flat inflection point ($c \approx 2.15$). It is trivial to show, that the optimal value is 0 and it is achieved for $x = \frac{\pi}{c} \approx 1.46$. The shape of the function around the optimum can be found on the left plot of Figure 6.3.

We used a 1-dimensional gradient descent method for the convergence and we reported on three quantities: The function suboptimality ($fx$), the magnitude of the gradient ($dfx$) and the rate predicted by the theory ($rate$). Theory states that either the function suboptimality or the magnitude of the gradient has to be below the predicted rate, which is what we observe in the right plot of Figure 6.3 as well. Note that the theory focuses on worst case bounds, which is possibly the case why the difference between the rate and the magnitudes is so huge.

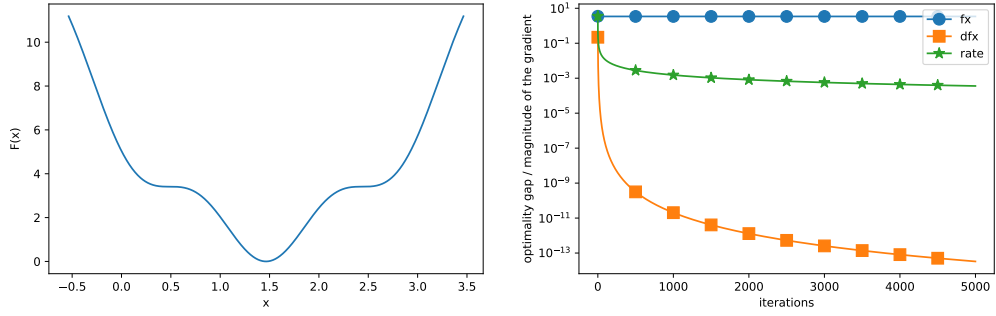

Figure 6.3: Plot of a one-dimensional function with a plateau being optimized (left) and the convergence of the magnitude of the gradient descent with its corresponding rate (right).

## 6.A  Some basic properties of WPL functions

Here we establish some basic properties of smooth WPL functions.

**Proposition 6.21.** *Assume $f \in \mathcal{W}_{PL}(\mu)$ for some $\mu > 0$. Then*

(i) $\mathcal{X}^* = \{\mathbf{x} \in \mathbb{R}^n \ : \ \nabla f(\mathbf{x}) = 0\}$

(ii) $\mathcal{X}^* = \{\mathbf{x} \in \mathbb{R}^n \ : \ \xi(\mathbf{x}) = 0\}$

(iii) *If $\nabla f(\mathbf{x})$ is continuous, then $\mathcal{X}^*$ is convex.*

(iv) *Assume $\nabla f(\mathbf{x})$ is continuous. Then $f$ is non-decreasing on all rays emanating from $\mathbf{x}^*$. That is, $\psi(t) := f(\mathbf{x}^* + t(\mathbf{x} - \mathbf{x}^*))$ is non-decreasing on $t \geq 0$ for all $\mathbf{x} \in \mathbb{R}^n$.*

**Proposition 6.22.** *Assume that $f$ has continuous gradient. If there exists a constant $c > 0$ such that*
$$\int_0^1 \|\nabla f(\mathbf{x}^* + t(\mathbf{x} - \mathbf{x}^*))\| \, dt \ \leq \ c \|\nabla f(\mathbf{x})\|, \qquad \mathbf{x} \in \mathbb{R}^n,$$
*then $f \in \mathcal{W}_{PL}(\frac{1}{c^2})$.*

*Proof.* Using the fundamental theorem of calculus, Cauchy-Schwartz inequality, and then applying the assumption, we get
$$\begin{aligned}
\xi(\mathbf{x}) = f(\mathbf{x}) - f(\mathbf{x}^*) \ &= \ \int_0^1 \langle \nabla f(\mathbf{x}^* + t(\mathbf{x} - \mathbf{x}^*)), \mathbf{x} - \mathbf{x}^* \rangle \, dt \\
&\leq \ \int_0^1 \|\nabla f(\mathbf{x}^* + t(\mathbf{x} - \mathbf{x}^*))\| \cdot \|\mathbf{x} - \mathbf{x}^*\| \, dt \\
&\leq \ c \|\nabla f(\mathbf{x})\| \|\mathbf{x} - \mathbf{x}^*\|.
\end{aligned}$$
□

Let us shed light on the above result. If $f$ is convex, then the directional derivative $\psi(t) := \langle \nabla f(\mathbf{x}^* + t(\mathbf{x} - \mathbf{x}^*)), \mathbf{x} - \mathbf{x}^* \rangle$ is an increasing function of $t$, and hence can be bounded above on $[0, 1]$ by $\psi(1) = \langle \nabla f(\mathbf{x}), \mathbf{x} - \mathbf{x}^* \rangle \leq \|\nabla f(\mathbf{x})\| \|\mathbf{x} - \mathbf{x}^*\|$. It follows that $f \in \mathcal{W}_{PL}(1)$, which we already know.

**Theorem 6.23.** *The following hold:*

1. *If $\mu_1 \geq \mu_2$, then $\mathcal{W}_{PL}(\mu_1) \subseteq \mathcal{W}_{PL}(\mu_2)$.*

2. *If $f$ is convex, then $f \in \mathcal{W}_{PL}(\mu)$ for all $\mu \leq 1$.*

3. *If $f \in \mathcal{W}_{PL}(\mu)$, then $af + b \in \mathcal{W}_{PL}(\mu)$ for all $a \geq 0$ and $b \in \mathbb{R}$.*

4. If $f \in \mathcal{C}^1(L)$, then $\mathcal{S}_{PL}(\mu) \subseteq \mathcal{W}_{PL}(4\mu/L)$.

5. Let $f \in \mathcal{C}^1(L)$, fix $\mathbf{x}^* \in \mathcal{X}^*$, and assume that there exists a constant $c > 0$ such that $\|\mathbf{x} - \mathbf{x}^*\| \leq c\|\nabla f(\mathbf{x})\|$ for all $\mathbf{x} \in \mathbb{R}^n$. Then $f \in \mathcal{W}_{PL}(\frac{4}{L^2 c^2})$.

6. Assume that there exists a constant $c > 0$ such that $\|\nabla f(\mathbf{x})\| \leq c$ for all $\mathbf{x} \in \mathbb{R}^n$. If $f \in \mathcal{W}_{PL}(\mu)$, then $\xi(\mathbf{x}) \leq \frac{c}{\sqrt{\mu}}\|\mathbf{x} - \mathbf{x}^*\|$ for all $\mathbf{x} \in \mathbb{R}^n$. That is, $f$ is Lipschitz on each ray emanating from $\mathbf{x}^*$ with Lipschitz constant $\frac{c}{\sqrt{\mu}}$.

7. Assume $f \in \mathcal{W}_{PL}(\mu)$ with $\mu > 0$. If
$$\frac{\langle \nabla f(\mathbf{x}), \mathbf{x} - \mathbf{x}^* \rangle}{\|\nabla f(\mathbf{x})\|\|\mathbf{x} - \mathbf{x}^*\|} \geq \frac{1}{\sqrt{\mu}},$$
then $f$ satisfies the restricted convexity property:
$$f(\mathbf{x}^*) \geq f(\mathbf{x}) + \langle \nabla f(\mathbf{x}), \mathbf{x}^* - \mathbf{x} \rangle, \qquad \mathbf{x} \in \mathbb{R}^n.$$

*Proof.* 1. Obvious.

2. This was established in the introduction for $\mu = 1$. It only remains to apply 1) to conclude that 2) holds for all $\mu \leq 1$.

3. Obvious.

4. If $f \in \mathcal{S}_{PL}(\mu)$, then
$$\frac{\|\nabla f(\mathbf{x})\|^2 \cdot \|\mathbf{x} - \mathbf{x}^*\|^2}{\xi^2(\mathbf{x})} \overset{(6.2)}{\geq} \frac{2\mu\|\mathbf{x} - \mathbf{x}^*\|^2}{\xi(\mathbf{x})} \overset{(6.4)}{\geq} \frac{4\mu}{L}.$$

5. We have $\xi(\mathbf{x}) \leq \frac{L}{2}\|\mathbf{x} - \mathbf{x}^*\|^2$ using smoothness. Combining this with the assumption $\|\mathbf{x} - \mathbf{x}^*\| \leq c\|\nabla f(\mathbf{x})\|$ from the claim we will prove $f \in \mathcal{W}_{PL}(\frac{4}{L^2 c^2})$ from the definition (6.7)
$$\frac{2}{Lc}\xi(\mathbf{x}) \leq \frac{1}{c}\|\mathbf{x} - \mathbf{x}^*\|^2 \leq \|\nabla f(\mathbf{x})\| \cdot \|\mathbf{x} - \mathbf{x}^*\|.$$

6. Directly using the definition of $f \in \mathcal{W}_{PL}(\mu)$ with the assumption $\|\nabla f(\mathbf{x})\| \leq c$ we get
$$\xi(\mathbf{x}) \overset{(6.7)}{\leq} \frac{1}{\sqrt{\mu}}\|\nabla f(\mathbf{x})\|\|\mathbf{x} - \mathbf{x}^*\| \leq \frac{c}{\sqrt{\mu}}\|\mathbf{x} - \mathbf{x}^*\|.$$

7. Combining the definition of $f \in \mathcal{W}_{PL}(\mu)$ (6.7) with the assumption from the claim we get that
$$\sqrt{\mu}\xi(\mathbf{x}) \overset{(6.7)}{\leq} \|\nabla f(\mathbf{x})\| \cdot \|\mathbf{x} - \mathbf{x}^*\| \leq \sqrt{\mu}\langle \nabla f(\mathbf{x}), \mathbf{x} - \mathbf{x}^* \rangle.$$
Dividing both sides by $\sqrt{\mu}$ and adding $f(\mathbf{x}^*)$ we get the restricted convexity. □

## 6.B Proofs

### 6.B.1 Proof of Lemma 6.13

*Proof.* From the definition of Algorithm 8 we have that $\mathbf{x}^+ = \mathbf{x} + \mathbf{u}^*_{[S]}$, where

$$\mathbf{u}^* \quad := \quad \underset{\mathbf{u} \in \mathbb{R}^n}{\arg\min} \left\{ \langle \nabla_{[S]} f(\mathbf{x}), \mathbf{u} \rangle + \frac{1}{2} \mathbf{u}^\top \mathbf{M}_{[S]} \mathbf{u} + \sum_{i \in S} [g_i(x_i + u_i) - g_i(x_i)] \right\} \quad (6.67)$$

It follows that

$$\begin{aligned}
\xi(\mathbf{x}^+) &= \xi(\mathbf{x} + \mathbf{u}^*) \\
&\stackrel{(6.19)}{=} F(\mathbf{x} + \mathbf{u}^*) - F(\mathbf{x}^*) \\
&\stackrel{(6.13)}{=} f(\mathbf{x} + \mathbf{u}^*_{[S]}) + g(\mathbf{x} + \mathbf{u}^*_{[S]}) - F(\mathbf{x}^*) \\
&\stackrel{(6.14)}{\leq} f(\mathbf{x}) + \langle \nabla f(\mathbf{x}), \mathbf{u}^*_{[S]} \rangle + \frac{1}{2} (\mathbf{u}^*_{[S]})^\top \mathbf{M} \mathbf{u}^*_{[S]} + g(\mathbf{x} + \mathbf{u}^*_{[S]}) - F(\mathbf{x}^*) \\
&\stackrel{(6.13)}{=} F(\mathbf{x}) + \langle \nabla f(\mathbf{x}), \mathbf{u}^*_{[S]} \rangle + \frac{1}{2} (\mathbf{u}^*_{[S]})^\top \mathbf{M} \mathbf{u}^*_{[S]} + g(\mathbf{x} + \mathbf{u}^*_{[S]}) - g(\mathbf{x}) - F(\mathbf{x}^*) \\
&\stackrel{(6.19)}{=} \xi(\mathbf{x}) + \langle \nabla f(\mathbf{x}), \mathbf{u}^*_{[S]} \rangle + \frac{1}{2} (\mathbf{u}^*_{[S]})^\top \mathbf{M} \mathbf{u}^*_{[S]} + g(\mathbf{x} + \mathbf{u}^*_{[S]}) - g(\mathbf{x}) \\
&\stackrel{(6.15)}{=} \xi(\mathbf{x}) + \langle \nabla f(\mathbf{x}), \mathbf{u}^*_{[S]} \rangle + \frac{1}{2} (\mathbf{u}^*_{[S]})^\top \mathbf{M} \mathbf{u}^*_{[S]} + \sum_{i \in S} [g_i(x_i + u_i^*) - g_i(x_i)] \\
&\stackrel{(6.67)}{=} \xi(\mathbf{x}) + \min_{\mathbf{u} \in \mathbb{R}^n} \left\{ \langle \nabla_{[S]} f(\mathbf{x}), \mathbf{u} \rangle + \frac{1}{2} \mathbf{u}^\top \mathbf{M}_{[S]} \mathbf{u} + \sum_{i \in S} [g_i(x_i + u_i) - g_i(x_i)] \right\} \\
&\stackrel{(6.24)}{\leq} \xi(\mathbf{x}) + \theta(S, \mathbf{x}) \cdot \min_{\mathbf{y} \in \mathbb{R}^n} \left\{ \langle \nabla f(\mathbf{x}), \mathbf{y} \rangle + \frac{L}{2} \|\mathbf{y}\|^2 + g(\mathbf{x} + \mathbf{y}) - g(\mathbf{x}) \right\} \\
&\stackrel{(6.21)}{=} [1 - \theta(S, \mathbf{x}) \cdot \mu(\mathbf{x})] \cdot \xi(\mathbf{x})
\end{aligned}$$

Note that the inequality in the one-to-last line might happen in the case, when $\mathbf{x} \notin \mathcal{X}$ in the definition of the proportion function in (6.24). Chaining up the resulting expressions from $\mathbf{x}^K$ all the way back to $\mathbf{x}^0$ proves the $K$-step bound. $\square$

### 6.B.2 Proof of Theorem 6.15

*Proof.* Using the result of Lemma 6.24, we have the bound (6.71). To get the result, we will use the bound $(a + b)^2 \geq 4ab$, which holds for any $a, b > 0$. Specifically, we will use it for $a = \xi(\mathbf{x})$ and $b = \frac{\lambda_F}{2} \|\mathbf{x} - \mathbf{x}^*\|^2$ in the expression in (6.71) to get

$$-L \cdot \min_{\mathbf{y} \in \mathbb{R}^n} \left\{ \langle \nabla f(\mathbf{x}), \mathbf{y} \rangle + \frac{L}{2} \|\mathbf{y}\|^2 + g(\mathbf{x} + \mathbf{y}) - g(\mathbf{x}) \right\}$$
$$\geq \xi(\mathbf{x}) \cdot \min \left\{ \frac{L}{2}, \frac{L \lambda_F}{\lambda_F - \lambda_f + L} \right\},$$

which is the claimed result. $\square$

### 6.B.3 Proof of Theorem 6.16

*Proof.* The one-step bound (6.26) combined with the definition of the class $\mathcal{S}^g_{PL}$ (6.7) gives

$$\xi(\mathbf{x}^{k+1}) \stackrel{(6.26)}{\leq} \left(1 - \theta(S_k, \mathbf{x}^k) \cdot \mu(\mathbf{x}^k)\right) \cdot \xi(\mathbf{x}^k) \stackrel{(6.7)}{\leq} \left(1 - \theta(S_k, \mathbf{x}^k) \cdot \mu\right) \cdot \xi(\mathbf{x}^k).$$

Now, by taking full expectation over the whole sampling procedure on both sides and using the definition of $\mu_k$ (6.53) we get

$$\begin{aligned} \mathbf{E}[\xi(\mathbf{x}^{k+1})] &\leq \mathbf{E}\left[(1 - \theta(S_k, \mathbf{x}^k)\mu) \cdot \xi(\mathbf{x}^k)\right] \\ &= \mathbf{E}[\xi(\mathbf{x}^k)] - \mu \mathbf{E}[\theta(S_k, \mathbf{x}^k) \cdot \xi(\mathbf{x}^k)] \\ &\stackrel{(6.53)}{=} (1 - \mu_k) \cdot \mathbf{E}[\xi(\mathbf{x}^k)]. \end{aligned}$$

To establish the convergence rate (6.54), we use the estimate $1 - s \leq e^{-s}$ to get $\mathbf{E}[\xi(\mathbf{x}^{k+1})] \leq e^{-\mu_k} \mathbf{E}[\xi(\mathbf{x}^k)]$, which we can simply chain together repeatedly to get

$$\mathbf{E}[\xi(\mathbf{x}^{k+1})] \leq e^{-\sum_{k=0}^{K-1} \mu_k} \xi(\mathbf{x}^0).$$

Setting the right-hand side less or equal to $\epsilon$ and rearranging we finally get (6.54). □

### 6.B.4 Proof of Theorem 6.18

*Proof.* Using the result of Lemma 6.24, we have the bound (6.71). Plugging $\lambda_F = \lambda_f = 0$ into the expression in (6.71) we get

$$-L \cdot \min_{\mathbf{y} \in \mathbb{R}^n} \left\{ \langle \nabla f(\mathbf{x}), \mathbf{y} \rangle + \frac{L}{2} \|\mathbf{y}\|^2 + g(\mathbf{x} + \mathbf{y}) - g(\mathbf{x}) \right\}$$
$$\geq \xi^2(\mathbf{x}) \cdot \min \left\{ \frac{L}{2\xi(\mathbf{x})}, \frac{1}{2\|\mathbf{x} - \mathbf{x}^*\|^2} \right\},$$

which is the first part of the claimed result. As for the second part, we can directly bound $\xi(\xi) \leq \xi(\mathbf{x}^0)$, as these are the only pairs of $\mathbf{x}, \mathbf{x}^0$ we need to consider according to Definition 6.17. Similarly, as the function value is bounded, the quantities $\|\mathbf{x} - \mathbf{x}^*\|$ can be upper bounded by the largest distance in the level set of $f(\mathbf{x}^0)$, which is given by $R$ in (6.60). Combining these arguments, we get that

$$\mu(\mathbf{x}) \geq \xi(\mathbf{x}) \cdot \min \left\{ \frac{L}{2\xi(\mathbf{x})}, \frac{1}{2\|\mathbf{x} - \mathbf{x}^*\|^2} \right\} \geq \xi(\mathbf{x}) \cdot \min \left\{ \frac{L}{2\xi(\mathbf{x}^0)}, \frac{1}{2R^2} \right\},$$

which is the claimed result. □

### 6.B.5 Proof of Theorem 6.19

*Proof.* Combining the one-step bound (6.26) with the definition of the class $\mathcal{W}_{PL}^g$ (6.57) we get that

$$\begin{aligned} \xi(\mathbf{x}^{k+1}) &\stackrel{(6.26)}{\leq} [1 - \theta(S_k, \mathbf{x}^k) \cdot \mu(\mathbf{x}^k)] \cdot \xi(\mathbf{x}^k) \\ &\stackrel{(6.57)}{\leq} [1 - \theta(S_k, \mathbf{x}^k) \cdot \rho(\mathbf{x}^0) \cdot \xi(\mathbf{x}^k)] \cdot \xi(\mathbf{x}^k) \end{aligned}$$

Now, taking full expectation over the whole sampling process on both sides and using the definition of $\mu_k$ (6.61) we get

$$\begin{aligned} \mathbf{E}[\xi(\mathbf{x}^{k+1})] &\leq \mathbf{E}\left[(1 - \theta(S_k, \mathbf{x}^k) \cdot \rho(\mathbf{x}^0) \cdot \xi(\mathbf{x}^k)) \cdot \xi(\mathbf{x}^k)\right] \\ &= \mathbf{E}[\xi(\mathbf{x}^k)] - \rho(\mathbf{x}^0) \mathbf{E}[\theta(S_k, \mathbf{x}^k) \cdot (\xi(\mathbf{x}^k))^2] \\ &\stackrel{(6.61)}{=} (1 - \mu_k \mathbf{E}[\xi(\mathbf{x}^k)]) \cdot \mathbf{E}[\xi(\mathbf{x}^k)]. \end{aligned}$$

Observe, that $\{\mu_k\}_{k=0}^K$ and $\{\mathbf{E}[\xi(\mathbf{x}^k)]\}_{k=0}^K$ are both positive scalars, therefore we can use Lemma 6.4 to get the bound

$$\mathbf{E}[\xi(\mathbf{x}^K)] \leq \frac{\xi(\mathbf{x}^0)}{1 + \xi(\mathbf{x}^0) \sum_{k=0}^{K-1} \mu_k}.$$

Putting the right-hand side less than $\epsilon$ and rearranging leads to the claimed result (6.62). $\square$

### 6.B.6 Proof of Theorem 6.20

*Proof.* If part $(i)$ from the claim holds, we are done. Now assume on the contrary, that $(i)$ does not hold, i.e.,

$$-L \cdot \min_{\mathbf{y} \in \mathbb{R}^n} \left\{ \langle \nabla f(\mathbf{x}^k), \mathbf{y} \rangle + \frac{L}{2} \|\mathbf{y}\|^2 + g(\mathbf{x}^k + \mathbf{y}) - g(\mathbf{x}^k) \right\} \geq \epsilon \qquad (6.68)$$

for all $k \in \{0, \ldots, K-1\}$. It follows, that

$$\mu(\mathbf{x}^k) \stackrel{(6.22)}{=} \frac{-L \cdot \min_{\mathbf{y} \in \mathbb{R}^n} \left\{ \langle \nabla f(\mathbf{x}^k), \mathbf{y} \rangle + \frac{L}{2} \|\mathbf{y}\|^2 + g(\mathbf{x}^k + \mathbf{y}) - g(\mathbf{x}^k) \right\}}{\xi(\mathbf{x}^k)} \stackrel{(6.68)}{\geq} \frac{\epsilon}{\xi(\mathbf{x}^k)}.$$

for all $k \in \{0, \ldots, K-1\}$. Using the result from Lemma 6.13, we have that $\xi(\mathbf{x}^k) \leq \xi(\mathbf{x}^0)$ for all $k$, which we can use to further bound

$$\mu(\mathbf{x}^k) \geq \frac{\epsilon}{\xi(\mathbf{x}^0)}. \qquad (6.69)$$

Using (6.26) combined with the above result (6.69) we get

$$\xi(\mathbf{x}^K) \stackrel{(6.26)}{\leq} \left[1 - \mu(\mathbf{x}^{K-1}) \cdot \theta(S_{K-1}, \mathbf{x}^{K-1})\right] \xi(\mathbf{x}^{K-1})$$
$$\stackrel{(6.69)}{\leq} \left[1 - \frac{\epsilon \theta(S_{K-1}, \mathbf{x}^{K-1})}{\xi(\mathbf{x}^0)}\right] \xi(\mathbf{x}^{K-1})$$

Taking the expectation over the whole sampling procedure on both sides and using the definition of $\mu_k$ in (6.63) we get

$$\mathbf{E}[\xi(\mathbf{x}^K)] \leq \mathbf{E}\left[\left(1 - \frac{\epsilon \theta(S_{K-1}, \mathbf{x}^{K-1})}{\xi(\mathbf{x}^0)}\right) \xi(\mathbf{x}^{K-1})\right]$$
$$= \mathbf{E}[\xi(\mathbf{x}^{K-1})] - \frac{\epsilon}{\xi(\mathbf{x}^0)} \mathbf{E}[\theta(S_{K-1}, \mathbf{x}^{K-1}) \cdot \xi(\mathbf{x}^{K-1})]$$
$$\stackrel{(6.63)}{=} \left(1 - \frac{\epsilon \mu_{K-1}}{\xi(\mathbf{x}^0)}\right) \mathbf{E}[\xi(\mathbf{x}^{K-1})] \qquad (6.70)$$

Combining the above inequality (6.70) with $(1-z) \leq \exp(-z)$ repeatedly, we get

$$\mathbf{E}[\xi(\mathbf{x}^K)] \stackrel{(6.70)}{\leq} \left(1 - \frac{\epsilon \mu_{K-1}}{\xi(\mathbf{x}^0)}\right) \mathbf{E}[\xi(\mathbf{x}^{K-1})]$$
$$\leq \exp\left(-\frac{\epsilon \mu_{K-1}}{\xi(\mathbf{x}^0)}\right) \mathbf{E}[\xi(\mathbf{x}^{K-1})]$$
$$\stackrel{(6.70)}{\leq} \ldots$$
$$\stackrel{(6.70)}{\leq} \exp\left(-\frac{\epsilon \sum_{k=0}^{K-1} \mu_k}{\xi(\mathbf{x}^0)}\right) \xi(\mathbf{x}^0)$$
$$\stackrel{(6.64)}{\leq} \epsilon,$$

where the last line follows from comparing the logarithms of both sides. This proves $(ii)$. $\square$

## 6.C Technical Lemmas

### 6.C.1 Lemma 6.24

**Lemma 6.24.** *Let $f$ be $\lambda_f$-strongly convex (6.51) with $\lambda_f \geq 0$ and $F$ be $\lambda_F$-strongly convex (6.50) with $\lambda_F \geq 0$. Then*

$$-\min_{\mathbf{y}\in\mathbb{R}^n}\left\{\langle\nabla f(\mathbf{x}),\mathbf{y}\rangle + \frac{L}{2}\|\mathbf{y}\|^2 + g(\mathbf{x}+\mathbf{y}) - g(\mathbf{x})\right\}$$
$$\geq \min\left\{\frac{1}{2}\xi(\mathbf{x}), \frac{\left(\xi(\mathbf{x}) + \frac{\lambda_F}{2}\|\mathbf{x}-\mathbf{x}^*\|^2\right)^2}{2(\lambda_F - \lambda_f + L)\|\mathbf{x}-\mathbf{x}^*\|^2}\right\}. \quad (6.71)$$

*Proof.* Let

$$\beta = \min\left\{1, \frac{\xi(\mathbf{x}) + \frac{\lambda_F}{2}\|\mathbf{x}-\mathbf{x}^*\|^2}{(\lambda_F - \lambda_f + L)\|\mathbf{x}-\mathbf{x}^*\|^2}\right\}. \quad (6.72)$$

Observe, that $\beta = 1$ implies, that

$$(\lambda_F - \lambda_f + L)\|\mathbf{x}-\mathbf{x}^*\|^2 \leq \xi(\mathbf{x}) + \frac{\lambda_F}{2}\|\mathbf{x}-\mathbf{x}^*\|^2,$$

from which it follows that

$$\frac{\lambda_f - L}{2}\|\mathbf{x}-\mathbf{x}^*\|^2 \geq -\frac{1}{2}\xi(\mathbf{x}) + \frac{\lambda_F}{4}\|\mathbf{x}-\mathbf{x}^*\|^2 \geq -\frac{1}{2}\xi(\mathbf{x}) \quad (6.73)$$

Now, it follows that

$$-\min_{\mathbf{y}\in\mathbb{R}^n}\left\{\langle\nabla f(\mathbf{x}),\mathbf{y}\rangle + \frac{L}{2}\|\mathbf{y}\|^2 + g(\mathbf{x}+\mathbf{y}) - g(\mathbf{x})\right\}$$
$$\stackrel{(6.13)}{=} F(\mathbf{x}) - \min_{\mathbf{y}\in\mathbb{R}^n}\left\{f(\mathbf{x}) + \langle\nabla f(\mathbf{x}),\mathbf{y}\rangle + \frac{L}{2}\|\mathbf{y}\|^2 + g(\mathbf{x}+\mathbf{y})\right\}$$
$$\stackrel{(6.51)}{\geq} F(\mathbf{x}) - \min_{\mathbf{y}\in\mathbb{R}^n}\left\{f(\mathbf{x}+\mathbf{y}) - \frac{\lambda_f}{2}\|\mathbf{y}\|^2 + \frac{L}{2}\|\mathbf{y}\|^2 + g(\mathbf{x}+\mathbf{y})\right\}$$
$$\stackrel{(6.13)}{=} F(\mathbf{x}) - \min_{\mathbf{y}\in\mathbb{R}^n}\left\{F(\mathbf{x}+\mathbf{y}) + \frac{L-\lambda_f}{2}\|\mathbf{y}\|^2\right\} \quad (\text{let } \mathbf{y} = \beta(\mathbf{x}^* - \mathbf{x}))$$
$$\geq F(\mathbf{x}) - F(\beta\mathbf{x}^* + (1-\beta)\mathbf{x}) - \frac{\beta^2(L-\lambda_f)}{2}\|\mathbf{x}-\mathbf{x}^*\|^2$$
$$\stackrel{(6.50)}{\geq} F(\mathbf{x}) - \beta F(\mathbf{x}^*) - (1-\beta)F(\mathbf{x}) + \left(\frac{\lambda_F\beta(1-\beta)}{2} - \frac{\beta^2(L-\lambda_f)}{2}\right)\|\mathbf{x}-\mathbf{x}^*\|^2$$
$$\stackrel{(6.19)}{=} \beta\left(\xi(\mathbf{x}) + \frac{\lambda_F}{2}\|\mathbf{x}-\mathbf{x}^*\|^2\right) - \beta^2\frac{(\lambda_F - \lambda_f + L)}{2}\|\mathbf{x}-\mathbf{x}^*\|^2$$
$$\stackrel{(6.72)}{=} \min\left\{\xi(\mathbf{x}) + \frac{\lambda_f - L}{2}\|\mathbf{x}-\mathbf{x}^*\|^2, \frac{\left(\xi(\mathbf{x}) + \frac{\lambda_F}{2}\|\mathbf{x}-\mathbf{x}^*\|^2\right)^2}{2(\lambda_F - \lambda_f + L)\|\mathbf{x}-\mathbf{x}^*\|^2}\right\}$$
$$\stackrel{(6.73)}{>} \min\left\{\frac{1}{2}\xi(\mathbf{x}), \frac{\left(\xi(\mathbf{x}) + \frac{\lambda_F}{2}\|\mathbf{x}-\mathbf{x}^*\|^2\right)^2}{2(\lambda_F - \lambda_f + L)\|\mathbf{x}-\mathbf{x}^*\|^2}\right\}.$$

□

### 6.C.2 Lemma 6.25

**Lemma 6.25.** *Let $f$ be a given function and let $X, Y$ be random variables for which it holds that $\mathbf{Prob}(f(X) \geq 0) = 1$ and $\mathbf{Prob}(\mathbf{E}[Y \mid X] \geq 0) = 1$. Additionally, let $c > 0$ be a scalar such that*

$$c \leq \mathbf{E}[Y \mid X] \quad (6.74)$$

*Then it holds that*
$$\mathbf{E}[f(X)Y] \geq c\mathbf{E}[f(X)]. \tag{6.75}$$

*Proof.* Using the tower property $\mathbf{E}[X] = \mathbf{E}[\mathbf{E}[X|Y]]$ combined with assumption (6.74) we get

$$\begin{aligned}
\mathbf{E}[f(X)Y] &= \mathbf{E}[\mathbf{E}[f(X)Y \mid X]] = \mathbf{E}[f(X)\mathbf{E}[Y \mid X]] \stackrel{(6.74)}{\geq} \mathbf{E}[f(X)c] \\
&= c\mathbf{E}[f(X)]
\end{aligned}$$

which concludes the proof. $\square$

## 6.D   Notation Glossary

| Notation | Description |
|---|---|
| $\mathbb{R}$ | the set of real numbers |
| $\mathbb{R}_+$ | the set of positive real numbers |
| $\bar{\mathbb{R}}$ | the set $\mathbb{R} \cup \{+\infty\}$ |
| $\mathbf{x}$ | a vector |
| $x_i$ | the $i$-th entry of the vector $\mathbf{x}$ |
| $\mathbf{X}$ | a matrix |
| $\mathbf{X}_{i:}$ | the $i$-th row of the matrix $\mathbf{X}$ |
| $\mathbf{X}_{:j}$ | the $j$-th column of the matrix $\mathbf{X}$ |
| $X_{ij}$ | the entry at the $i$-th row and $j$-th column of the matrix $\mathbf{X}$ |
| $[n]$ | a shorthand for the set $\{1, \ldots, n\}$ |
| $\mathbf{x}_S$ | a $|S|$-dimensional vector containing entries of $\mathbf{x}$ with indices in $S$ |
| $\mathbf{x}_{[S]}$ | the vector $\mathbf{x}$ with the entries with indices outside $S$ zeroed out |
| $\mathbf{X}_S$ | a $|S| \times |S|$ submatrix of $\mathbf{X}$ containing only rows and columns with indices in $S$ |
| $\mathbf{X}_{[S]}$ | the matrix $\mathbf{X}$ with all entries on columns or rows outside of $S$ zeroed out |
| $\mathbf{X}_{[S]}^{-1}$ | the $n \times n$ matrix containing $(\mathbf{X}_S)^{-1}$ at the rows and columns indicated by $S$. |
| $f$ | a smooth function from $\mathbb{R}^n$ to $\mathbb{R}$ (see Def. 6.8) |
| $g$ | a separable (6.9) and possibly non-smooth function from $\mathbb{R}^n$ to $\mathbb{R}$ |
| $F$ | the objective function from $\mathbb{R}^n$ to $\mathbb{R}$, defined as $f + g$ (6.13) |
| $\nabla_S f(\mathbf{x})$ | a shorthand for $(\nabla f(\mathbf{x}))_S$ |
| $\nabla_{[S]} f(\mathbf{x})$ | a shorthand for $(\nabla f(\mathbf{x}))_{[S]}$ |
| $\xi(\mathbf{x})$ | the optimality gap defined in (6.19) |
| $\lambda(\mathbf{x})$ | the auxilary function defined in (6.20) |
| $\mu(\mathbf{x})$ | the forcing function defined in (6.21) |
| $\theta(S, \mathbf{x})$ | the proportion function defined in (6.24) |
| $\mathcal{X}^*$ | the set of global minimizers of $F$ |
| $\mathcal{X}$ | the set of vectors with nonzero $\lambda(\mathbf{x})$ (6.23) |
| $\lambda_{\min}(\mathbf{M})$ | the smallest eigenvalue of a square matrix $\mathbf{M}$ |
| $\lambda_{\max}(\mathbf{M})$ | the largest eigenvalue of a square matrix $\mathbf{M}$ |
| $\mathbf{M}, L\mathbf{I}$ | the $\mathbf{M}$-smoothness matrices of $f$ (6.14) for smooth and non-smooth cases |
| $L_\tau$ | the smoothness parameter defined as $\max_{S:|S|=\tau}\{\lambda_{\max}(\mathbf{M}_S)\}$ |
| $\lambda_F, \lambda_f$ | strong convexity parameters of $F$ and $f$, respectively |
| $U_S(\mathbf{x}, \mathbf{u})$ | the quadratic upper bound function at a given point $\mathbf{x}$, defined in (6.16) |

Table 6.3: Notation Glossary.

# Conclusion and Extensions

In this work we tackled the fundamental problem of empirical risk minimization

$$\mathbf{w}^* \;=\; \arg\min_{\mathbf{w}\in\mathcal{W}}\left\{\frac{1}{n}\sum_{j=1}^{n} f_j(\mathbf{w}) + r(\mathbf{w})\right\},$$

in a few different formulations. In each chapter, we proposed a new approach for solving this problem using stochastic methods, with a focus on the sampling strategies.

In Chapter 1 we introduced the emprical risk minimization problem in detail together with some baseline methods for solving it.

In Chapter 2 we introduced a new adaptive sampling strategy for stochastic dual coordinate ascent. The adaptivity comes from the fact that the method changes its sampling distribution over the course of the iterative process. We show that the proposed sampling is theoretically superior to all previously known samplings, but it is inpractical as it has high implementation cost. For this reason, we include a practical variant, which empirically outperforms other existing approaches.

In Chapter 3 we proposed a new variance-reduced primal method, which we analyzed for arbitrary sampling schemes. The method can be viewed as a dual-free version of SDCA, as it acts on the primal formulation of empirical risk minimization, but has an update reminiscent of SDCA. We also proposed a new sampling method, which aims to help with the problems with synchronizing parallel coordinate descent updates. The main idea of the sampling is to make sure that every core has the same amount of computation to do on each iteration in order to minimize the idle time of the cores.

In Chapter 4 we proposed a novel sampling scheme. The main idea of the new sampling is to combine parallel sampling strategies with importance sampling – giving rise to the first ever importance sampling for minibatches. We give theoretical rates for our new sampling and we also empirically test it against uniform parallel sampling strategies. The empirical results for the parallel setting are very pleasing, as the observed speedup against uniform strategies corresponds to the speedup observed for the importance sampling in the serial case.

In Chapter 5 we compared the performance of coordinate descent and stochastic dual coordinate ascent on the problem of linear ERM. In the case that the data is dense, we confirmed the community belief that primal coordinate descent is advantegous if $d \geq n$ and vice-versa – dual coordinate ascent is better if $n \geq d$. However, we showed that in the case of sparse data, the situation gets more complicated. First, we proved that serial importance sampling is optimal even if we take into account the iteration costs. Second, we showed that the superiority of one of the approaches is connected to a certain data-dependent quantity, which can be computed for each dataset. Additionally, we gave examples of datasets where the dual approach outperforms the primal even if $d \gg n$ and vice-versa.

Finally, in Chapter 6 we introduced the notion of Weak Polyak-Łojasiewicz condition, which can be seen as the Polyak-Łojasiewicz condition for convex functions. Using this notion, we established global convergence guarantees for proximal block descent methods for a large class of previously unconsidered non-convex objectives. Additionally, our analysis allows for greedy, cylic and other deterministic sampling strategies on top of traditional randomized samplings, in which we also include adaptive strategies. We also showed how we recover the state-of-the art guarantees for a large class of functions and sampling strategies as special cases of our analysis.

To conclude, we will now briefly outline possible future extensions to our work. First and foremost, we did not include any accelerated methods in our approaches. This is not a major issue, as acceleration can be achieved by using frameworks as [36, 65, 66]. However, by directly considering accelerated methods in our analysis, we believe that we will get better practical speedup. Some of the parts of this work can be extended to incorporate accelerated methods in a straightforward fashion (e.g. Chapter 5), while others might prove more challenging (e.g. Chapter 2). There is definitely a lot of research to be done in this direction.

Secondly, the challenge of inventing a sampling which theoretically beats importance sampling in serial coordinate descent / serial dual coordinate ascent is still on. Indeed, all the methods which offer better bounds in theory are inpractical, due to their high iteration cost, and only their heuristic versions are outperforming importance sampling in wall clock time. We believe that our contribution in Chapter 2 lays solid groundwork for devising such a sampling strategy.

Lastly, due to the popularity of non-convex objectives appearing in empirical risk minimization and machine learning in general, we find it important to extend the current analysis of gradient descent type methods to an even wider class of objectives. In Chapter 6 we managed to extend the analysis of block descent methods to a new class of objectives, but it is definitely not the final station on this journey. The optimization community is well aware of this research direction and its utmost importance, as there is currently an increasing amount of papers being published on non-convex optimization every year, largely due to the recent popularity of deep learning methods.